\documentclass[preprint,12pt,3p]{elsarticle}

\usepackage{framed,multirow}

\usepackage[utf8]{inputenc}

\usepackage{bm}
\usepackage{amsmath}
\usepackage{mathrsfs}
\usepackage{graphicx}
\usepackage{epsfig, setspace}
\usepackage{url}
\usepackage{graphicx, transparent, color}

\usepackage{etex}
\reserveinserts{28}

\usepackage[bookmarks=false]{hyperref}
\usepackage{amsmath}
\usepackage{xcolor}
\usepackage{rotating} % Rotating the figures but keeps the page in portrait mode
\usepackage{pdflscape} % Rotates single page into Landscape mode
\usepackage{silence}
\usepackage{ulem,cancel}
\newtheorem{remark}{\bf Remark}[section]
\newtheorem{example}{Example}[section]
\usepackage{graphicx}
\usepackage{float}
\usepackage{subfigure}% subcaptions for subfigures
\usepackage{subfigmat}% matrices of similar subfigures,

\newcommand{\mathd}{\mathrm{d}}

\usepackage{setspace,lipsum}

\begin{document}

\begin{frontmatter}

\title{{\color{black}High-Order Central-Upwind shock capturing scheme using a \textcolor{black}{Boundary} Variation Diminishing (BVD) Algorithm}}

% Group authors per affiliation:
\author[AA_address]{Amareshwara Sainadh Chamarthi \cortext[cor1]{Corresponding author. \\ 
E-mail address: s.chamarthi@al.t.u-tokyo.ac.jp (Amareshwara Sainadh  Ch.).}}
\author[AA_address]{Steven H.\ Frankel}
\address[AA_address]{Faculty of Mechanical Engineering, Technion - Israel Institute of Technology, Haifa, Israel}
\begin{abstract}
In this paper, we present a novel hybrid nonlinear explicit-compact scheme for shock-capturing based on a boundary variation diminishing (BVD) reconstruction. In our approach, we combine a non-dissipative sixth-order central compact interpolation and a fifth-order monotonicity preserving scheme (MP5) through the BVD algorithm. For a smooth solution, the BVD reconstruction chooses the highest order possible interpolation, which is central, i.e. non-dissipative in the current approach and for the discontinuities, the algorithm selects the monotone scheme. \textcolor{black}{This method provides an alternative to the existing adaptive upwind-central schemes in the literature.} Several numerical examples are conducted with the present approach, which suggests that the current method is capable of resolving small scale flow features and has the same ability to capture sharp discontinuities as the MP5 scheme. 
\end{abstract}
\begin{keyword}
Compact reconstruction, Boundary variation diminishing, Finite volume method, Low dissipation, Shock capturing
\end{keyword}
\end{frontmatter}

\section{Introduction}\label{sec-1}

Numerical simulations of high-speed turbulent flows involving discontinuities such as shockwaves and material interfaces represents a significant computational challenge. High-order central schemes, which are non-dissipative, can resolve smooth and small-scale features but often introduce high-frequency oscillations (Gibbs phenomenon) near discontinuities. In contrast, the nonlinear limiting typically used for capturing discontinuities  \textcolor{black}{is} too dissipative in relatively smooth turbulent regions. This poses a significant numerical challenge  \textcolor{black}{in designing} a scheme that can effectively treat these contradictory requirements. Considerable efforts have been devoted to the development of numerical schemes that are essentially non-oscillatory near the discontinuities and low-dissipative in the smooth regions.

One of the first attempts to develop such methods is the piecewise parabolic method (PPM), proposed by Woodward and Collela \cite{woodward1984numerical}, which uses a four-point centred polynomial to approximate the interface values. Limiters then correct the interface values to achieve non-oscillatory results. A well-known drawback of these slope limiters is that they tend to clip the smooth extrema of the flow and the accuracy degenerates to first-order. Collela and Sekora have proposed an improvement of PPM that can preserve the smooth extrema in Ref. \cite{colella2008limiter}. Suresh and Hyunh \cite{suresh1997accurate} argued that the three-point stencil typically used in total variation diminishing (TVD) methods is not able to distinguish between local extrema and a genuine discontinuity. They have proposed a monotonicity-preserving (MP) scheme that uses a five (or more) point stencil to make this distinction.

In contrast, Harten et al. \cite{harten1987uniformly} have proposed an essentially non-oscillatory (ENO) approach that will adaptively choose the appropriate stencil containing the smoothest data, thereby avoiding interpolations across discontinuities.  \textcolor{black}{ENO schemes require so many logical conditions for their stencil choosing procedure which make them ill-suited for faster computations.} The weighted essentially non-oscillatory (WENO) schemes, originated by Liu et al. \cite{liu1994weighted} overcame the problems associated with the ENO scheme and were further significantly improved by Jiang and Shu \cite{jiang1995}. WENO schemes achieve high-order accuracy by analyzing all possible stencils through a convex combination enabling the scheme to choose the highest order interpolation for smooth flows and one-sided interpolations in the vicinity of discontinuities. Despite their success, WENO schemes are overly dissipative with regard to resolving  small scale solution features. Various modifications of WENO schemes have been proposed over the years to improve their accuracy and reduce their dissipative nature. Henrick et al. \cite{Henrick2005} have noticed that the WENO scheme proposed by Jiang and Shu \cite{jiang1995} losses its accuracy at critical points and proposed the WENO-M scheme as a fix. Borges et al. \cite{Borges2008} proposed another alternative approach for improving the accuracy at these critical points by introducing a global smoothness indicator for the stencils. 

Pirozzoli \cite{Pirozzoli2002} has proposed a hybrid compact-WENO scheme which combines the conservative compact scheme and the WENO scheme. Compact schemes are implicit in space and have superior dissipation and dispersive properties over explicit schemes of the same order, which is beneficial in smooth regions \cite{lele1992compact}.  Ren et al. \cite{ren2003characteristic} have improved the hybrid scheme of Pirozzoli through a continuous weighing function which assisted in a smooth transition from one sub-scheme to the other. Kim and Kwon \cite{kim2005high} proposed a hybrid central-WENO approach through the weighing function of Ren et al..  Their scheme combined a central scheme and the numerical dissipation of the upwind WENO scheme controlled by a weighing function. Hu et al. \cite{Hu2010} proposed a sixth-order adaptive central upwind WENO scheme called WENO-CU6. The scheme adapts between central and upwind schemes smoothly by blending the smoothness indicators of the optimal high-order central stencil and the lower order upwind stencils thereby reducing the inherent numerical dissipation associated with the upwind schemes in the smooth regions. Fu et al. \cite{fu2018new} have developed a new family of schemes called Targeted ENO (TENO) schemes which further improved the capabilities of WENO-CU6, while maintaining the accuracy at first and second-order critical points along with low numerical dissipation. {\color{black} Ghosh and Baeder \cite{ghosh2012compact} developed a class of upwind biased compact-reconstruction finite difference WENO schemes called CRWENO, which combined the compact upwind schemes and WENO schemes. CRWENO method is purely compact, but the scheme is upwind biased and excessively damps the fine-scale structures of turbulence and is computationally expensive as it involves inverting a block-tridiagonal matrix.

Weighted Compact Nonlinear Schemes (WCNS) developed by Deng et al. \cite{Deng2000}  \textcolor{black}{have} similar discontinuity capturing abilities as that of WENO but are much more flexible than WENO schemes as one can interpolate not only fluxes \cite{Zhang2008}, but also conservative variables \cite{nonomura2012}, primitive or characteristic variables \cite{Wong2017} and still maintain high-order accuracy. Recently, Liu et al. \cite{Liu2015} developed a new class of nonlinear compact schemes with a weighted hybrid interpolation of an upwind and a central interpolation. It was shown that the proposed WCNS with hybrid weighted interpolation displays a more localized dissipation than the classical WENO schemes. Finally, Subramaniam et al. \cite{subramaniam2019high} have proposed an explicit-compact interpolation along with compact finite differences that provides higher resolution and more localized dissipation compared to that of any WCNS methods in the current literature. Despite the apparent advantages of the WENO schemes, the smoothness indicators used for discontinuity detection are expensive to compute, and high-order WENO schemes are not robust enough for effective suppression of numerical oscillations in the presence of discontinuities.

Another relevant study by Sun et al. \cite{sun2016boundary}  proposed a novel approach for constructing high-fidelity shock-capturing schemes with small numerical dissipation called Boundary Variation Diminishing (BVD) algorithm in a general finite-volume framework. The BVD approach adaptively chooses appropriate reconstruction polynomials from a given set of polynomials to minimize the jumps at the cell interfaces, effectively reducing  numerical dissipation in the Riemann solvers. Sun et al. \cite{sun2016boundary} designed a BVD algorithm by hybridizing WENO and THINC (Tangent of Hyperbola for INterface Capturing) schemes, and the resulting scheme maintained the accuracy of WENO in smooth regions and substantially improved solution quality near discontinuities via THINC. Subsequently, various BVD algorithms have been proposed in the literature. Deng et al. \cite{deng2018high} have extended the approach to multi-component flows by hybridizing TVD and THINC schemes leading to significant scheme improvements. In a series of papers, Deng et al. \cite{deng2020implicit,deng2019fifth,deng2020constructing} have also proposed a new BVD algorithm that combines a high-order unlimited linear reconstruction polynomial and THINC schemes. The resulting schemes called $P_{n}T_{m}$, where $n$ is the degree of the unlimited polynomial and $m$ is the number of stages in the THINC scheme.  They showed that the BVD principle retrieves the underlying linear scheme for smooth regions while still effectively capturing discontinuities.

 The $P_{n}T_{m}$ schemes still use upwind biased interpolation, which is a source of numerical dissipation. To address this issue, we present a new algorithm, named HOCUS (High-Order Central Upwind Scheme), which combines the MP scheme and a linear-compact scheme using the BVD principle. Unlike the earlier studies where one of the reconstruction candidates is a non-polynomial function, such as THINC, we consider two fifth-order polynomials. The proposed method has the following advantages:
\begin{description}
	\item[(a)] Unlike the $P_{n}T_{m}$ schemes, which use multiple stages of evaluation of the BVD algorithm, the current approach requires a \textit{single-stage evaluation}
	\item[(b)] The unlimited linear scheme is evaluated through \textit{compact reconstruction} which has superior dissipation and dispersive properties compared to explicit schemes (see \cite{ghosh2012compact, subramaniam2019high, Pirozzoli2002,lele1992compact, adams1996high} 
	\item[(c)] Finally, since the underlying linear scheme is central, which is non-dissipative, the \textit{inherent numerical dissipation in Riemann solver is reduced in the smooth regions of the flow}. 
\end{description}

The rest of the paper is organized as follows. The new HOCUS with BVD algorithm along with the various reconstruction procedures is presented in Section \ref{sec-2}. Details of time advancement are given in Section \ref{sec3.1}.  \textcolor{black}{Several one- and two-dimensional test cases} for linear advection and Euler equations are presented in Section \ref{sec-3}. The numerical experiments clearly demonstrate the new numerical scheme can provide high-order oscillation-free results. Finally, Section \ref{sec-4} summarizes our findings.

\section{Numerical method: Spatial discretization}\label{sec-2}	

Here, we present upwind flux reconstruction using the BVD algorithm in the context of the finite-volume method for  conservation laws. For simplicity, we first consider a scalar hyperbolic conservation law represented by the following partial differential equation in one-dimension:

\begin{equation}
\frac{{\partial U}}{\partial t} + \frac{{\partial F(U)}}{\partial x} = 0 \label{eq:scalar_conservation},
\end{equation}
where, $U(x,t)$ is the solution function and $F(U)$ is the physical flux function.  

\subsection{Finite-volume method}

Equation (\ref{eq:scalar_conservation}) is discretized on a uniform grid with $N$ cells on a spatial domain spanning $x \in \left[x_a, x_b \right]$. The cell center locations are at $x_j = x_a + (j - 1/2) \Delta x$, $\forall j \in \{1, \: 2, \: \dots, \: N\}$, where $\Delta x = (x_b - x_a)/N$. The cell interfaces, indexed by half integer values, are at $x_{j+\frac{1}{2}}$, $\forall j \in \{ 0, \: 1, \: 2, \: \dots, \: N \}$. Let $I_j = [x_{j - 1 / 2}, x_{j + 1 / 2}]$ be a control volume (a computational cell) of width $\Delta x_j = x_{j + 1 / 2} - x_{j - 1 / 2}$. Integrating Equation (\ref{eq:scalar_conservation}) over  $I_j$, we obtain the following semi-discrete relation, expressed as an ordinary differential equation:

\begin{equation}\label{eqn-differencing}
 \frac{\mathd}{\mathd t} \hat{U}_j(t) = - \frac{1}{\Delta x_j} [\hat{F}_{j+ 1 / 2} - \hat{F}_{j - 1 / 2}], 
 \end{equation}
where $\hat{U}_j(t)$ is the cell average of the solution in cell $I_j$ at time $t$ and $\hat{F}_{j + 1 / 2}$ is the numerical flux over cell interface, respectively:

\begin{equation}
\hat{U}_j(t) = \frac{1}{\Delta x} \int_{x_{j - 1 / 2}}^{x_{j + 1 / 2}} U(x,t) \mathd x
\end{equation}

\begin{equation}\label{Numerical Flux}
\hat{F}_{j+\frac{1}{2}}={F}^{\rm Riemann}_{j+\frac{1}{2}}(U_{j+\frac{1}{2}}^{L},U_{j+\frac{1}{2}}^{R}), 
\end{equation}
where $L$ and $R$ are adjacent values of a cell interface as show in Fig. \ref{fig:interface}. The numerical fluxes at the cell boundaries can be computed by a variety of Riemann solvers which can be written in a canonical form as
\begin{equation}
{F}^{\rm Riemann}_{j+\frac{1}{2}} = \frac{1}{2}\underbrace{({F^L_{j+\frac{1}{2}}} + {F^R}_{j+\frac{1}{2}})}_{\text{Physical central flux}} -
 \frac{1}{2} | {A_{j+\frac{1}{2}}}|\underbrace{({U^R_{j+\frac{1}{2}}}-{U^L_{j+\frac{1}{2}}})}_{\text{Numerical dissipation}},
\label{eqn:Riemann}
\end{equation}
\begin{figure}[H]
\centering
 \includegraphics[width=0.6\textwidth]{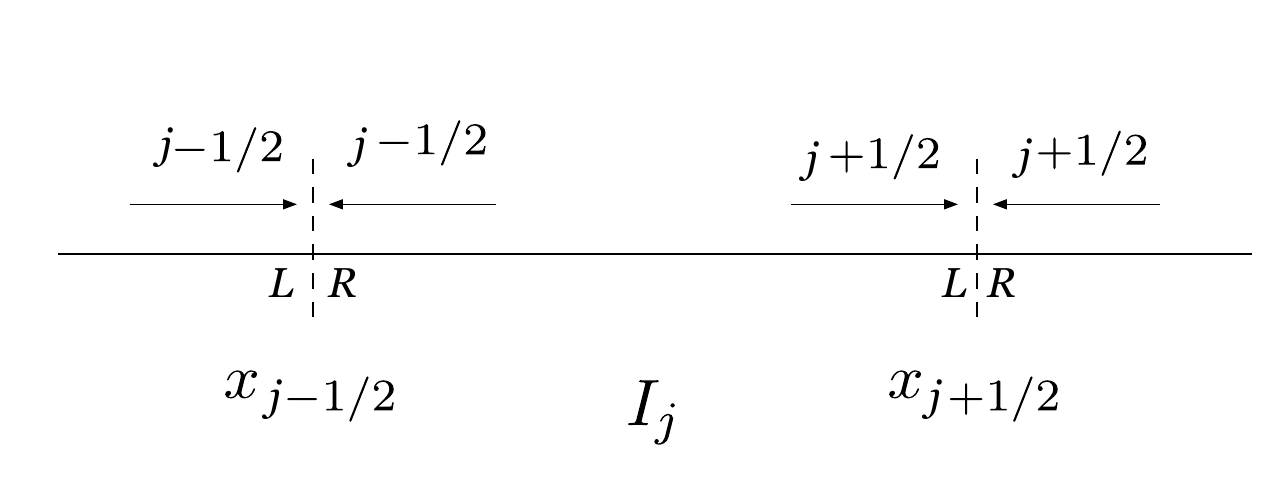}
 \caption{Illustration of cell $I_j$ and location of the cell interfaces L, R at $x_{j - 1 / 2}$ and $x_{j + 1 / 2}$.}
 \label{fig:interface}
\end{figure}
where $|{A_{j+\frac{1}{2}}}|$ denotes the characteristic signal velocity evaluated at the cell interface, and the Equation (\ref{eqn:Riemann}) clearly indicates the physical central flux and the numerical dissipation terms. The procedure of obtaining the values at the interface, $U_{j + 1 / 2}$, from cell center variables, $\hat{U}_j$, is called reconstruction and then the numerical fluxes at cell interfaces are computed by the Riemann solver. It is obvious from the Equation (\ref{eqn:Riemann}) that our objective is to obtain the values at the left and right interfaces, $U_{i+\frac{1}{2}}^{L}$ and $U_{i+\frac{1}{2}}^{R}$, at cell interfaces $x_{j+\frac{1}{2}}$, $\forall j \in \{ 0, \: 1, \: 2, \: \dots, \: N \}$. We compute these interface values by using two different reconstruction procedures, candidate polynomials, and the BVD algorithm will choose the appropriate value such that the numerical dissipation is minimized. In the following subsections, we provide the details of the calculations of candidate polynomials used for the BVD algorithm.

\subsubsection{Linear upwind compact reconstruction}
The first candidate considered for the evaluation of the interface states is the linear upwind compact reconstruction scheme.  \textcolor{black}{The advantage of the compact schemes is that they yield better accuracy and resolution than the non-compact or explicit schemes despite having a smaller stencil}. These advantages of compact schemes are well known and have been discussed in detail in Ref. \cite{ lele1992compact,adams1996high,Pirozzoli2002,ghosh2012compact,subramaniam2019high}. In the present paper, we consider the fifth-order compact reconstruction, denoted as C5, to compute the interface values at $j+1/2$. The left-biased states are obtained using cell-averaged values, $\hat U$, of the 3-point stencil $(j - 1, j, j+1)$, and right biased states uses $(j, j+1, j+2)$ respectively and the corresponding equations are shown below,

\begin{subequations}
     \begin{alignat}{1}
\frac{1}{2}  U^{L, C5}_{j-\frac{1}{2}}+  U^{L, C5}_{j+\frac{1}{2}} + \frac{1}{6}  U^{L, C5}_{j+\frac{3}{2}}&= \frac{1}{18}\hat  U_{j-1}+  \frac{19}{18}\hat U_{j} + \frac{5}{9}\hat  U_{j+1} \label{eq:left}\\
\frac{1}{6}  U^{R, C5}_{j-\frac{1}{2}}+ U^{R, C5}_{j+\frac{1}{2}} + \frac{1}{2}  U^{R, C5}_{j+\frac{3}{2}}&= \frac{5}{9}\hat  U_{j}+ \frac{19}{18}\hat U_{j+1} + \frac{1}{18}\hat  U_{j+2} \label{eq:CR}
     \end{alignat}
     \label{eqn:upwind-compact}
   \end{subequations}

The solution of the above equations requires an inversion of a tridiagonal system of equations due to their inherent, implicit nature. However, the computational complexity of a tridiagonal solution scales linearly with the number of grid points and the inversion can be efficiently carried out by Thomas algorithm. For the implementation of the boundary fluxes, i.e., the first and last interfaces along each grid line, we use the MP5 scheme presented in the next subsection. The resulting tridiagonal system of equations can be represented as

\begin{equation}
\begin{bmatrix}
1 & 0 & 0 & \cdots & &0 \\
\\
a & b & c & & & & \\

\vdots &  & \ddots & \ddots & \ddots & \vdots \\

& & & a & b &  c\\
\\
0 & &\dots &0 &0 & 1
\end{bmatrix}
\begin{bmatrix}
U^{C5}_{\frac{1}{2}} \\

\vdots \\
U^{C5}_{j+\frac{1}{2}} \\
\vdots \\

U^{C5}_{N+\frac{1}{2}}
\end{bmatrix}
=
\begin{bmatrix}
d^{MP5}_ {\frac{1}{2}}\\

\vdots \\
d^{C5}_ {j+\frac{1}{2}}\\
\vdots \\

d^{MP5}_ {N+\frac{1}{2}}
\end{bmatrix},
\end{equation}
where, $N$ is the number of cells, $a$, $b$, and $c$ are the coefficients of the left-hand side and $d$ represents the right-hand side of the system of the Equations (\ref{eqn:upwind-compact}). We use ghost points for the data required by the MP5 scheme at the boundaries. We use the same approach for both periodic and non-periodic test cases, and the numerical tests do not indicate any loss of accuracy or restrictions due to this approach. This scheme is intended to be applied in the smooth regions of the flow field as it is a linear scheme and thereby produces oscillations near discontinuities, and to achieve a shock-capturing capability, and it will be coupled with a shock-capturing scheme through BVD algorithm. In the next subsection, the candidate polynomials used for shock-capturing is explained.

\subsubsection{Monotonicity-Preserving scheme}
  The second candidate considered for the evaluation of the solution vector at the interface is the Monotonicity-Preserving fifth-order (MP5) reconstruction of Suresh and Hyunh \cite{suresh1997accurate}, which uses geometric based approach for shock-capturing. The key advantage of MP5 reconstruction is that it not only captures discontinuities but also preserves the extrema and maintains high-order accuracy in smooth regions. MP5 reconstruction is carried out in two steps. For brevity, we only explain the procedure for the left interface values, $U^{L,MP5}_{j+1/2}$, since the right interface values, $U_{j+1/2}^{R,\rm{MP5}}$ can be obtained via symmetry. In the first step, a fifth-order polynomial is constructed using cells $j-2, j-1, j, j+1,j+2$, to interpolate cell-average values $\hat U_j$ needed to obtain the cell interface values $U^{L}_{j+1/2}$: 
 \begin{equation}\label{eqn:linear5}
U^{L,P5}_{j+1/2}=(2\hat U_{j-2} - 13\hat U_{j-1} + 47\hat U_j + 27\hat U_{j+1} - 3\hat U_{j+2})/60
\end{equation}
In the second step, interpolated values $U^{L,P5}_{j+1/2}$ are limited based on the following condition:
\begin{equation} \label{eq:mp5-cond}
(U^{L,P5}_{j+1/2}- \hat U_i)(U^{L,P5}_{j+1/2}-U^{MP}) \leq \epsilon, \\
\end{equation}

\begin{equation} \label{eqn:alpha}
\begin{aligned}
\text{where,} &\ U^{MP} = \hat U_j + \operatorname{minmod}(\hat U_{j+1}-\hat U_{i}, \tilde{\alpha}\,(\hat U_j - \hat U_{j-1}))\,\,,\\
\text{and,} &\operatorname{minmod}(a,b) = \frac{1}{2} \left(\operatorname{sign}(a)+\operatorname{sign}(b)\right)\min(|a|,|b|)\,\,.
\end{aligned}
\end{equation}
 The parameter ${\alpha}$ is a constant which will be explained further in Section \ref{sec2.2} - Remark \ref{remark-5} and $\epsilon$ is a small constant set as $\epsilon=10^{-20}$. If the condition given in Equation (\ref{eq:mp5-cond}) is violated then the following algorithm is used. First, we compute the second derivatives:
\begin{eqnarray}
D_j^{-} &=& \hat U_{j-2} - 2\hat U_{j-1} + \hat U_j\,\,, \\
D_j^{0} &=& \hat U_{j-1} - 2\hat U_{j} + \hat U_{j+1}\,\,, \\
D_j^{+} &=& \hat U_{j} - 2\hat U_{j+1} + \hat U_{j+2}\,\,.
\end{eqnarray}
Next, we compute:
\begin{eqnarray}
D^{M4}_{j+1/2} &=& \operatorname{minmod4}(4D_j^0-D_j^+, 4D_j^+-D_j^0, D_j^0, D_j^+)\,\,, \\
D^{M4}_{j-1/2} &=& \operatorname{minmod4}(4D_j^0-D_j^-, 4D_j^--D_j^0, D_j^0, D_j^-)\,\,,
\end{eqnarray}
where the function $\operatorname{minmod4}$ is given by:
\begin{eqnarray}
\operatorname{minmod4}(w,a,b,c) &=& 0.125(\operatorname{sign}(w)+\operatorname{sign}(a)) \times
\\\nonumber 
&&|(\operatorname{sign}(w)+\operatorname{sign}(b))(\operatorname{sign(w)}+\operatorname{sign}(c))| \times
\\\nonumber
&&\min(|w|,|a|,|b|,|c|))\,\,.
\end{eqnarray}
We then compute
\begin{eqnarray}
U^{UL} &=& \hat U_j + \alpha(\hat U_j- \hat U_{j-1})\,\,, \\
U^{AV} &=& 0.5(\hat U_j+ \hat U_{j+1})\,\,, \\
U^{MD} &=& U^{AV} - 0.5 D^{M4}_{j+1/2}\,\,,\\
U^{LC} &=& \hat U_j + 0.5(\hat U_j-\hat U_{j-1}) + \frac{4}{3}D^{M4}_{j-1/2}\,\,.
\end{eqnarray}
Using these expressions, we compute:
\begin{eqnarray}
U_{\min} = \max(\min(\hat U_j, \hat U_{j+1}, U^{MD}), \min(\hat U_j, U^{UL}, U^{LC}))\,\,,\\
U_{\max} = \min(\max(\hat U_j,\hat U_{j+1}, U^{MD}), \max(\hat U_j, U^{UL}, U^{LC}))\,\,.
\end{eqnarray}
Finally, a new limited value for the cell interface $U^L_{j+1/2}$ is obtained via
\begin{equation}\label{eqn:sc}
U_{j+1/2}^{L,\rm{MP5}} = U_{j+1/2}^L + \operatorname{minmod}(U_{\min}-U_{j+1/2}^L, U_{\max}-U_{j+1/2}^L)\,\,.
\end{equation}

Recently, Zhao et al. \cite{zhao2019general} have evaluated various shock-capturing schemes and found that MP5 scheme is an excellent choice for wave propagation and is the most efficient of all the evaluated fifth-order accurate schemes. They also suggested combining WENO schemes with MP limiters, as in Balsara and Shu \cite{balsara2000monotonicity}, as an effective strategy for shock capturing since the monotonicity preserving WENO schemes gave the lowest overall numerical error.

\subsection{Central-upwind scheme with BVD algorithm}\label{sec2.2}

In this section, we describe the central-upwind scheme using BVD algorithm. Similar to earlier studies using the BVD algorithm, we consider two-different reconstruction polynomials C5 or C6 (details of C6 scheme are explained below) and MP5 schemes as the possible candidates for the evaluation of the interface states. The BVD algorithm selects the reconstruction polynomial with minimal numerical dissipation from the two candidate reconstructions in a given cell by evaluating the Total Boundary Variation (TBV) at a cell interface, shown in Fig. \ref{fig:BVD_scheme}, given by the following equation:
\begin{equation}\label{Eq:TBV}
TBV_{j}=\big|U_{j-\frac{1}{2}}^{L}-U_{j-\frac{1}{2}}^{R}\big|+\big|U_{j+\frac{1}{2}}^{L}-U_{j+\frac{1}{2}}^{R} \big|
\end{equation} 

\begin{figure}[H]
\centering
\includegraphics[width=1.0\textwidth]{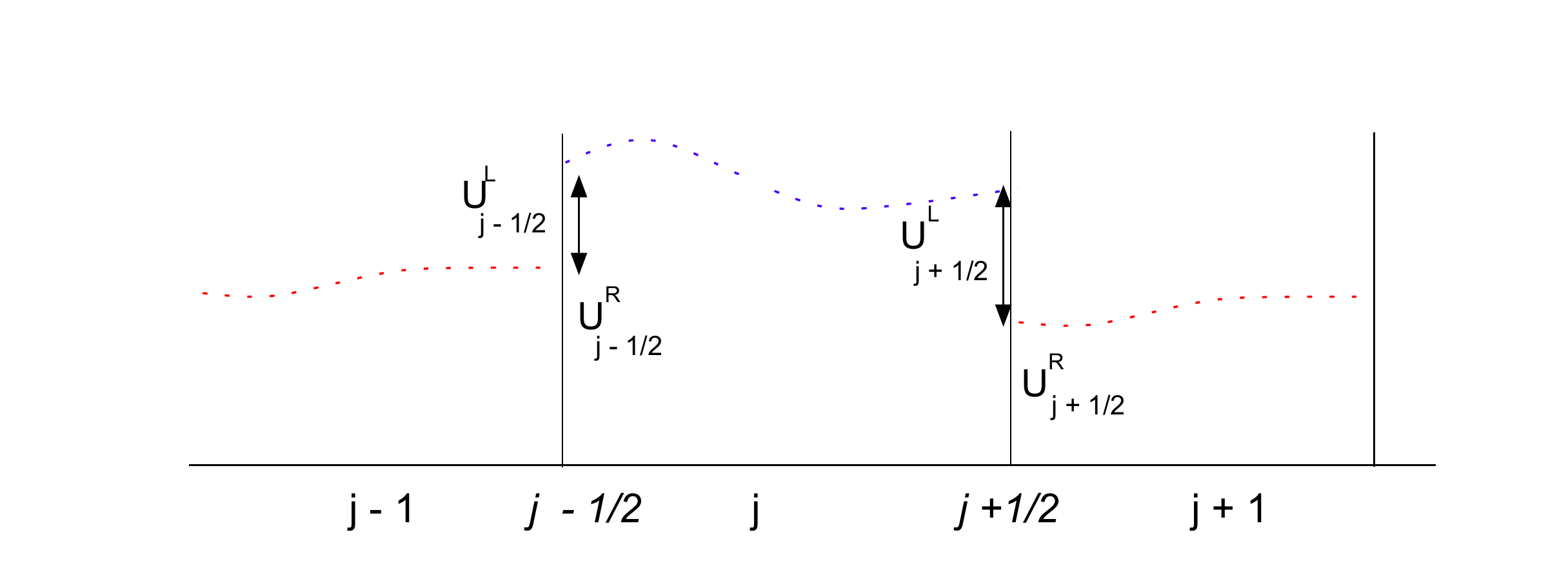}
\caption{Boundary variations at the cell interfaces.}  
\label{fig:BVD_scheme}
\end{figure}

The terms on the right-hand side of the Equation (\ref{Eq:TBV}) represent the amount of numerical dissipation introduced in the numerical flux in Equation (\ref{eqn:Riemann}) at each interface of a given cell. When two reconstruction functions of the same data are available at an interface, the BVD algorithm compares the TBVs of the concerned polynomials and chooses the least dissipative one in a given cell. The linear schemes C5 or C6  will be used in the smooth region and the BVD algorithm will switch to the MP5 scheme in the presence of discontinuities. The complete procedure for the new scheme presented in this work is summarized as follows:

\begin{description}
\item[Step 1.] Evaluate the interface values by using two different reconstruction procedures:
\begin{enumerate}[(a)]
\item Linear upwind compact reconstruction given by Equation (\ref{eqn:upwind-compact}) and
\item MP5 scheme given by Equations (\ref{eqn:linear5}) - (\ref{eqn:sc}).
\end{enumerate}
\item[Step 2.] Calculate the TBV values for each cell $I_{j}$ by using the compact reconstruction:

\begin{equation}\label{Eq:TBVC5}
TBV_{j}^{C5}=\big|U_{j-\frac{1}{2}}^{L,C5}-U_{j-\frac{1}{2}}^{R,C5}\big|+\big|U_{j+\frac{1}{2}}^{L,C5}-U_{j+\frac{1}{2}}^{R,C5} \big|
\end{equation} 
and  MP5 scheme:
\begin{equation}\label{Eq:TBVSC}
TBV_{j}^{MP5}=\big|U_{j-\frac{1}{2}}^{L,MP5}-U_{j-\frac{1}{2}}^{R,MP5}\big|+\big|U_{j+\frac{1}{2}}^{L,MP5}-U_{j+\frac{1}{2}}^{R,MP5} \big|.
\end{equation} 

\item[Step 3.] By averaging the left and right interface values computed by the C5 scheme we obtain the central non-dissipative linear scheme. This step is one of the key contributions of this paper in the development of the non-dissipative central scheme, denoted hereafter as HOCUS6. 

\begin{equation}\label{eqn-add}
\begin{cases} 
U_{{j+\frac{1}{2}}}^{L} = U_{{j+\frac{1}{2}}}^{L, C6} = \frac{1}{2}(U_{{j+\frac{1}{2}}}^{L, C5} +U_{{j+\frac{1}{2}}}^{R, C5})\\
\\
U_{{j+\frac{1}{2}}}^{R}=U_{{j+\frac{1}{2}}}^{R, C6} = \frac{1}{2}(U_{{j+\frac{1}{2}}}^{L, C5} +U_{{j+\frac{1}{2}}}^{R, C5})\\
\end{cases}
\rightarrow \textbf{HOCUS6}
\end{equation} 

This step is equivalent to computing the interface fluxes by the following sixth order linear compact reconstruction scheme, denoted as C6 in this paper:
\begin{equation}\label{eqn-C6}
\frac{1}{3} U^{C6}_{j-\frac{1}{2}}+ U^{C6}_{j+\frac{1}{2}} + \frac{1}{3} U^{C6}_{j+\frac{3}{2}} = \frac{29}{36} ( \hat U_{j+1}+\hat U_{j}) + \frac{1}{36} (\hat U_{j-1}+ \hat U_{j+2})
\end{equation}

\item[Step 4.] Now,  \textcolor{black}{all the interface values of $j-\frac{3}{2}$, $j-\frac{1}{2}$, $j+\frac{1}{2}$ and $j+\frac{3}{2}$, both $L$ and $R$} are modified according to the following algorithm to obtain non-oscillatory results:
\begin{equation}\label{eqn:BVDstep}
\begin{aligned} 
{\rm if } \ \ TBV_{j}^{MP5} < TBV_{j}^{C5} 
\end{aligned} 
\end{equation}

\item[Step 5.] Finally, evaluate the interface flux $\hat{F}_{j + \frac{1}{2}}$ from ${U}_{j+\frac{1}{2}}^{L,R}$ using Equation (\ref{eqn:Riemann})
\end{description}

\begin{remark}
\normalfont It is important to note that reconstruction polynomials at the cell interface are first computed by the compact upwind scheme, given by Equations (\ref{eqn:upwind-compact}), and averaged later to obtain the central scheme. Direct evaluation by the central scheme given by Equation (\ref{eqn-C6}) \textcolor{black}{cannot be done as the values of left and right interface values are the same for the central scheme. It will lead to oscillations and failure of the simulation. The total boundary variations of MP5 and C5 schemes should be compared in Equation (\ref{eqn:BVDstep})}. 
\end{remark}

\begin{remark}
\normalfont   \textcolor{black}{If the underlying solution is smooth on this stencil it results in $U_{{j+\frac{1}{2}}}^{L}=U_{{j+\frac{1}{2}}}^{R}$ and through this step we can \textit{reduce} the numerical dissipation in the Riemann solver shown in Equation (\ref{eqn:Riemann}).} 

\begin{equation}
\text{if} \ U_{{j+\frac{1}{2}}}^{L}=U_{{j+\frac{1}{2}}}^{R} \rightarrow  \underbrace{({U^R_{j+\frac{1}{2}}}-{U^L_{j+\frac{1}{2}}})}_{\text{Numerical dissipation}} = 0.
\end{equation}
\end{remark}

\begin{figure}[H]
\centering
\includegraphics[width=0.6\textwidth]{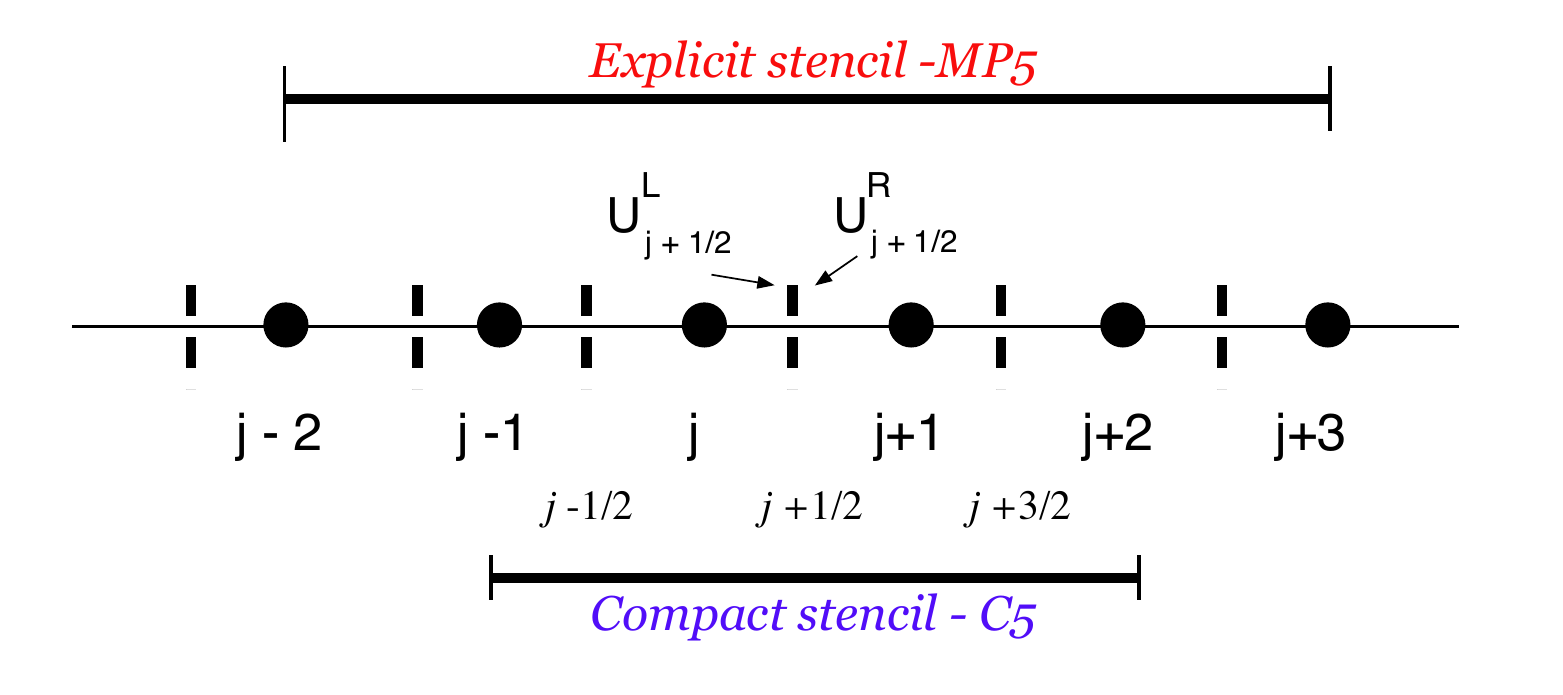}
\caption{Computational stencils of MP5 and C5 schemes.}   
        		\label{fig:scheme}
\end{figure}

%\begin{remark}
%\normalfont Another important note is that 
%\end{remark}

\begin{remark}
\normalfont We can also use the left and right interface values computed by the linear fifth order upwind compact scheme, which is denoted as HOCUS5 scheme in this paper, i.e. ${U^{(L,R) C5}}_{j+\frac{1}{2}}$ in Equation (\ref{eqn:BVDstep}), which will also give superior results but the numerical dissipation will still be present as $(U_R - U_L)_{j+\frac{1}{2}} \neq 0$ due to upwinding. The implementation of HOCUS5 scheme is as follows:
\begin{equation}\label{eqn-hocus5}
\begin{cases} 
 U_{{j+\frac{1}{2}}}^{L} =  U_{{j+\frac{1}{2}}}^{L, C5}\\
\\											
 U_{{j+\frac{1}{2}}}^{R}=  U_{{j+\frac{1}{2}}}^{R, C5}
\end{cases}
\rightarrow \textbf{HOCUS5},
\end{equation}
\begin{equation}\label{eqn-blah5}
 \text{if} ~~ TBV_{j}^{MP5} < TBV_{j}^{C5} ~~~  {U^{L,R}}_{j+\frac{1}{2}}={U^{(L,R)MP5}}_{j+\frac{1}{2}}.
\end{equation}
\end{remark}

\begin{remark}
\normalfont  \textcolor{black}{The objective of the current paper is not to study various combinations. However, at the same time, we want to present the reasons that motivated us to arrive at the current method, which is a combination of MP5 and C5/C6. In the development of our new scheme, several different combinations of linear and nonlinear schemes were considered, which are given below:} 
\begin{enumerate}
\item For example, it is also possible to use a third-order reconstruction with the \textcolor{black}{minmod limiter \cite{van1979towards} as the shock-capturing scheme in combination with the compact schemes with the present algorithm. The results are shown in Appendix \ref{sec-appb}. }
\item  \textcolor{black}{We have also considered the THINC scheme in combination with the compact schemes and the results are shown in Appendix \ref{sec-appc}.}
\item Finally, by combining the WENO-Z scheme with the compact scheme C6 by using the BVD algorithm, Equation (\ref{eqn:BVDstep}), we observed that the results are too diffusive which motivated us to use monotone schemes as the base shock-capturing scheme. We have also studied a different BVD algorithm proposed in \cite{xie2017hybrid, ivan2014high}, and the results are shown in Appendix C.
\end{enumerate}
\end{remark}

\begin{remark}
\normalfont \textcolor{black}{In P4T2 scheme (\cite{deng2019fifth}) the interfaces of the cells $j-1$, $j$ and $j+1$ are corrected by using BVD algorithm and in the current approach the interfaces values of $j-\frac{3}{2}$, $j-\frac{1}{2}$, $j+\frac{1}{2}$ and $j+\frac{3}{2}$ are modified.}
\end{remark}

\begin{remark}\label{remark-5}
\normalfont The value of ${\alpha}$ in Equation (\ref{eqn:alpha}) puts a restriction on the CFL (Courant-Friedrichs-Lewy) number such that CFL $\leq$ 1/(1+${\alpha}$) (see Suresh and Hyunh \cite{suresh1997accurate} for further details). The values considered various schemes for parameter ${\alpha}$ are shown in Table \ref{tab:alpha}. It is the only heuristic parameter in the present approach. For all the simulations in this paper, we considered CFL to be less than or equal to 0.2.
\end{remark}

\begin{table}[H]
 \centering
 \footnotesize
 \caption{Parameter ${\alpha}$ for MP5, HOCUS5 and HOCUS6 schemes}
  \begin{tabular}{ccccccc}
 \hline
     & MP5         & HOCUS5 &  HOCUS6 & \\
 \hline
   ${\alpha}$        & 4 &  7 &   7 \\
   \hline
  \end{tabular}%
 \label{tab:alpha}%
\end{table}%
The numerical method described above can be easily extended to multi-dimensional (2D and 3D) problems via dimension by dimension approach using the method of lines. The extension of the scalar conservation equation to Euler equations is discussed in the following section.

\subsection{Extension to Euler equations}\label{sec-2.1}	

The governing equations for the two-dimensional Euler equations are:
\begin{equation}\label{eq:2d}
\frac{\partial{\mathbf{Q}}}{\partial{t}}+\frac{\partial{\mathbf{F}}}{\partial{x}}+\frac{\partial{\mathbf{G}}}{\partial{y}}=0,
\end{equation}
\begin{equation}
\begin{split}
&\mathbf{Q}=\left[\begin{matrix}
\rho\\
\rho u\\
\rho v\\
E\\
\end{matrix}
\right],\
\mathbf{F}=\left[\begin{matrix}
\rho u\\
\rho u^2+p\\
\rho uv\\
(E+p)u\\
\end{matrix}
\right],\
\mathbf{G}=\left[\begin{matrix}
\rho v\\
\rho uv\\
\rho v^2+p\\
(E+p)v\\
\end{matrix}
\right]\textcolor{black}{,}
\mathbf{U}=\left[\begin{matrix}
\rho\\
u\\
v\\
p\\
\end{matrix}
\right],
\end{split}
\end{equation}
where $\mathbf{Q}$, $\mathbf{F}$ and $\mathbf{G}$ are the conservative variables, inviscid flux vectors in $x$ and $y$ directions, respectively. The primitive variables $\rho,u,v\ \text{and} \ p$ represent the density, the $x-$ and $y-$ components of the velocity, and the pressure are stored in the vector $\mathbf{U}$. The total specific energy is defined:
\begin{equation}
E=\frac{p}{\gamma-1}+\frac{\rho}{2}(u^2+v^2),
\end{equation}
where $\gamma$ is the constant specific heat ratio. The left and right eigenvectors of the two-dimensional Euler equations, denoted by $\bm{L_n}$ and $\bm{R_n}$, are used for characteristic variable projection are as given as:
\begin{align}\label{eqn:leftright}
	\bm{R_n} = \begin{bmatrix} 
   		 1 & 1 & 1 & 0 \\
		 \\
   		 -\frac{n_xc}{\rho} & 0 & \frac{n_xc}{\rho} & \frac{l_x}{\rho} \\
		 \\
   		 -\frac{n_yc}{\rho} & 0 & \frac{n_yc}{\rho} & \frac{l_y}{\rho} \\
		 \\
     c^2 & 0 & c^2 & 0 \\
   			 \end{bmatrix}, \enskip
	\bm{L_n} = \begin{bmatrix}
   		 0 & -\frac{n_x\rho}{2c} &-\frac{n_y\rho}{2c} & \frac{1}{2c^2} \\
		 \\
   		 1 & 0 & 0 & -\frac{1}{c^2} \\
		 \\
   		 0 & -\frac{n_x\rho}{2c} & \frac{n_y\rho}{2c} & \frac{1}{2c^2} \\
		 \\
     0 & \rho l_x & \rho l_y & 0 \\
   			 \end{bmatrix}, &&\\\nonumber
\end{align}
where $\bm{n}$ = $[n_x \ n_y]^t$ and $[l_x \ l_y]^t$ is a tangent vector (perpendicular to $\bm{n}$) such as $[l_x \ l_y]^t$ = $[-n_y \ n_x]^t$. By taking $\bm{n}$ = $[1, 0]^t$ and $[0, 1]^t$ we obtain the corresponding eigenvectors in $x$ and $y$ directions. In the present paper, we use simple dimension-wise implementation for two-dimensional Euler equations which is akin to evaluation of integrals by mid-point rule rather than using high-order quadratures, which is as follows:
\begin{equation}
 \frac{\mathbf{d\hat Q}_{j,i}}{\mathbf{dt}} = -\frac{1}{\Delta x}\left(\mathbf{\hat F_{j + \frac{1}{2}}} - \mathbf{\hat F_{j- \frac{1}{2}}}\right) - \frac{1}{\Delta y} \left(\mathbf{\hat G_{i + \frac{1}{2}}} - \mathbf{\hat G_{i - \frac{1}{2}}}\right), \label{eqn:e6m}
 \end{equation}
where, $\mathbf{\hat Q_{j,i}}$ are cell averaged conservative variables and  $\mathbf{\hat{F}}_{}$ and $\mathbf{\hat{G}}_{}$ are interface numerical fluxes computed by an approximate Riemann solver. Here we consider the HLLC (Harten, Leer and Lax with Contact) approximate Riemann solver of Harten et al. \cite{harten1983upstream} and Toro et al. \cite{toro1994restoration}. Considering $x-$direction as an example, the interface flux computed by the HLLC Riemann solver is:

%\noindent This results in the following definition of flux depending
%on the wave speeds,
\begin{equation}
\mathbf{\hat{F}}_{}=\begin{cases}
\mathbf{\hat{F}}_{L} & S_{L}\geq0\\
\mathbf{\hat{F}}_{*L} & S_{L}\leq0\leq S_{*}\\
\mathbf{\hat{F}}_{*R} & S_{*}\leq0\leq S_{R}\\
\mathbf{\hat{F}}_{R} & S_{R}\leq0,
\end{cases}\label{eqn:hllcformulation}
\end{equation}
where, $\mathbf{\hat{F}}_{L}=\mathbf{F}\left(\mathbf{\hat{Q}}_{L}\right)\ $ and $\quad\mathbf{\hat{F}}_{R}=\mathbf{F}\left(\mathbf{\hat{Q}}_{R}\right)\ $ are the exact local fluxes at either side of the interface and using the Rankine-Hugoniot jump conditions across the left and
right wave:
\begin{equation}
\mathbf{\hat{F}}_{*L}=\mathbf{\hat{F}}_{L}+S_{L}\left(\mathbf{\hat{Q}}_{*L}-\mathbf{\hat{Q}}_{L}\right)\ \label{eq:FStarL}
\end{equation}
\begin{equation}
\mathbf{\hat{F}}_{*R}=\mathbf{\hat{F}}_{R}+S_{R}\left(\mathbf{\hat{Q}}_{*R}-\mathbf{\hat{Q}}_{R}\right)\,\label{eq:FStarR}
\end{equation}
where, $S_L$ and $S_R$ are numerical approximations to the speeds of the left most and right most running characteristics that emerge as the solution of the Riemann problem at an interface. It has been shown that under appropriate choice of wave speeds $S_L$ and $S_R$, the HLLC scheme is both positivity preserving and entropy satisfying \cite{batten1997choice}. This choice of wave speeds are  \textcolor{black}{as follows}: 

\begin{align}
\nonumber
 S_L = min(u_{nL}-c_L, \tilde{u}_n-\tilde{c})\\
 S_R = max(u_{nR}+c_R, \tilde{u}_n +\tilde{c})
 \label{eqn:HLLEwavespeedestimate}
\end{align}
where $u_{nL,R}$ are the normal velocities across an interface, $c_{L,R}$ are the respective sonic speeds and $\tilde{u}_n,\tilde{c}$ are the standard Roe averaged quantities at the interface. Using the integral form of the conservation laws, closed form expressions for the conserved quantities in the states in the star region $S_*$ can be derived as: 
\begin{equation}  \label{eqn:HLLCstate} 
\begin{aligned}
  \mathbf{\hat Q}_{*K}^{HLLC} &= \rho_{K} \left(\frac{S_{K}-u_{nK}}{S_{K}-S_*}\right) \left( \begin{array}{c}
                                       1\\
                                       S_*\\
                                       u_{tK}\\
                                       \frac{(\rho E)_{K}}{\rho_{K}} + (S_*-u_{nK})(S_* + \frac{p_{K}}{\rho_{K}(S_{K}-u_{nK})})                             
                                       \end{array} \right)\\
                                       \end{aligned}
\end{equation}
where K =L, R and $u_{tK}$ denote the tangential velocities across an interface. In the above expressions, $S_L$ and $S_R$ can be obtained using Eq.(\ref{eqn:HLLEwavespeedestimate}) and Batten et al. \cite{batten1997choice} provided a closed form expression for $S_*$ which is as follows:
\begin{align}
 \label{eqn:HLLCmiddlewaveestimate}
 S_* = \frac{p_R - p_L + \rho_Lu_{nL}(S_L - u_{nL}) - \rho_Ru_{nR}(S_R - u_{nR})}{\rho_L(S_L - u_{nL}) -\rho_R(S_R - u_{nR})}.
\end{align}

Unlike scalar advection, shock-capturing should be carried out using characteristic variables for coupled hyperbolic equations like the Euler equations for optimal \textit{cleanest} results \cite{van2006upwind}. Direct reconstruction of the interfaces values using primitive variables resulted in small oscillations. The complete numerical algorithm for the Euler equations is summarized below, which includes the transforming of primitive variables into characteristics variables necessary for capturing discontinuities.

\begin{description}

\item[Step 1.] To evaluate the RHS contribution of the advective fluxes in the $x-$ direction from $\mathbf{\hat Q_{j,i}}$ to build $\mathbf{F_{j + \frac{1}{2}}}$ in Equation \ref{eqn:e6m} \\

\begin{enumerate}
\item Compute primitive variables $\mathbf{U_{j,i}}$	from $\mathbf{\hat Q_{j,i}}$. \\

\item Interpolate the primitive variables,$\mathbf{U_{j,i}}$, from cell-center to cell interface by compact scheme C5, Equations (\ref{eqn:upwind-compact}), and obtain $\mathbf U_{j+\frac{1}{2},i}^{L,C5}$ and $\mathbf U_{j+\frac{1}{2},i}^{R,C5}$. \\

\item Average the interface values $\mathbf U_{j+\frac{1}{2},i}^{L,C5}$ and $\mathbf U_{j+\frac{1}{2},i}^{R,C5}$ and obtain $\mathbf U_{j+\frac{1}{2},i}^{L}$ and $\mathbf U_{j+\frac{1}{2},i}^{R}$ as shown in Equation (\ref{eqn-add}). \\

\item Interpolate the characteristic variables, $\bm{W}$, projected from the primitive variables by the MP5 scheme, through Equations (\ref{eqn:linear5} - \ref{eqn:sc}). The procedure is as follows, \\

\begin{enumerate}

\item Compute the arithmetic averages at the interface $(x_{j+\frac{1}{2}}, y_i)$ by using neighbouring cells, $(x_j, y_i)$ and $(x_{j+1},y_i)$. Compute the left and right eigenvectors $\bm{L_{n}}$ and $\bm{R_{n}}$, Equation (\ref{eqn:leftright}), at the interface using the arithmetic averages, with $\bm{n}$ = $[1, 0]^t$ which corresponds to $x-$ direction. The primitive variables are then transformed to characteristic variables by the following equation: \\
\begin{equation}
	\bm{W}_{j,i} = \bm{L}_{\bm{n}_{j+\frac{1}{2},i}} \bm{U}_{j,i}.
\end{equation}
\item Carry out the reconstruction of characteristic variables by MP5 scheme, through Equations (\ref{eqn:linear5} - \ref{eqn:sc}), and obtain left- and right- interface values denoted by $\tilde{\bm{W}}_{j+\frac{1}{2},i}^L$ and $\tilde{\bm{W}}_{j+\frac{1}{2},i}^R$. \\
\item After obtaining $\tilde{\bm{W}}_{j+\frac{1}{2},i}^L$ and $\tilde{\bm{W}}_{j+\frac{1}{2},i}^R$ from the MP5 reconstruction the primitive variables are then recovered by projecting the characteristic variables back to physical fields: \\
\begin{equation}
\begin{aligned}
	{\bm{U}}_{j+\frac{1}{2},i}^{L,MP5} &= \bm{R}_{\bm{n}_{j+\frac{1}{2},i}} \tilde{\bm{W}}_{j+\frac{1}{2},i}^L, \\
  {\bm{U}}_{j+\frac{1}{2},i}^{R,MP5} &= \bm{R}_{\bm{n}_{j+\frac{1}{2},i}} \tilde{\bm{W}}_{j+\frac{1}{2},i}^R.
\end{aligned}
\end{equation}

\end{enumerate}
\item Modify the interface values by using BVD algorithm, Equation (\ref{eqn:BVDstep}). Primitive variables are used for the BVD algorithm.

\item Compute the upwind flux $\mathbf{F_{j + \frac{1}{2}}}$ at the cell interface by using Riemann solver, through Equation (\ref{eqn:hllcformulation}). \\
\end{enumerate} 
\item[Step 2.] To evaluate the RHS contribution of the advective fluxes in the $y-$ direction to build $\mathbf{G_{i + \frac{1}{2}}}$ repeat steps \textbf{(a)-(f)} in Step 1 with $\bm{n}$ = $[0,1]^t$.\\
\item[Step 3.] Evaluate the right-hand side of the finite-volume equation, Equation (\ref{eqn:e6m}) and perform time integration as described in the next section. 
\end{description}

\section{Numerical method: Temporal discretization}\label{sec3.1}
For time integration we use the explicit third-order TVD Runge-Kutta method \cite{jiang1995}:
\begin{eqnarray} \label{rk}
\mathbf{Q^{(1)}} & = & \mathbf{Q^n} + \Delta t \mathbf{Res}(\mathbf{Q^n}), \nonumber \\
\mathbf{Q^{(2)}} & = & \frac{3}{4} \mathbf{Q^n} + \frac{1}{4} \mathbf{Q^{(1)}} + \frac{1}{4}\Delta t \mathbf{Res}(\mathbf{Q^{(1)}}) ,\\
\mathbf{Q^{n+1}} & = & \frac{1}{3} \mathbf{Q^n} + \frac{2}{3} \mathbf{Q^{(2)}} + \frac{2}{3}\Delta t \mathbf{Res}(\mathbf{Q^{(2)}}), \nonumber
\end{eqnarray}
where $\mathbf{Res}$ is the residual that is evaluated on the right hand side of the Equation (\ref{eqn:e6m}), the superscripts $\mathbf{n}$ and $\mathbf{n+1}$ denote the current and the subsequent time-steps, and superscripts $\mathbf{(1)-(2)}$ corresponds to intermediate steps respectively. The time step $\Delta t$ is taken as suggested by Titarev and Toro \cite{titarev2004finite}:

\begin{equation}
\Delta{t}=\text{CFL} \times \text{$min_{j,i}$} (\frac{\Delta{x}}{\left(|u_{j,i}+c_{j,i}|\right)}, \frac{\Delta{y}}{\left(|v_{j,i}+c_{j,i}|\right)}),
\end{equation}
where $c$ is the speed of sound and given by $c=\sqrt{\gamma{p}/\rho}$. Time integration is performed with a CFL = 0.2 for all the problems unless otherwise stated. Simulations are carried out by MP5, HOCUS5, and HOCUS6 schemes presented in this paper on a uniform grid. We also carried out simulations with a fifth-order WENO spatial reconstruction scheme of Borges et al. \cite{Borges2008}, denoted as WENO-Z, for comparison. A comparison of HOCUS5 and C5 schemes is carried out for the one-dimensional test cases for Euler equations in Section \ref{sec-3.2.1}.

\section{Results}\label{sec-3}	
									
In this section, the proposed numerical scheme is tested for both the linear advection equation and Euler equations. 

\subsection{Linear advection equation}
\begin{example}\label{ex:adv_comp} Advection of complex waveforms
\end{example}

\noindent First we consider the linear advection equation:
\begin{equation}\label{eq:lin}
\dfrac{\partial{u}}{\partial{t}}+\dfrac{\partial{u}}{\partial{x}}=0
\end{equation}
with the following initial condition

\begin{equation}\label{eq:1st_case}
u_0(x)=\left\{
\begin{array}{ll}
exp(-log(2)(x+0.7)^2/0.0009 & -0.8\le x\le -0.6,\\
1 & -0.4\le x \le -0.2,\\
1-|10(x-0.1)| &\ 0\le x \le 0.2,\\
(1-100(x-0.5))^{0.5} &\ 0.4\le x\le 0.6,\\
0,& \ otherwise
\end{array}
\right .
\end{equation}
The initial condition contains a combination of a square wave, a Gaussian, a sharp triangle and a semi-ellipse. Two sets of results are presented.  The first uses a grid size of $200$ cells and results are presented at $t=2$.  The second uses a grid size of $400$ cells and results are presented at $t=500$.  A CFL of $0.1$ is used for all the simulations.  The results are presented in Figs. \ref{fig_advectioncomplex1} and \ref{fig_advectioncomplex2}.  Both the HOCUS5 and HOCUS6 schemes produce non-oscillatory results and preserve the initial wave profiles better than the MP5 and WENO-Z schemes. The difference between HOCUS5 and HOCUS6 schemes is marginal. For the results at $t=500$, shown in the second column of Fig. \ref{fig_advectioncomplex1}, HOCUS6 can preserve the waves even after $1 \times 10^6$ time steps similar to that of the $P4T2$ scheme proposed by Deng et al. in \cite{deng2019fifth} (see their Fig. 5). On the other hand, we can see a squaring effect for the MP5 and WENO scheme is completely distorted in Fig. \ref{fig_advectioncomplex2}.

\begin{figure}[H]
\begin{onehalfspacing}
\centering
\subfigure[HOCUS5, t=2]{\includegraphics[width=0.48\textwidth]{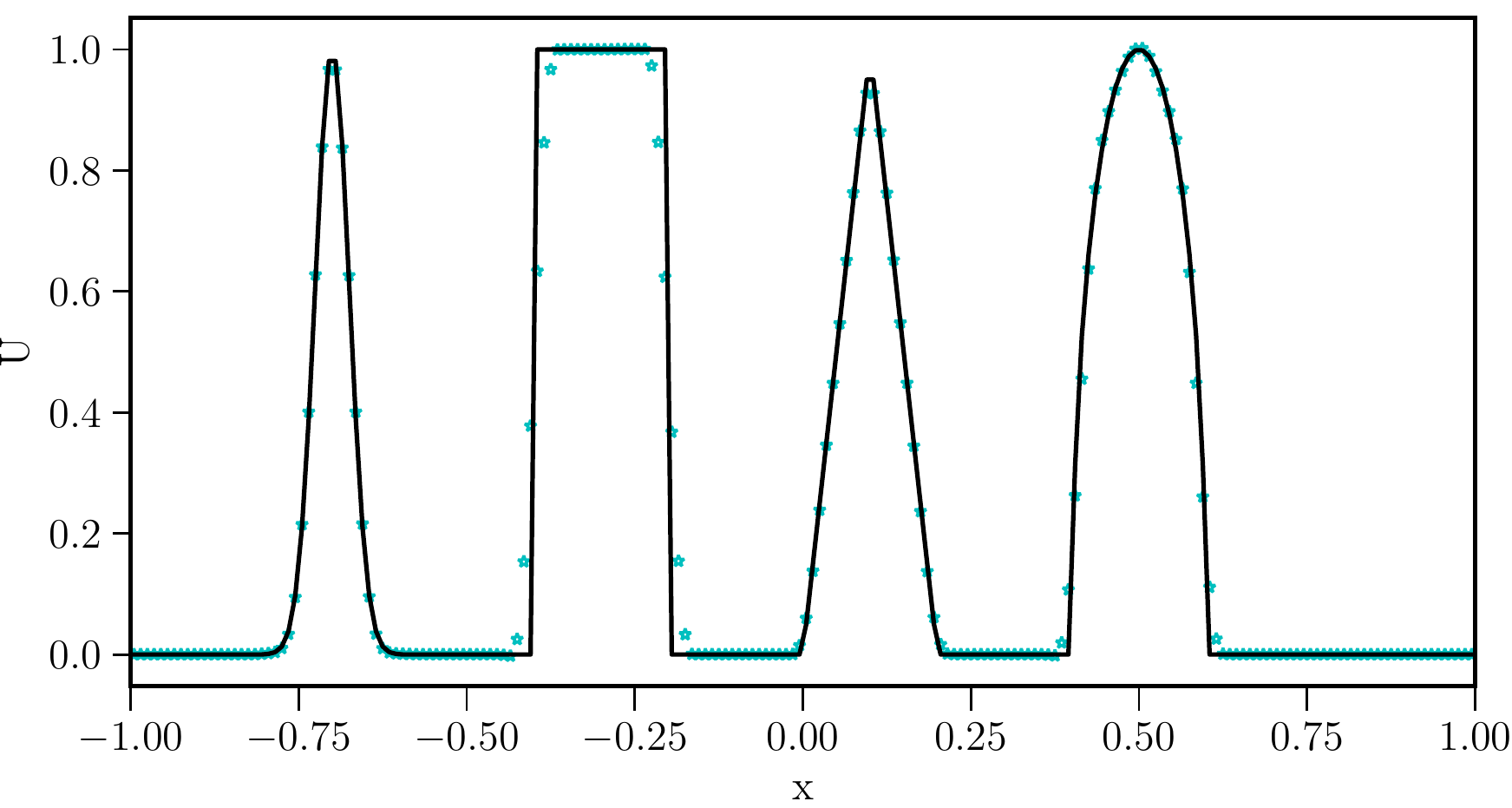}
\label{fig:HOCUS5_AC}}
\subfigure[HOCUS5, t=500]{\includegraphics[width=0.48\textwidth]{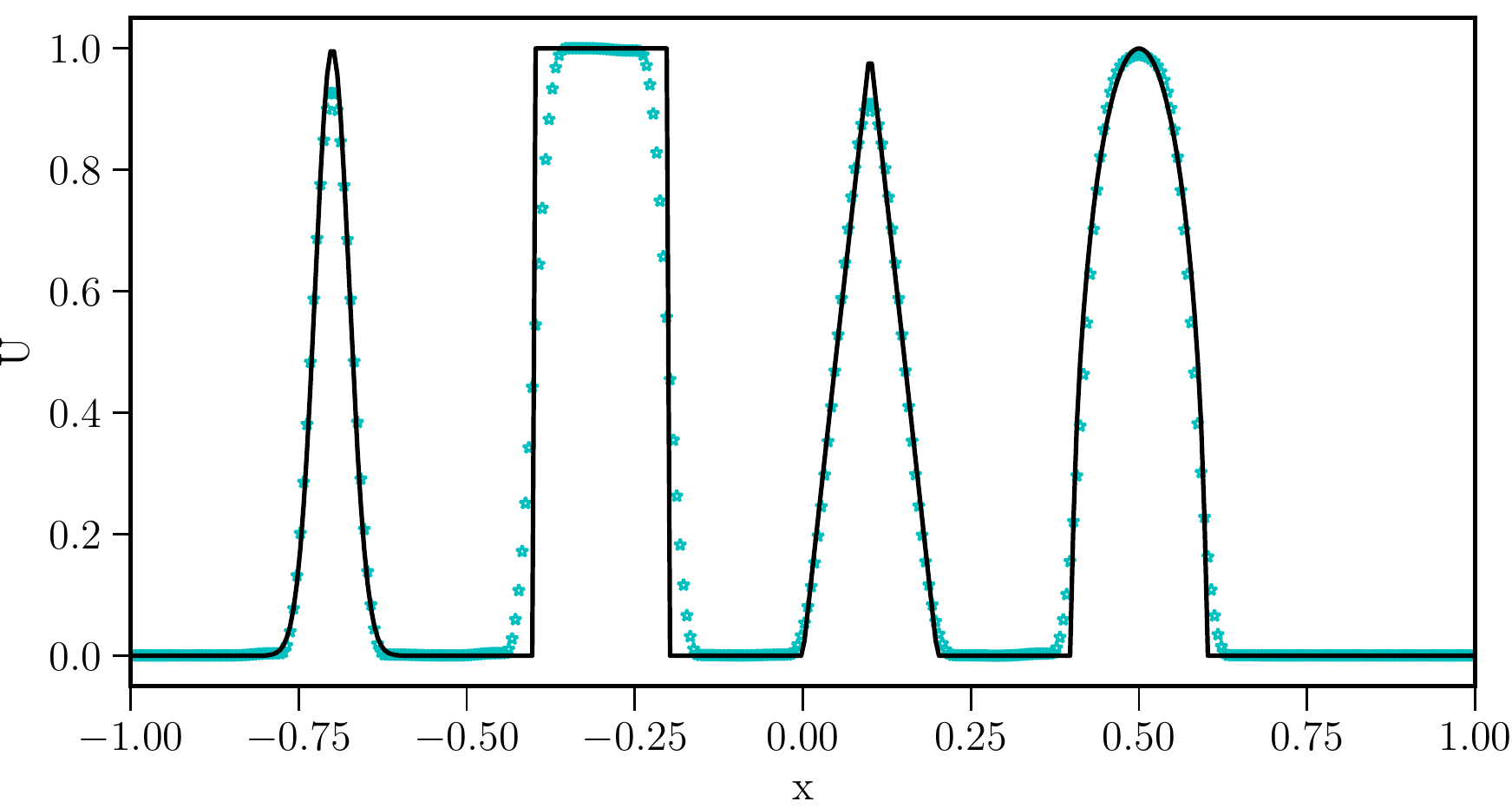}
\label{fig:HOCUS5_AC1}}
\subfigure[HOCUS6, t=2]{\includegraphics[width=0.48\textwidth]{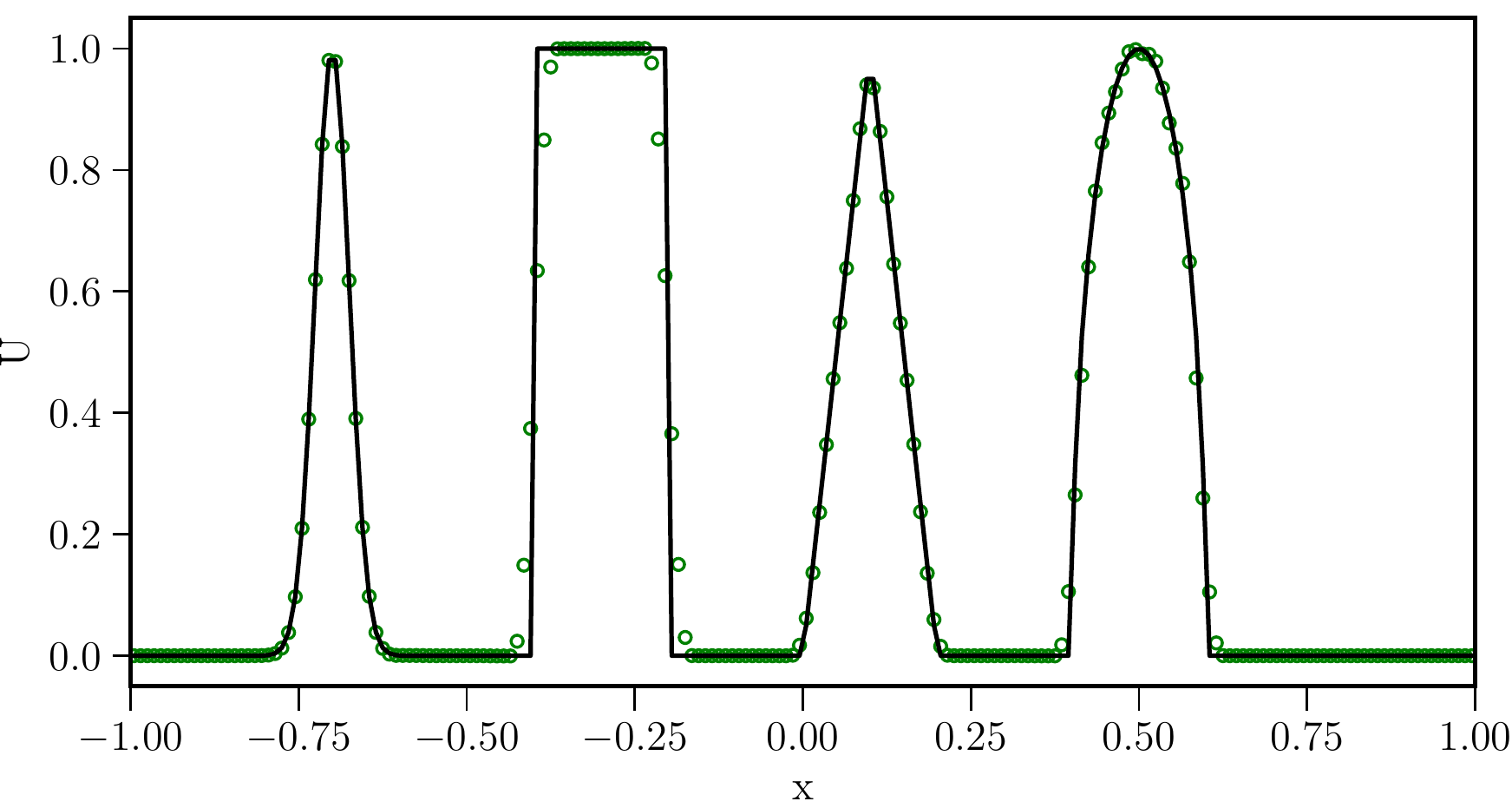}
\label{fig:HOCUS6_AC}}
\subfigure[HOCUS6, t=500]{\includegraphics[width=0.48\textwidth]{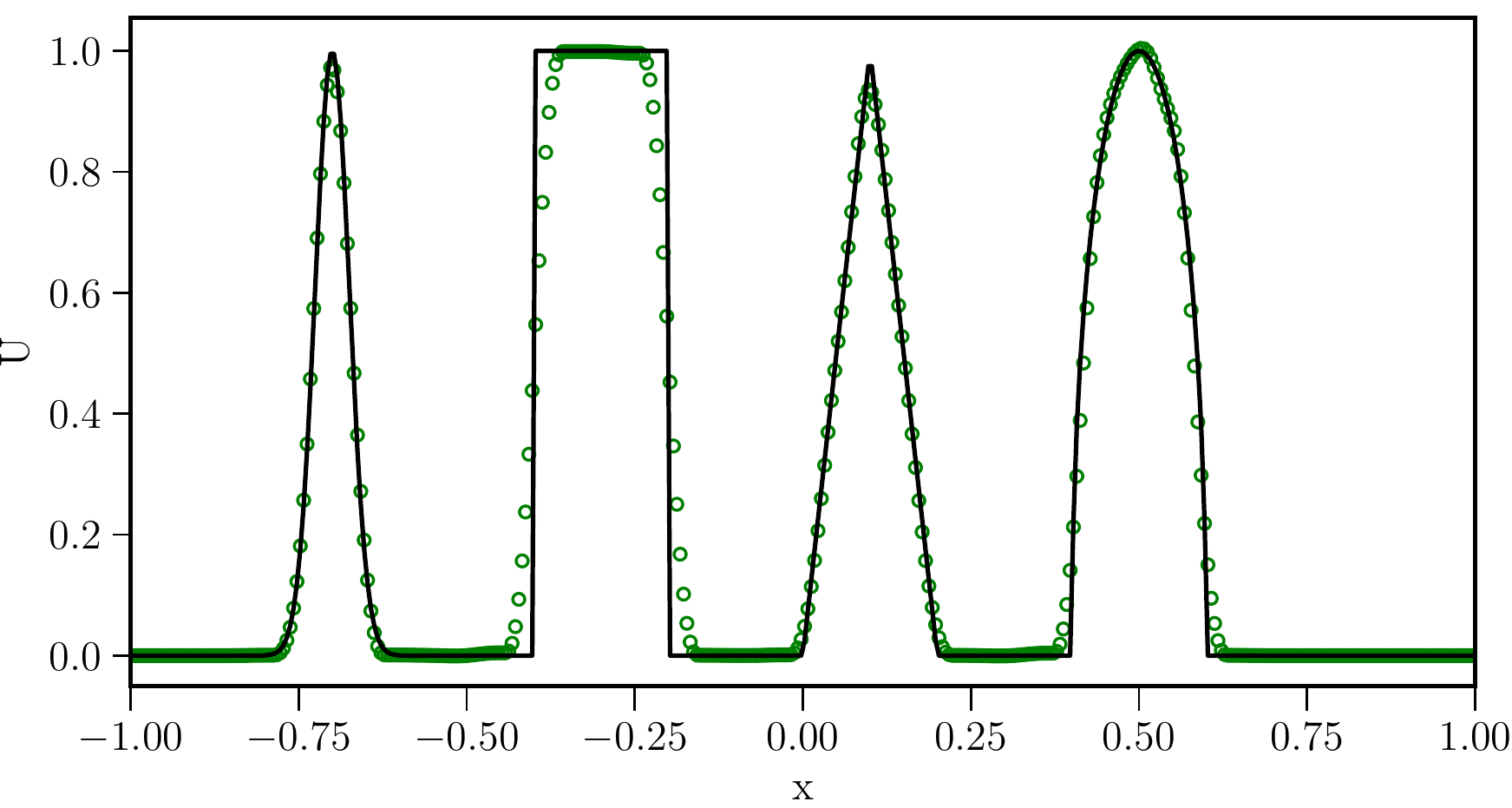}
\label{fig:HOCUS6_AC1}}
\caption{Results for HOCUS5 and HOCUS6 schemes for Example \ref{ex:adv_comp}, where solid line: reference solution; \textcolor{black}{cyan stars: HOCUS5}, green circles: HOCUS6.}
\label{fig_advectioncomplex1}
\end{onehalfspacing}
\end{figure}

\begin{figure}[H]
\begin{onehalfspacing}
\centering
\subfigure[\textcolor{black}{WENOZ, t=2}]{\includegraphics[width=0.48\textwidth]{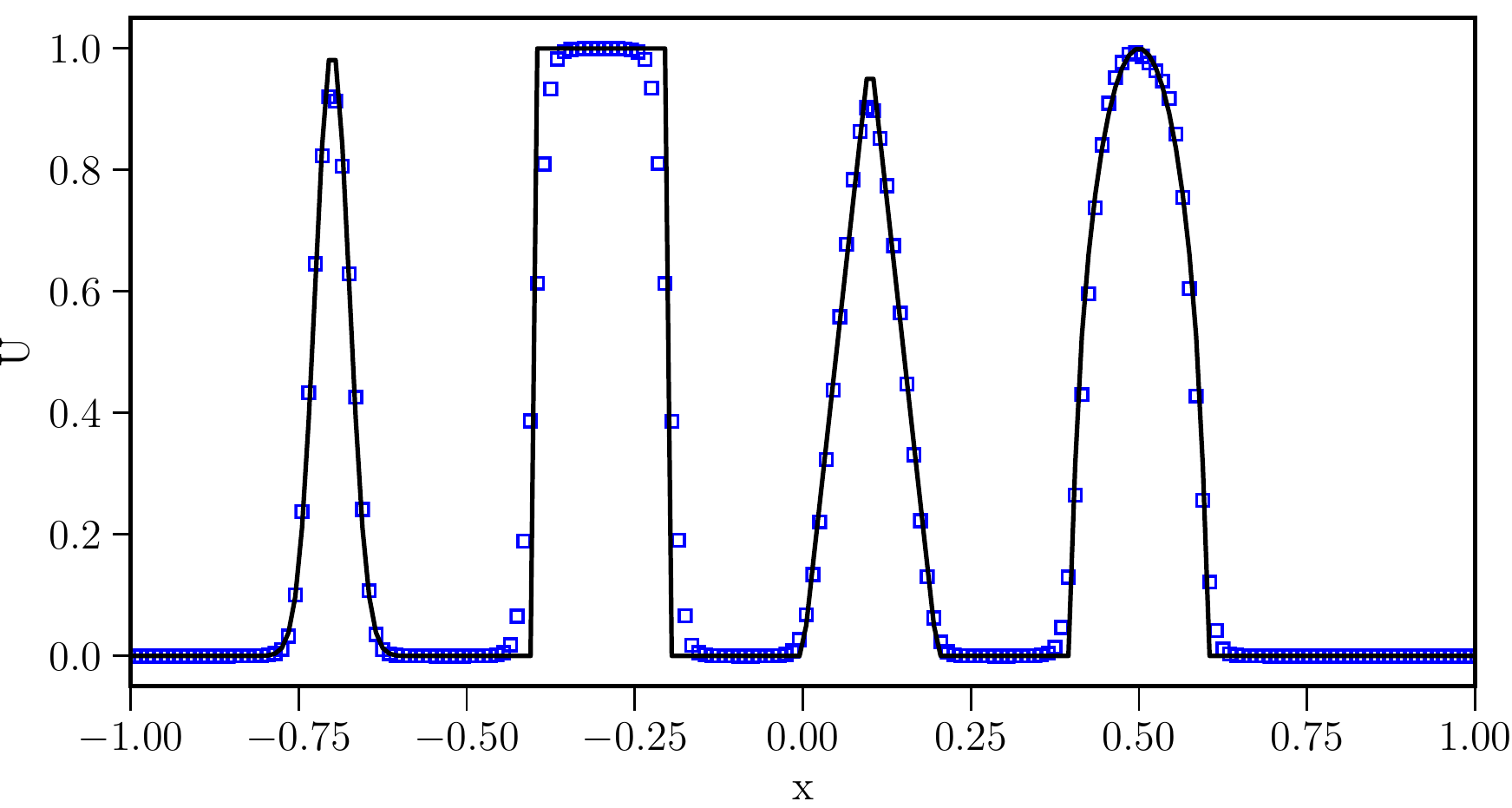}
\label{fig:WENO-Z_AC}}
\subfigure[\textcolor{black}{WENOZ, t=500}]{\includegraphics[width=0.48\textwidth]{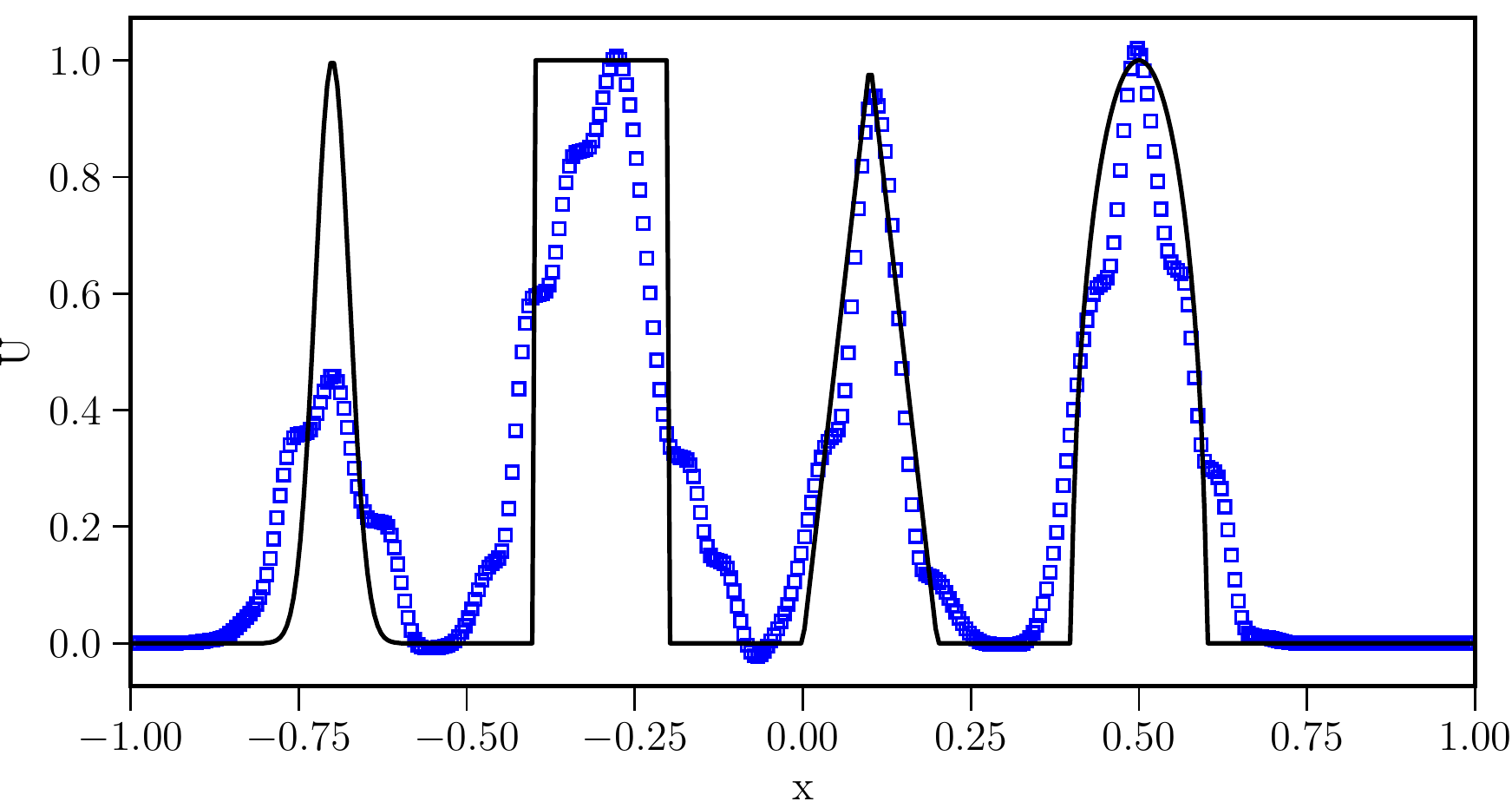}
\label{fig:WENO-Z_AC1}}
\subfigure[MP5, t=2]{\includegraphics[width=0.48\textwidth]{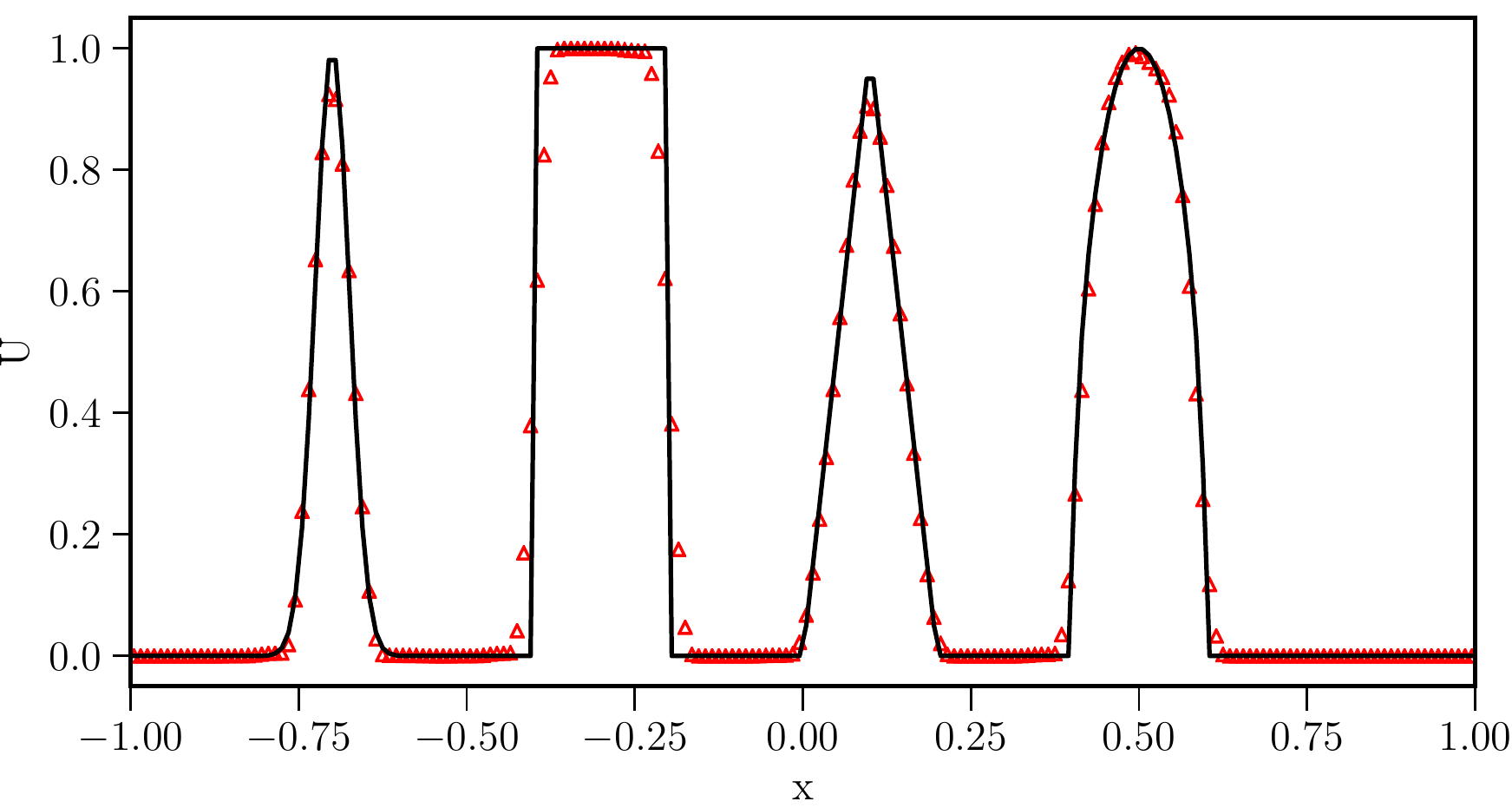}
\label{fig:MP5_AC}}
\subfigure[MP5, t=500]{\includegraphics[width=0.48\textwidth]{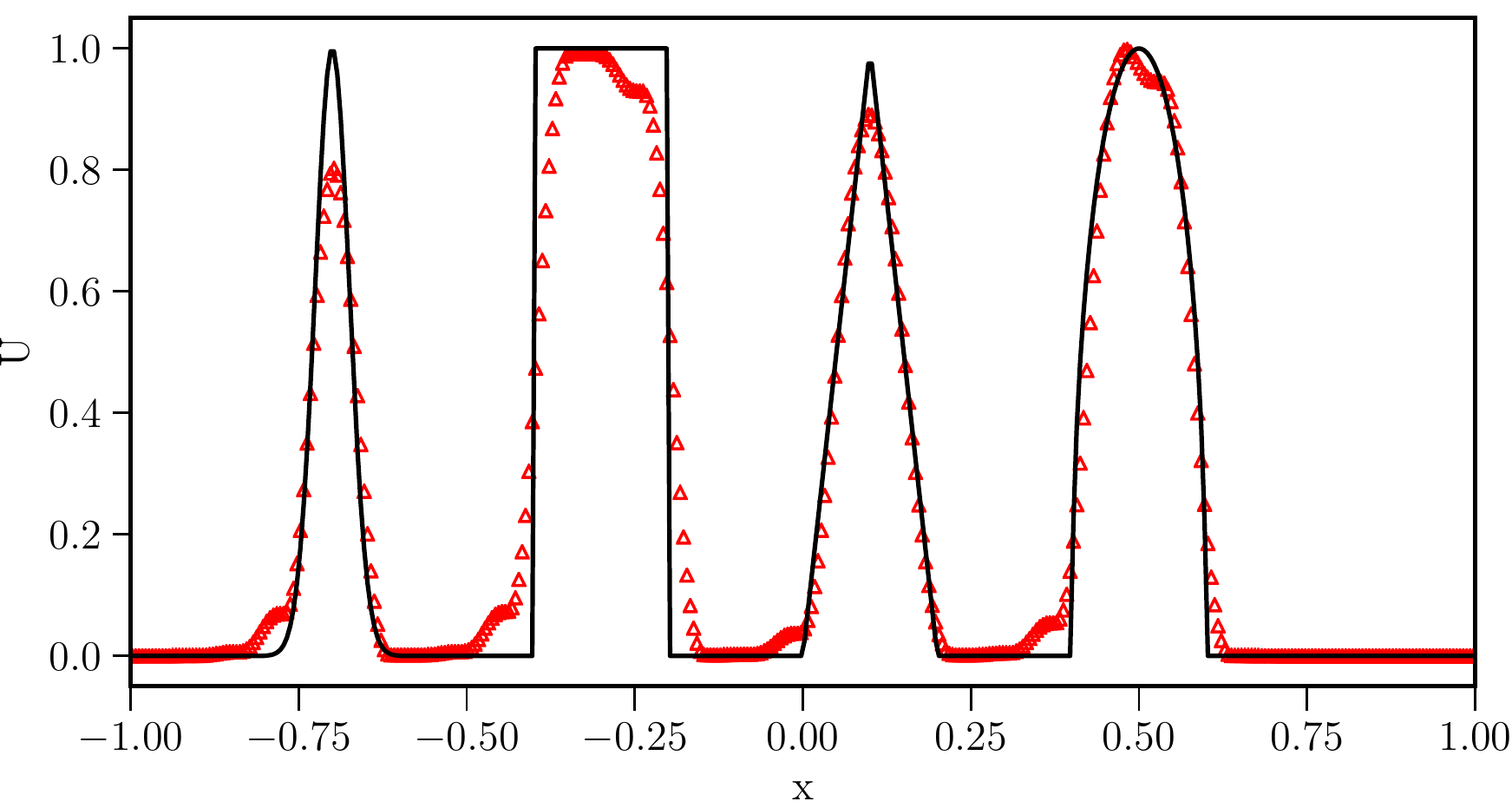}
\label{fig:MP5_AC1}}
\caption{Results for WENO-Z and MP5 schemes for Example \ref{ex:adv_comp}, where solid line: reference solution; blue squares: WENO-Z; red triangles: MP5.}
\label{fig_advectioncomplex2}
\end{onehalfspacing}
\end{figure}

\begin{example}\label{ex:accuracy}{Order of accuracy tests}

\end{example}
First, we consider the one-dimensional Gaussian pulse advection problem \cite{Yamaleev2009} to verify our scheme's order of accuracy.  Here we consider the following initial condition: 
\begin{equation}\label{eq:26}
u_0(x) = exp(-300(x-0.5)^2)
\end{equation}
on a computational domain of $[0,1]$ subject to periodic boundary conditions. Simulations are performed until $t=1$. To ensure that the effect of time-discretization is negligible a time step of \textcolor{black}{$\Delta t$ = 0.1$\Delta x^{2.0}$} is selected. As can be seen from Table \ref{tab:case1a} the $L_1$ convergence histories and order of accuracies for both HOCUS5 and HOCUS6 are the same as the corresponding linear schemes, C5 and C6, which shows the BVD algorithm is choosing appropriate smoother polynomials resulting in formal order of accuracies. \textcolor{black}{Comparison has also been made with the sixth order explicit central scheme, denoted as E6, which indicates the current approach, HOCUS6, is one order better than the E6 scheme.}

% Table generated by Excel2LaTeX from sheet 'Sheet1'
\begin{table}[htbp]
  \centering
 \footnotesize
\caption{\textcolor{black}{$L_1$ errors and numerical orders of accuracy on $u_t$ + $u_x$ = 0 with $u_0$(x) = $exp(-300(x-0.5)^2)$. $N$ is the total number of cells on a uniform mesh and $t = 1$.}}
    \begin{tabular}{rcccccccr}
\hline
\hline
    \multicolumn{1}{c}{N} & MP5   &       & WENO-Z &       & HOCUS5 &       & C5    &  \\
    \hline
          &       & order &       & order &       & order &       & \multicolumn{1}{c}{order} \\
          \hline
    \multicolumn{1}{c}{40} & 1.04E-02 &       & 1.00E-02 &       & 5.94E-03 &       & 4.95E-03 &  \\
\hline
    \multicolumn{1}{c}{80} & 1.30E-03 & 3.00  & 1.29E-03 & 2.96  & 1.82E-04 & 5.03  & 1.78E-04 & \multicolumn{1}{c}{4.80} \\
\hline
    \multicolumn{1}{c}{160} & 4.71E-05 & 4.78  & 4.71E-05 & 4.78  & 5.16E-06 & 5.14  & 5.13E-06 & \multicolumn{1}{c}{5.12} \\
\hline
    \multicolumn{1}{c}{320} & 1.51E-06 & 4.97  & 1.51E-06 & 4.97  & 1.55E-07 & 5.05  & 1.55E-07 & \multicolumn{1}{c}{5.05} \\
\hline
    \multicolumn{1}{c}{640} & 4.73E-08 & 4.99  & 4.73E-08 & 4.99  & 4.79E-09 & 5.02  & 4.79E-09 & \multicolumn{1}{c}{5.02} \\
\hline
\hline
          & N     & HOCUS6 &       & C6    &       & E6    &       &  \\
\hline
          &       &       & order &       & order &       & Order &  \\
\hline
          & 40    & 3.28E-03 &       & 2.45E-03 &       & 2.02E-02 &       &  \\
\hline
          & 80    & 3.99E-05 & 6.36  & 3.49E-05 & 6.13  & 4.68E-04 & 5.43  &  \\
\hline
          & 160   & 5.98E-07 & 6.06  & 4.94E-07 & 6.14  & 7.82E-06 & 5.90  &  \\
\hline
          & 320   & 1.05E-08 & 5.83  & 7.19E-09 & 6.10  & 1.26E-07 & 5.96  &  \\
\hline
          & 640   & 1.68E-10 & 5.97  & 1.11E-10 & 6.01  & 1.97E-09 & 5.99  &  \\
\hline
    \end{tabular}%
 \label{tab:case1a}%
\end{table}%

\begin{figure}[H]
\centering
\includegraphics[width=0.45\textwidth]{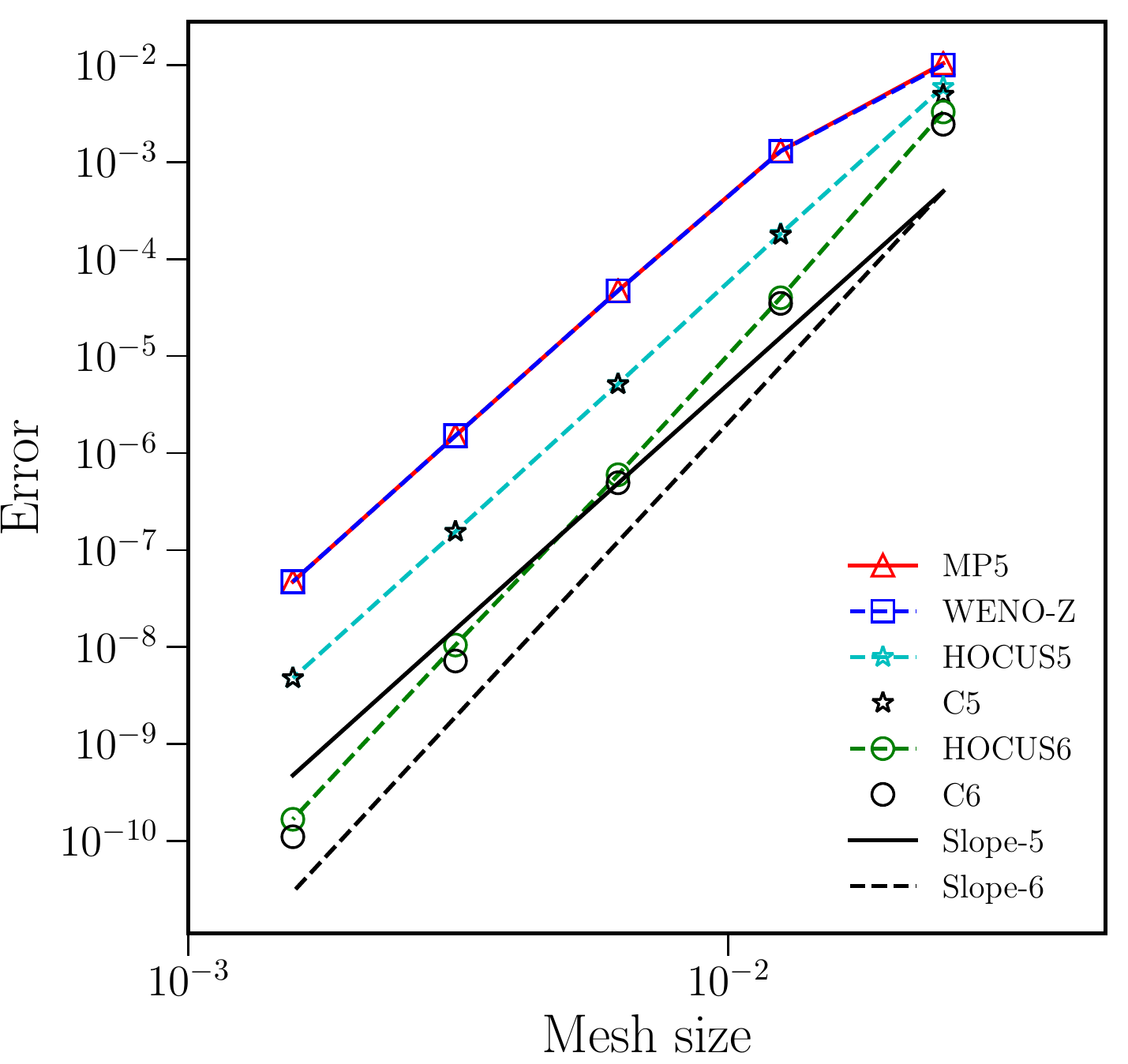}
\caption{\textcolor{black}{Accuracy plot}}
\label{fig_acc}
\end{figure}

Second, we consider the test case with critical points given in Henrick et al. \cite{Henrick2005}. It is known from the literature that WENO schemes lose their accuracy at critical points where the first and second derivatives go to zero. The initial condition is as follows:
\begin{equation}\label{eq:25}
u_0(x) = sin(\pi x - sin(\pi x)/\pi ), 
\end{equation}
on a computational domain of $[-1,1]$ subject to periodic boundary conditions. Simulations are performed until \textcolor{black}{$t=8$}. This initial condition has two critical points \cite{Henrick2005}. We chose the same time step as the previous case.  Table \ref{tab:case2a} shows the results for the HOCUS5 scheme in comparison with the corresponding linear compact scheme C5. We can see that the HOCUS5 scheme exactly retrieves the linear scheme with fifth-order accuracy, even in the presence of critical points. In Table \ref{tab:case2a} a comparison between HOCUS6 and C6 is shown, where we can see that there is a small difference between the two schemes with HOCUS6 being less accurate. Incorporating additional constraints in the BVD algorithm given by Equation (\ref{eqn:BVDstep}) seems to retrieve the accuracy but these constraints could not be generalized, especially for Euler equations, and are not considered further. \textcolor{black}{Results considering additional constraints are presented in Appendix \ref{sec-appd}}. Nevertheless, the HOCUS6 scheme shows superior absolute accuracy compared to the HOCUS5 scheme. 

% Table generated by Excel2LaTeX from sheet 'Sheet2'
\begin{table}[H]
  \centering
 \footnotesize
 \caption{\textcolor{black} {$L_1$ errors and numerical orders of accuracy on $u_t$ + $u_x$ = 0 with $u_0$(x) = sin($\pi$x - sin($\pi$x)/$\pi$ ). $N$ is the total number of cells on a uniform mesh and $t = 8$.}}
    \begin{tabular}{rcccccccr}
    \hline
    \hline
    \multicolumn{1}{c}{N} & WENOZ &       & MP5   &       & HOCUS5 &       & \multicolumn{1}{l}{C5} &  \\
    \hline
          &       & Order &       & Order &       & Order &       & \multicolumn{1}{c}{Order} \\
   \hline
    \multicolumn{1}{c}{20} & \multicolumn{1}{r}{7.89E-03} &       & \multicolumn{1}{r}{7.89E-03} &       & 1.40E-03 &       & 1.40E-03 &  \\
\hline
    \multicolumn{1}{c}{40} & \multicolumn{1}{r}{2.89E-04} & 4.77  & \multicolumn{1}{r}{2.89E-04} & 4.77  & 3.62E-05 & 5.27  & 3.62E-05 & \multicolumn{1}{c}{5.27} \\
\hline
    \multicolumn{1}{c}{80} & \multicolumn{1}{r}{9.28E-06} & 4.96  & \multicolumn{1}{r}{9.28E-06} & 4.96  & 1.01E-06 & 5.16  & 1.01E-06 & \multicolumn{1}{c}{5.16} \\
\hline
    \multicolumn{1}{c}{160} & \multicolumn{1}{r}{2.91E-07} & 4.99  & \multicolumn{1}{r}{2.91E-07} & 4.99  & 3.04E-08 & 5.06  & 3.04E-08 & \multicolumn{1}{c}{5.06} \\
\hline
\hline
          &  N     & HOCUS6 &       & C6    &       & E6    &       &  \\
\hline
          &       &       & Order &       & Order &       & Order &  \\
\hline
          & 20    & 7.47E-04 &       & 1.30E-04 &       & 2.90E-03 &       &  \\
\hline
          & 40    & 1.61E-05 & 5.53  & 1.72E-06 & 6.23  & 4.84E-05 & 5.91  &  \\
\hline
          & 80    & 3.13E-07 & 5.69  & 2.74E-08 & 5.97  & 7.84E-07 & 5.95  &  \\
\hline
          & 160   & 4.81E-09 & 6.03  & 5.35E-10 & 5.68  & 1.23E-08 & 5.99  &  \\
\hline
    \end{tabular}%
 \label{tab:case2a}%
\end{table}%

%----------------------------------------%----------------------------------------%----------------------------------------%----------------------------------------%----------------------------------------
%----------------------------------------%----------------------------------------%----------------------------------------%----------------------------------------%----------------------------------------
%----------------------------------------%----------------------------------------%----------------------------------------%----------------------------------------%----------------------------------------
%----------------------------------------%----------------------------------------%----------------------------------------%----------------------------------------%----------------------------------------
%----------------------------------------%----------------------------------------%----------------------------------------%----------------------------------------%----------------------------------------
%----------------------------------------%----------------------------------------%----------------------------------------%----------------------------------------%----------------------------------------
%\subsection{Euler equations}
\subsection{One-dimensional Euler equations}

In this subsection, we consider the test cases for the one-dimensional Euler equations, as outlined in Section 2.

%----------------------------------------
\begin{example}\label{Titarev-Toro}{ Titarev-Toro problem}
\end{example}
In this test case, we consider the shock-entropy wave problem of Titarev-Toro \cite{titarev2004finite}, where a high frequency oscillating sinusoidal wave interacts with a shock wave. It tests the ability of the high-order scheme to capture the extremely high-frequency waves. The initial conditions are given as:

\begin{align}\label{shock_tita}
(\rho,u,p)=
\begin{cases}
&(1.515695,\ \ 0.523326,\ \ 1.805),\quad x<-4.5,\\
&(1+0.1\sin(20x\pi),\ \ 0,\ \ 1),\quad x>-4.5, \text{Case-1}\\ 
&(1+0.1\sin(10x\pi),\ \ 0,\ \ 1),\quad x>-4.5, \text{Case-2},
\end{cases}
\end{align}

Case-1 corresponds to the original test case proposed in \cite{titarev2004finite} and Case-2 is the modified version with a smaller frequency considered in \cite{deng2019fifth}. Simulations are carried out on a grid with $N=1000$ cells for Case-1 and $N=400$ cells for Case-2. The reference solution is obtained by the WENO-Z scheme on the grid size of $N=3000$ and $N=1600$ for test Case-1 and Case-2, respectively. From Figs.\ \ref{fig_titazoom} and \ref{fig_titaglobal}, we can observe that the HOCUS6 scheme yields the best resolution and captures the fine-scale structures in the solution much better than other schemes. For Case-1, results obtained by the HOCUS6 scheme are similar to that of the those obtained by adaptive upwind-central sixth-order adaptive scheme proposed by Wong and Lele in \cite{Wong2017} (see their Fig. 13) on the same grid size using the same Riemann solver. For Case-2, we obtained superior results in comparison with P4T2 scheme proposed by Deng et al. in \cite{deng2019fifth} (see their Fig. 11) which required $800$ cells to capture the high-frequency waves whereas the current approach requires only $400$ cells. This test shows the capabilities of the adaptive-central scheme as it can capture the high-frequency density fluctuations and demonstrates the capabilities of the proposed HOCUS6 scheme over the upwind HOCUS5 scheme, MP5 and WENO-Z schemes. Also, it can be observed that the compact schemes, HOCUS6 and HOCUS5, are superior to the explicit, MP5 and WENO-Z, schemes for both cases. 

\begin{figure}[H]
\centering
\subfigure[Case-1]{\includegraphics[width=0.48\textwidth]{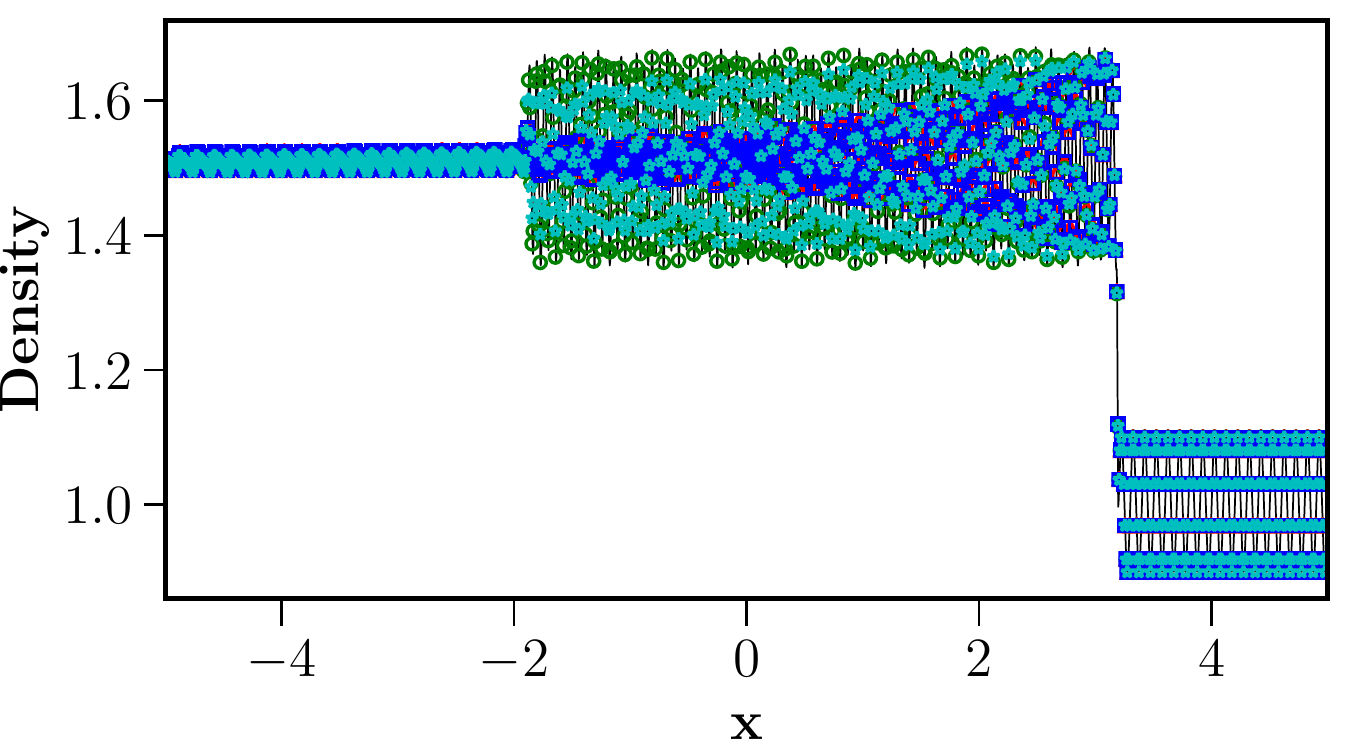}
\label{fig:tita1}}
\subfigure[Case-2]{\includegraphics[width=0.48\textwidth]{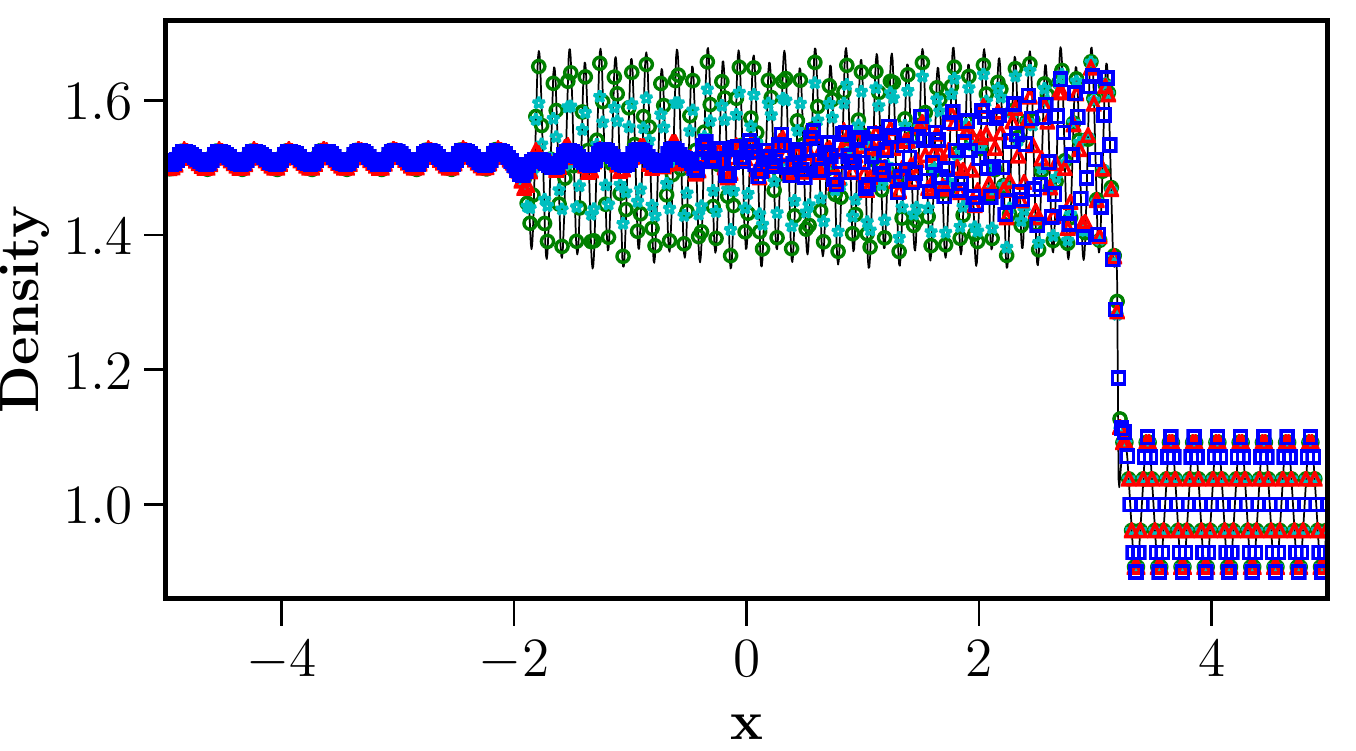}
\label{fig:tita2}}
\caption{Titarev-Toro problem, Example \ref{Titarev-Toro}, at $t=5.0$ by using WENO-Z, MP5, HOCUS5 and HOCUS6 schemes. Fig. (a) Case-1 with N=1000 and Fig. (b) Case-2 with $N=400$, where solid line: reference solution; cyan stars: HOCUS5, green circles: HOCUS6; blue squares: WENO-Z; red triangles: MP5.}
\label{fig_titaglobal}
\end{figure}

\begin{figure}[H]
\centering
\subfigure[Case-1, local profile]{\includegraphics[width=0.48\textwidth]{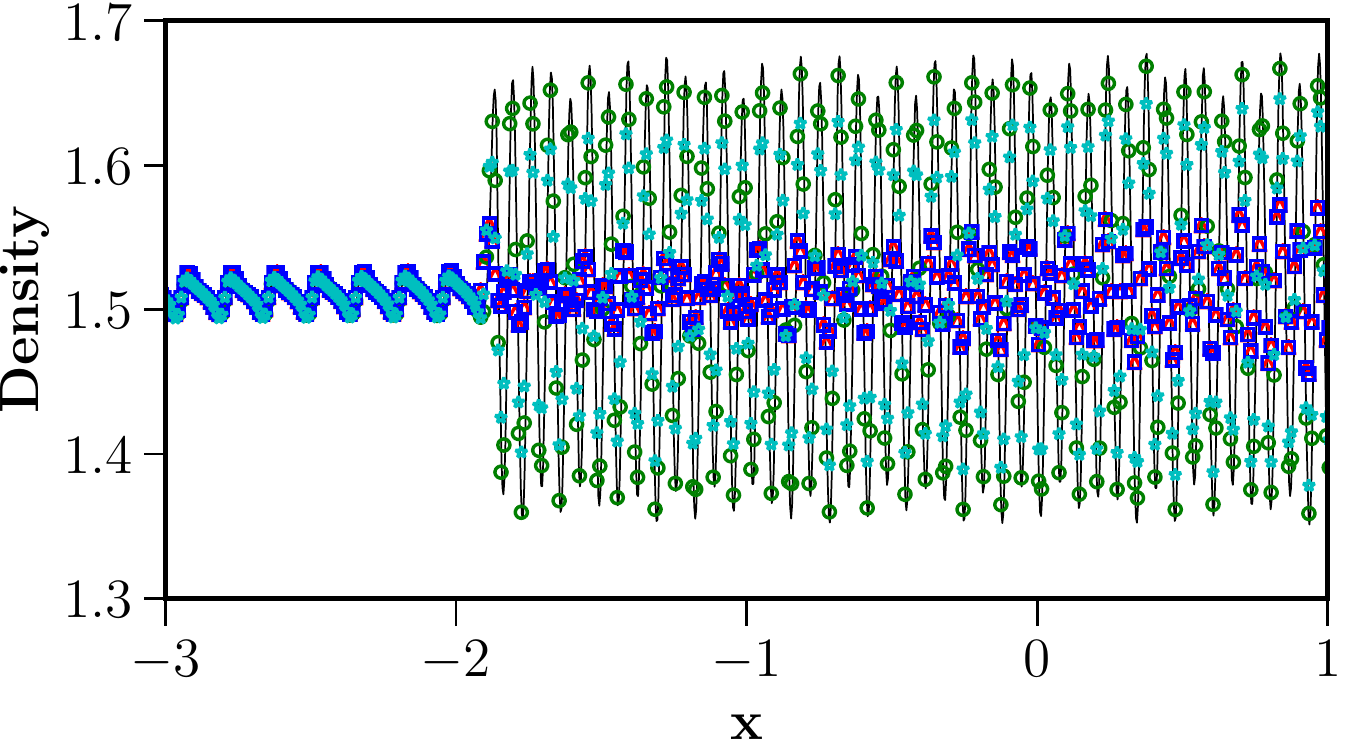}
\label{fig:tita1local}}
\subfigure[Case-2, local profile]{\includegraphics[width=0.48\textwidth]{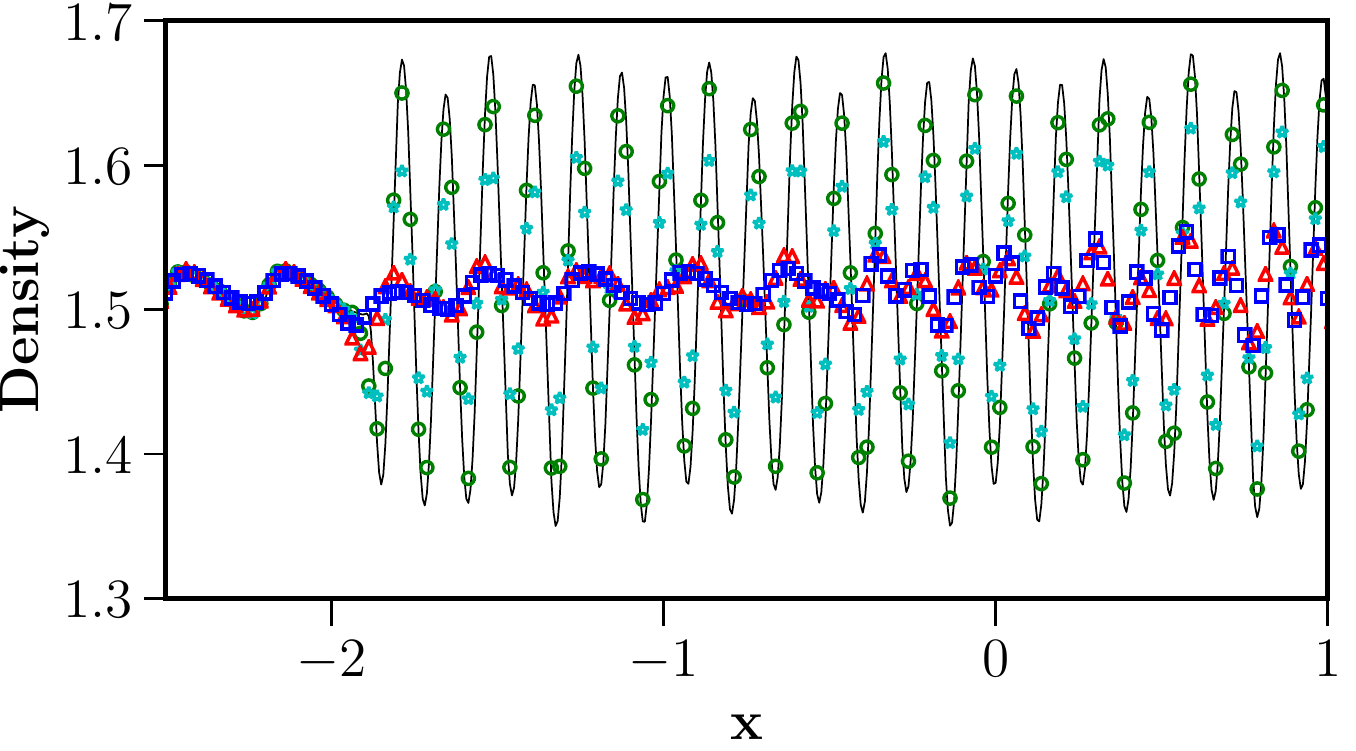}
\label{fig:tita2local}}
\caption{Local density profile of Fig. \ref{fig_titaglobal} for the Titarev-Toro problem, Example \ref{Titarev-Toro}, at $t=5.0$ by using WENO-Z, MP5, HOCUS5 and HOCUS6 schemes. Fig. (a) Case-1 with $N=1000$ and Fig. (b) Case-2 with $N=400$, where solid line: reference solution; cyan stars: HOCUS5, green circles: HOCUS6; blue squares: WENO-Z; red triangles: MP5.}
\label{fig_titazoom}
\end{figure}
%----------------------------------------%----------------------------------------%----------------------------------------%----------------------------------------%----------------------------------------
%----------------------------------------%----------------------------------------%----------------------------------------%----------------------------------------%----------------------------------------
%----------------------------------------%----------------------------------------%----------------------------------------%----------------------------------------%----------------------------------------
%----------------------------------------%----------------------------------------%----------------------------------------%----------------------------------------%----------------------------------------
%----------------------------------------%----------------------------------------%----------------------------------------%----------------------------------------%----------------------------------------
%----------------------------------------%----------------------------------------%----------------------------------------%----------------------------------------%----------------------------------------
\begin{example}\label{sod}{Shock tube problems}
\end{example}
 Here, we present the simulation results for the new algorithm for \textcolor{black}{ three shock-tube problems, namely  Sod's problem \cite{sod1978survey}, Lax's problem \cite{lax1954weak} and Le Blanc problem \cite{loubere2005subcell}.  For Sod and Lax problems} the specific heat ratio is $\gamma=1.4$ and the computational domain is $[0,1]$. Exact solutions for both cases are presented on a grid size of $1000$ points by an exact Riemann solver \cite{toro2009riemann}. 

First, Sod's  problem is calculated to verify the shock-capturing ability of HOCUS6 scheme. The initial conditions for the Sod problem are given below:

\begin{align}\label{sod_prob}
(\rho,u,p)=
\begin{cases}
&(0.125,\ \ 0,\ \ 0.1),\quad 0<x<0.5,\\
&(1,\ \ 0,\ \ 1),\quad 0.5<x<1.
\end{cases}
\end{align}

In the test case, grid points are set to be 100. All results are obtained on a $N=100$ grid and numerical predictions of density and pressure obtained from HOCUS6 scheme at time $t=0.2$ are shown in Fig. \ref{fig_sod}. No overshoots at the discontinuities are observed and the results are in good agreement with the exact solution. 

Second, Lax's problem is considered with the following initial conditions:

\begin{align}\label{lax_prob}
(\rho,u,p)=
\begin{cases}
&(0.445,\ \ 0.698,\ \ 3.528),\quad 0<x<0.5,\\
&(0.5,\ \ 0,\ \ 0.571),\quad 0.5<x<1,
\end{cases}
\end{align}
Numerical results on a $N=200$ grid for density and velocity obtained from HOCUS6 scheme at time $t=0.14$ are shown in Fig.\ \ref{fig_lax}. The HOCUS6 scheme accurately resolves the shock and contact discontinuities without over- and under-shoots. \textcolor{black}{Numerical tests for these test cases using HOCUS5 scheme are shown in \ref{sec-3.2.1}}

\textcolor{black}{Third, we consider the Le Blanc problem \cite{loubere2005subcell}, which is an extreme shock-tube problem, with the following initial conditions:}

\textcolor{black}{\begin{align}\label{blanc_prob}
(\rho,u,p)=
\begin{cases}
&(1.0,\ \ 0,\ \ \frac{2}{3}\times10^{-1}),\quad 0<x<3.0,\\
&(10^{-3},\ \ 0,\ \ \frac{2}{3}\times10^{-10}),\quad 3.0<x<9,
\end{cases}
\end{align}}

\textcolor{black}{The specific heat ratio for this test case is $\frac{5}{3}$ and the final time is $t=6$ and the exact solution is computed by an exact Riemann solver \cite{toro2009riemann}. Numerical results on a $N=200$ grid for density and pressure obtained from both HOCUS5 and HOCUS6 schemes are shown in Fig.\ \ref{fig_blanc}. There is a small error in the shock position which decreased when the grid resolution is increased (not shown here).}

\begin{figure}[H]
\centering
\subfigure[Density]{\includegraphics[width=0.35\textwidth]{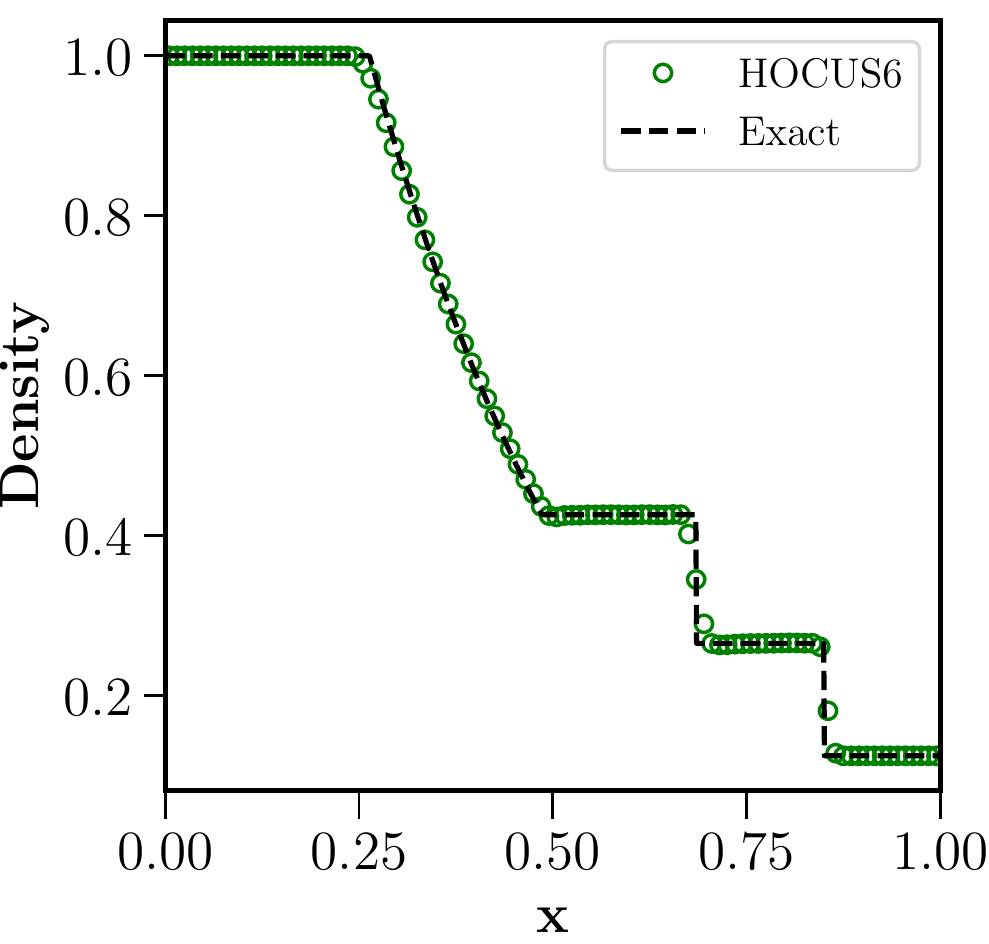}
\label{fig:sod-den}}
\subfigure[Pressure]{\includegraphics[width=0.35\textwidth]{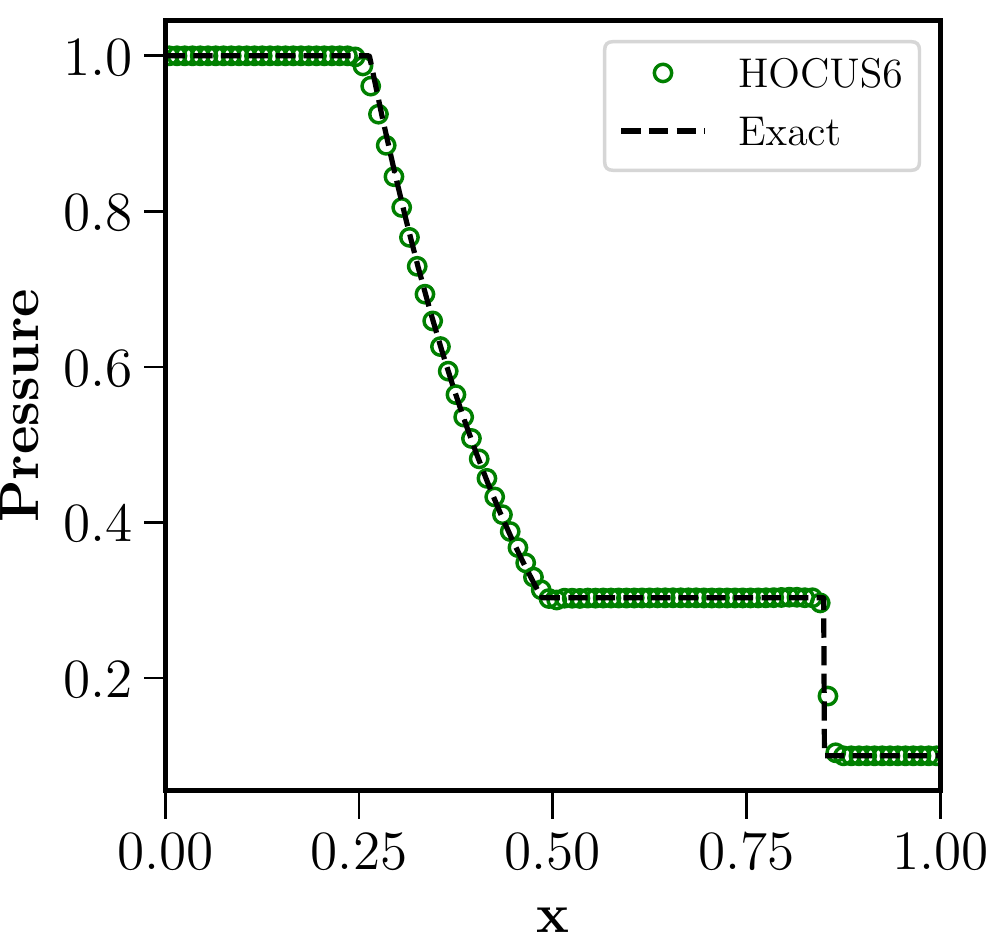}
\label{fig:sod-pres}}
\caption{Numerical solution for Sod problem in Example \ref{sod}}
\label{fig_sod}
\end{figure}

\begin{figure}[H]
\centering
\subfigure[Density]{\includegraphics[width=0.35\textwidth]{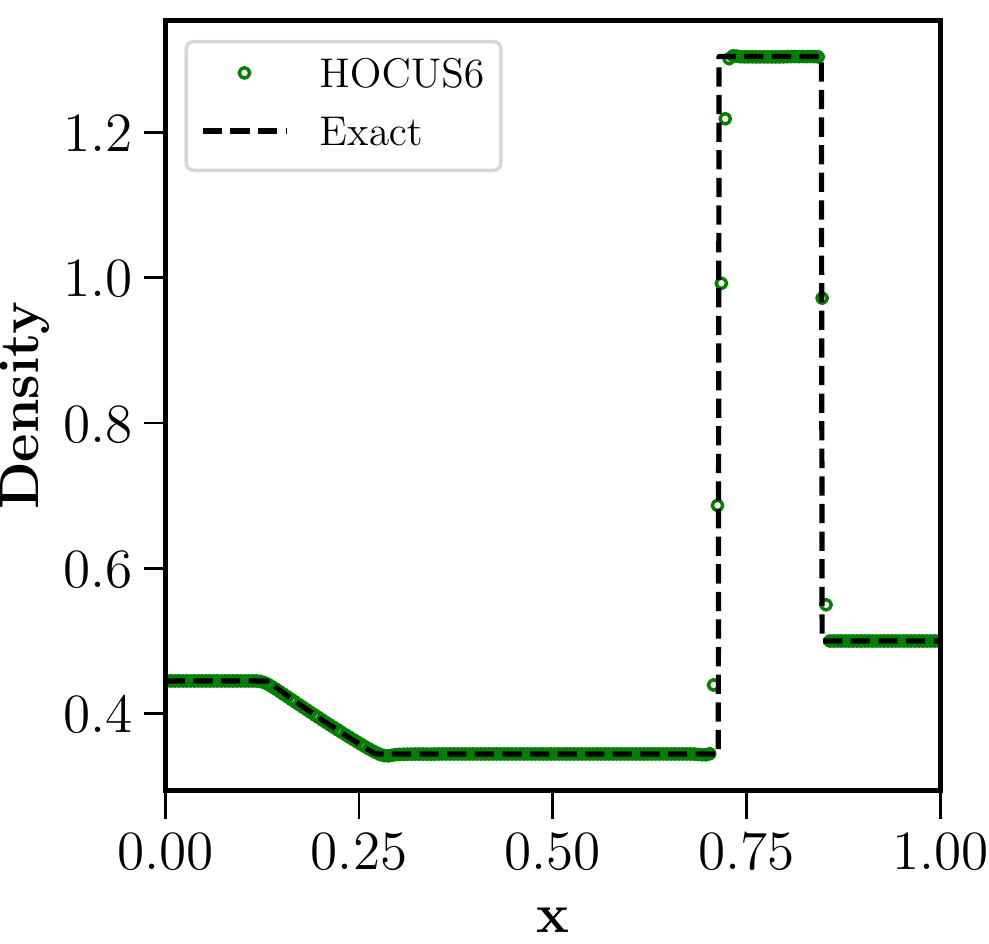}
\label{fig:lax-den}}
\subfigure[Velocity]{\includegraphics[width=0.35\textwidth]{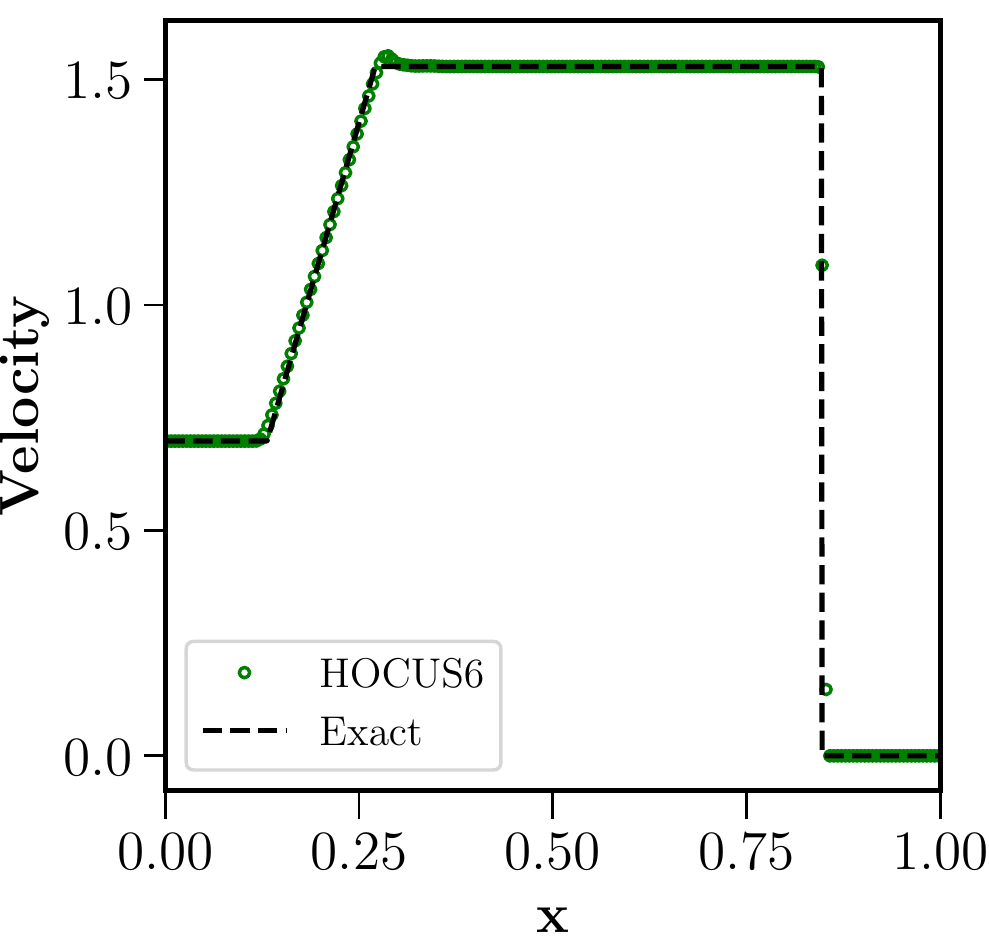}
\label{fig:lax-pres}}
\caption{Numerical solution for Lax problem in Example \ref{sod}}
\label{fig_lax}
\end{figure}

\begin{figure}[H]
\centering
\subfigure[Density]{\includegraphics[width=0.35\textwidth]{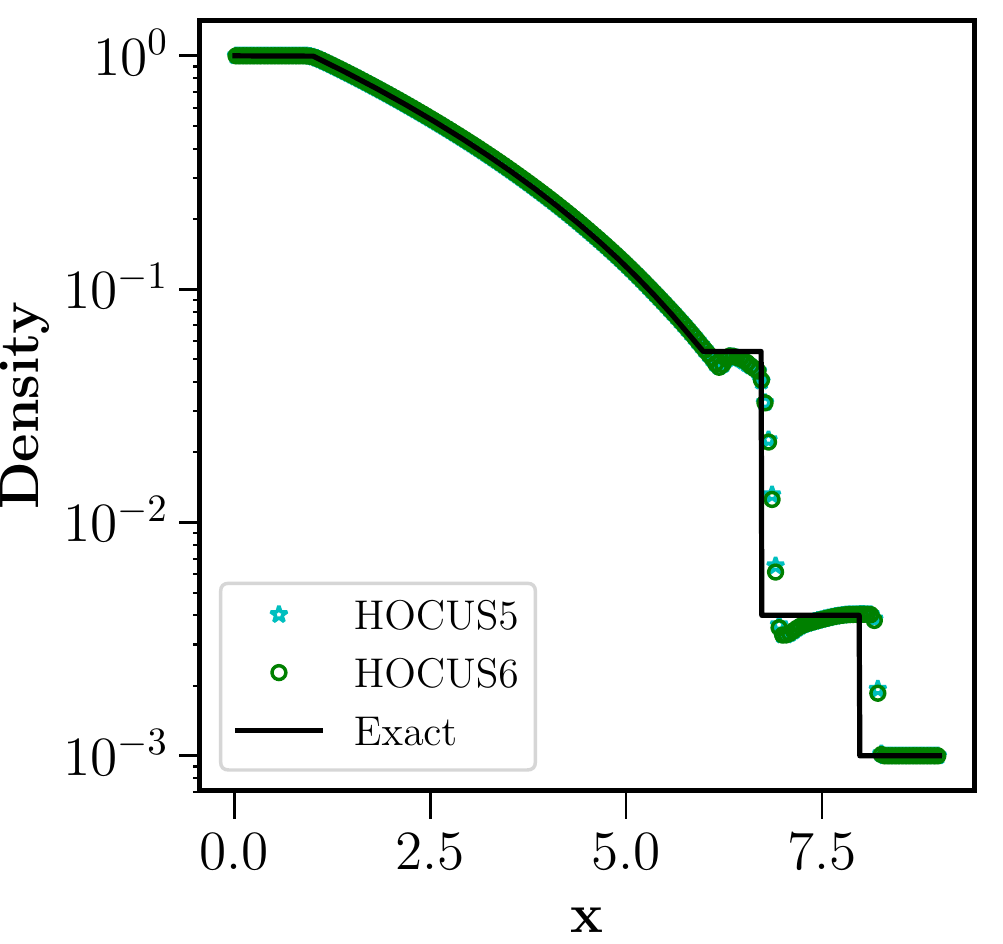}
\label{fig:blanc-den}}
\subfigure[Pressure]{\includegraphics[width=0.35\textwidth]{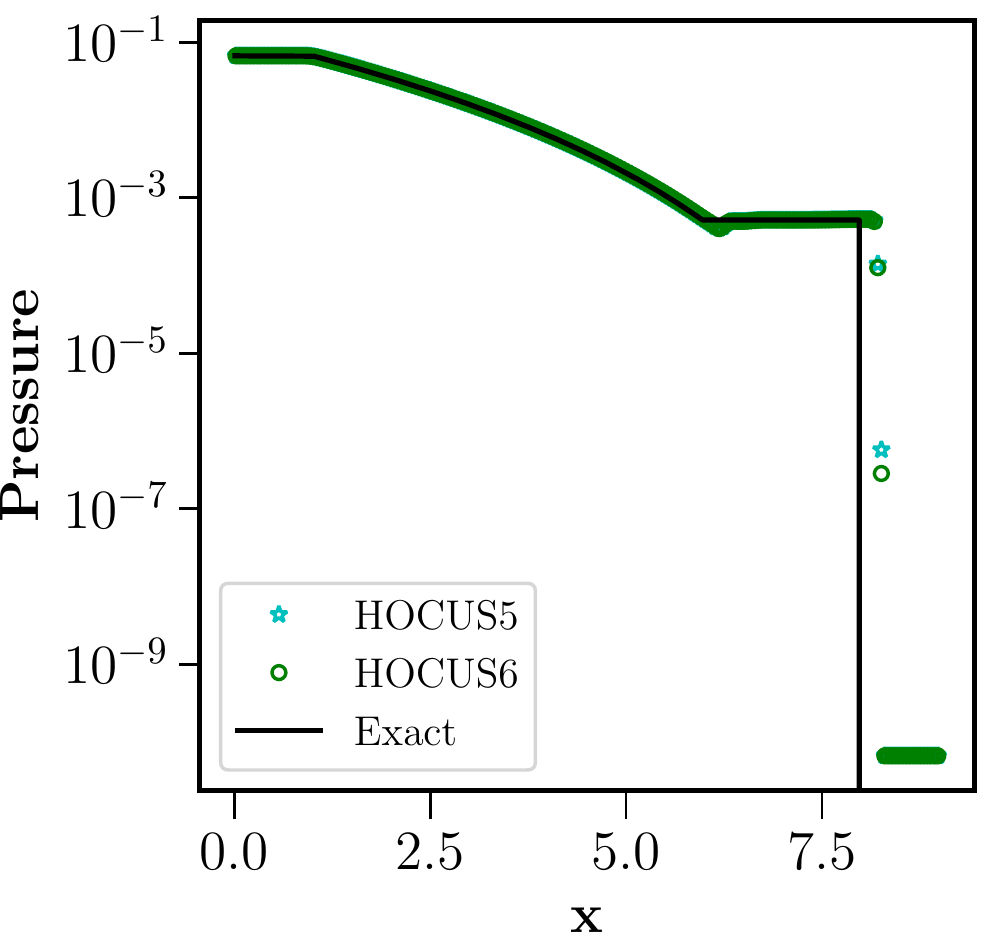}
\label{fig:blanc-pres}}
\caption{\textcolor{black}{Numerical solution for Le Blanc problem in Example \ref{sod}}}
\label{fig_blanc}
\end{figure}

%----------------------------------------%----------------------------------------%----------------------------------------%----------------------------------------%----------------------------------------
%----------------------------------------%----------------------------------------%----------------------------------------%----------------------------------------%----------------------------------------
%----------------------------------------%----------------------------------------%----------------------------------------%----------------------------------------%----------------------------------------
%----------------------------------------%----------------------------------------%----------------------------------------%----------------------------------------%----------------------------------------
%----------------------------------------%----------------------------------------%----------------------------------------%----------------------------------------%----------------------------------------
%----------------------------------------%----------------------------------------%----------------------------------------%----------------------------------------%----------------------------------------

\begin{example}\label{Shu-Osher}{Shu-Osher problem}
\end{example}
In this test case, the Shu-Osher problem \cite{Shu1988} is computed. This problem simulates the interaction of a right-moving Mach $3$ shock with a perturbed density field. The following initial data is considered:
\begin{align}\label{shock_den_prob}
(\rho,u,p)=
\begin{cases}
&(3.857143,\ \ 2.629369,\ \ 10.3333),\quad -5<x<-4,\\
&(1+0.2\sin(5x),\ \ 0,\ \ 1),\quad -4<x<5,
\end{cases}
\end{align}
subject to zero-gradient boundary conditions. Fig. \ref{fig:1d-SO} shows numerical results for the MP5, HOCUS6 and WCNS-Z schemes at time $t=1.8$ on a grid of $N=300$. The \textit{exact} solution is computed using WENO-Z on a much finer grid of $N=1600$. It is shown that all schemes give satisfactory approximations of the wavelike structures behind the shock. But the HOCUS6 scheme displays less dissipation and does a better job resolving the wave phenomenon compared to the WENO-Z and MP5 schemes.

\begin{figure}[H]
\centering\offinterlineskip
\subfigure[Global profile]{\includegraphics[width=0.35\textwidth]{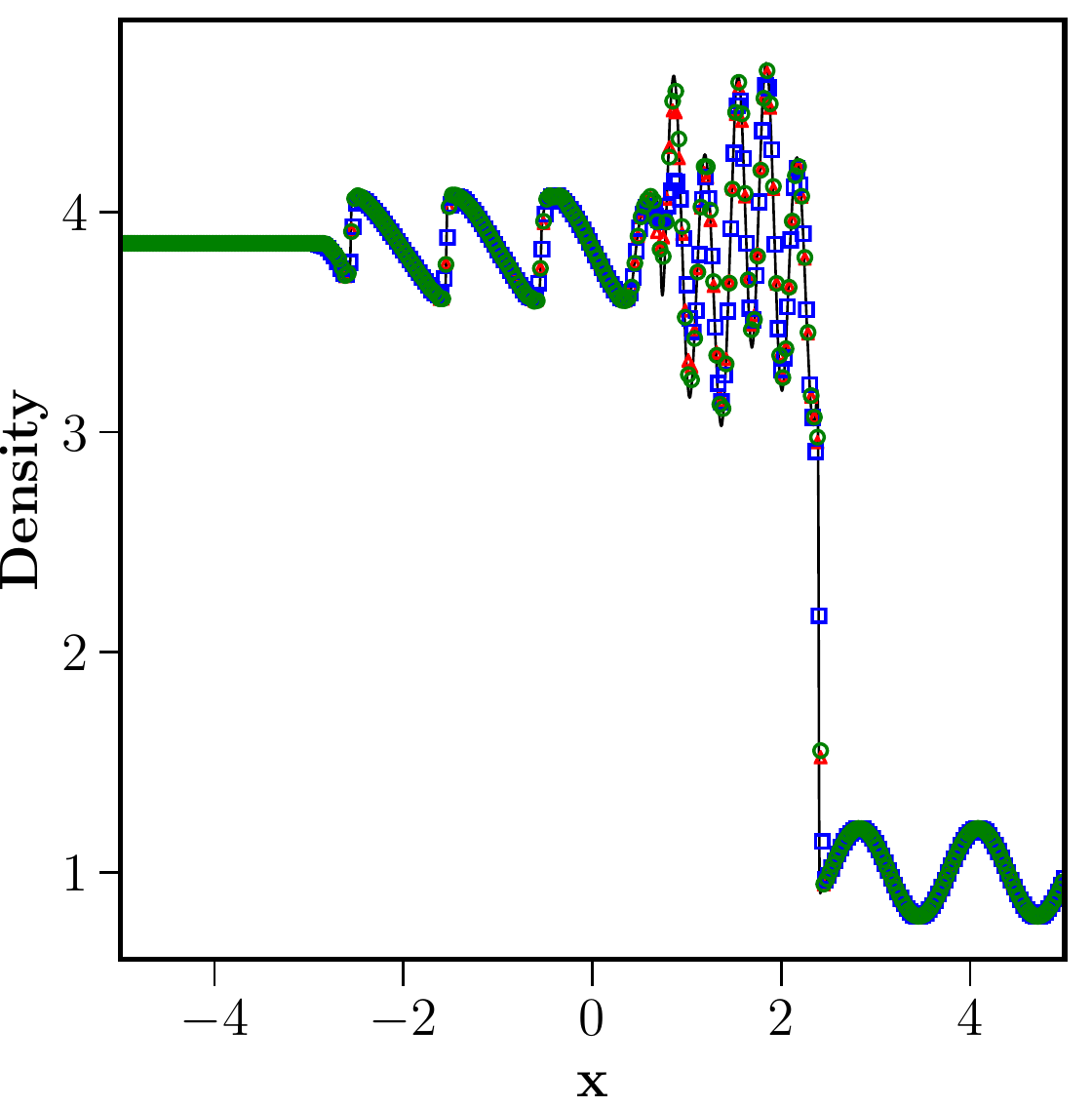}
\label{fig:shu-global}}
\subfigure[Local profile]{\includegraphics[width=0.375\textwidth]{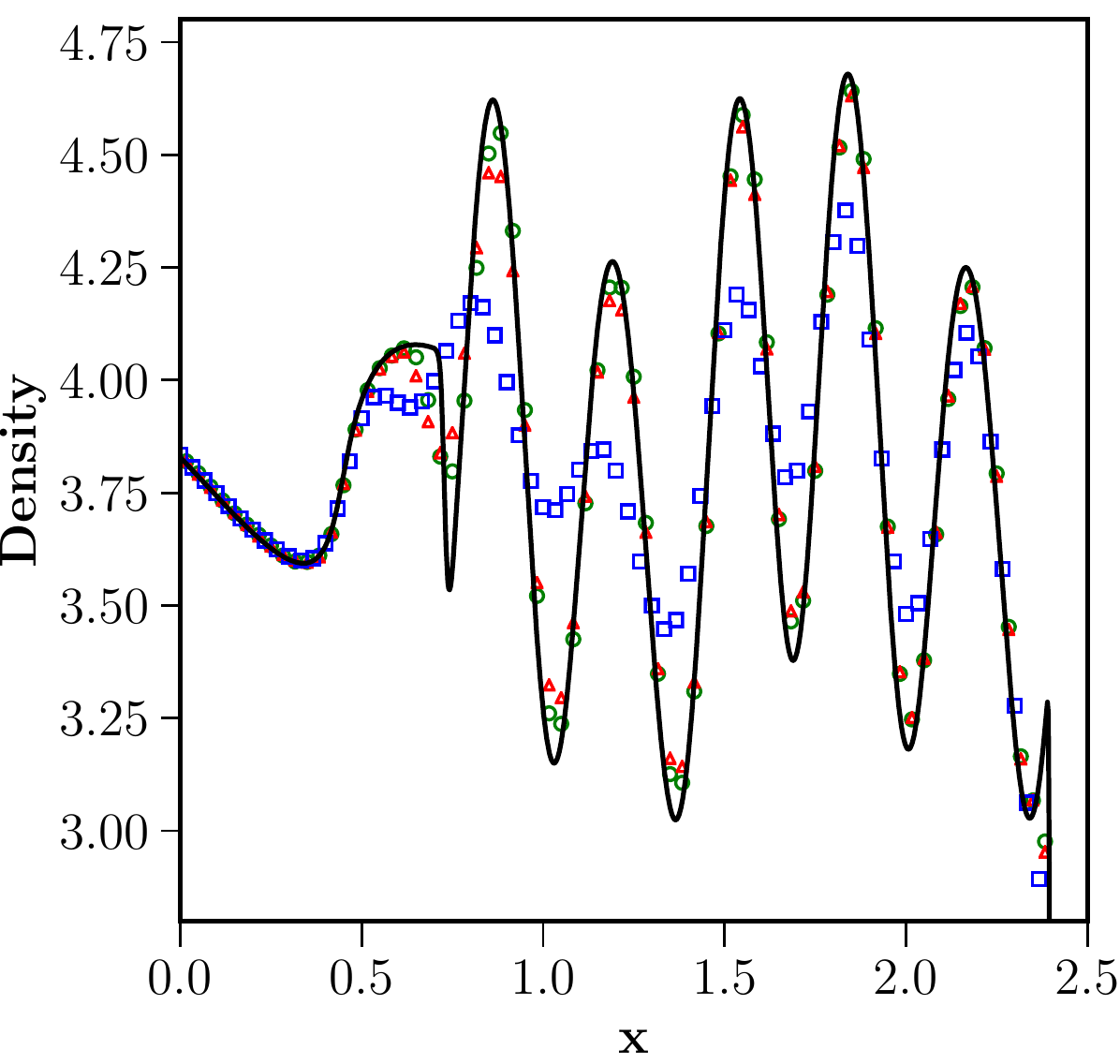}
\label{fig:shu-local}}
\caption{Density profile for Shu-Osher problem, Example \ref{Shu-Osher}, on a $N=300$ points, where solid line: reference solution; green circles: HOCUS6; blue squares: WENO-Z; red triangles: MP5.}
\label{fig:1d-SO}
\end{figure}
%----------------------------------------%----------------------------------------%----------------------------------------%----------------------------------------%----------------------------------------
%----------------------------------------%----------------------------------------%----------------------------------------%----------------------------------------%----------------------------------------
%----------------------------------------%----------------------------------------%----------------------------------------%----------------------------------------%----------------------------------------
%----------------------------------------%----------------------------------------%----------------------------------------%----------------------------------------%----------------------------------------
%----------------------------------------%----------------------------------------%----------------------------------------%----------------------------------------%----------------------------------------
%----------------------------------------%----------------------------------------%----------------------------------------%----------------------------------------%----------------------------------------
\begin{example}\label{blast}{Blast wave problem}
\end{example}
In this last one-dimensional test case, we consider the problem taken from \cite{woodward1984numerical} which involves the interaction of two blast waves with the following initial condition. 
 \begin{equation}
(\rho, u, p)=\left\{\begin{array}{ll}
        (1.0,0.0,1000) ~~~~  0.0<x<0.1, \\
        (1.0,0.0,0.01) ~ ~~~~0.1<x<0.9.\\
        (1.0,0.0,100) ~ ~~~~~ 0.9<x<1.0.
        \end{array}\right.
\end{equation}
and 
Reflective boundary conditions at both ends are used for this simulation. The numerical solutions are computed for time $t=0.038$ over a domain of $[0,1]$ with $N=400$ and $N=800$ cells and depicted in Fig. \ref{fig_blast}. Once again, we have compared the numerical solutions with the reference solution, denoted as \textit{Exact}, which is obtained using WENO-Z scheme on a grid of $1600$ points. It can be seen from Fig. \ref{fig:blast-800} that good agreement is observed compared with the reference solution. HOCUS6 shows an improved solution near the right peak, see the zoomed inset, at about $x = 0.78$, compared to the results of HOCUS5, MP5 and WENO-Z due to reduced numerical dissipation. We also observe no oscillations with the present scheme despite using non-dissipative interpolation for this test case as opposed to the P4T2 scheme proposed by Deng et al. in \cite{deng2019fifth} (see their Fig. 13) which shows small oscillations in density right after the shock. 

\begin{figure}[H]
\centering
\subfigure[400-cells]{\includegraphics[width=0.35\textwidth]{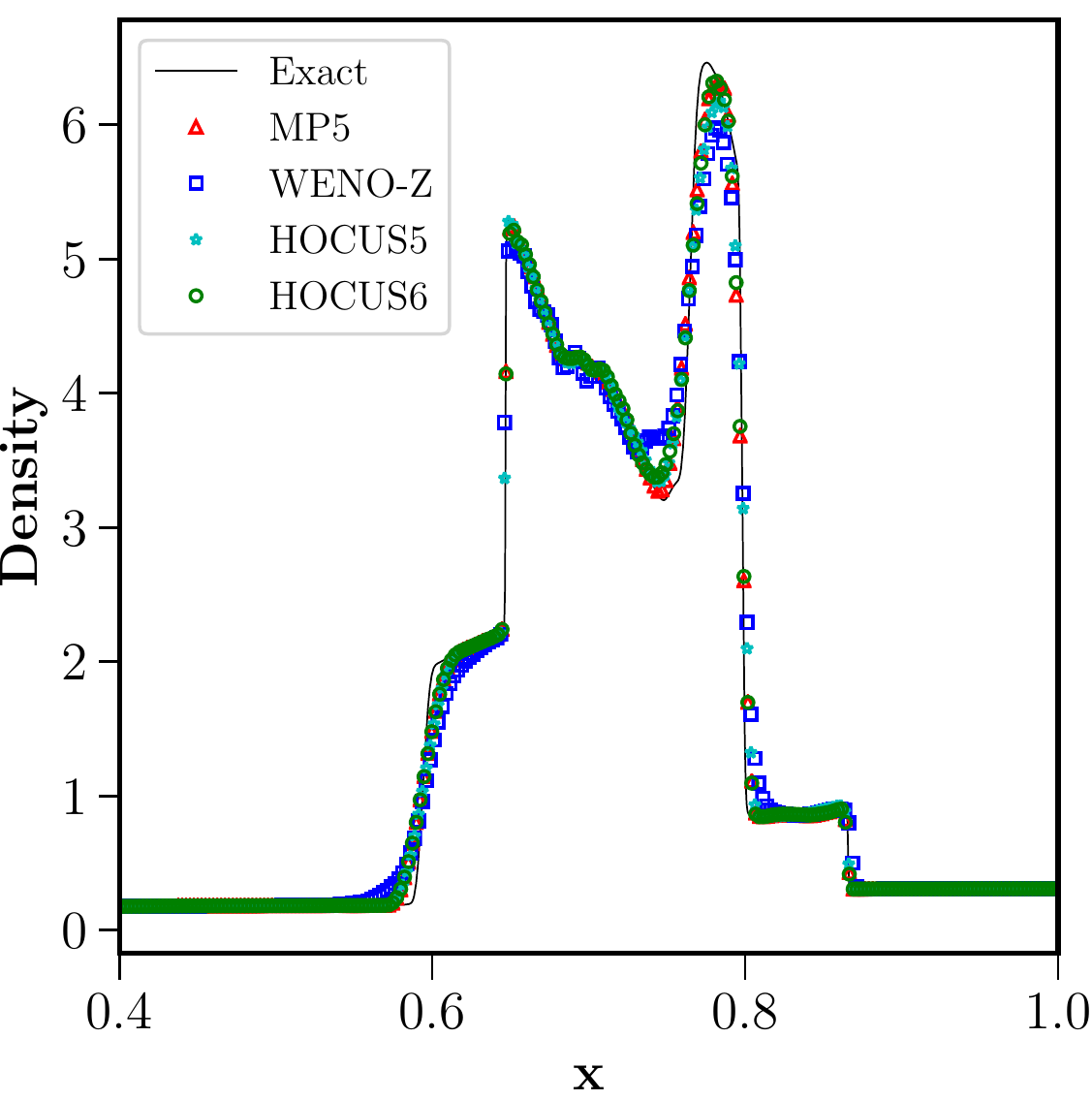}
\label{fig:blast-400}}
\subfigure[800-cells]{\includegraphics[width=0.35\textwidth]{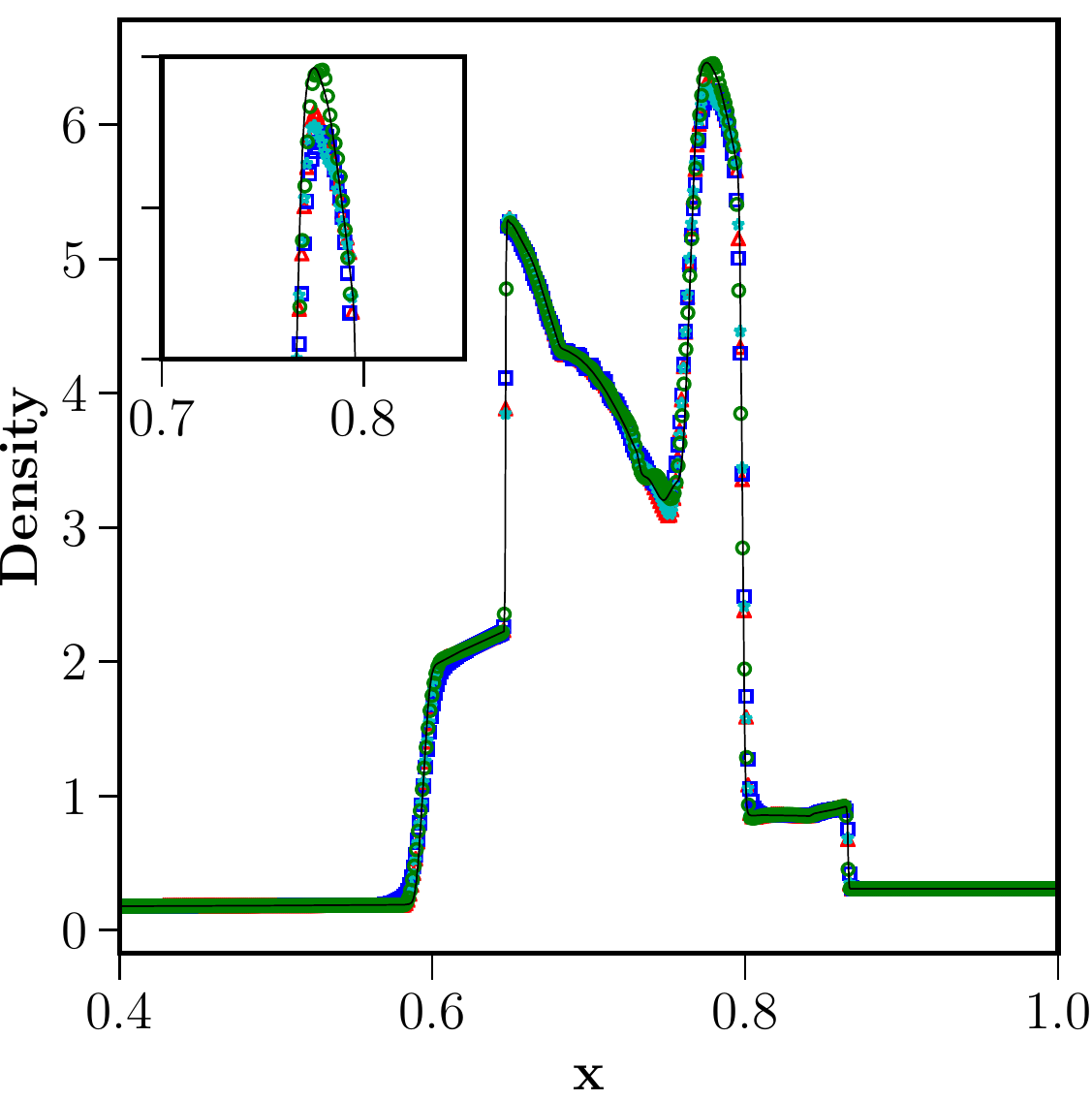}
\label{fig:blast-800}}
\caption{Density profiles obtained by various schemes for Example \ref{blast}. Solution obtained on a grid on $N=400$ and $N=800$ cells are shown in left and right panels, respectively.}
\label{fig_blast}
\end{figure}
%----------------------------------------%----------------------------------------%----------------------------------------%----------------------------------------%----------------------------------------
%----------------------------------------%----------------------------------------%----------------------------------------%----------------------------------------%----------------------------------------
%----------------------------------------%----------------------------------------%----------------------------------------%----------------------------------------%----------------------------------------
%----------------------------------------%----------------------------------------%----------------------------------------%----------------------------------------%----------------------------------------
%----------------------------------------%----------------------------------------%----------------------------------------%----------------------------------------%----------------------------------------
%----------------------------------------%----------------------------------------%----------------------------------------%----------------------------------------%----------------------------------------
\subsubsection{Analysis of BVD algorithm through HOCUS5 and C5 schemes}\label{sec-3.2.1}	

We carried out further numerical experiments for the test cases \ref{Titarev-Toro}, \ref{sod} and \ref{Shu-Osher} using the linear fifth-order compact scheme C5 given by Equation (\ref{eqn:upwind-compact}), i.e. without any nonlinear shock-capturing mechanism and compared with the HOCUS5 scheme for better understanding of the BVD algorithm. The results are shown in Fig.\ \ref{fig_unlim}. All the numerical simulations  \textcolor{black}{are carried out} with a CFL of $0.2$. Primitive variables are directly interpolated to the cell interfaces for the C5 scheme. For Titarev-Toro problem in Example \ref{Titarev-Toro} grid size of $N = 1000$ is used, for Sod and Lax problems in Example \ref{sod} we considered $N=200$ and finally for Shu-Osher problem in Example \ref{Shu-Osher} we carried out the simulations on the grid size of $N=300$. We make the following observations:

\begin{itemize}
\item In Figs. \ref{fig:unlim-sod} and \ref{fig:unlim-lax} we can see the density profiles for the Sod and Lax problems. We can see that the oscillations are only observed near the discontinuities for the C5 scheme, and it is sufficient to \textit{correct} the cells in that particular region as the linear scheme along with the Riemann solver tolerates the oscillations. We also carried out the same tests with the explicit linear schemes, for example, Equation (\ref{eqn:linear5}), and observed similar results, albeit with a significant numerical dissipation.% Corresponding figures are shown in Appendix B, Fig.\ \ref{fig_unlime5}.
\item Figs. \ref{fig:unlim-shu} and \ref{fig:unlim-tita}  show the density profile for the Shu-Osher problem and Titarev-Toro problems. Even though there are small oscillations near the shock-lets for the C5 scheme, the overall wave-like structures are well-preserved and are significantly better than the results obtained by the MP5 and WENO-Z schemes shown previously. The superior dispersion and dissipation properties of the compact schemes are helpful in this regard. Local density profile for the Titarev-Toro test case, shown in Fig.\ref{fig:unlim-tita-compare}, indicates that the results obtained by HOCUS5 scheme and the compact scheme C5 are almost identical. \textcolor{black}{As we can observe from the test cases presented in this section, the HOCUS6 is superior to HOCUS5 and also captures the discontinuities without oscillations.}
\begin{figure}[H]
\centering\offinterlineskip
\subfigure[]{\includegraphics[width=0.33\textwidth]{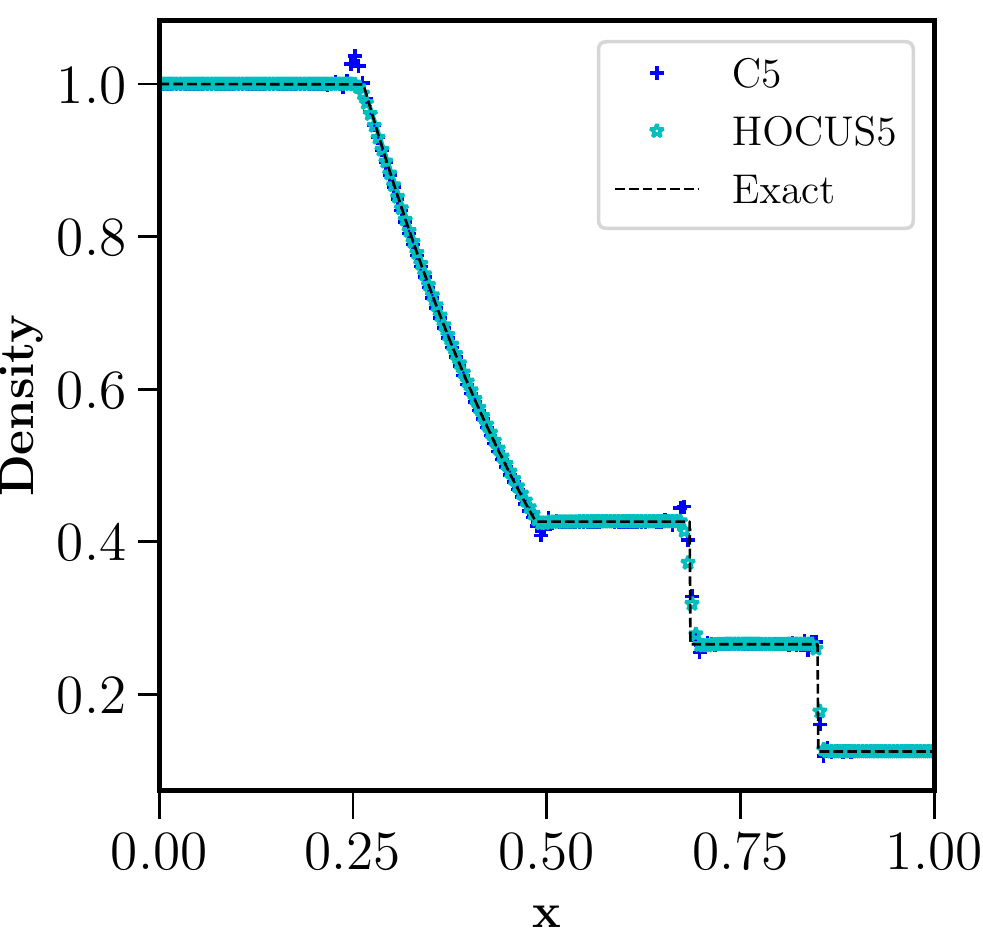}
\label{fig:unlim-sod}}
\subfigure[]{\includegraphics[width=0.33\textwidth]{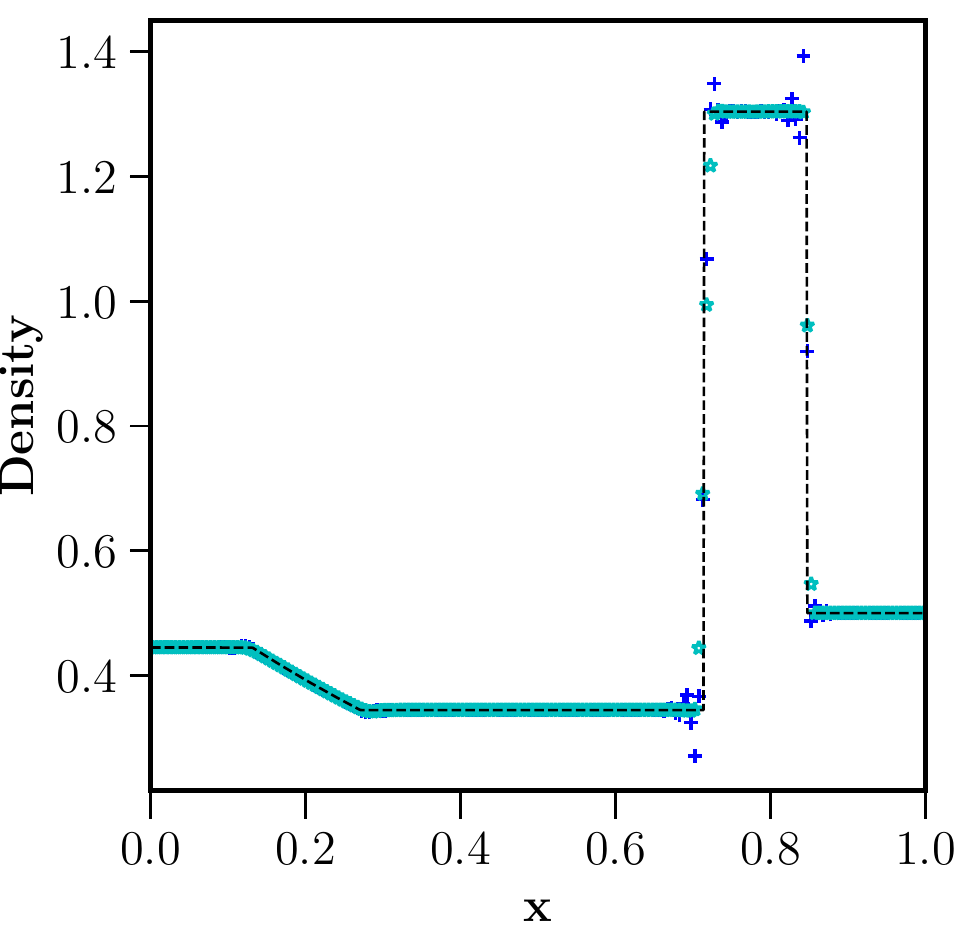}
\label{fig:unlim-lax}}
\subfigure[]{\includegraphics[width=0.325\textwidth]{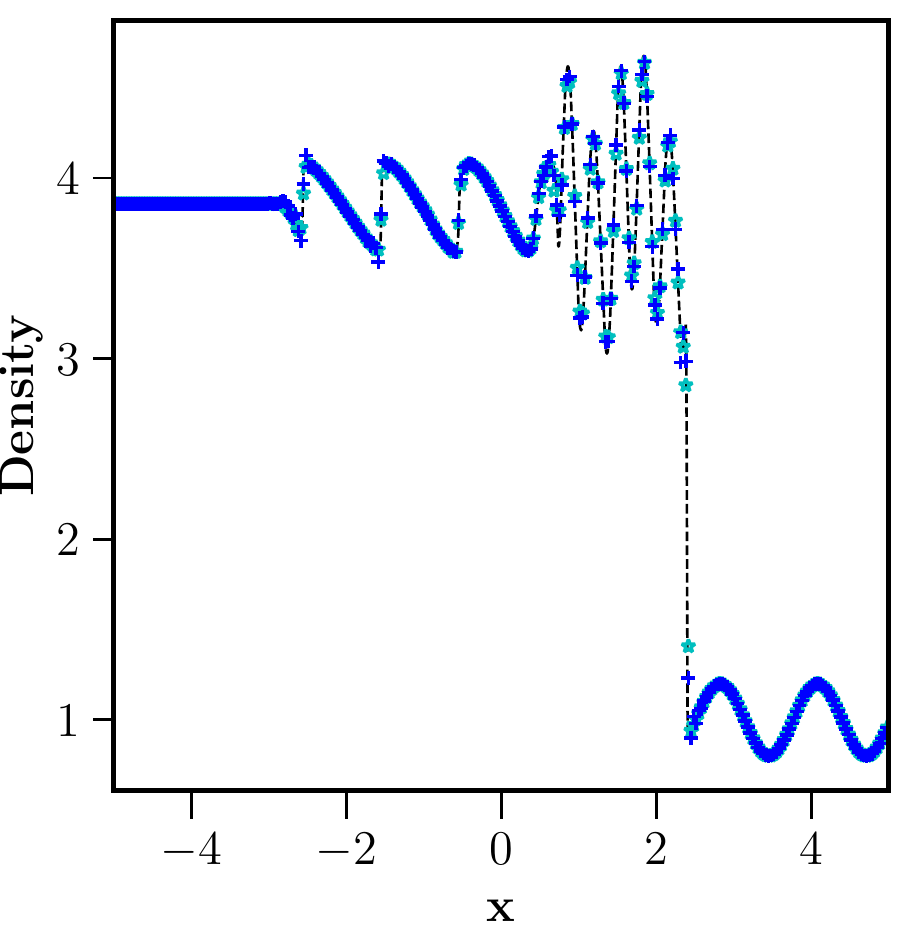}
\label{fig:unlim-shu}}
\subfigure[]{\includegraphics[width=0.325\textwidth]{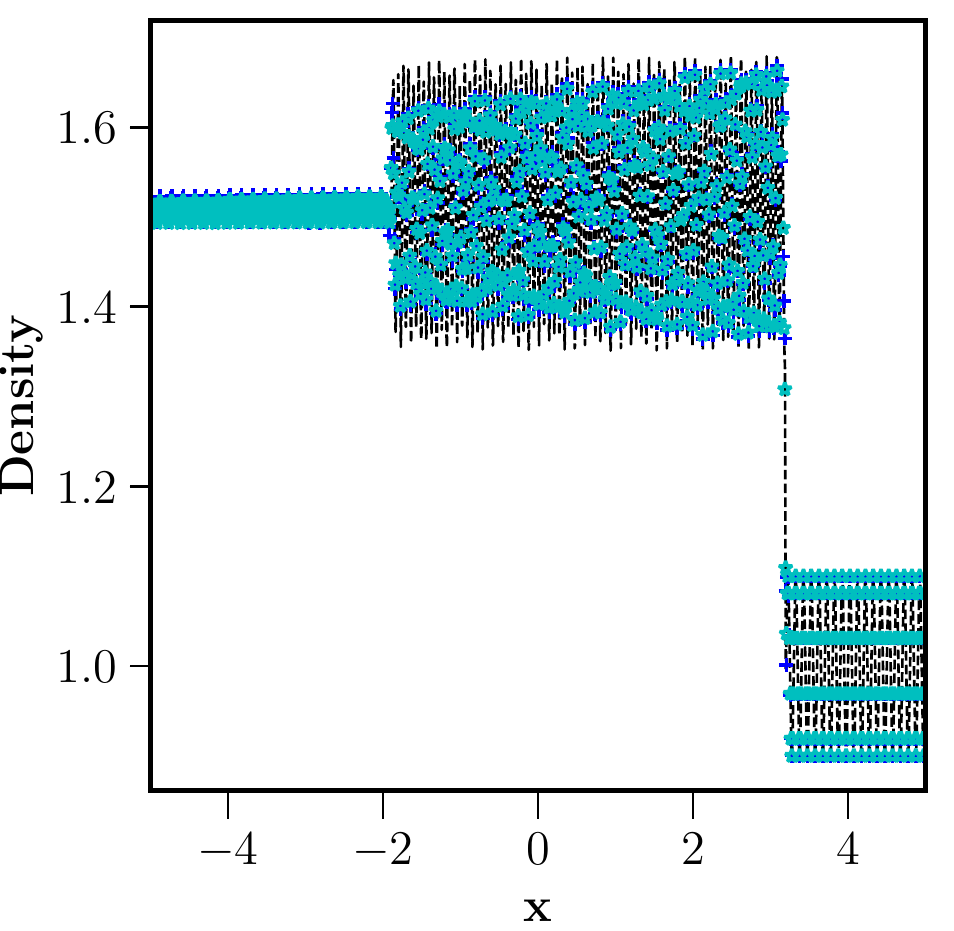}
\label{fig:unlim-tita}}
\subfigure[]{\includegraphics[width=0.475\textwidth]{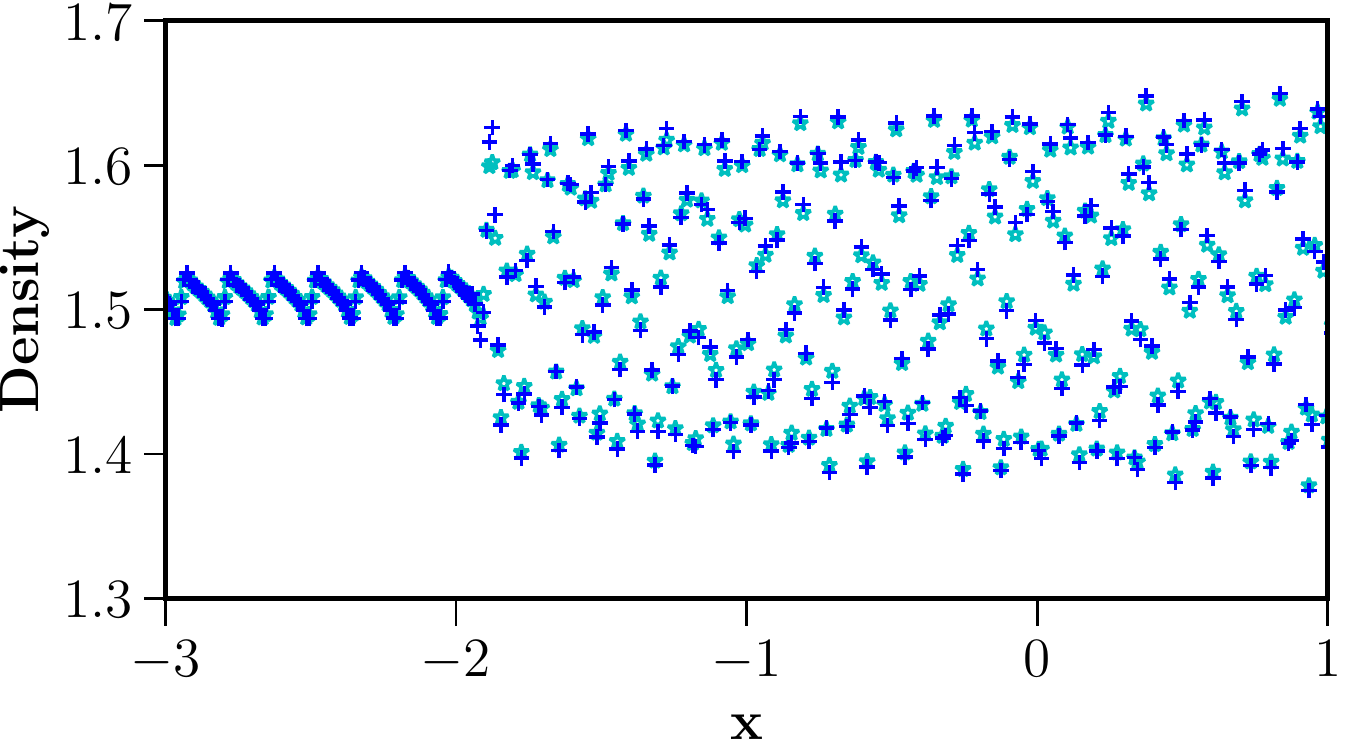}
\label{fig:unlim-tita-compare}}
\caption{Density profiles by C5 and HOCUS5 schemes for various test cases. Figs. (a) and (b) Sod and Lax problem in \ref{sod}, Fig. (b) corresponds to Example \ref{Shu-Osher} and Fig. (d) Titarev-Toro problem case-1 of \ref{Titarev-Toro}. }
\label{fig_unlim}
\end{figure}
\item Nishikawa has shown in \cite{nishikawa2018green} that for a discontinuous solution, linear compact finite-volume schemes converged even without a gradient limiter with mild oscillations even on an unstructured grid, where the other solvers diverged which motivated us to test linear-compact reconstruction scheme for these test cases. Hu et al. also showed it in Ref \cite{hu2015efficient}, although with a different shock-detection algorithm for shock-capturing, that the fifth-order linear explicit upwind scheme tolerates a small overshoot or undershoot.
\item These results indicate that the BVD algorithm automatically chooses the highest polynomial for the smooth solution, acts as a capable discontinuity detector, and in turn preserves the wave-like phenomenon. 
\end{itemize}

%----------------------------------------%----------------------------------------%----------------------------------------%----------------------------------------%----------------------------------------
%----------------------------------------%----------------------------------------%----------------------------------------%----------------------------------------%----------------------------------------
%----------------------------------------%----------------------------------------%----------------------------------------%----------------------------------------%----------------------------------------
%----------------------------------------%----------------------------------------%----------------------------------------%----------------------------------------%----------------------------------------
%----------------------------------------%----------------------------------------%----------------------------------------%----------------------------------------%----------------------------------------
%----------------------------------------%----------------------------------------%----------------------------------------%----------------------------------------%----------------------------------------
\subsection{Two-dimensional test cases for Euler equations}

In this section we show numerical results of the proposed method for two-dimensional Euler equations.

%\subsection{Case 2: Shock-turbulence interaction}
%------------------------------------------------------------------------------
\begin{example}\label{Explosion}{Explosion problem}
\end{example}
From this test case onwards we consider the two-dimensional Euler equations. In this example, the initial condition consists of two constant states for the flow variables, a circular region of radius $r=0.4$ centered at $(1,1)$ and the region outside of it as mentioned in Toro \cite{toro2009riemann}. The initial conditions are given as:
\begin{equation}\label{explosion}
(\rho, u,v, p)=\left\{\begin{array}{ll}
        (1, 0, 0, 1), ~~~~~~~~~ \mbox{If} ~~{x^2+y^2} <r^2, \\
        (0.125, 0, 0, 0.1), ~~~~~~~ \mbox{else}. \\
        \end{array}\right.
\end{equation}
In the present case, numerical simulations are carried out over a square domain of size $[0,2]\times [0,2]$ until a final time $t=0.25$ on a uniform grid of resolution $400\times 400$. In Fig.\ \ref{fig:RM-birb}, we can see the density distribution of the simulation computed by the proposed scheme. Cross-sectional slices of density and pressure along the plane $y=0$ computed using the present scheme are shown in Figs.\ \ref{fig:RM-density} and \ref{fig:RM-pressure}, respectively. In Fig.\ \ref{fig:RM-Compare} we compare the density profile obtained by various schemes where the numerical results for the reference solution are computed using WENO-Z scheme over a uniform mesh of resolution $1000\times 1000$. We can observe that the present approach resolves the discontinuities without oscillations.

\begin{figure}[H]
\centering\offinterlineskip
\subfigure[]{\includegraphics[width=0.35\textwidth]{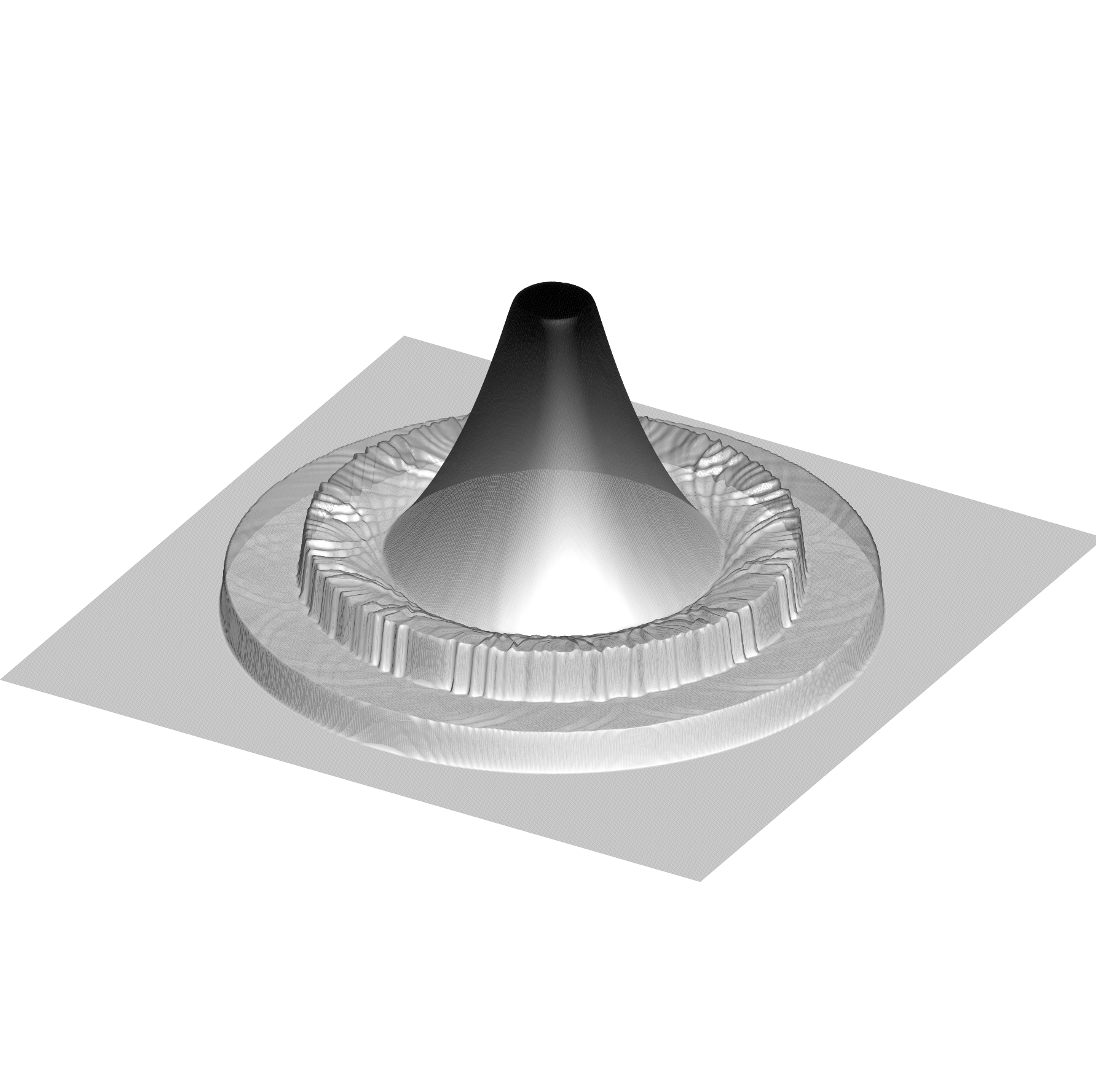}
\label{fig:RM-birb}}
\subfigure[]{\includegraphics[width=0.35\textwidth]{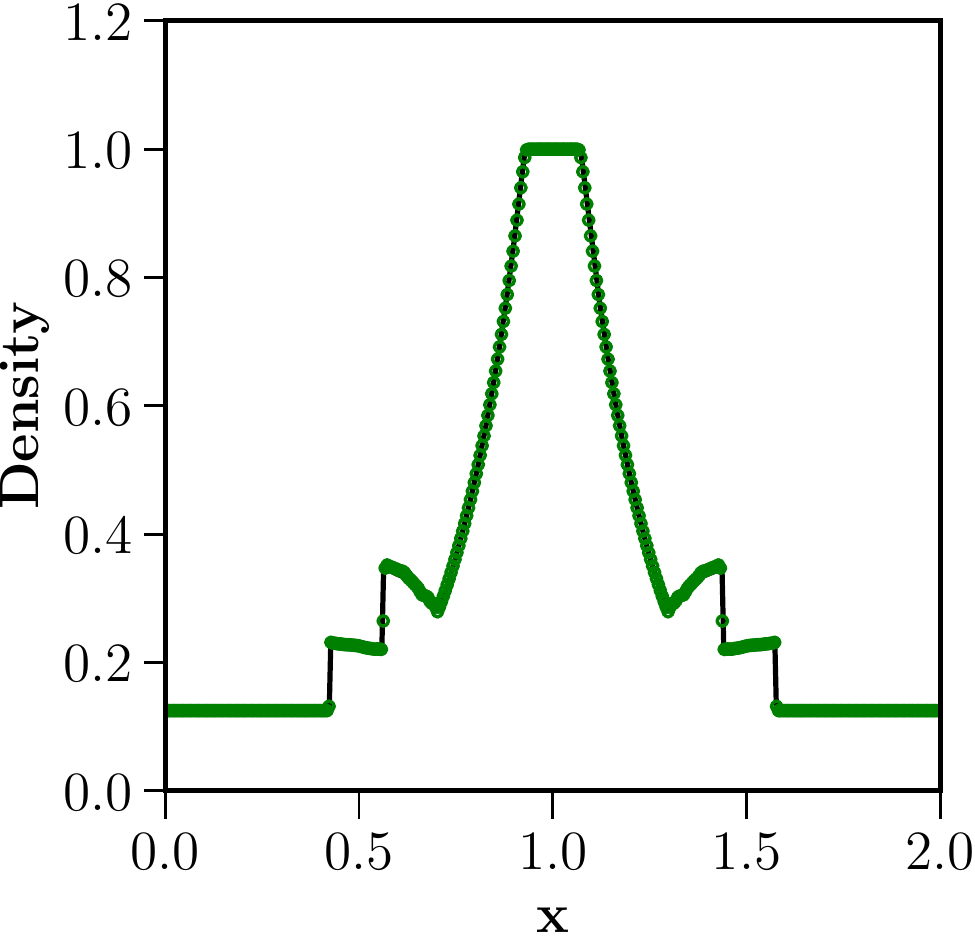}
\label{fig:RM-density}}
\subfigure[]{\includegraphics[width=0.35\textwidth]{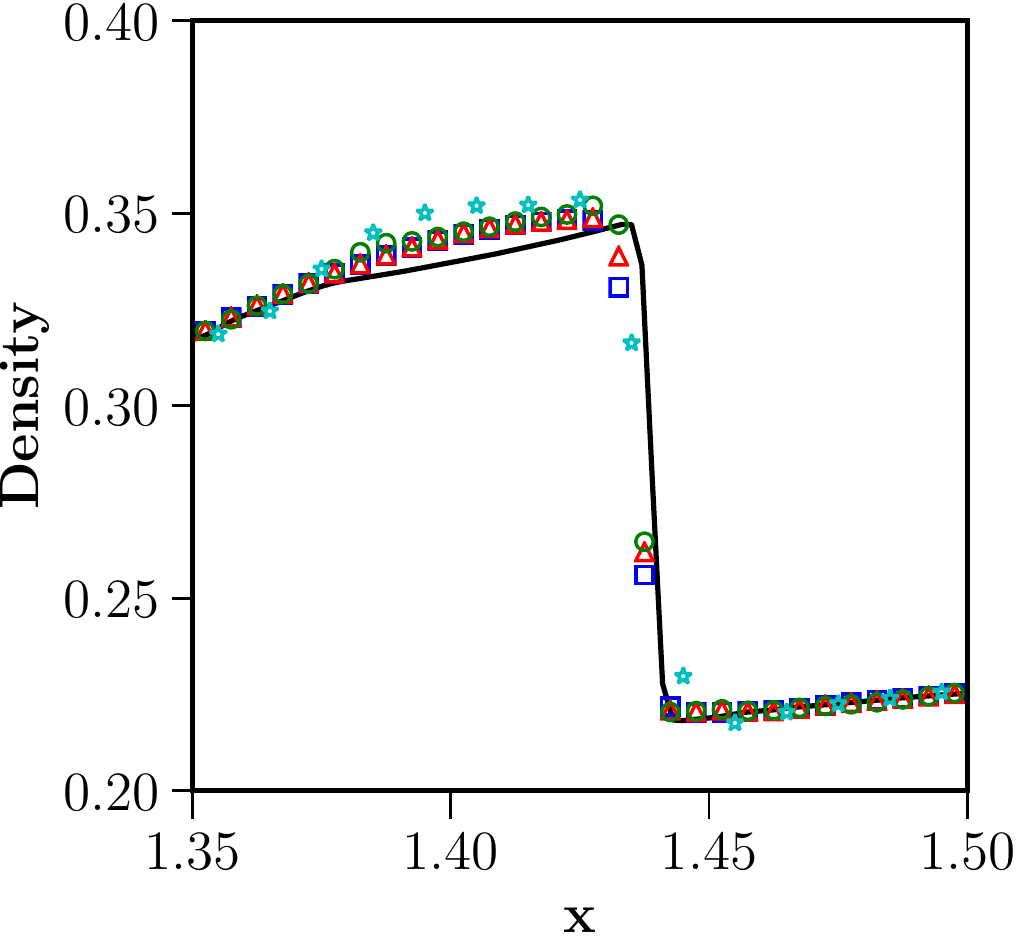}
\label{fig:RM-Compare}}
\subfigure[]{\includegraphics[width=0.35\textwidth]{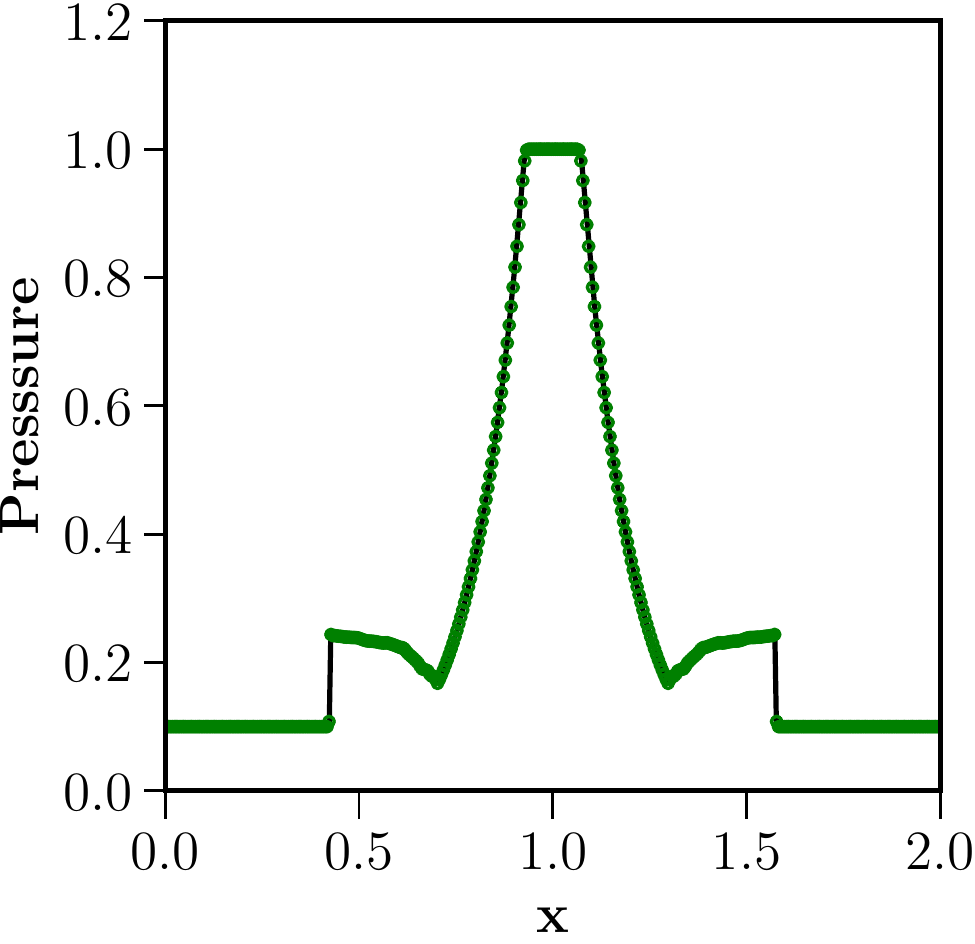}
\label{fig:RM-pressure}}
\caption{\textcolor{black}{Solution for the Example \ref{Explosion}. Fig.\ (a) Bird's eye view of the density distribution, Figs.\ (b) and (d) shows the density and pressure along the plane $y=0$ for HOCUS6 and Fig.\ (c) comparison of all the schemes for the density profile at a local region, where solid line: reference solution; green circles: HOCUS6; cyan stars: HOCUS5; blue squares: WENO-Z; red triangles: MP5.}}
\label{fig_explosion}
\end{figure}
%----------------------------------------%----------------------------------------%----------------------------------------%----------------------------------------%----------------------------------------
%----------------------------------------%----------------------------------------%----------------------------------------%----------------------------------------%----------------------------------------
%----------------------------------------%----------------------------------------%----------------------------------------%----------------------------------------%----------------------------------------
%----------------------------------------%----------------------------------------%----------------------------------------%----------------------------------------%----------------------------------------
%----------------------------------------%----------------------------------------%----------------------------------------%----------------------------------------%----------------------------------------
%----------------------------------------%----------------------------------------%----------------------------------------%----------------------------------------%----------------------------------------
\begin{example}\label{euler-accuracy}{Accuracy of 2D Euler equations}
\end{example}
In this test case we again evaluate the accuracy of the proposed schemes for the compressible Euler equations considering the following initial conditions:
\begin{align}\label{accu-euler}
(\rho,u,v,p)=
\begin{cases}
&(1+0.5 \sin(x+y),\ \ 1.0,\ \, 1.0 \ \ 1.0)\\
\end{cases}
\end{align}
 over the domain $[-1, 1]$ $\times$ $[-1, 1]$.
 Numerical solutions are computed until a final time $t=2$. The time step is taken as $\Delta t=0.1\Delta x^{2.0}$. Table \ref{tab:case3a} shows the $L_1$-errors and convergence rates of the numerical schemes. \textcolor{black}{We can observe from Table \ref{tab:case3a} that the both HOCUS5 and HOCUS6 schemes have achieved design order of accuracy and convergence rate.} 

\begin{table}[H]
 \centering
 \footnotesize
 \caption{\textcolor{black}{$L_1$ errors and numerical orders of accuracy for test case in Example \ref{euler-accuracy}. $N$ is the total number of cells on a uniform mesh. t = 2.}}
    \begin{tabular}{ccccccccc}
 \hline
          & WENOZ &       & MP5   &       & HOCUS5 &       & HOCUS6 &  \\
 \hline
          &       & order &       & order &       & order &       &  \\
 \hline
     $20^2$    & 2.27E-04 &       & 2.12E-04 &       & 5.38E-05 &       & 4.20E-05 &  \\
 \hline
    $40^2$    & 6.84E-06 & 5.05  & 6.81E-06 & 4.96  & 1.19E-06 & 5.50  & 7.24E-07 & 5.86 \\
 \hline
     $80^2$    & 2.20E-07 & 4.96  & 2.20E-07 & 4.96  & 2.97E-08 & 5.32  & 1.04E-08 & 6.12 \\
 \hline
    $160^2$   & 7.30E-09 & 4.91  & 7.13E-09 & 4.94  & 8.13E-10 & 5.19  & 1.77E-10 & 5.87 \\
 \hline
    \end{tabular}%
 \label{tab:case3a}%
\end{table}%

The following points are important to emphasize:
\begin{itemize}
\item  \textcolor{black}{In general, flux integrals should be computed using a high-order Gaussian quadrature with suitable Gaussian integration points over the faces of the control volume to achieve high-order accuracy, third-order or more,} see \cite{titarev2004finite,coralic2014finite,haimovich2017numerical}. Since the flux Jacobian $\partial F/ \partial Q$ is constant for a linear equation, the dimension-by-dimension variable extrapolation also achieves formal order of accuracy without quadratures for a linear problem. The above test case is a linear system, and therefore we can obtain the high-order of accuracy.
\item In order to obtain high-order accuracy for a nonlinear problem, a high-order Gaussian quadrature integration or any other suitable approach, has to be used, which is not considered this paper. We are of aware of this as the first author has previously used the same finite-volume schemes in \cite{chamarthi2018high} for the computation of \textit{linear} diffusion and advection-diffusion equations. For nonlinear problems, a finite-difference scheme \cite{chamarthi2019first} has been used to obtain high-order accuracy.
\item It is also demonstrated by Zhang et al. in Ref. \cite{zhang2011order} that the finite-volume WENO method with mid-point rule is only second-order accurate for nonlinear systems and Gaussian integral rule is necessary for high-order accuracy. However, they also noted that for flows involving shock waves, which are non-smooth, the resolution characteristics are often comparable despite the obvious difference in the order of accuracy. In the present paper, we are also interested in flows involving discontinuities and we have used the mid-point rule for our computations. In the original BVD algorithm of Sun et al. \cite{sun2016boundary} and the subsequent papers that improved the original scheme, the computations are also based on the finite-volume method with mid-point rule.
\item It is important to note that the present scheme interpolates the primitive or characteristic variables to the interfaces, and the BVD algorithm given by Equation (\ref{eqn:BVDstep}) will not work with conservative variables or fluxes. It is also possible to attain high-order accuracy for smooth problems without quadratures for a finite volume approach using the BVD algorithm and is currently being investigated and will be presented elsewhere.
\end{itemize}
%----------------------------------------%----------------------------------------%----------------------------------------%----------------------------------------%----------------------------------------
%----------------------------------------%----------------------------------------%----------------------------------------%----------------------------------------%----------------------------------------
%----------------------------------------%----------------------------------------%----------------------------------------%----------------------------------------%----------------------------------------
%----------------------------------------%----------------------------------------%----------------------------------------%----------------------------------------%----------------------------------------
%----------------------------------------%----------------------------------------%----------------------------------------%----------------------------------------%----------------------------------------
%----------------------------------------%----------------------------------------%----------------------------------------%----------------------------------------%----------------------------------------

\begin{example}\label{ex:rp}{Riemann Problem}
\end{example}
The third 2D problem is Riemann problem taken from Schulz-Rinne et al. \cite{schulz1993numerical}, described as configuration 3. The simulation is carried out over a unit square domain $[0,1]\times [0,1]$, with the following initial data:

 \begin{equation}
(\rho, u,v, p)=\left\{\begin{array}{ll}
        (1.5, 0, 0, 1.5), ~~~~~~~~~~~~~~~~~~~~~~~~~\mbox{if} ~~x > 0.8, ~~y > 0.8, \\
        (33/62, 4/\sqrt{11}, 0, 0.3), ~~~~~~~~~~~~~~\mbox{if} ~~ x\leq 0.8, ~~y > 0.8, \\
        (77/558, 4/\sqrt{11}, 4/\sqrt{11}, 9/310), ~~\mbox{if} ~~ x \leq 0.8, ~~y\leq 0.8, \\
        (33/62, 0, 4/\sqrt{11}, 0.3), ~~~~~~~~~~~~~~\mbox{if} ~~x > 0.8, ~~ y\leq 0.8.
        \end{array}\right.
\end{equation}
This initial condition, with constant states of primitive variables along the lines $x=0.8$ and $y=0.8$, produces four shocks at the interfaces of the four quadrants. The small-scale complex structures generated along the slip lines due to Kelvin-Helmholtz instability serves as a benchmark to test numerical dissipation of the scheme. The numerical solutions are computed for time $t=0.8$ on a grid of size $400\times 400$. Non-reflective boundary conditions are applied at all the boundaries for this test case. Figure \ref{fig_riemann} shows the density contours obtained by the various schemes. A closer look at Fig. \ref{fig:HOCUS6_RM} indicates that the roll-up of the slip line forms more vortices for the HOCUS6 scheme in comparison with the WENO-Z, \textcolor{black}{MP5 and HCOUS5 schemes}, and also the various reflected shocks are well captured despite using the non-dissipative scheme. This test case can also be called as double Mach reflection and additional comments regarding roll-up vortices are given in Example \ref{ex:dmr}.

\begin{figure}[H]
\begin{onehalfspacing}
\centering\offinterlineskip
\subfigure[WENO-Z]{\includegraphics[width=0.4\textwidth]{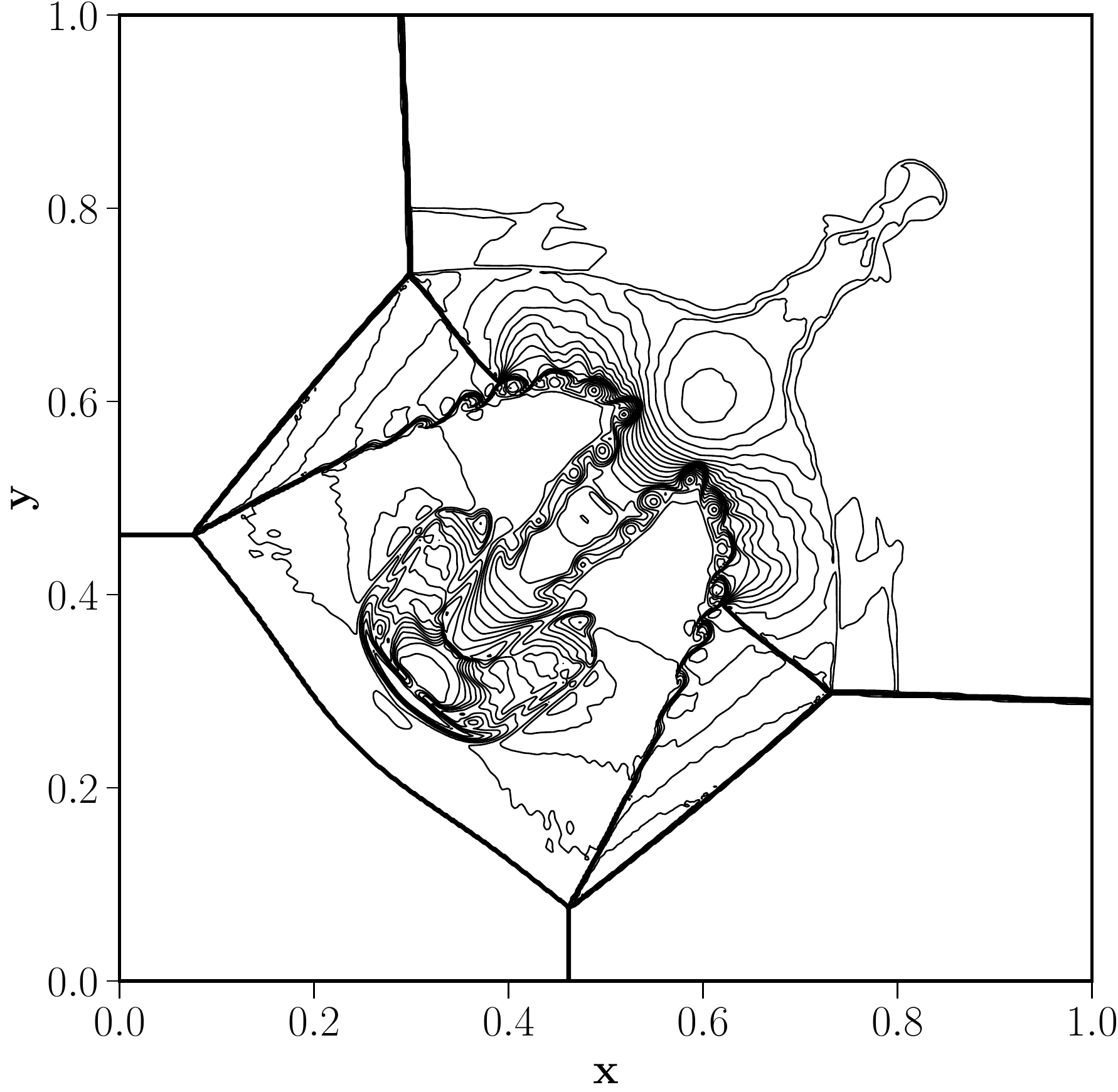}
\label{fig:WENO-Z_RM}}
\subfigure[MP5]{\includegraphics[width=0.4\textwidth]{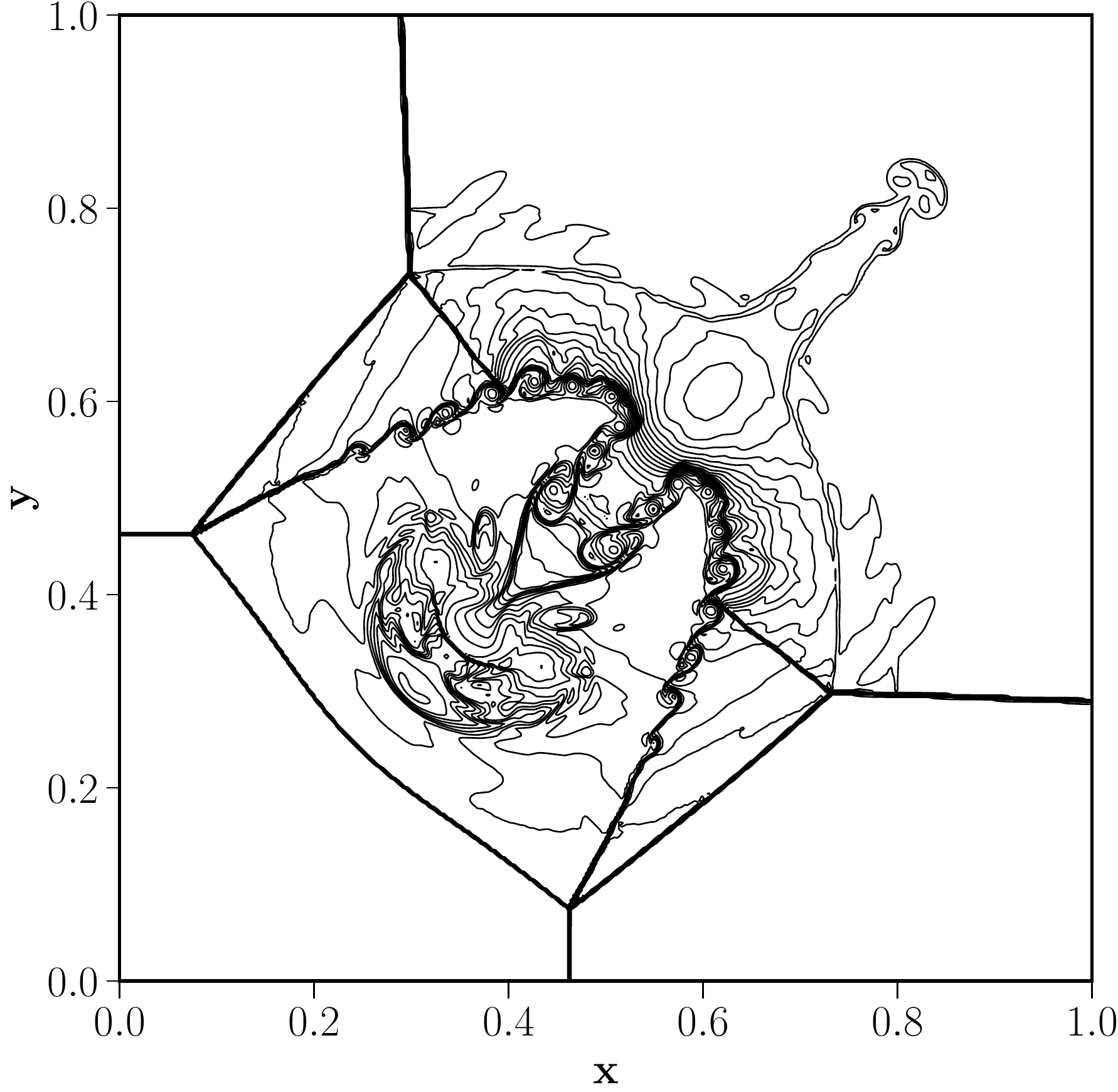}
\label{fig:MP5_RM}}
\subfigure[\textcolor{black}{HOCUS5}]{\includegraphics[width=0.4\textwidth]{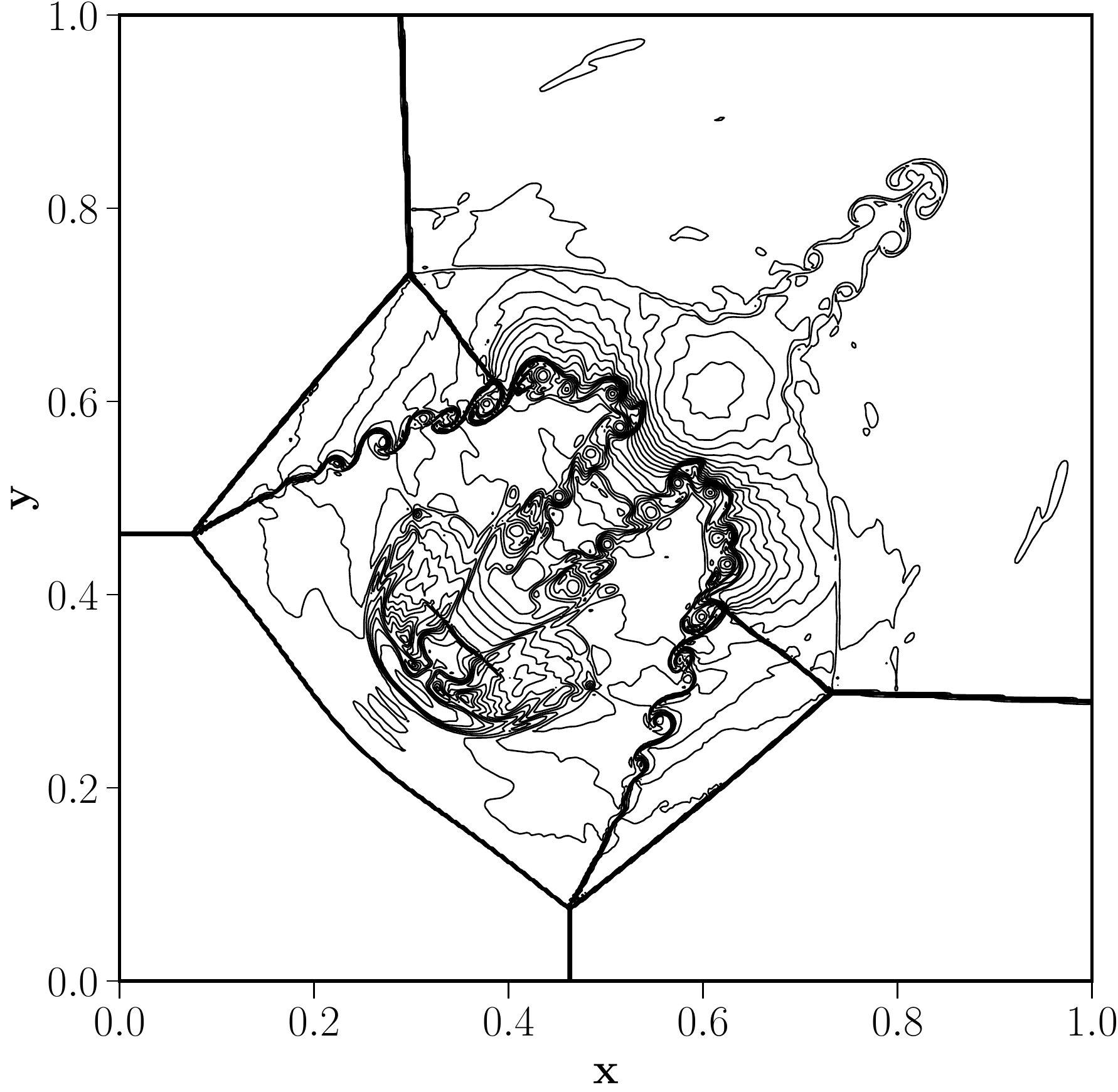}
\label{fig:HOCUS5_RM}}
\subfigure[HOCUS6]{\includegraphics[width=0.4\textwidth]{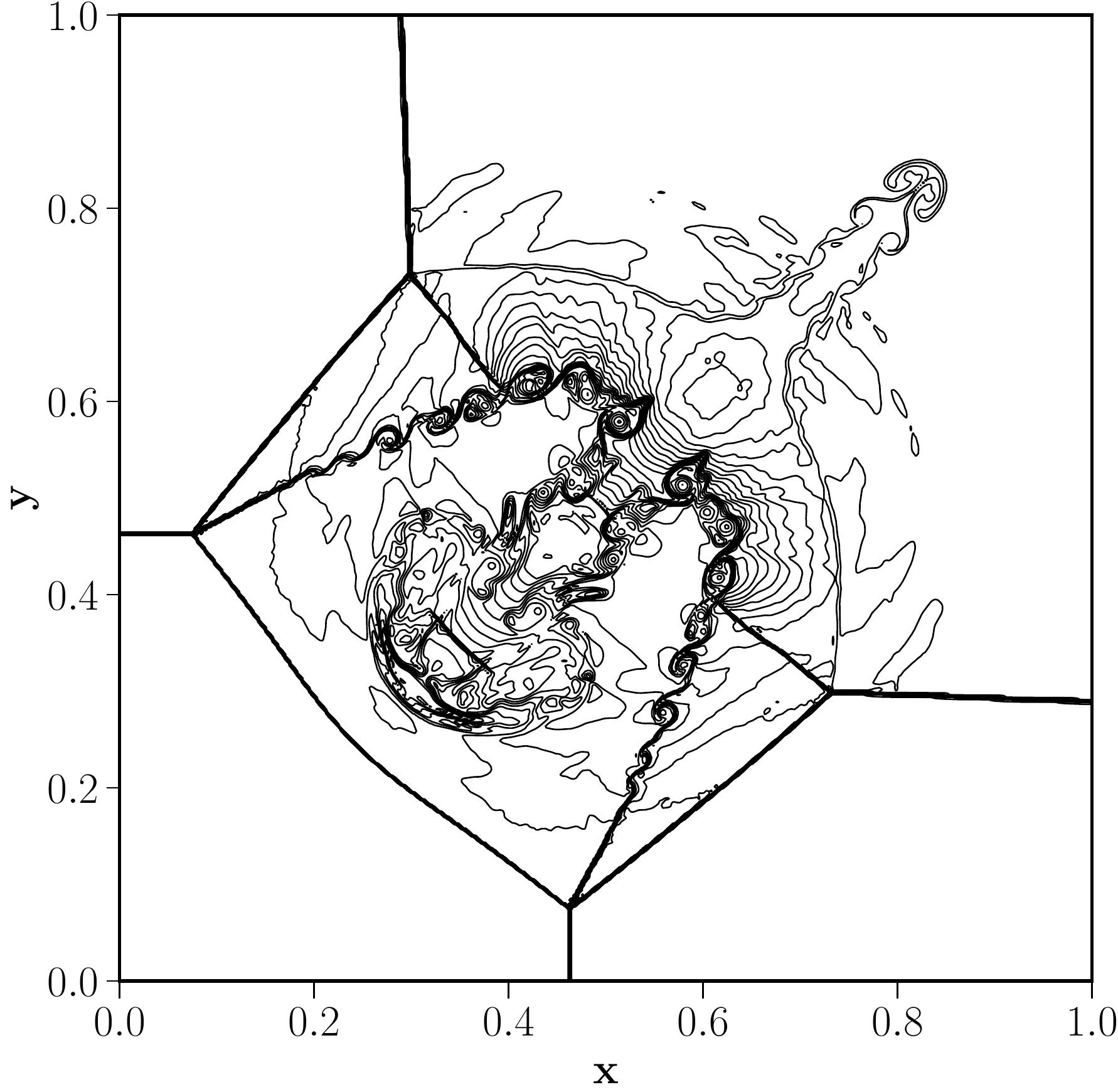}
\label{fig:HOCUS6_RM}}
\caption{\textcolor{black}{Solution of the Riemann problem, Example \ref{ex:rp}, indicating density contours obtained for WENO-Z, MP5, HOCUS5 and HOCUS6 schemes.}}
\label{fig_riemann}
\end{onehalfspacing}
\end{figure}

%----------------------------------------%----------------------------------------%----------------------------------------%----------------------------------------%----------------------------------------
%----------------------------------------%----------------------------------------%----------------------------------------%----------------------------------------%----------------------------------------
%----------------------------------------%----------------------------------------%----------------------------------------%----------------------------------------%----------------------------------------
%----------------------------------------%----------------------------------------%----------------------------------------%----------------------------------------%----------------------------------------
%----------------------------------------%----------------------------------------%----------------------------------------%----------------------------------------%----------------------------------------
%----------------------------------------%----------------------------------------%----------------------------------------%----------------------------------------%----------------------------------------
\begin{example}\label{shock-entropy}{2D shock-entropy wave test}
\end{example}
In this test case we consider the two-dimensional shock-entropy wave interaction problem proposed in \cite{acker2016improved}. The initial conditions for the test case are as follows,
\begin{align}\label{shock_entropy}
(\rho,u,p)=
\begin{cases}
&(3.857143, \ \ 2.629369,\ \ 10.3333),\quad x<-4,\\
&(1+0.2\sin(10x \cos\theta+10y\sin\theta),\ \ 0,\ \ 1),\quad otherwise,
\end{cases}
\end{align}
with $\theta$ = $\pi/6$ over a domain of $[-5,5]\times [-1,1]$. The initial sine waves make an angle of $\theta$ radians with the $x$ axis. Initial conditions are modified as in \cite{deng2019fifth} with a higher frquency for the initial sine waves compared to that of \cite{acker2016improved} to show the benefits of the proposed method. Mesh size of $400\times 80$, which corresponds to $\Delta x= \Delta y =1/40$, is used for the computations. Numerical results in Fig. \ref{fig_shock_entropy} indicate that the resolution of the flow structures is significantly improved by HOCUS6 scheme, Fig. \ref{fig:SE-BVD}, in comparison with the \textcolor{black}{HOCUS5}, WENO-Z and MP5, shown in Figs. \ref{fig:SE-MP5} and \ref{fig:SE-HOCUS5}, respectively. The local density profile in the high-frequency region along $y=0$, shown in Fig. \ref{fig:SSE-Compare}, indicates that the HOCUS6 approach predicts the density amplitudes with lower numerical dissipation compared to the other methods. The reference solution in Fig. \ref{fig:SSE-Compare} is computed on a very fine mesh of $1600 \times 320$ by the WENO-Z scheme. The current approach is also significantly better than the P4T2 scheme proposed by Deng et al. in \cite{deng2019fifth} (see their Fig. 21(d)) which improves only marginally over WENO-Z. These results indicate that the present approach is well suited for compressible shock-turbulence interactions.\\

\begin{figure}[H]
\centering\offinterlineskip
\subfigure[MP5]{\includegraphics[width=0.48\textwidth]{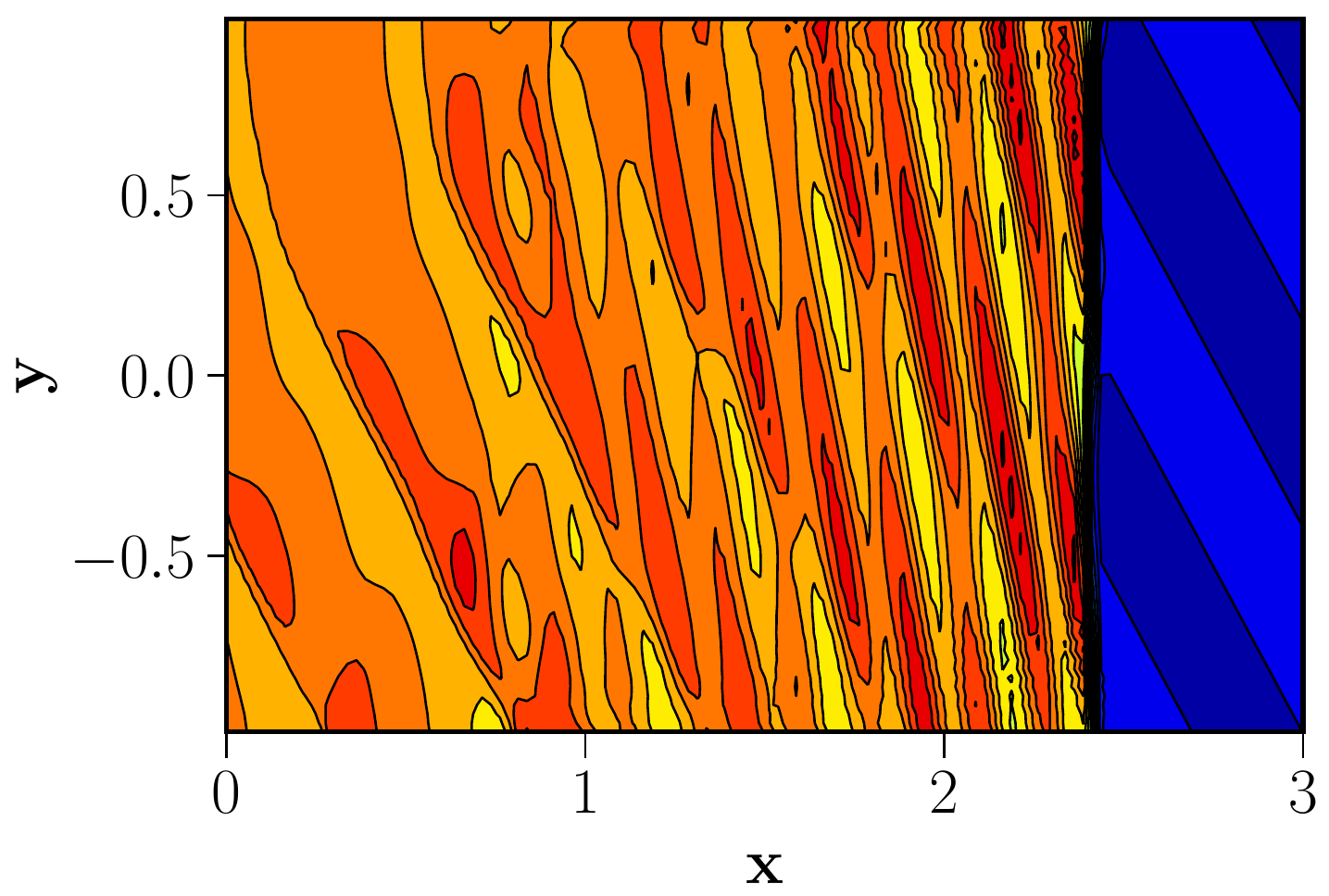}
\label{fig:SE-MP5}}
\subfigure[\textcolor{black}{HOCUS5}]{\includegraphics[width=0.48\textwidth]{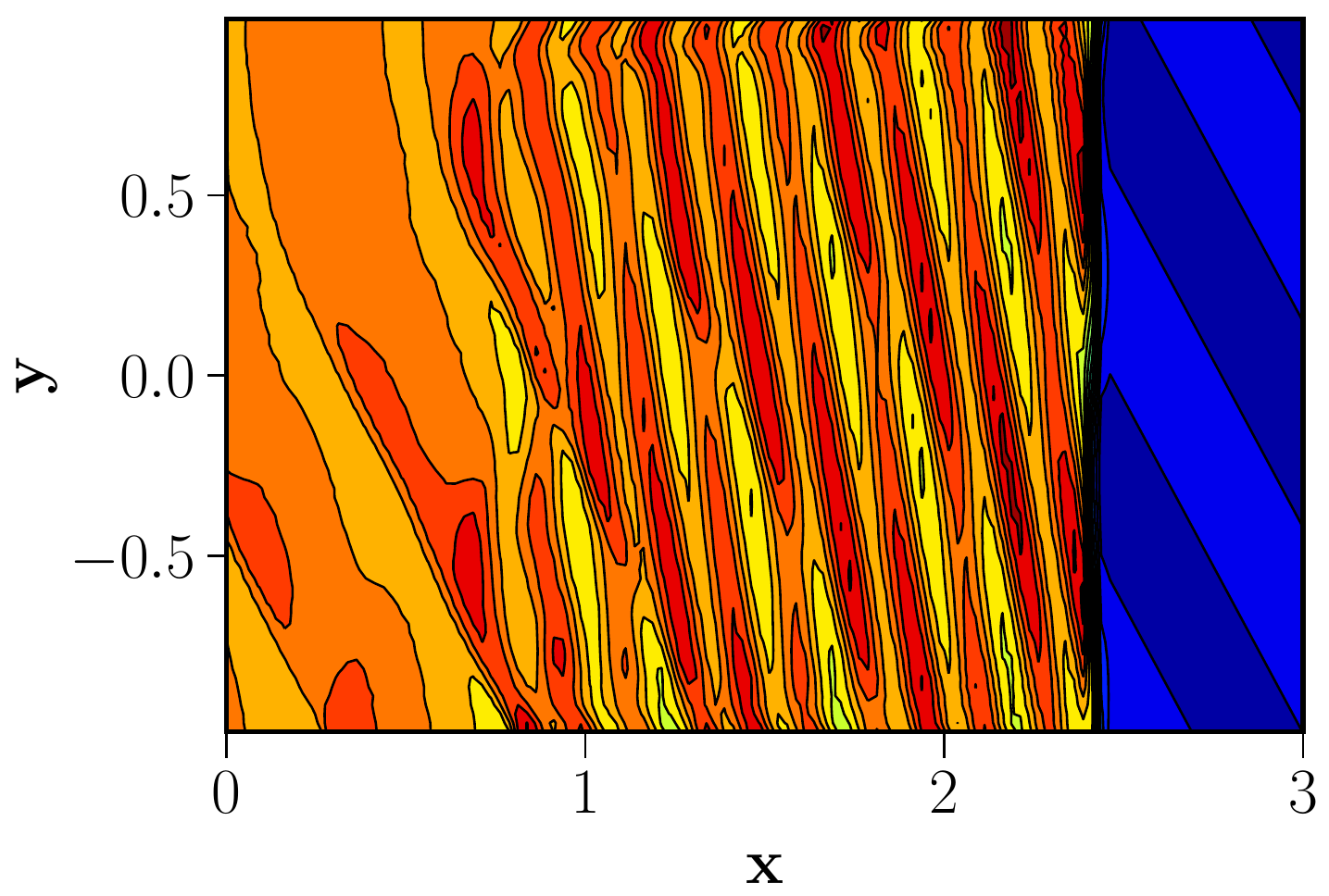}
\label{fig:SE-HOCUS5}}
\subfigure[HOCUS6]{\includegraphics[width=0.48\textwidth]{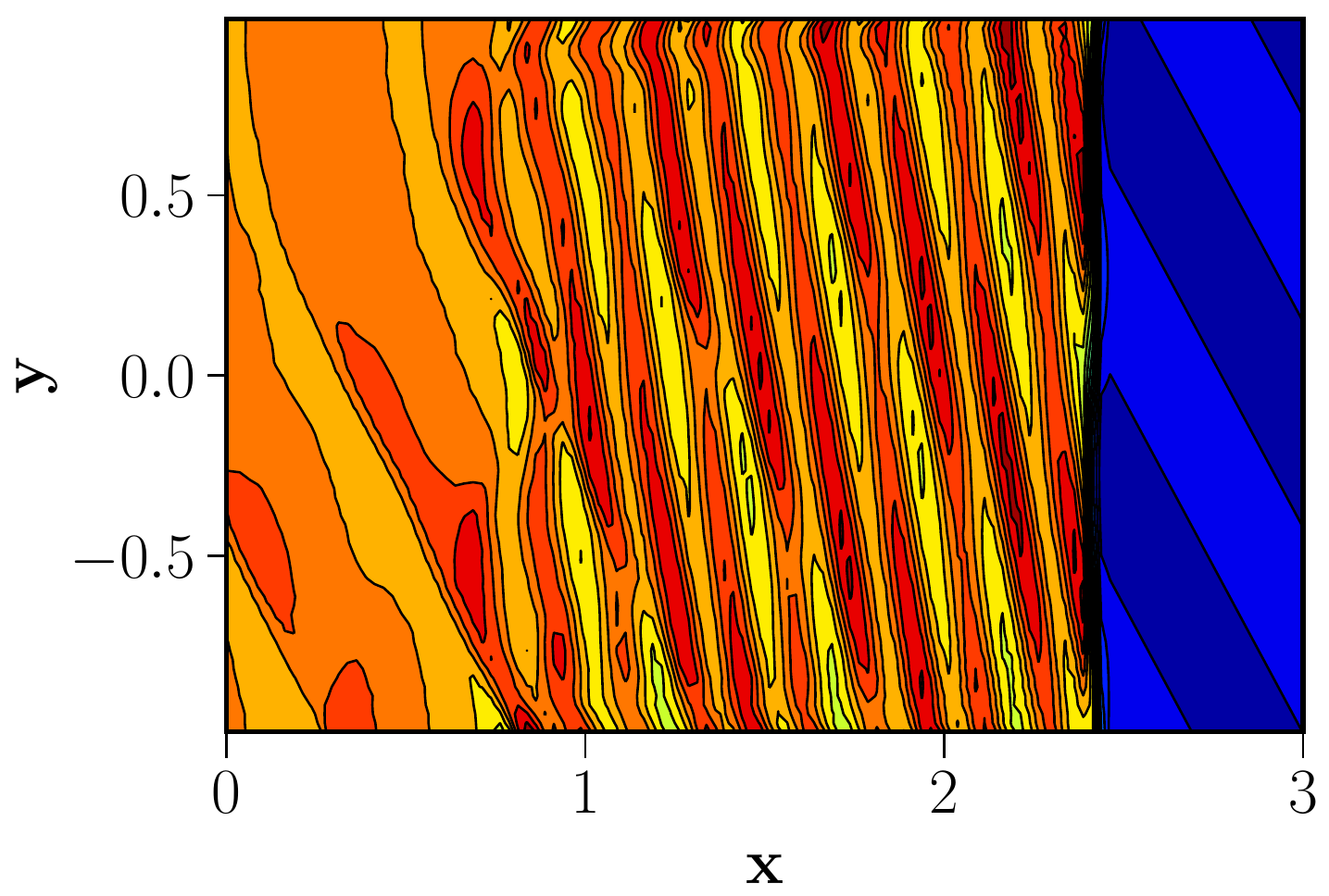}
\label{fig:SE-BVD}}
\subfigure[Local profile]{\includegraphics[width=0.495\textwidth]{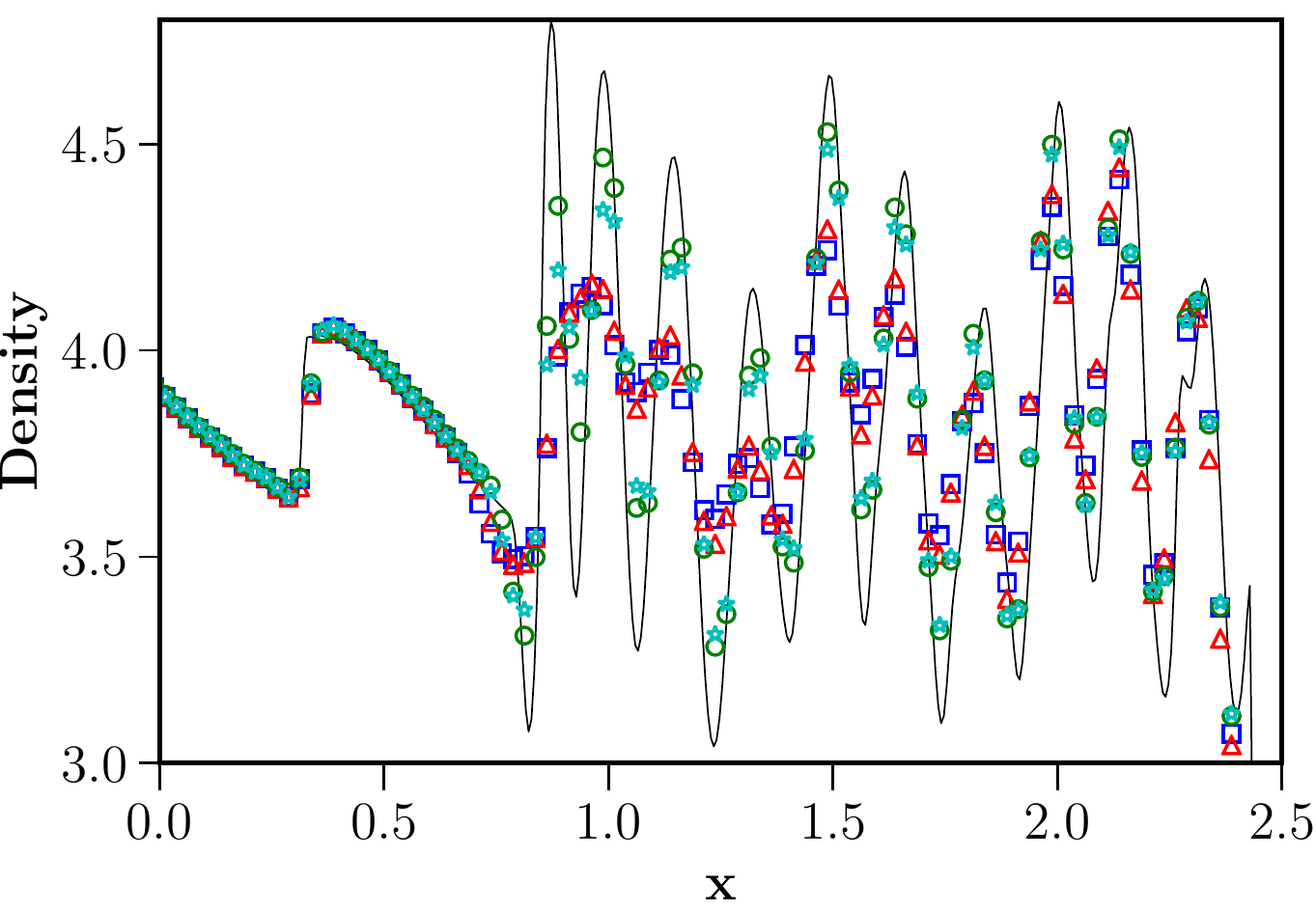}
\label{fig:SSE-Compare}}
\caption{\textcolor{black}{15 density contours for the 2D shock-entropy wave test at $t=1.8$, Example \ref{shock-entropy}, for various schemes are shown in Figs. (a), (b) and (c). Fig. (d) shows the local density in the region with high-frequency waves for all the schemes where solid line: reference solution; green circles: HOCUS6; cyan stars: HOCUS5; blue squares: WENO-Z; red triangles: MP5.}}
\label{fig_shock_entropy}
\end{figure}

%----------------------------------------%----------------------------------------%----------------------------------------%----------------------------------------%----------------------------------------
%----------------------------------------%----------------------------------------%----------------------------------------%----------------------------------------%----------------------------------------
%----------------------------------------%----------------------------------------%----------------------------------------%----------------------------------------%----------------------------------------
%----------------------------------------%----------------------------------------%----------------------------------------%----------------------------------------%----------------------------------------
%----------------------------------------%----------------------------------------%----------------------------------------%----------------------------------------%----------------------------------------
%----------------------------------------%----------------------------------------%----------------------------------------%----------------------------------------%----------------------------------------

\begin{example}\label{ex:rm} {Richtmeyer - Meshkov instability}
\end{example}
In this example, we consider the two-dimensional single-mode \textcolor{black}{Richtmeyer--Meshkov} instability problem \cite{terashima2009front}. This phenomenon occurs when a shock wave approaches a perturbed interface separating two fluids of different densities. To simplify the physical model we consider that the two different fluids have same specific heat ratio of $\gamma$ =1.4 as is considered in \cite{deng2019fifth}. The computational domain has size of $\left[0, 4 \right] \times \left[ 0, 1 \right]$ and the following initial conditions are used:
\begin{equation*}
	\left( \rho, u, v, p \right)
 =
 \begin{cases}
 	\left(1, 0, 0, 1 \right), &\mbox{if x $<$ 3.2}, \\
  \left(1.4112, -665/1556, 0, 1.628 \right), &\mbox{if x $<$ 2.9 sin(2$\pi$(y+0.25), perturbed interface}, \\
  \left(5.04, 0, 0, 1 \right), &\mbox{for $\mathrm{otherwise}$}. \\
 \end{cases}
\end{equation*}
This initial condition indicates that a shock wave is propagating from the left to the interface imposed minimal disturbance. The computation is conducted until $t=9$, on the uniform mesh size of $320 \times 80$. The upper and lower boundaries were treated as periodic boundary conditions through the use of ghost cells and for the left, and right boundary values are fixed to initial conditions. Fig. \ref{fig_rm} shows the density distribution obtained by various schemes, and we can notice that the HOCUS6 scheme has lower numerical dissipation and produces small scale roll-up vortices due to the instability compared to MP5 and WENO-Z. \textcolor{black}{The difference between HOCUSS5 and HOCUS6 is minimal for this test case.}

\begin{figure}[H]
\begin{onehalfspacing}
\centering\offinterlineskip
\subfigure[WENO-Z]{\includegraphics[width=0.6\textwidth]{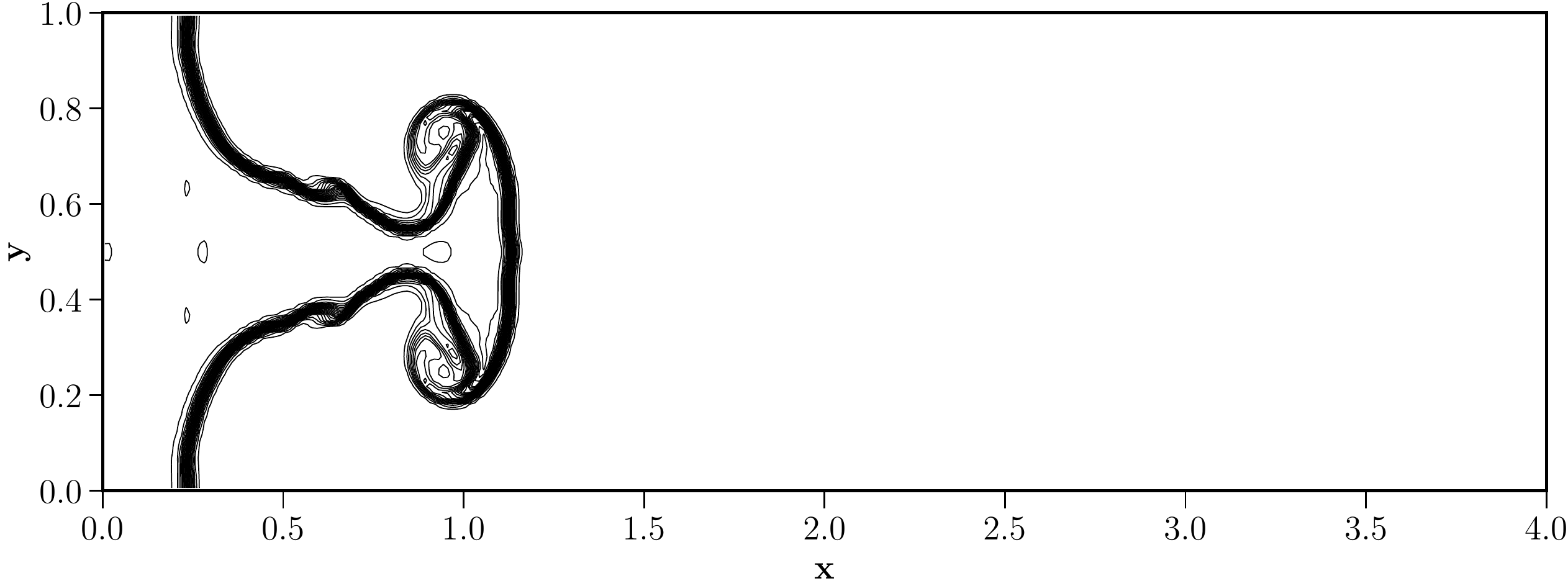}
\label{fig:WENO-Z_RT}}
\subfigure[MP5]{\includegraphics[width=0.6\textwidth]{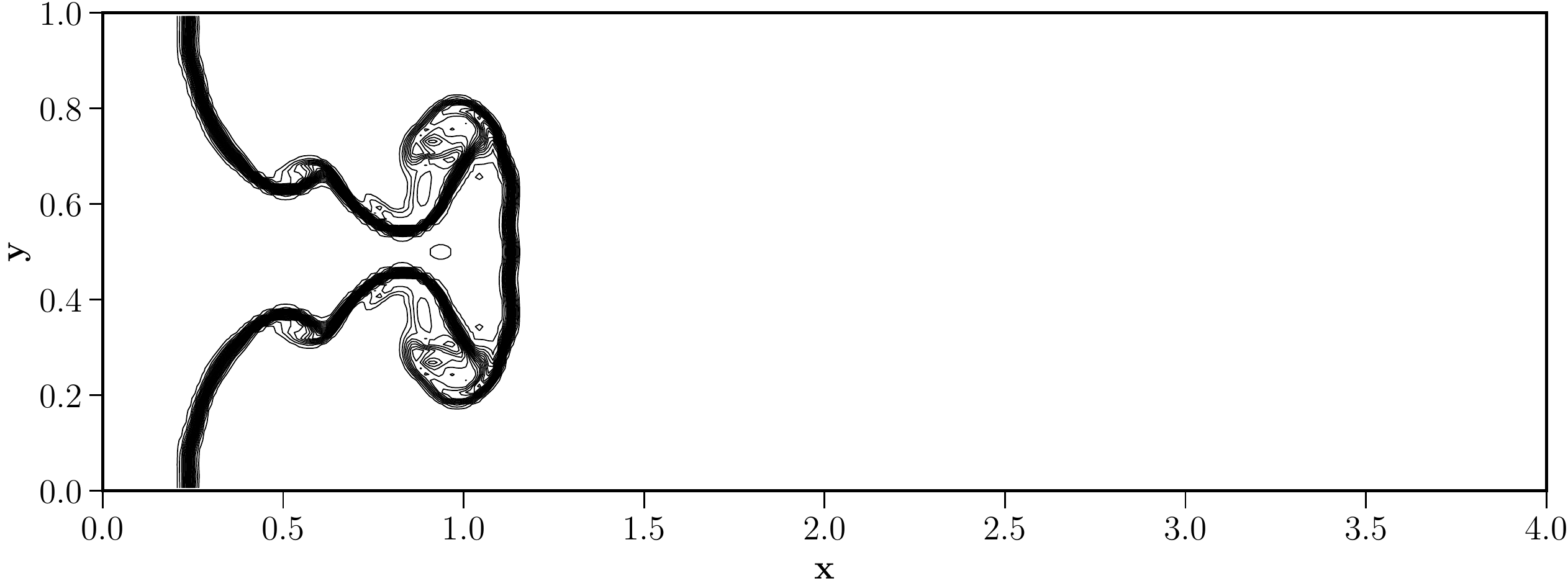}
\label{fig:MP5_RT}}
\subfigure[HOCUS5]{\includegraphics[width=0.6\textwidth]{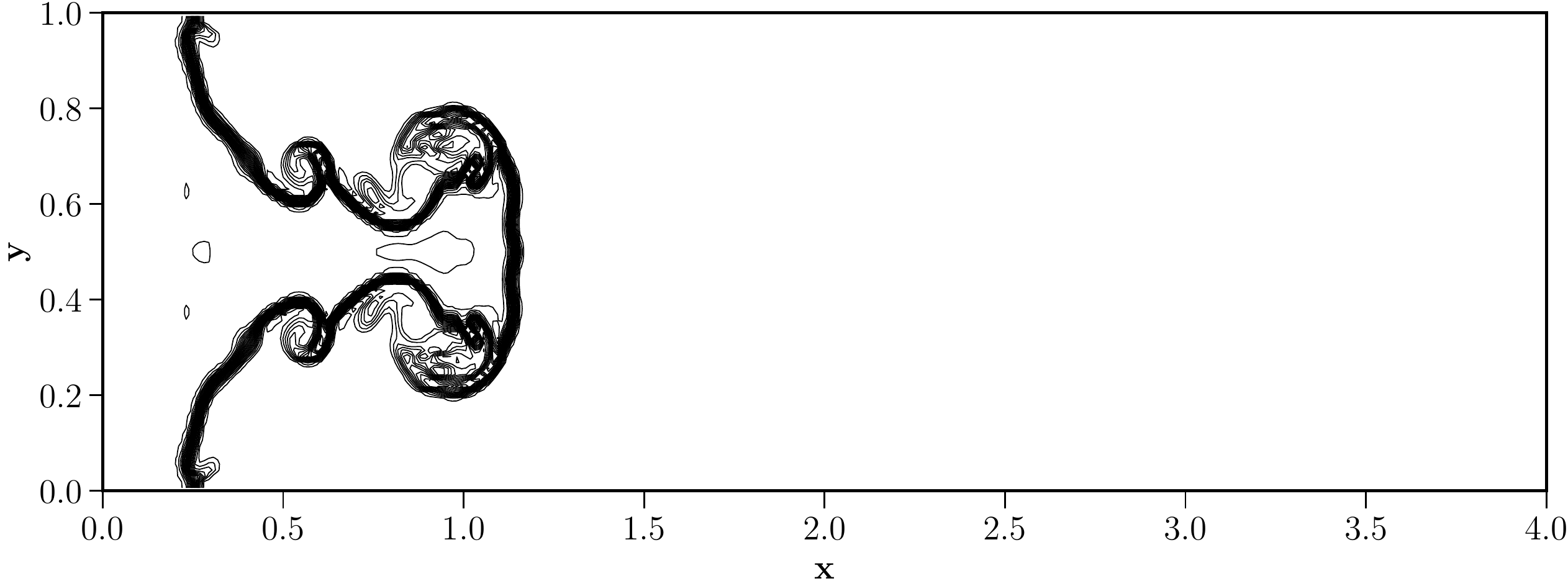}
\label{fig:HOCUS5_RT}}
\subfigure[HOCUS6]{\includegraphics[width=0.6\textwidth]{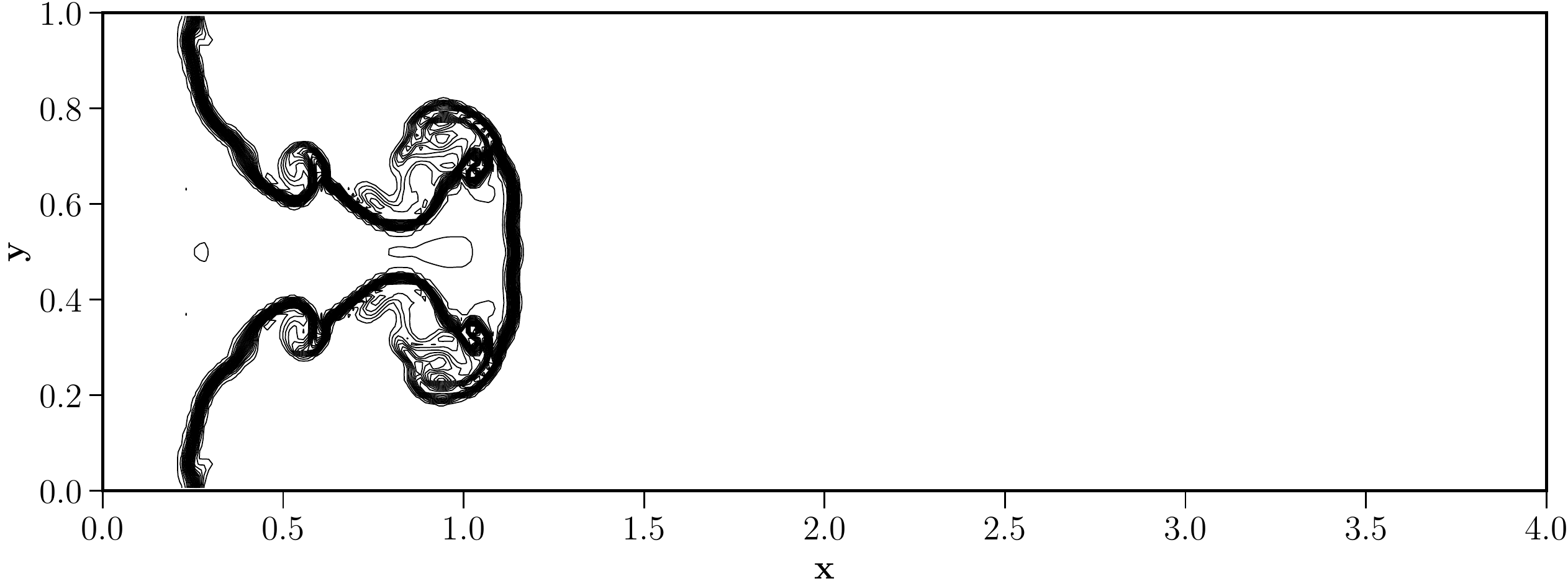}
\label{fig:HOCUS6_RT}}
\caption{\textcolor{black}{24 equally spaced contours of density at t = 9 for various schemes for the Richtmeyer-Meshkov instability test case, Example \ref{ex:rm}}}
\label{fig_rm}
\end{onehalfspacing}
\end{figure}

%----------------------------------------%----------------------------------------%----------------------------------------%----------------------------------------%----------------------------------------
%----------------------------------------%----------------------------------------%----------------------------------------%----------------------------------------%----------------------------------------
%----------------------------------------%----------------------------------------%----------------------------------------%----------------------------------------%----------------------------------------
%----------------------------------------%----------------------------------------%----------------------------------------%----------------------------------------%----------------------------------------

\begin{example}\label{ex:rt} {Rayleigh-Taylor instability}
\end{example}
This type of instability occurs at an interface between fluids of different densities when an acceleration is directed from heavier fluid to lighter fluid. The objective of this test case is to show the feature of dissipation, with bubbles of light fluid rising into the ambient heavy fluid and spikes of heavy fluid falling into the light fluid and resulting in fingering like structure, in two dimensions. The initial conditions of Rayleigh-Taylor instability are \cite{xu2005anti}
 \begin{equation}
\begin{aligned}
(\rho,u,v,p)=
\begin{cases}
&(2.0,\ 0,\ -0.025\sqrt{\frac{5p}{3\rho}\cos(8\pi x)},\ 2y+1.00),\quad 0\leq y< 0.5,\\
&(1.0,\ 0,\ -0.025\sqrt{\frac{5p}{3\rho}\cos(8\pi x)},\ 1y+0.75),\quad 0\leq y\leq 0.5,
\end{cases}
\end{aligned}
\label{eu2D_RT}
\end{equation}
over the computational domain $[0, 1/4]\times [0,1]$. Reflective boundary conditions are imposed on the right and left boundaries via ghost cells. The flow conditions are set to $\rho=1$, $p=2.5$, and $u=v=$0 on top boundary and $\rho=2$, $p=1.0$, and $u=v=0$ on bottom boundary. The source term $S=(0,0,\rho, \rho v)$ is added to the Euler equations. We performed simulations on a uniform mesh of resolutions $80 \times 320$ and \textcolor{black}{$120 \times 480$} and the computations are conducted until $t = 1.95$. The value of adiabatic constant
 $\gamma$ is taken to be $5/3$. Figs. \ref{fig:2d-RT} and \ref{fig:2d-RT1} indicates the density distribution of the Rayleigh-Taylor instability problem and it can be seen that HOCUS6 resolved the finer structures in comparison with WENO-Z, \textcolor{black}{MP5 and HOCUS5}.

\begin{figure}[H]
\begin{onehalfspacing}
\centering\offinterlineskip
\subfigure[WENO-Z]{%
\includegraphics[width=0.20\textwidth]{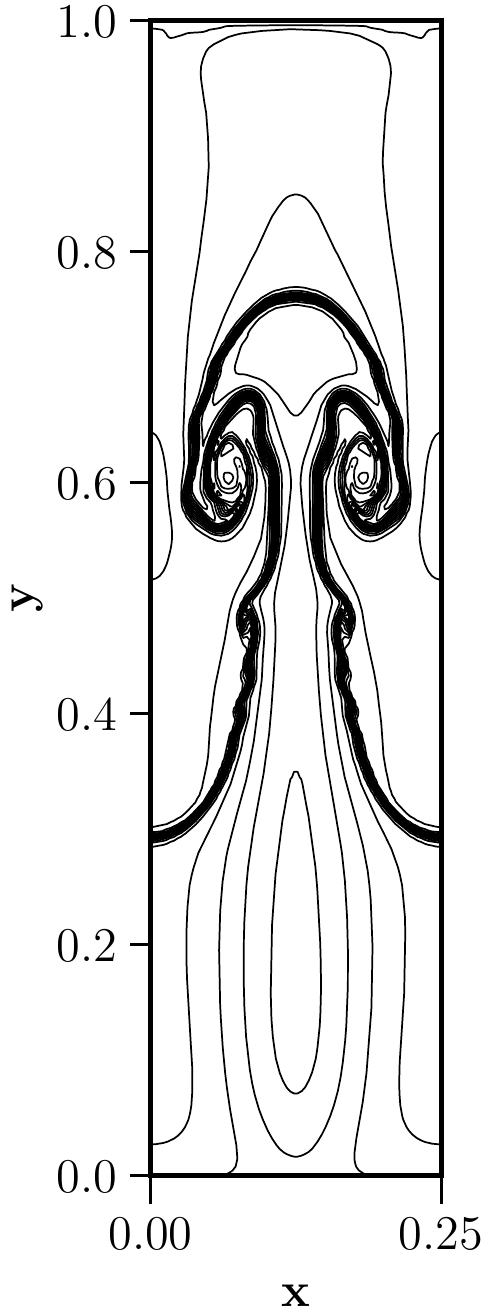}
\label{fig:RT_Z}}
\subfigure[MP5]{%
\includegraphics[width=0.20\textwidth]{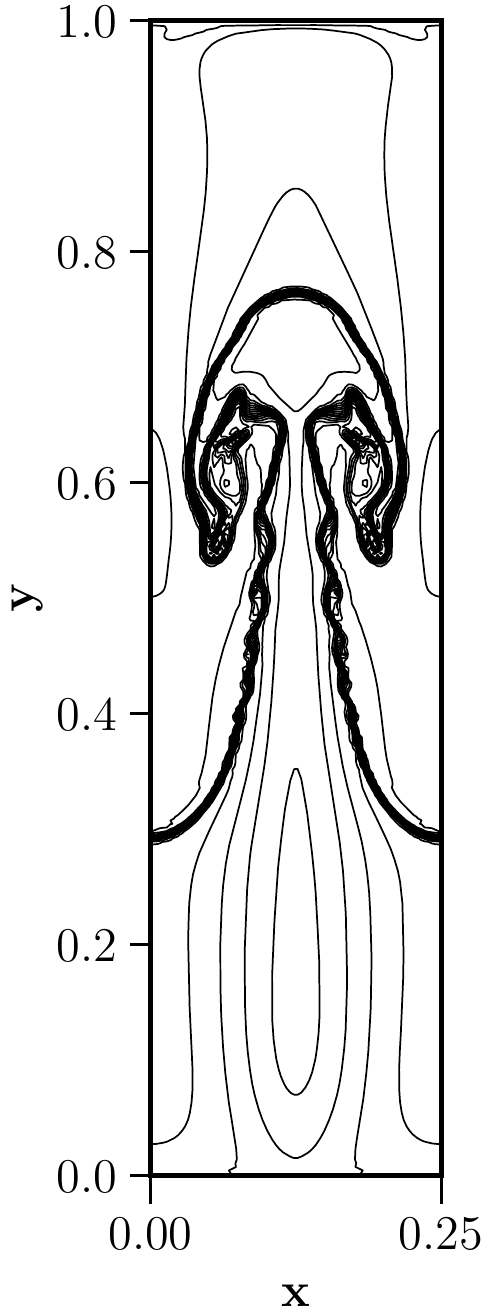}
\label{fig:RT_MP5}}
\subfigure[HOCUS5]{%
\includegraphics[width=0.20\textwidth]{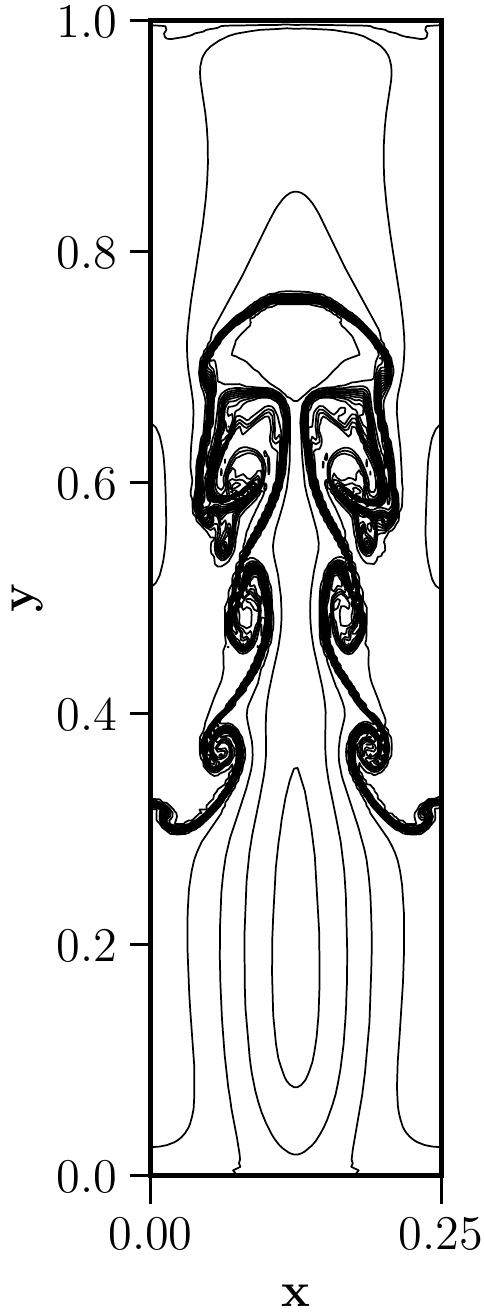}
\label{fig:RT_HOCUS5}}
\subfigure[HOCUS6]{%
\includegraphics[width=0.20\textwidth]{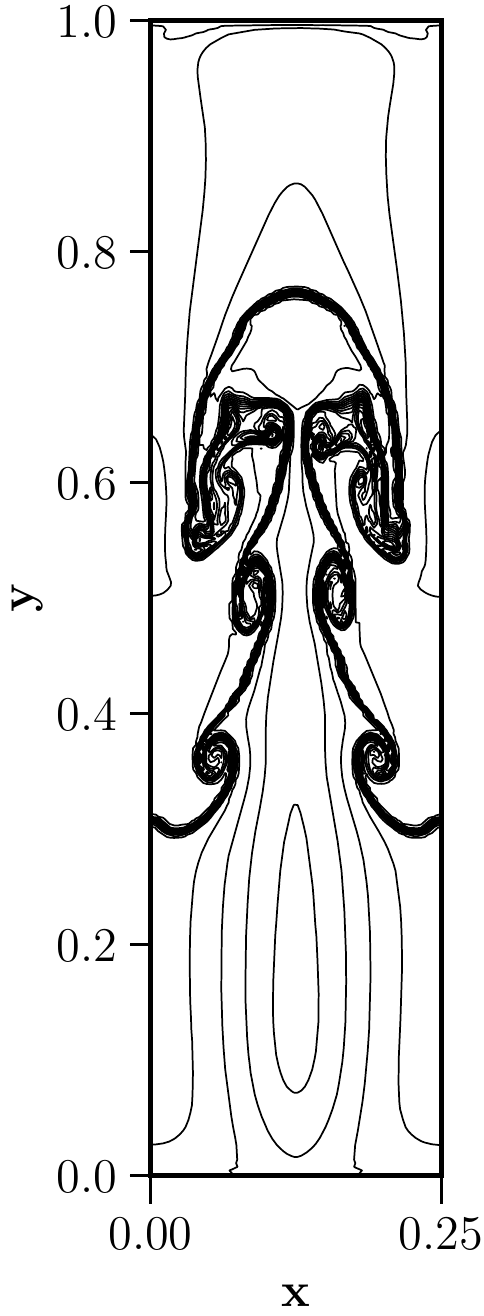}
\label{fig:RT_HOCUS6}}
\caption{\textcolor{black}{Comparison of density contours obtained by WENO-Z, MP5, HOCUS5 and HOCUS6 on a grid size of 80 $\times$ 320 for the test case in Example \ref{ex:rt}.}}
\label{fig:2d-RT}
\end{onehalfspacing}
\end{figure}

\begin{figure}[H]
\begin{onehalfspacing}
\centering\offinterlineskip
\subfigure[WENO-Z]{%
\includegraphics[width=0.20\textwidth]{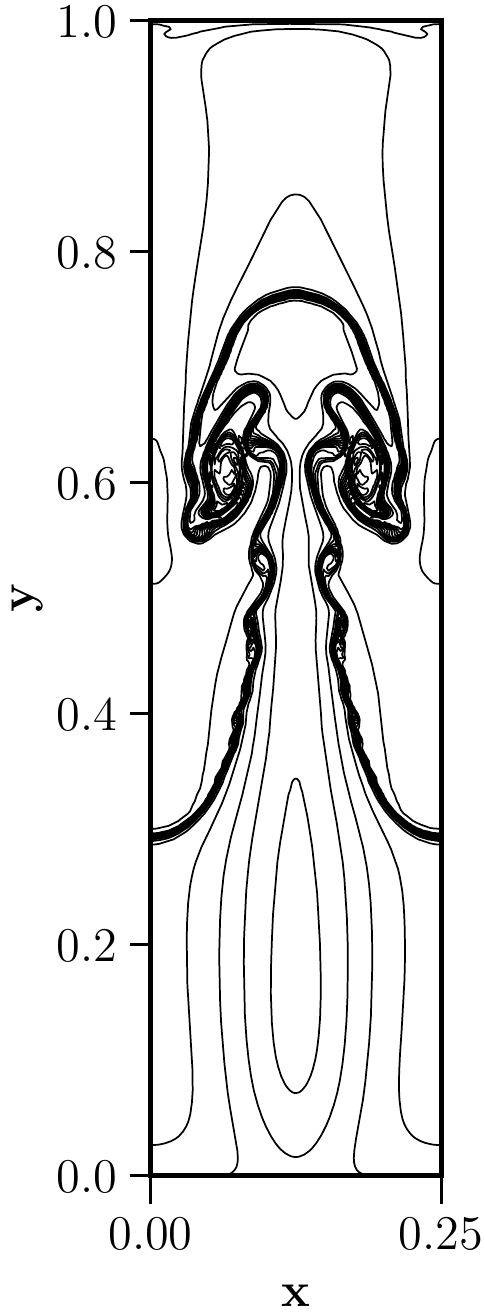}
\label{fig:RT_Z1}}
\subfigure[MP5]{%
\includegraphics[width=0.20\textwidth]{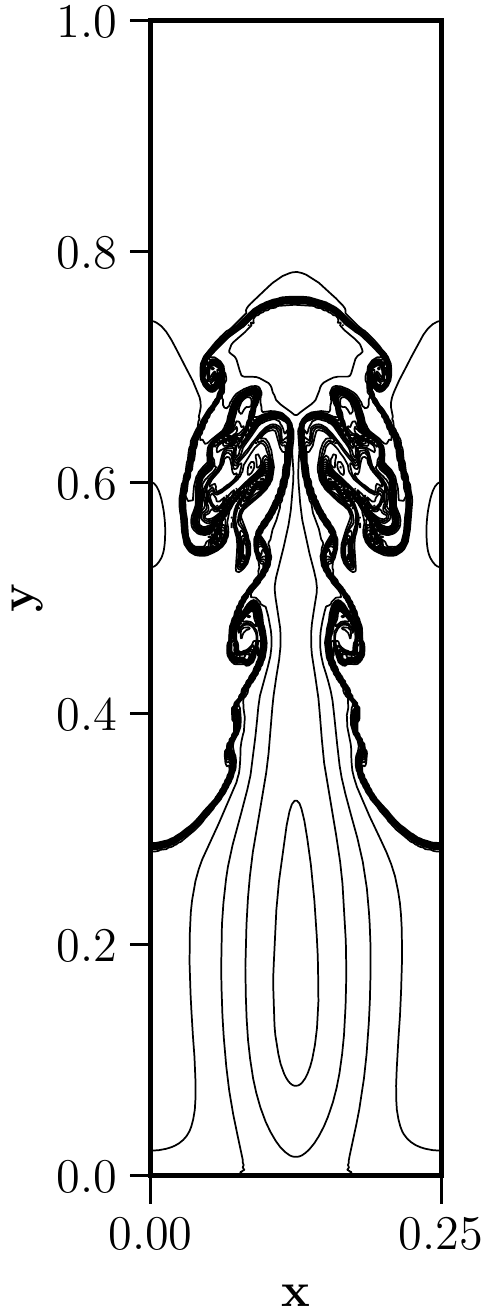}
\label{fig:RT_MP51}}
\subfigure[HOCUS5]{%
\includegraphics[width=0.20\textwidth]{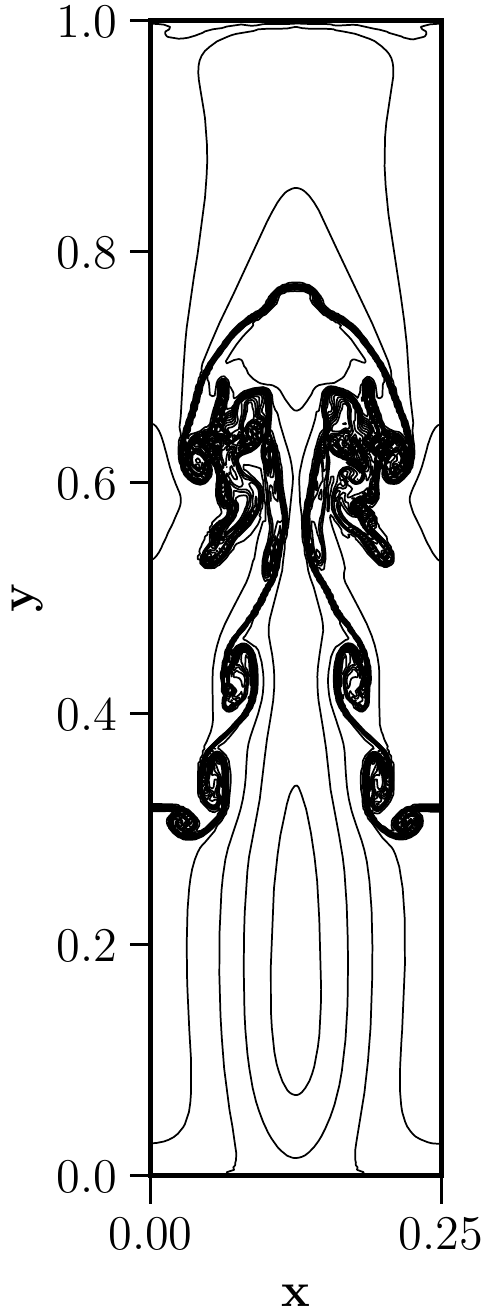}
\label{fig:RT_HOCUS51}}
\subfigure[HOCUS6]{%
\includegraphics[width=0.20\textwidth]{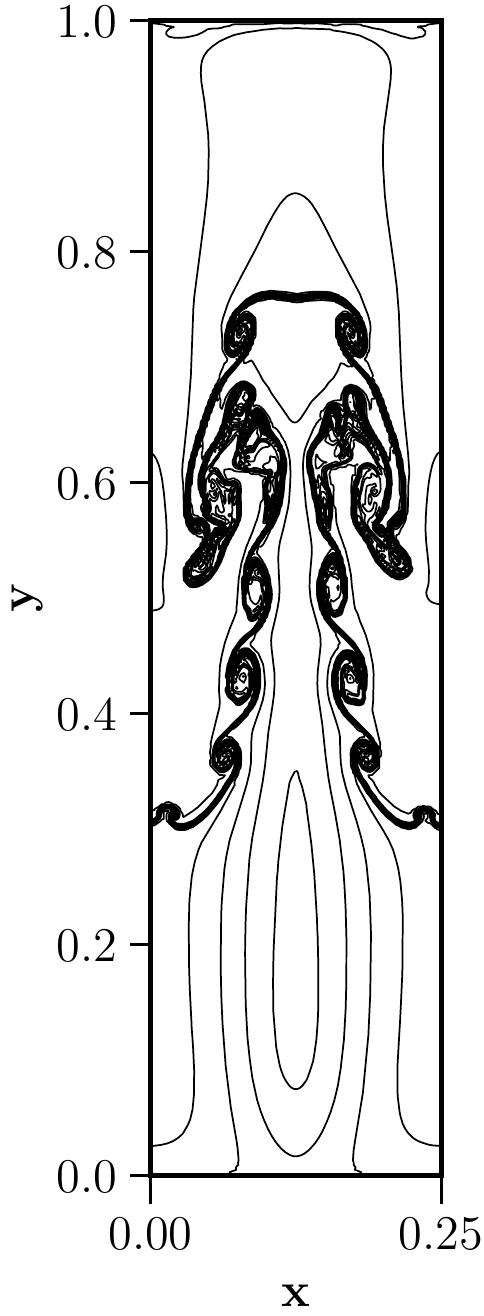}
\label{fig:RT_HOCUS61}}
\caption{\textcolor{black}{Comparison of density contours obtained by WENO-Z, MP5, HOCUS5 and HOCUS6 on a grid size of $120 \times 480$ for the test case in Example \ref{ex:rt}.}}
\label{fig:2d-RT1}
\end{onehalfspacing}
\end{figure}

%----------------------------------------%----------------------------------------%----------------------------------------%----------------------------------------%----------------------------------------
%----------------------------------------%----------------------------------------%----------------------------------------%----------------------------------------%----------------------------------------
%----------------------------------------%----------------------------------------%----------------------------------------%----------------------------------------%----------------------------------------
%----------------------------------------%----------------------------------------%----------------------------------------%----------------------------------------%----------------------------------------
%----------------------------------------%----------------------------------------%----------------------------------------%----------------------------------------%----------------------------------------
%----------------------------------------%----------------------------------------%----------------------------------------%----------------------------------------%----------------------------------------

\begin{example}\label{ex:sb} {Shock-Bubble interaction}
\end{example}
In this test case, we consider shock-bubble interaction, where a Mach 1.22 shock wave impacts a helium bubble \cite{shankar2010numerical}. For simplicity, both Air (Shock wave) and Helium are treated as an ideal gas. The helium bubble is placed at $x = 3.5$, $y = 0.89$ within in a domain of size $[0, 6.5] \times [0, 1.78]$. The initial radius of the bubble is taken as $0.5$. The shock front is initially placed at $x = 4.5$. The initial condition of the test case is shown in Fig. \ref{fig_initibubble}. Inflow and outflow conditions are applied at the right and left boundary, respectively. Slip-wall boundary conditions are set at the remaining boundaries. The initial conditions for the problem are as follows, which are computed using exact Rankine-Hugoniot jump conditions,

 \begin{equation}
\begin{aligned}
(\rho,u,v,p)=
\begin{cases}
&(1.0,\ 0,\ 0,\ 1),~~~~~~~~~~~~~~~~~~~~~~ \rm{pre-shocked \ air},\\
&(1.3764,\ -0.3947,\ 0,\ 1.5698) ~~~ \rm{post-shocked \ air},\\
&(0.1819,\ 0,\ 0,\ 1)~~~~~~~~~~~~~~~~~~~ \rm{helium \ bubble}.
\end{cases}
\end{aligned}
\label{sb}
\end{equation}
\textcolor{black}{Simulations are carried out on a mesh resolution of $400 \times 400$, which corresponds to $\Delta x= \Delta y = 0.0025$ until a time of $t=3.25$. Once again, we can see that HOCUS6 resolves more small scale vortices than HOCUS5, WENO-Z and MP5 in Fig.\ref{fig_shockbubble}.}
%%---------------------------------------------------------------

\begin{figure}[H]
\centering\offinterlineskip
\includegraphics[width=0.7\textwidth]{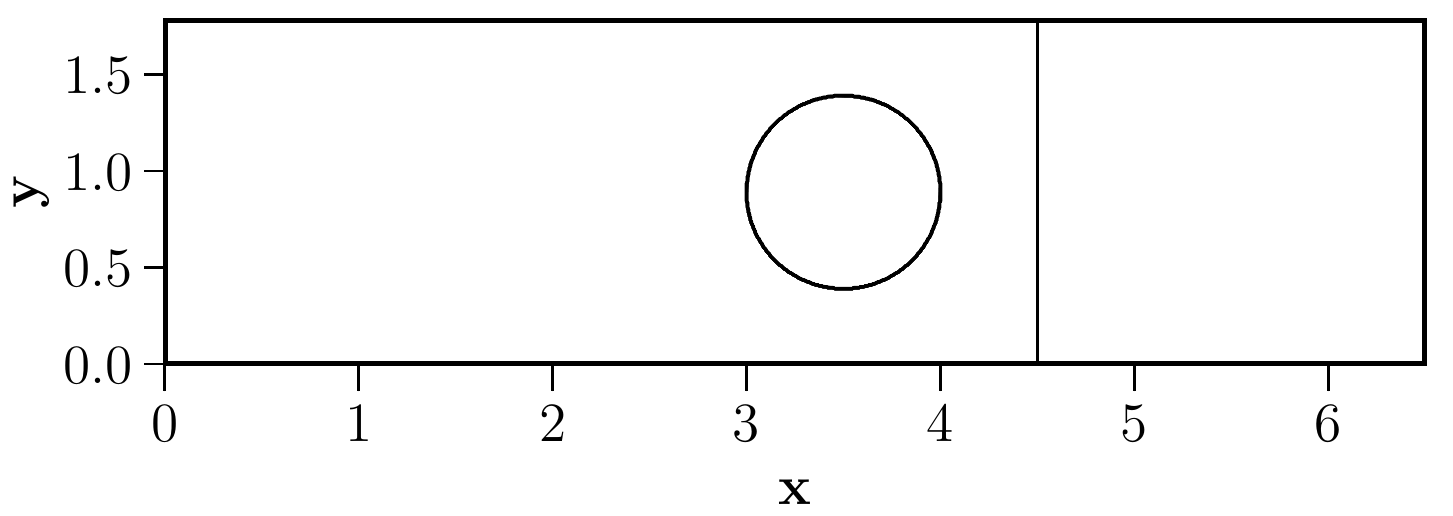}
\caption{Initial condition for the shock-bubble interaction}
\label{fig_initibubble}
\end{figure}

\begin{figure}[H]
\begin{onehalfspacing}
\centering\offinterlineskip
\subfigure[WENO-Z]{\includegraphics[width=0.35\textwidth]{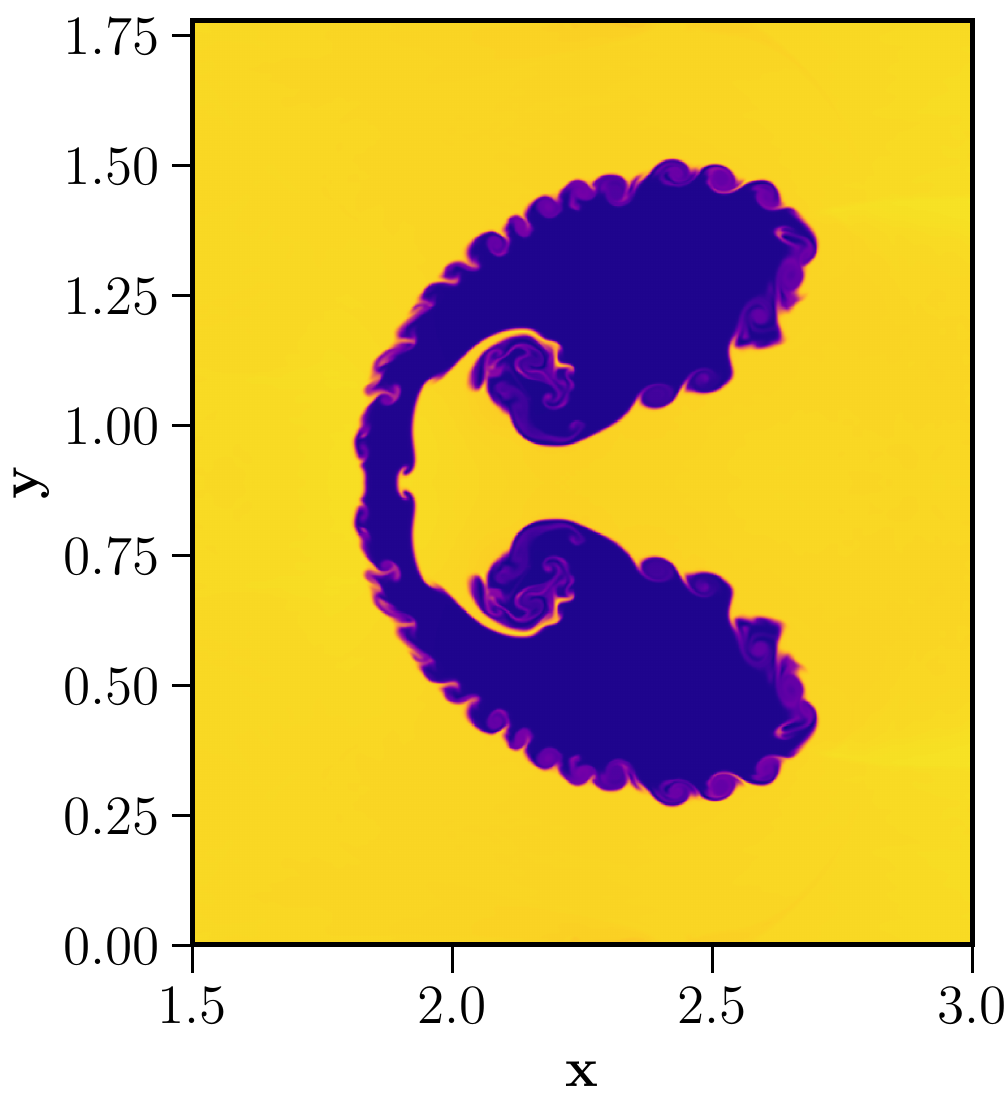}
\label{fig:WENO-Z_BB1}}
\subfigure[MP5]{\includegraphics[width=0.35\textwidth]{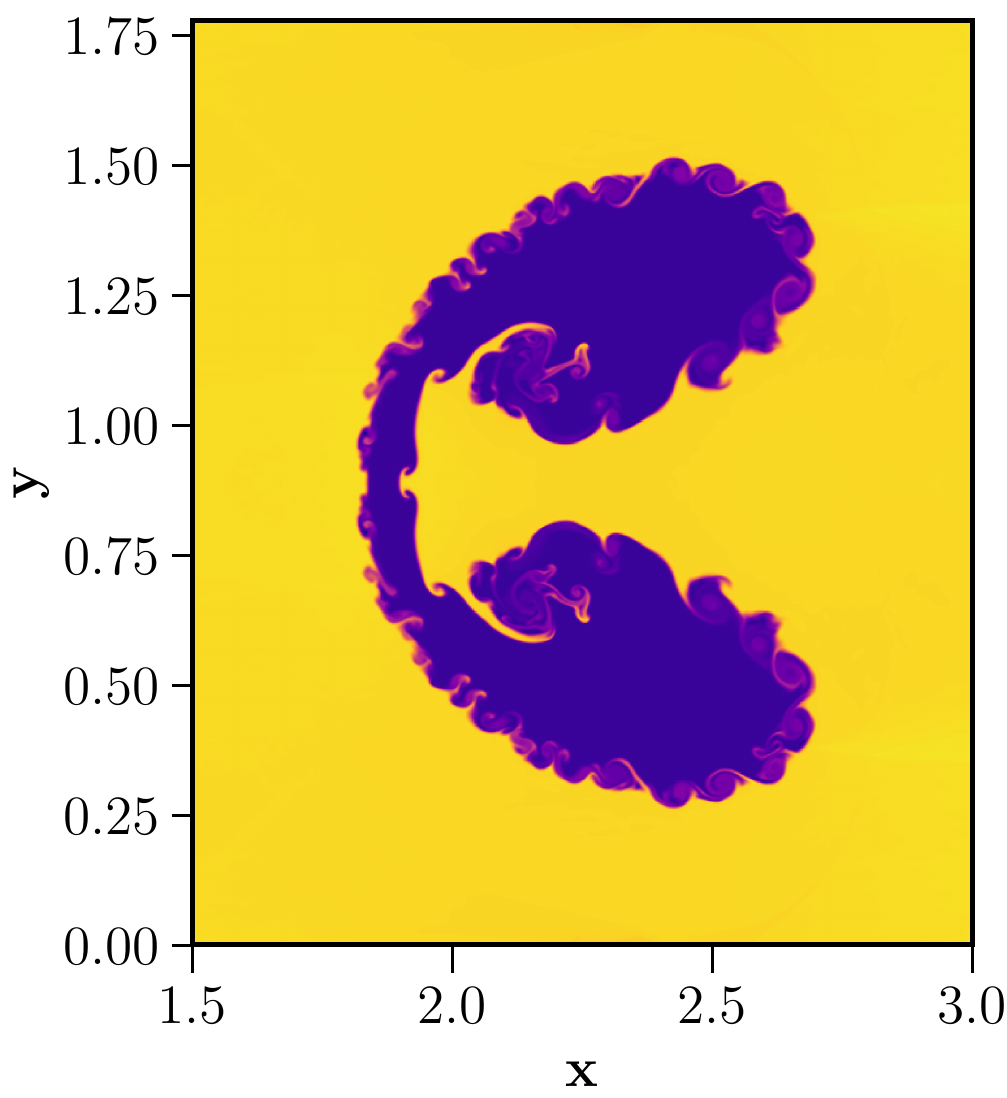}
\label{fig:MP5_BB1}}
\subfigure[HOCUS5]{\includegraphics[width=0.35\textwidth]{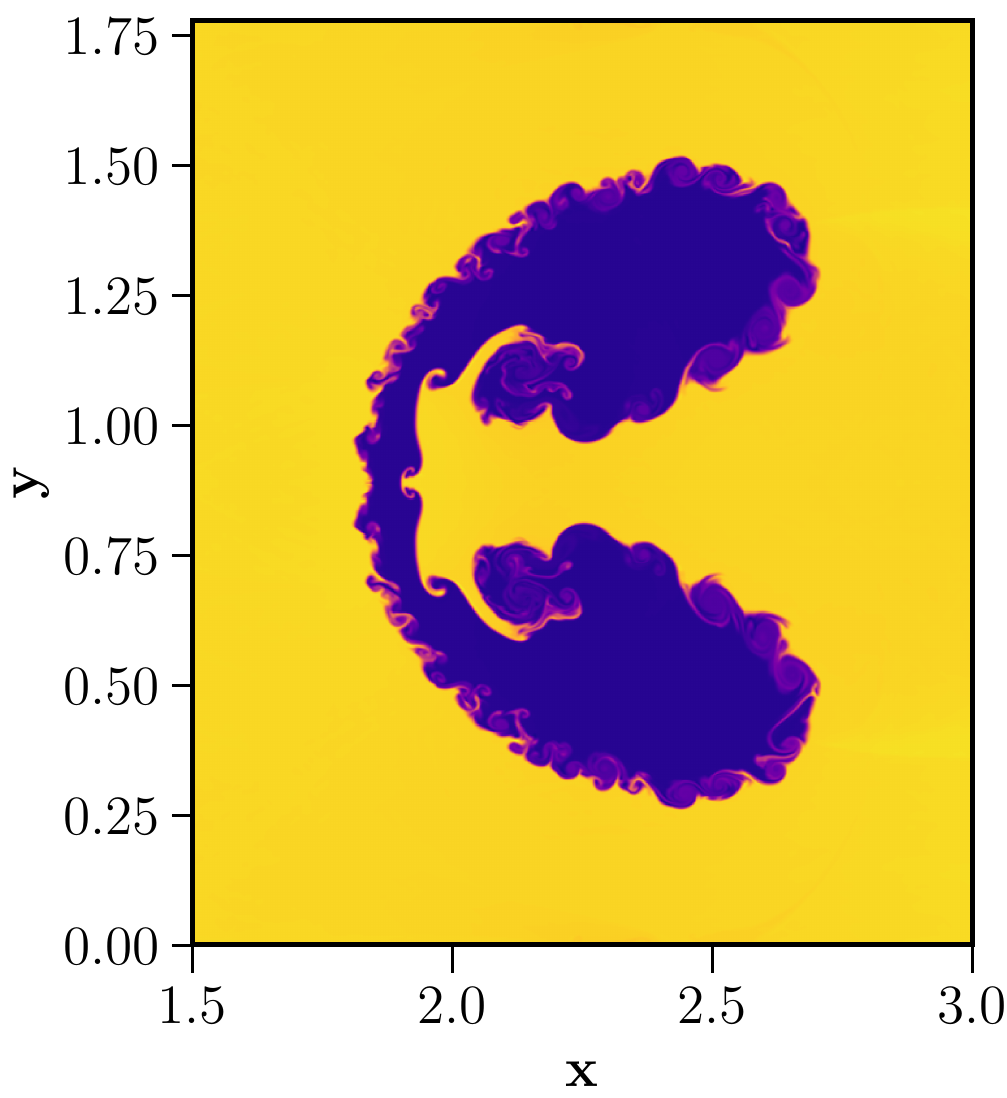}
\label{fig:HOCUS5_BB1}}
\subfigure[HOCUS6]{\includegraphics[width=0.35\textwidth]{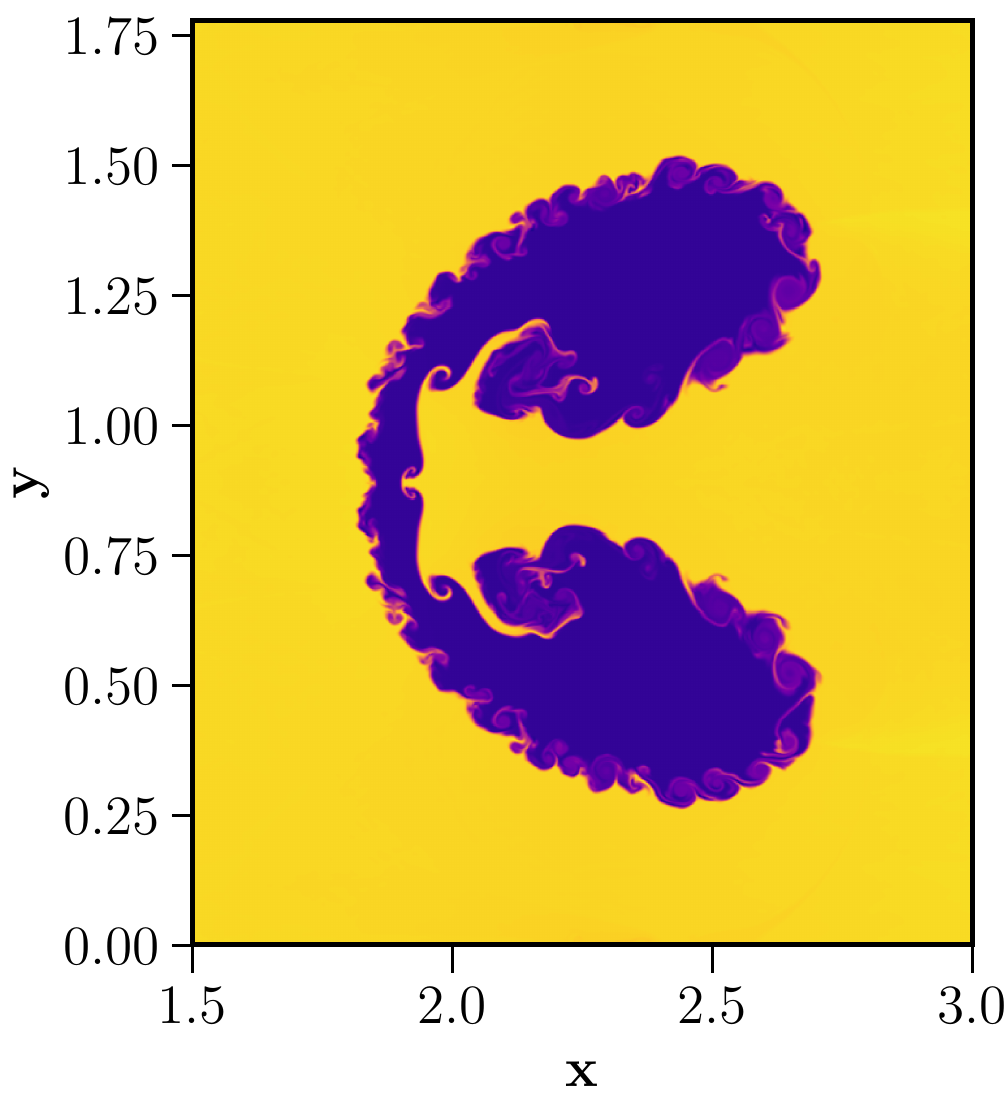}
\label{fig:HOCUS6_BB1}}
\caption{\textcolor{black}{Density plots for various schemes, WENO-5Z, MP5, HOCUS5 and HOCUS6, at $t$ =3.25 for the Example \ref{ex:sb} on a grid spacing of 0.0025.}}
\label{fig_shockbubble}
\end{onehalfspacing}
\end{figure}

%----------------------------------------%----------------------------------------%----------------------------------------%----------------------------------------%----------------------------------------
%----------------------------------------%----------------------------------------%----------------------------------------%----------------------------------------%----------------------------------------
%----------------------------------------%----------------------------------------%----------------------------------------%----------------------------------------%----------------------------------------
%----------------------------------------%----------------------------------------%----------------------------------------%----------------------------------------%----------------------------------------
%----------------------------------------%----------------------------------------%----------------------------------------%----------------------------------------%----------------------------------------
%----------------------------------------%----------------------------------------%----------------------------------------%----------------------------------------%----------------------------------------

\begin{example}\label{ex:dmr}{Double Mach reflection}
\end{example}
In this test case, introduced by Woodward and Collela \cite{woodward1984numerical}, we simulate an unsteady planar shock wave of Mach $10$ impinging on an inclined surface which is at an angle of $30$ degrees. As typical to the setup of this problem, we tilt the incident shock wave rather than the wall to avoid modelling oblique physical wall boundary. The initial conditions for the problem are as follows:
\begin{equation}
\begin{aligned}
(\rho,u,v,p)=
\begin{cases}
&(8,\ 8.25 \cos 30^{0},\ -8.25 \sin 30^{0},\
116.5),\quad x<1/6+\frac{y}{\tan 60^{0}},\\
&(1.4,\ 0,\ 0,\
1),\quad\quad\quad\quad\quad\quad\quad\quad\quad\quad\quad\ x>
1/6+\frac{y}{\tan 60^{0}}.
\end{cases}
\end{aligned}
\label{eu2D_mach}
\end{equation}
Post-shock flow conditions are set for the left boundary, and zero-gradient conditions are applied at the right boundary. At the bottom boundary, reflecting boundary conditions are used from $x =1/6 $ to $x=4.0$ and post-shock conditions for $x \in \left[ 0, 1/6 \right]$. The exact solution of the Mach 10 moving oblique shock is imposed at the upper boundary $y = 1$, which is time-dependent. Simulation is carried over the domain $[0,4]\times [0,1]$ and solutions are computed until time $t=0.2$, with a constant CFL number of $0.2$ for MP5, \textcolor{black}{HOCUS5} and HOCUS6 and $0.4$ for WENO-Z scheme respectively, on a grid size of $1024 \times 256$. As there is no physical viscosity in the Euler equations, vortices generated along the slip line due to Kelvin Helmholtz instability and the near-wall jet will indicate the numerical dissipation of a scheme. We can notice from Fig.\ref{fig_doublemach}, the results obtained by the HOCUS6  \textcolor{black}{are better than} WENO-Z and MP5 as it captures more number of vortices along the slip lines and a strong near-wall jet. \textcolor{black}{For this test case the differences between HOCUS5 and HOCUS6 are marginal} Additional simulations are also conducted for this test case, and the observations are as follows:

\begin{itemize}
\item Simulations are also carried out with the dissipative Global Lax-Friedrich (GLF) Riemann solver, shown in Fig.\ref{fig_doublemach_LF}. It can be seen from Figs. \ref{fig_doublemach} and \ref{fig_doublemach_LF} that WENO-Z is dissipative regardless of the Riemann solver and whereas the adaptive scheme, HOCUS6, has the small scale vortices along the slip line which indicates the advantage of the adaptive non-dissipative scheme.
\item Both MP5 and HOCUS6 schemes are found to be preserving positivity of density and pressure for both HLLC and global Lax-Friedrich Riemann solvers for a constant CFL of $0.2$, and stable computations are carried out by WENO-Z by using a CFL of $0.4$ similar to that of Ref. \cite{fu2019low}. Positivity preserving method developed in \cite{hu2013positivity} may be considered to improve robustness.
\item HLLC Riemann solver used in this paper will lead to \textit{carbuncle} phenomenon, grid-aligned shock instability, and it can be avoided by using a recently proposed low-dissipation modified component-wise Lax-Friedrichs flux by Fleischmann et al. \cite{fleischmann2020low}. Simulations are also carried out with the proposed method, results not shown here, which further improved the simulations.
\item Even though it is not related to the current paper but it is observed that the THINC method used in the original BVD algorithm \cite{sun2016boundary} and the subsequent papers \cite{deng2019fifth, deng2018limiter} with values $\beta=1.2$ or $1.3$ can have compressive or anti-diffusion effect in the numerical results. This anti-diffusion effect improves the roll-up vortices along the slip lines for this test case, see Fig. 11 in Deng et al. \cite{deng2018limiter} and the Riemann problem in Example \ref{ex:rp}. Similar results are also observed for WENO schemes coupled with anti-diffusive flux corrections, see Fig.\ 17 in Xu and Shu \cite{xu2005anti}. A similar approach may be considered for the present HOCUS6 scheme in future.
%\item Deng et al. \cite{deng2019fifth} reported smaller numerical dissipation when primitive variables are reconstructed rather than characteristic variables for this test case. Effect of the numerical dissipation on the choice of the variables used for reconstruction will be investigated in future.
\end{itemize}

\begin{figure}[H]
%\begin{halfspacing}
\centering\offinterlineskip
\subfigure[WENO-Z - HLLC]{\includegraphics[width=0.45\textwidth]{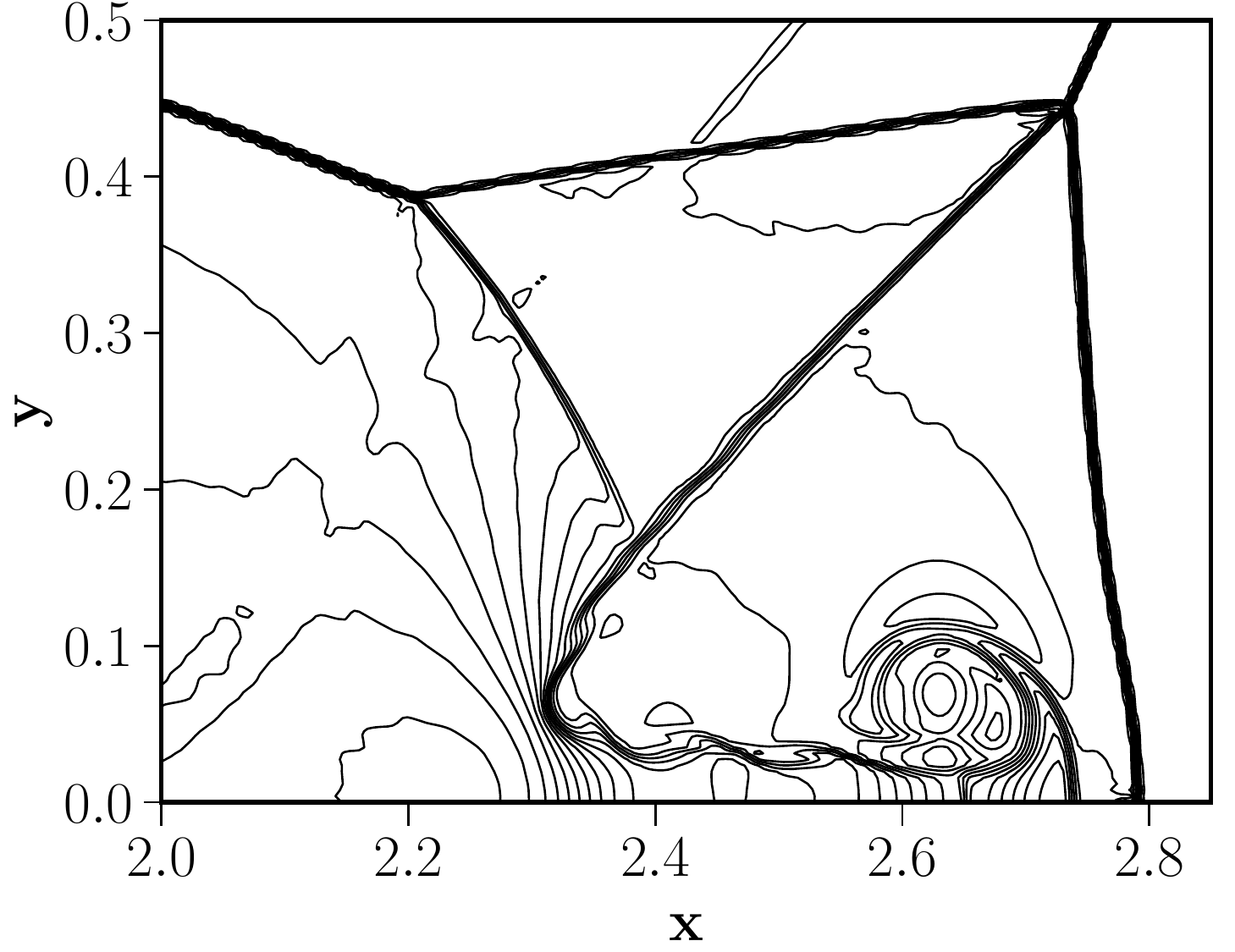}
\label{fig:WENO-Z_DB}}
\subfigure[MP5 - HLLC]{\includegraphics[width=0.45\textwidth]{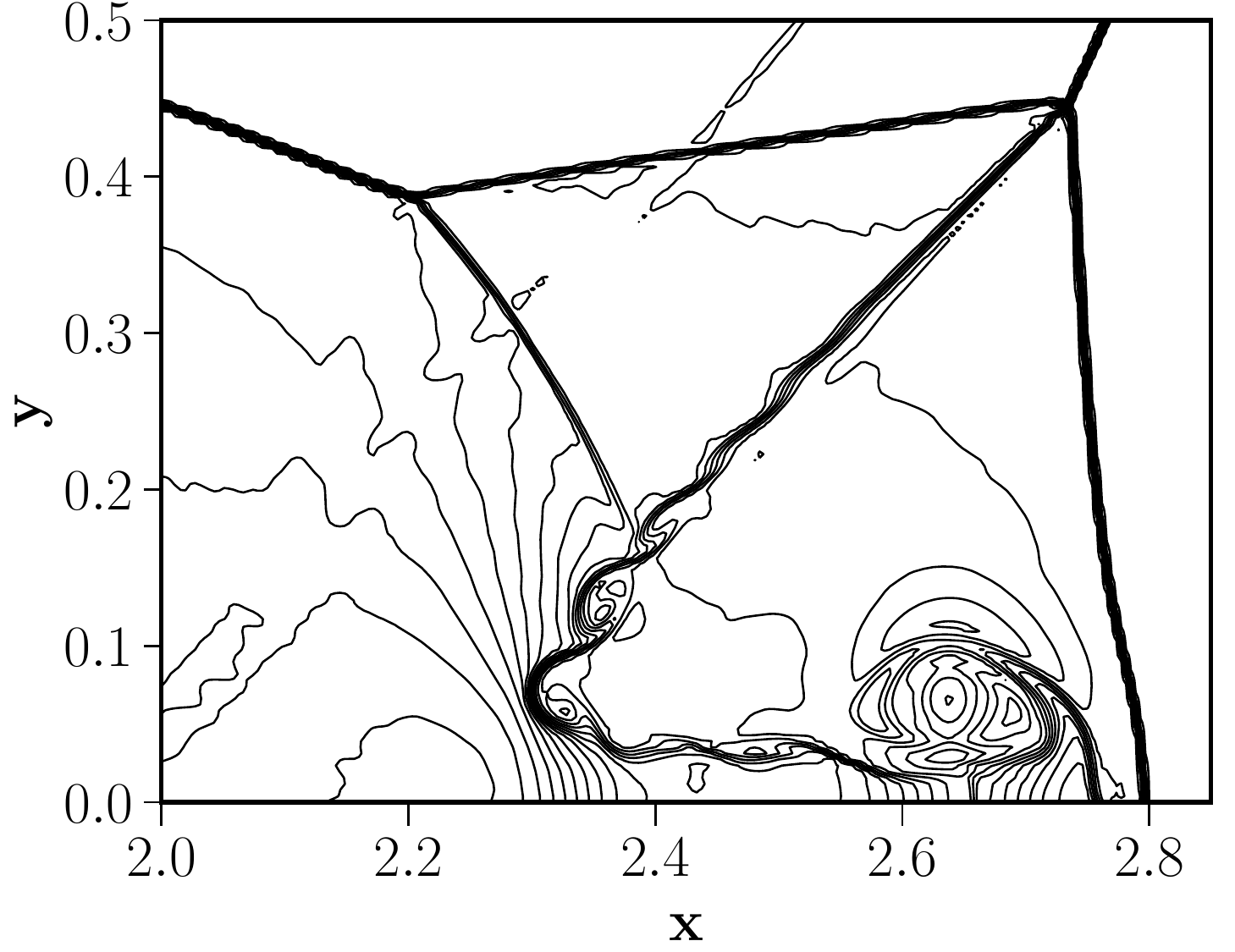}
\label{fig:MP5_DB}}
\subfigure[HOCUS5 - HLLC]{\includegraphics[width=0.45\textwidth]{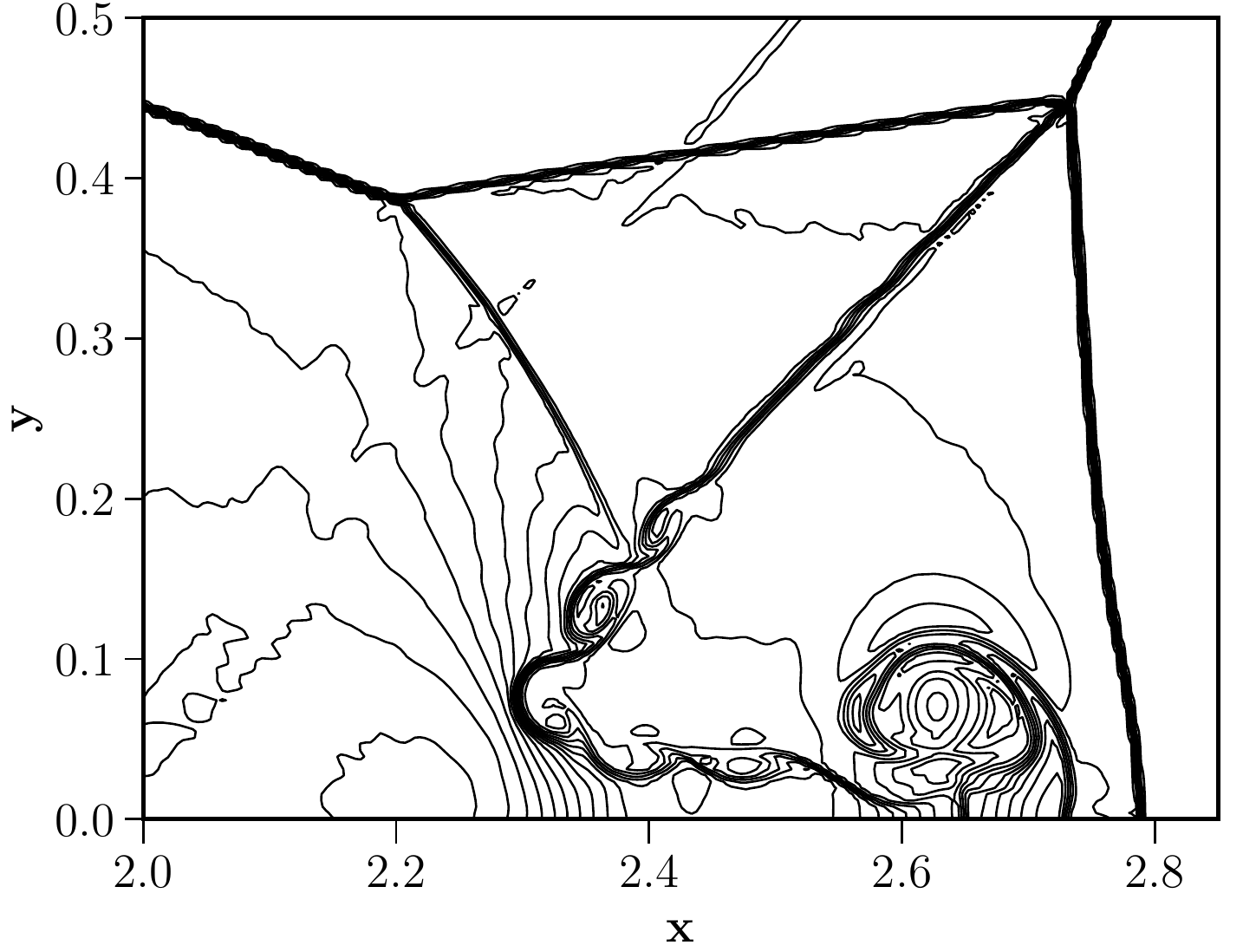}
\label{fig:HOCUS5_DB}}
\subfigure[HOCUS6 - HLLC]{\includegraphics[width=0.45\textwidth]{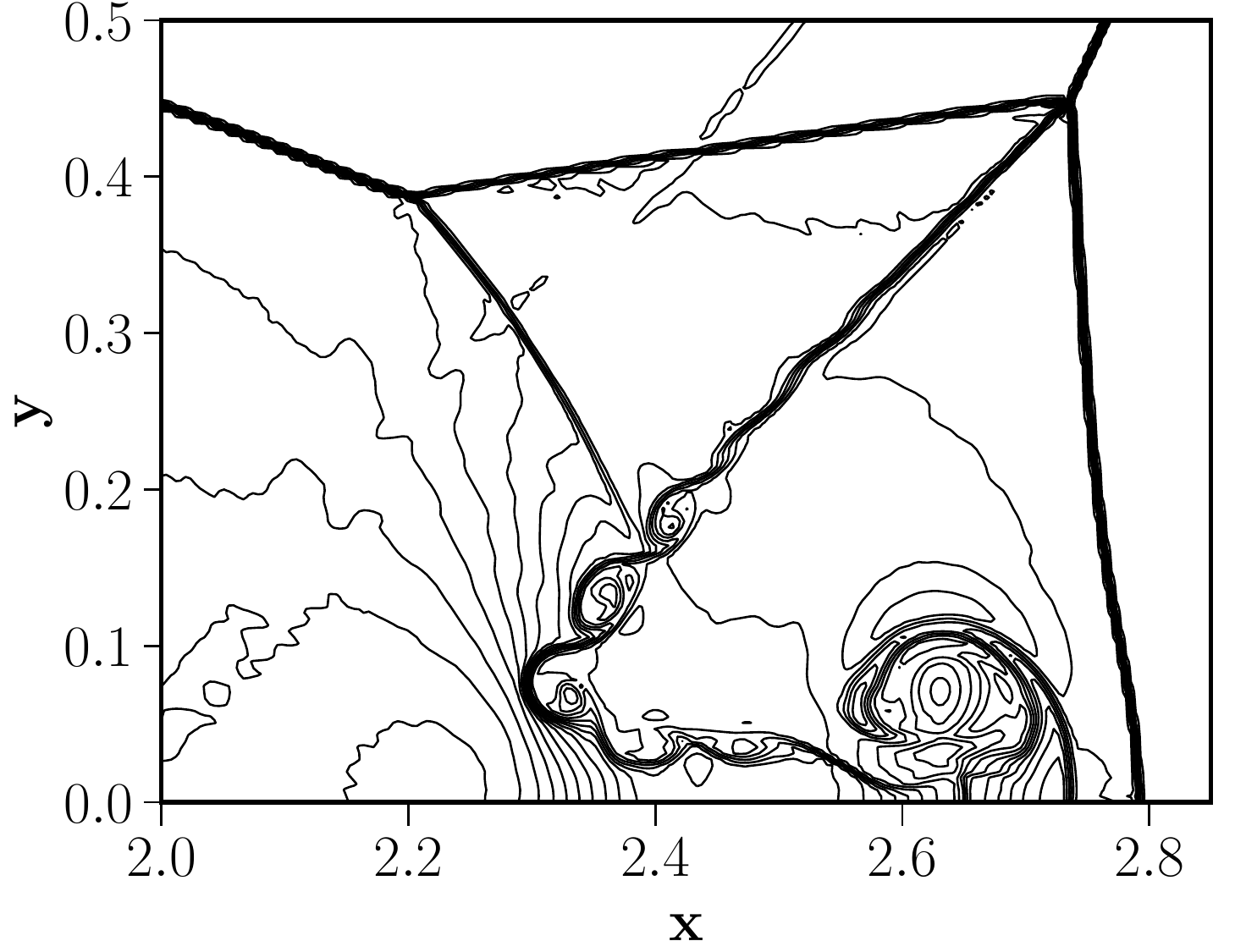}
\label{fig:HOCUS6_DB}}
\caption{\textcolor{black}{Density contours at $t = 0.2$ for the three numerical schemes in the blown-up region around the Mach stem for Example \ref{ex:dmr} using HLLC Riemann solver on a grid size of $1024 \times 256$.}}
\label{fig_doublemach}
%\end{halfspacing}
\end{figure}

\begin{figure}[H]
\centering
\subfigure[WENO-Z - GLF]{\includegraphics[width=0.4\textwidth]{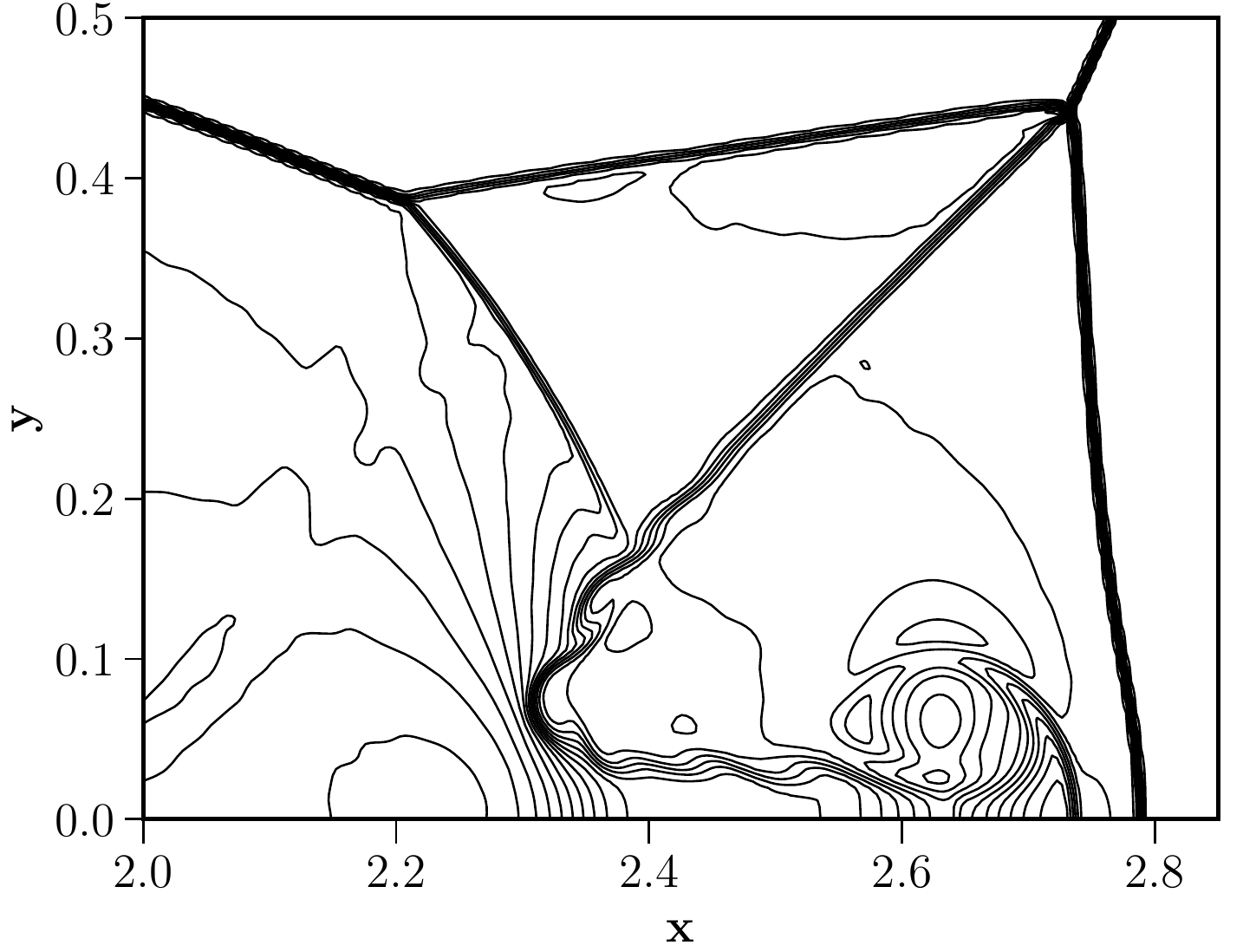}
\label{fig:WENO-Z_DB-LF}}
\subfigure[MP5 - GLF]{\includegraphics[width=0.4\textwidth]{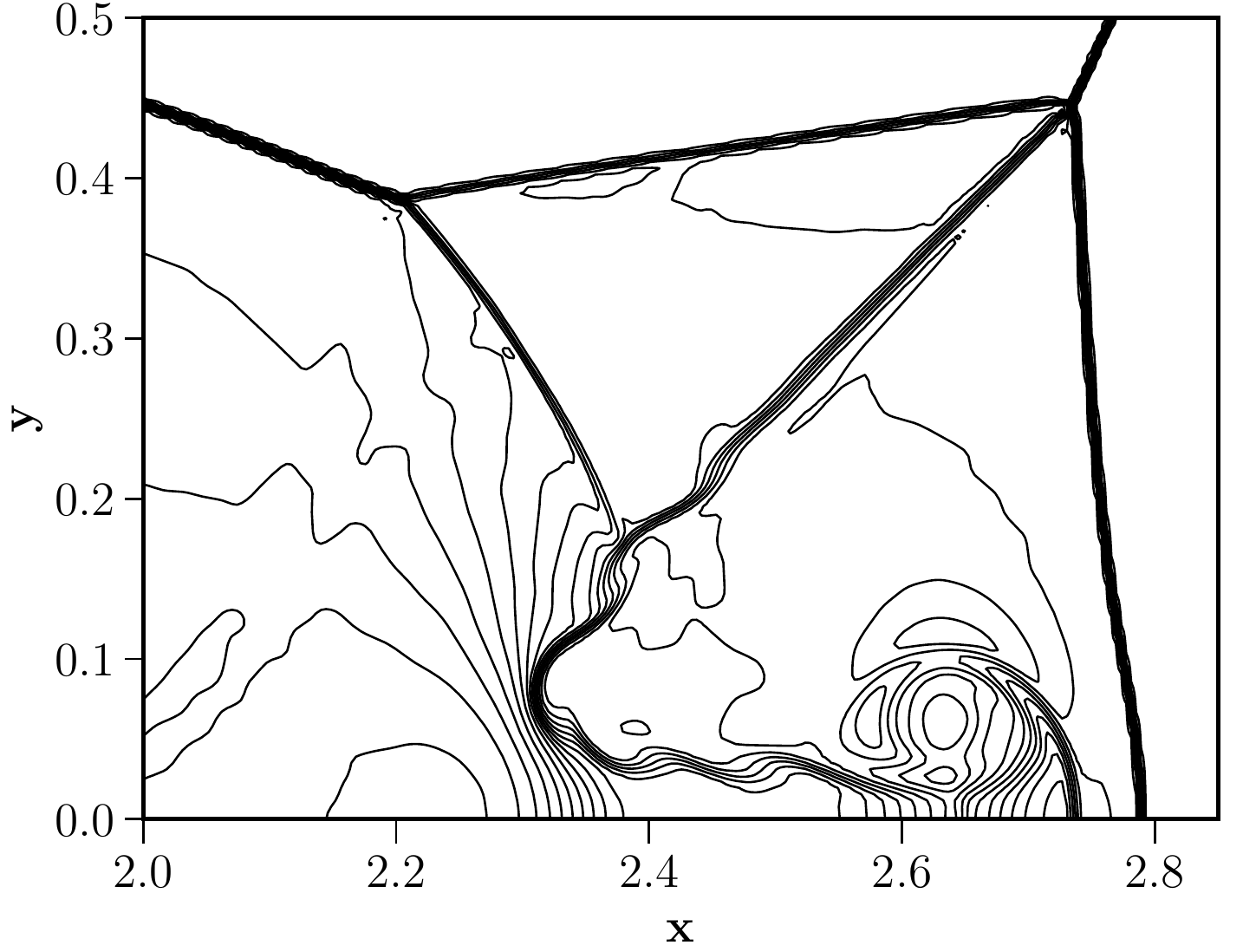}
\label{fig:MP5_DB-LF}}
\subfigure[HOCUS5 - GLF]{\includegraphics[width=0.4\textwidth]{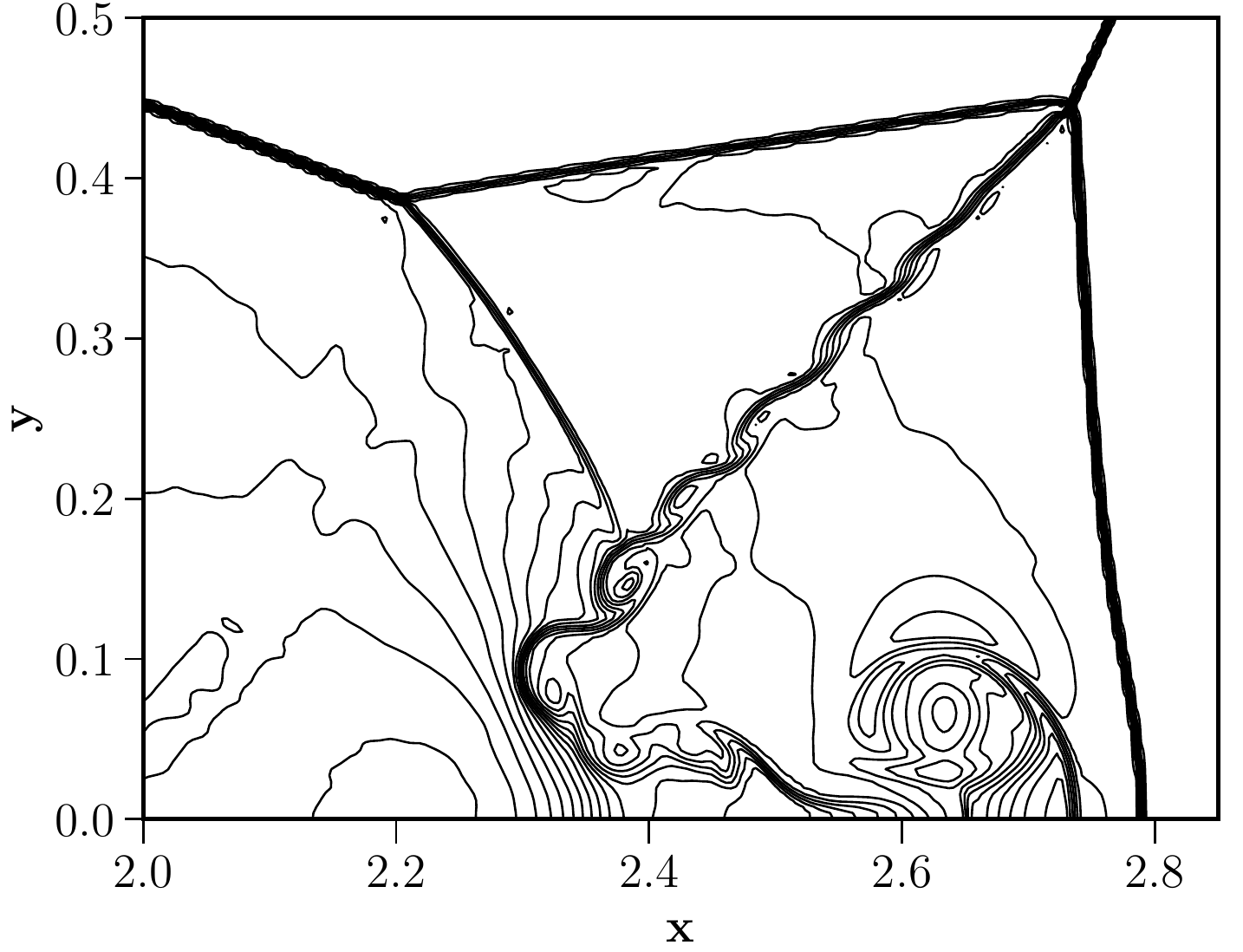}
\label{fig:HOCUS5_DB-LF}}
\subfigure[HOCUS6 - GLF]{\includegraphics[width=0.4\textwidth]{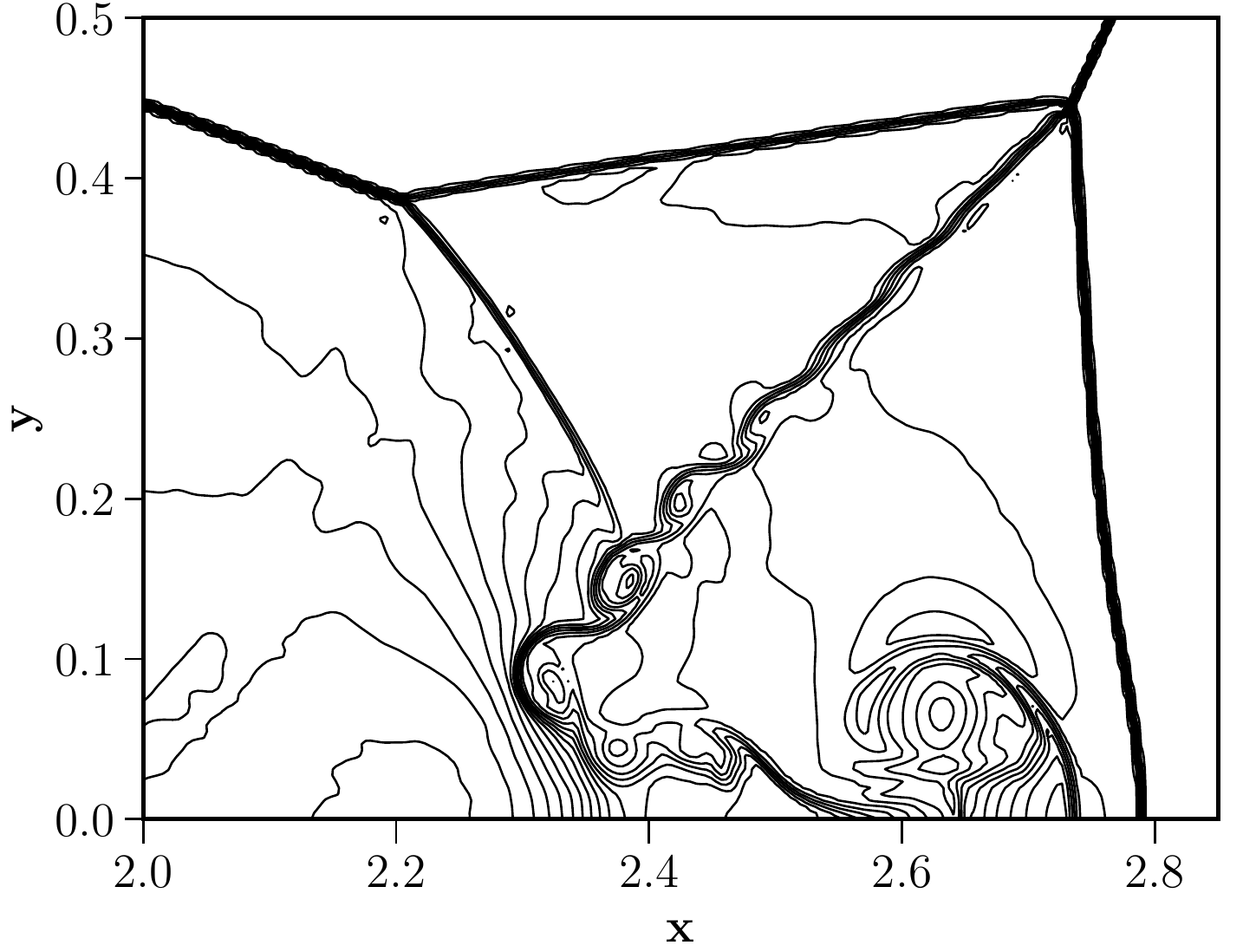}
\label{fig:HOCUS6_DB-LF}}
\caption{\textcolor{black}{Density contours at $t = 0.2$ for the three numerical schemes in the blown-up region around the Mach stem for Example \ref{ex:dmr} using using GLF Riemann solver on a grid size of $1024 \times 256$.}}
\label{fig_doublemach_LF}
\end{figure}

\begin{figure}[H]
\begin{onehalfspacing}
\centering\offinterlineskip
\subfigure[WENO-Z-HLLC]{\includegraphics[width=0.7\textwidth]{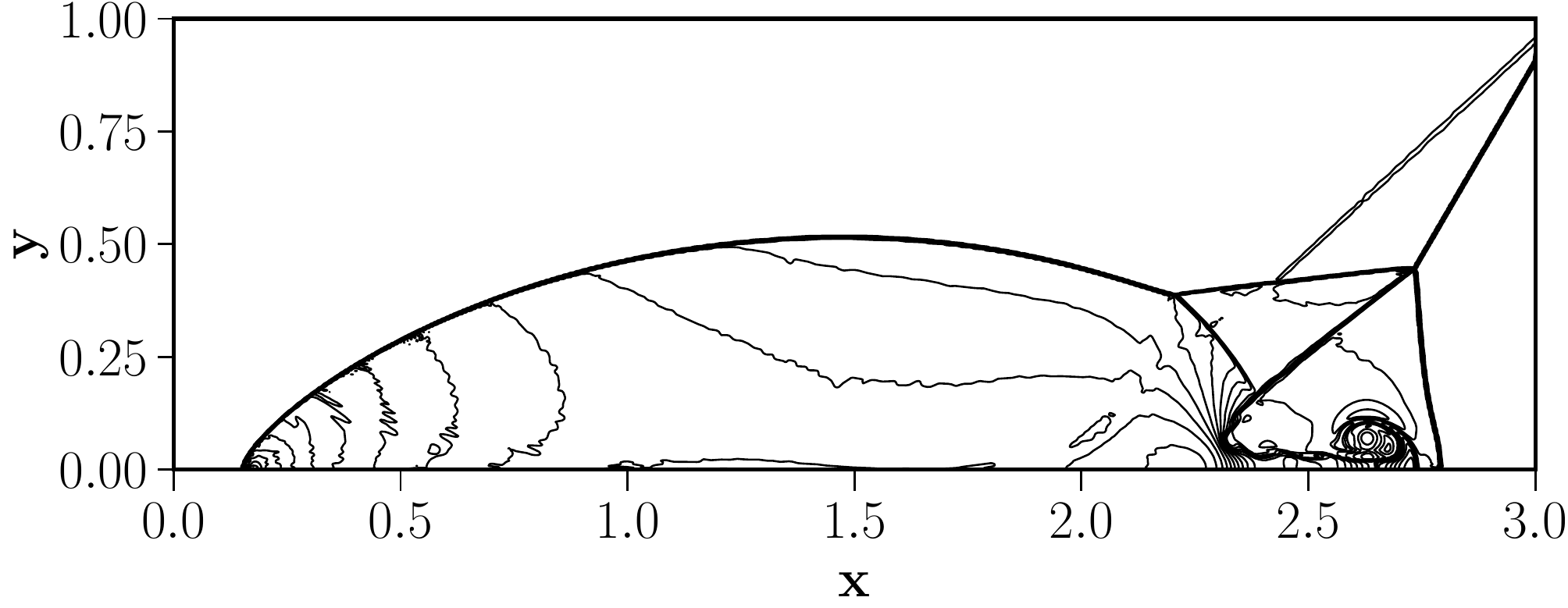}
\label{fig:WENO-Z_DB-full}}
\subfigure[MP5-HLLC]{\includegraphics[width=0.7\textwidth]{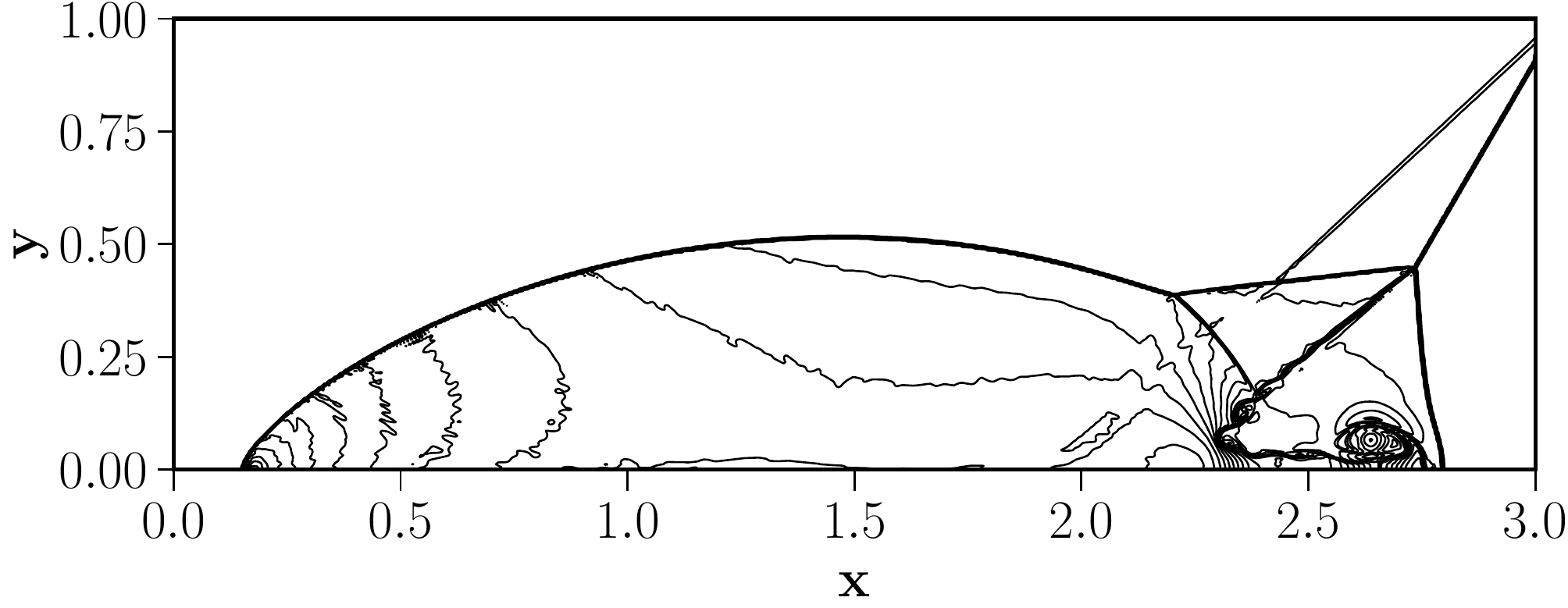}
\label{fig:MP5_DB-full}}
\subfigure[HOCUS5-HLLC]{\includegraphics[width=0.7\textwidth]{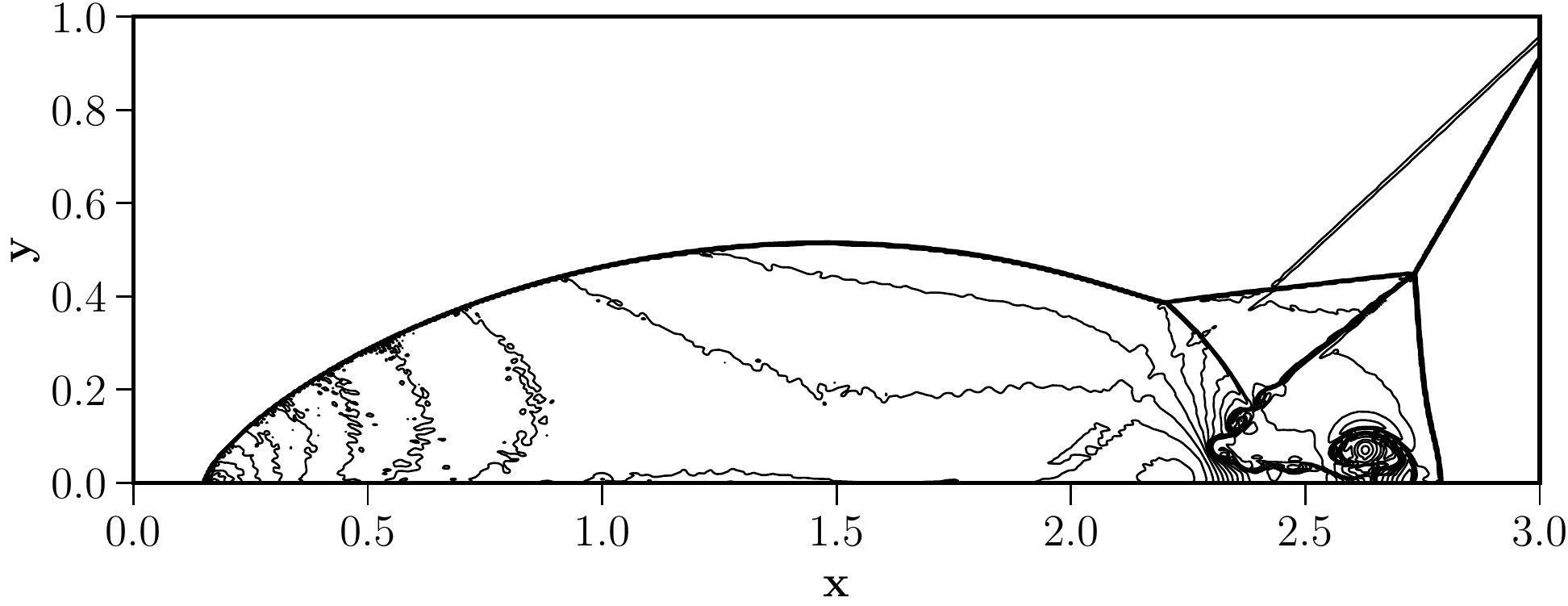}
\label{fig:HOCUS5_DB-full}}
\subfigure[HOCUS6-HLLC]{\includegraphics[width=0.7\textwidth]{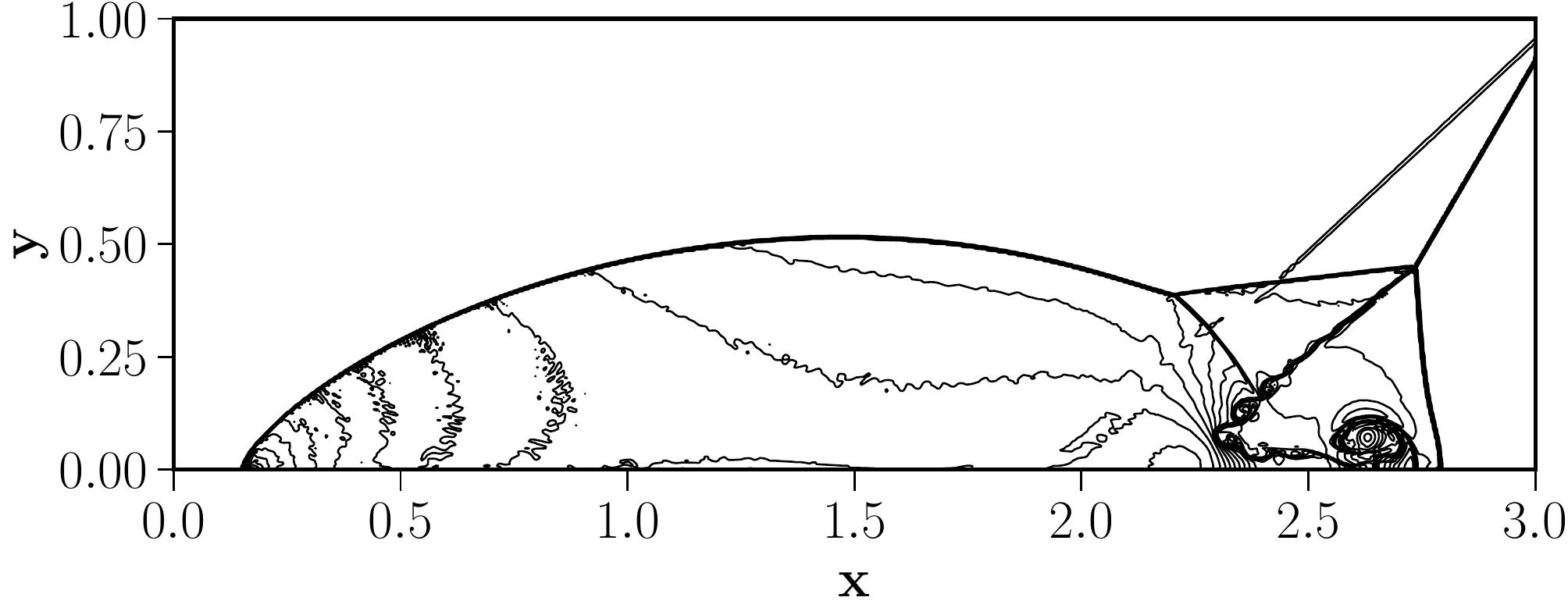}
\label{fig:HOCUS6_DB-full}}
\caption{\textcolor{black}{Double Mach reflection, Example \ref{ex:dmr}, $40$ equally spaced contours for the density on a grid size of $1024 \times 256$, $t = 0.2$ by HLLC Riemann solver with various schemes.}}
\label{fig_doublemach-full}
\end{onehalfspacing}
\end{figure}

\begin{example}\label{ex:vs}{Viscous Shock tube}
\end{example}

\textcolor{black}{In this last test case, we consider the viscous shock tube test case to demonstrate the resolution ability of the proposed schemes. This test case simulates the two-dimensional shock-wave and boundary layer interaction. A shock wave is initialized at the mid-location in the domain which interacts with a developing boundary layer, as it reflects from the slip wall at the right boundary. Also, the contact discontinuity and the reflected shock wave interact with each other. The domain for this test case is taken as $x \in [0,1], y \in [0,0.5]$. Viscosity is assumed to be constant, and the Prandtl number is taken as 0.73. A symmetry boundary condition is imposed on the top boundary, and no-slip wall boundary condition is used for other boundaries. The initial conditions are given by equation \eqref{vst} keeping the Mach number of the shock wave at 2.37.}

\textcolor{black}{
\begin{equation}\label{vst}
\left( {\rho , u,v, p} \right) = \left\{ \begin{array}{l}
\left( {120, 0 ,0,120/\gamma } \right),     0 < x < 0.5,\\
\left( {1.2, 0,0, 1.2/\gamma } \right),       0.5 \le x < 1,
\end{array} \right.
\end{equation}}

\textcolor{black}{The problem is solved for a Reynolds number of $Re=500$ on a grid size of $500 \times 200$ cells for time $t=1$. The sixth-order compact scheme is used to compute the velocity and temperature gradients, and the viscous stresses and heat-fluxes are calculated by 6th-order interpolation scheme \cite{lele1992compact,fang2014investigation}. As shown in Fig. \ref{fig_VST}, the results of WENO-Z and MP5 schemes show significant differences from the ``converged'' results in Ref. \cite{daru2009numerical} (see Fig 8b). Both HOCUS5 and HOCUS6 are closer to the reference result despite using a coarse mesh. The primary vortex of HOCUS6 matches well with the results of the reference, obtained on a grid size of 1500 $\times$ 750, whereas that from HOCUS5 is slightly distorted.}

 \textcolor{black}{As shown in Table \ref{tab:comptime}, it is observed that the HOCUS6 scheme takes more computational time than the WENOZ and MP5 schemes for a given grid size. Concerning the performance, the WENOZ and MP5 schemes need twice the resolution (results not shown here) in order to obtain comparable results as HOCUS6, which leads to two-times higher computing times. It shows that HOCUS schemes do not significantly increase the computational cost.}

% Table generated by Excel2LaTeX from sheet 'Sheet2'
\begin{table}[H]
  \centering
 \footnotesize
  \caption{\textcolor{black}{Computing time for Viscous shock tube test case, Example \ref{ex:vs}}}
    \begin{tabular}{lcccc}
  \hline
          &       & \multicolumn{2}{c}{Numerical method} &  \\
  \hline
          & WENOZ & MP5   & HOCUS5 & HOCUS6 \\
  \hline
      Grid  &       &       &       &  \\
  \hline
    \multicolumn{1}{c}{500 $\times$ 250} & 1240s & 1269s & 1718s & 1704s \\
  \hline
    750 $\times$ 375 & 4200s & 4312s & 5841s & 5823s \\
  \hline
    \end{tabular}%
  \label{tab:comptime}%
\end{table}%

\begin{figure}[H]
\centering
\subfigure[WENO-Z]{\includegraphics[width=0.45\textwidth]{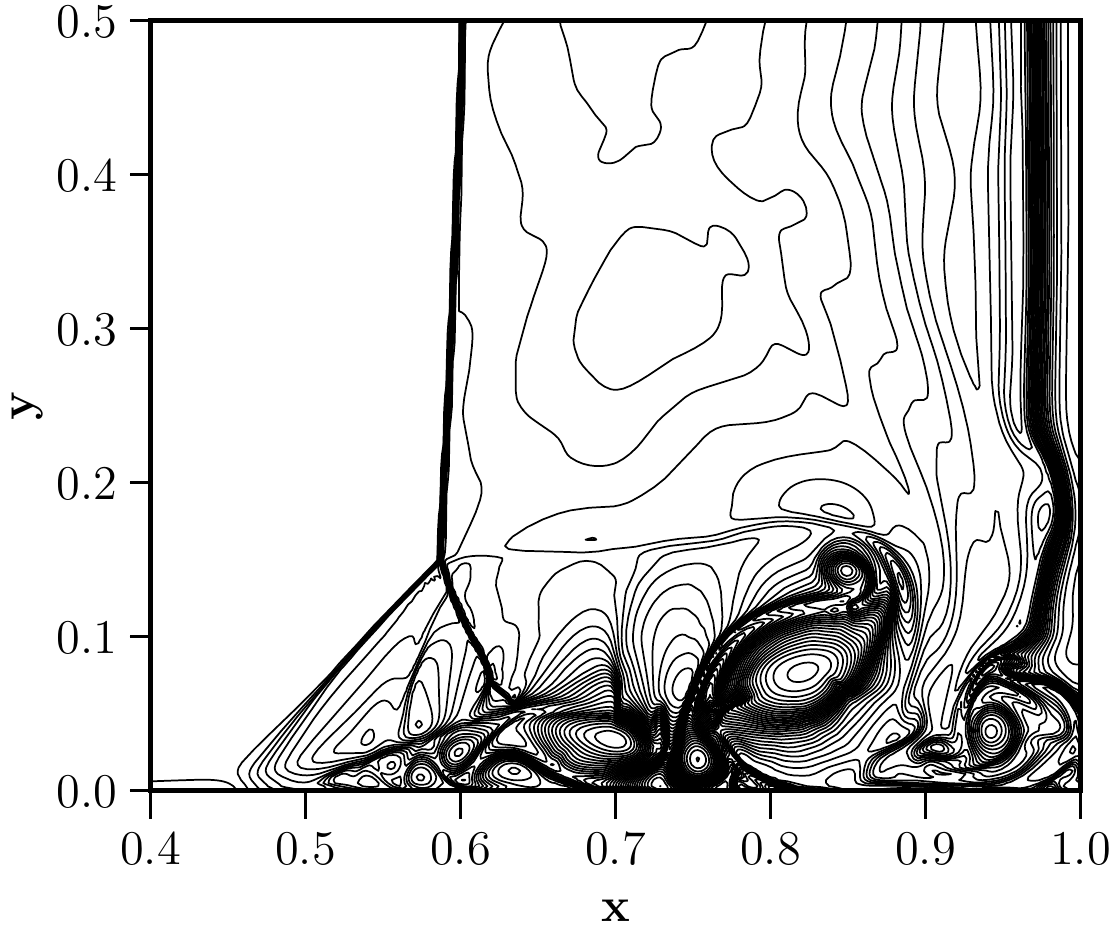}
\label{fig:WENO-Z_VST}}
\subfigure[MP5]{\includegraphics[width=0.45\textwidth]{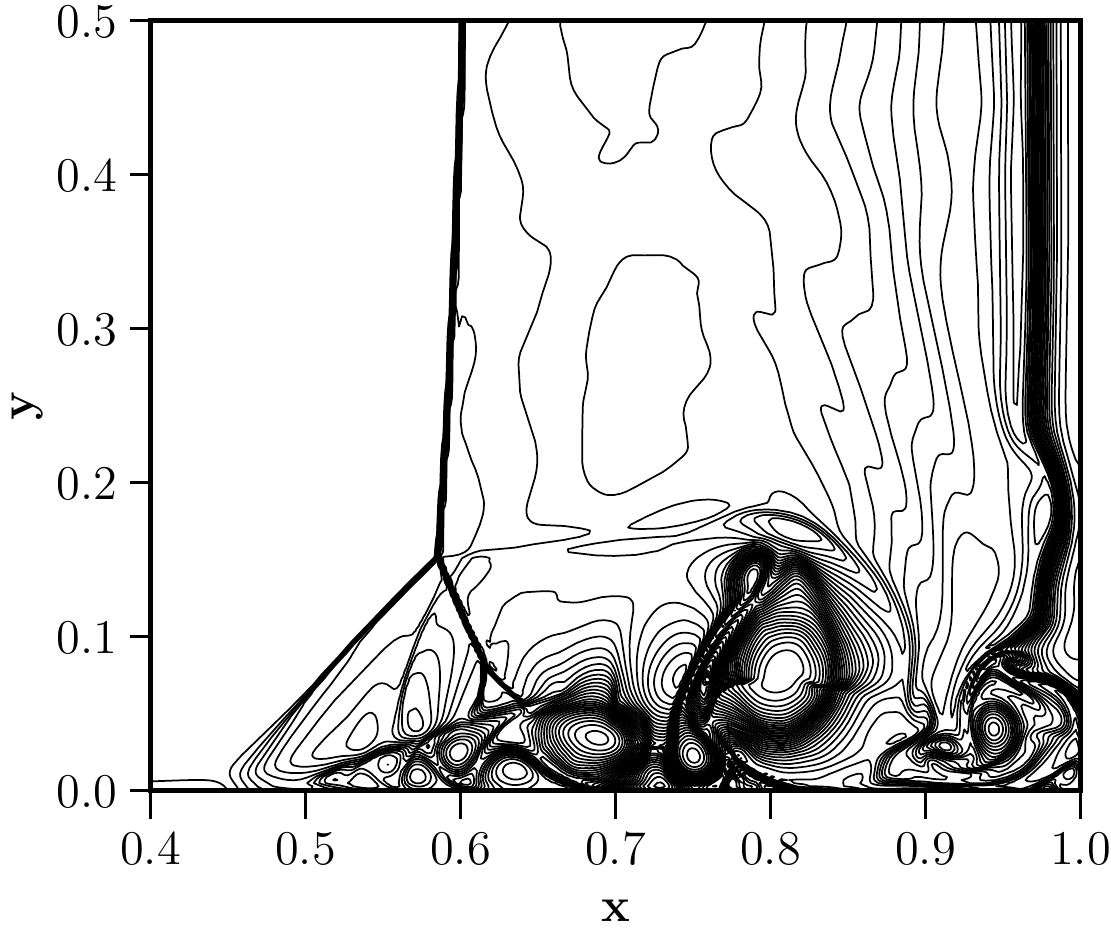}
\label{fig:MP5_VSST}}
\subfigure[HOCUS5]{\includegraphics[width=0.45\textwidth]{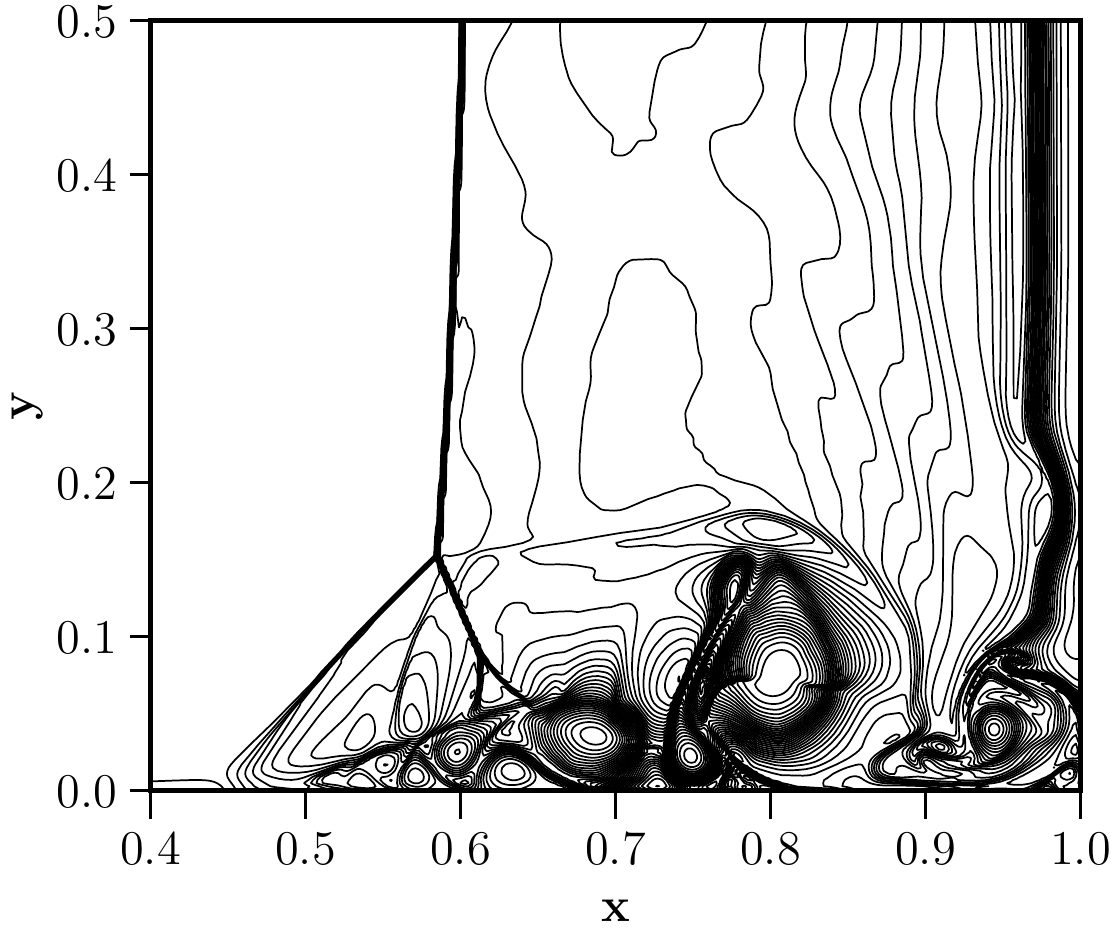}
\label{fig:HOCUS5_VST}}
\subfigure[HOCUS6]{\includegraphics[width=0.45\textwidth]{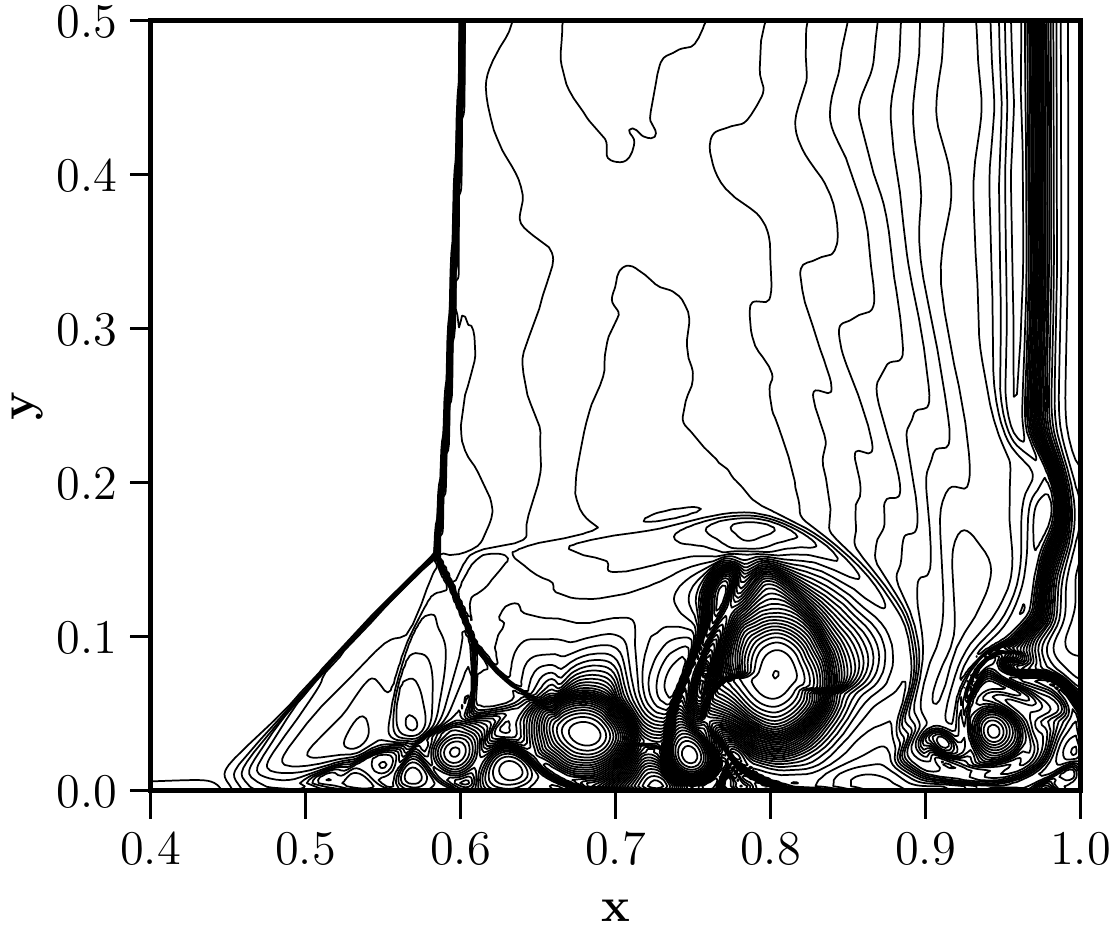}
\label{fig:HOCUS6_VST}}
\caption{\textcolor{black}{Density contours for different schemes for the Example \ref{ex:vs} with $Re=500$ on a grid size of $500 \times 250$.}}
\label{fig_VST}
\end{figure}

%----------------------------------------%----------------------------------------%----------------------------------------%----------------------------------------%----------------------------------------
%----------------------------------------%----------------------------------------%----------------------------------------%----------------------------------------%----------------------------------------
%----------------------------------------%----------------------------------------%----------------------------------------%----------------------------------------%----------------------------------------
%----------------------------------------%----------------------------------------%----------------------------------------%----------------------------------------%----------------------------------------
%----------------------------------------%----------------------------------------%----------------------------------------%----------------------------------------%----------------------------------------
%----------------------------------------%----------------------------------------%----------------------------------------%----------------------------------------%----------------------------------------
\section{Conclusions}\label{sec-4}	
We have developed a new hybrid scheme that combines an explicit MP5 scheme and a linear upwind compact reconstruction scheme using a BVD algorithm to enable high-resolution of smooth flow features and low-dissipation non-oscillatory shock capturing.  The BVD algorithm, which serves as a discontinuity detector, chooses the highest possible interpolation for the smooth solution, which is central and non-dissipative in the current approach and for the discontinuities, the algorithm prefers the monotone scheme. Primitive variables are interpolated for the compact scheme, and characteristic variable interpolation is  \textcolor{black}{carried out for the monotone scheme} for cleaner results. It, therefore, reduces the numerical dissipation inherently present in the upwind schemes and Riemann solvers, which \textcolor{black}{lead} to significant improvements in the numerical solutions for all the test cases considered in this paper. This method provides an alternative to the existing adaptive upwind-central schemes in the literature. It is possible to improve the current approach by optimizing dissipation properties of the MP5 scheme as in \cite{fang2013optimized}. In future, we will extend the methodology to the compressible multi-component flows and unstructured grids. 

%----------------------------------------%----------------------------------------%----------------------------------------%----------------------------------------%----------------------------------------
%----------------------------------------%----------------------------------------%----------------------------------------%----------------------------------------%----------------------------------------
%----------------------------------------%----------------------------------------%----------------------------------------%----------------------------------------%----------------------------------------
%----------------------------------------%----------------------------------------%----------------------------------------%----------------------------------------%----------------------------------------
%----------------------------------------%----------------------------------------%----------------------------------------%----------------------------------------%----------------------------------------
%----------------------------------------%----------------------------------------%----------------------------------------%----------------------------------------%----------------------------------------

\section*{Acknowledgements}
A.S. is supported by Technion fellowship during this work.
%----------------------------------------%----------------------------------------%----------------------------------------%----------------------------------------%----------------------------------------
%----------------------------------------%----------------------------------------%----------------------------------------%----------------------------------------%----------------------------------------
%----------------------------------------%----------------------------------------%----------------------------------------%----------------------------------------%----------------------------------------
%----------------------------------------%----------------------------------------%----------------------------------------%----------------------------------------%----------------------------------------
%----------------------------------------%----------------------------------------%----------------------------------------%----------------------------------------%----------------------------------------
%----------------------------------------%----------------------------------------%----------------------------------------%----------------------------------------%----------------------------------------
\section*{Appendix}
\renewcommand{\thesubsection}{\Alph{subsection}}
\textcolor{black}{In this appendix we present the various combinations of the numerical schemes that are also considered as part of this research work.}
%----------------------------------------%----------------------------------------%----------------------------------------%----------------------------------------%----------------------------------------
%----------------------------------------%----------------------------------------%----------------------------------------%----------------------------------------%----------------------------------------
%----------------------------------------%----------------------------------------%----------------------------------------%----------------------------------------%----------------------------------------
%----------------------------------------%----------------------------------------%----------------------------------------%----------------------------------------%----------------------------------------
%----------------------------------------%----------------------------------------%----------------------------------------%----------------------------------------%----------------------------------------
%----------------------------------------%----------------------------------------%----------------------------------------%----------------------------------------%----------------------------------------

\subsection{\textcolor{black}{HOCUS-TVD}} \label{sec-appb}

{\textcolor{black}{In this subsection, we present the hybrid approach by combining the compact reconstruction and the MUSCL schemes, named as HOCUS-TVD, with minmod limiter through boundary variation diminishing algorithm. The algorithm is the same as in Section \ref{sec2.2} except for the modification where the MUSCL scheme replaces the MP5 scheme in Equation \ref{Eq:TBVSC}.}

{\textcolor{black}{\begin{equation}\label{Eq:TBV-SC}
TBV_{j}^{MUSCL}=\big|U_{j-\frac{1}{2}}^{L,MUSCL}-U_{j-\frac{1}{2}}^{R,MUSCL}\big|+\big|U_{j+\frac{1}{2}}^{L,MUSCL}-U_{j+\frac{1}{2}}^{R,MUSCL} \big|.
\end{equation}
   Third order MUSCL scheme is briefly outlined below.}
\begin{eqnarray}
u_{i+1 / 2}^{R}&=u_{i+1}-\frac{1}{4}\left[(1-\bar{\eta}) \tilde \Delta_{i+3 / 2} u+(1+\bar{\eta}) \ddot \Delta_{i+1 / 2} u\right] \\
u_{i+1 / 2}^{L}&=u_{i+0}+\frac{1}{4}\left[(1-\bar{\eta})  \ddot \Delta_{i-1 / 2} u+(1+\bar{\eta}) \tilde \Delta_{i+1 / 2} u\right]
\end{eqnarray}

\begin{eqnarray}
\tilde\Delta_{i+1 / 2} u &=\operatorname{minmod}\left(\Delta_{i+1 / 2} u, \omega \Delta_{i-1 / 2} u\right) \\
\ddot\Delta_{i+1 / 2} u&=\operatorname{minmod}\left(\Delta_{i+1 / 2} u, \omega \Delta_{i+3 / 2} u\right)
\end{eqnarray}
{\textcolor{black}{ where, $1 \leq \omega \leq \frac{3-\eta}{1-\eta}$ and minmod function defined as }
\[ 
\mbox{minmod}(a,b) = \left\{ \begin{array}{cl}
     a  & \mbox{if $|a|<|b|$ and $ab>0$} \\
     b  & \mbox{if $|b|<|a|$ and $ab>0$} \\
     0  & \mbox{if $ab\leq 0$}
                  \end{array}
                  \right. 
\]

{\textcolor{black}{The coefficient $\bar{\eta}$ determines the order of accuracy. A second-order accuracy is achieved by using $\bar{\eta}=-1$ and a third-order accuracy is achieved by using $\bar{\eta}=\frac{1}{3}$. In this paper we considered third order accurate scheme for all our computations. We carried out numerical experiments for the test cases \ref{ex:adv_comp}, \ref{Titarev-Toro}, \ref{Shu-Osher}, \ref{blast}, \ref{ex:dmr} and \ref{ex:vs} using the HOCUS-TVD algorithm. The parameter $\omega$ is taken as 4 for all the simulations. The results are shown in Fig.\ \ref{fig_tvd-1D}. All the numerical simulations \textcolor{black}{are carried out} with a CFL of $0.2$ and the grid sizes are same as that of HOCUS6. We make the following observations:}

\begin{itemize}
\item \textcolor{black}{Figure \ref{fig:Jiang-tvd} shows the results for the advection of complex waves test case, \ref{ex:adv_comp}. We can see that the HOCUS-TVD is slightly more dissipative than the HCOUS6 method for the semi-elliptical profile. Figures \ref{fig:tita-tvd}, \ref{fig:shu-tvd} and \ref{fig:tita-blast} show density profile for the one-dimensional Euler test cases \ref{Titarev-Toro}, \ref{Shu-Osher}, and \ref{blast} respectively. We can see that in all the cases HCOUS-TVD is more dissipative than the HOCUS6.}

\item \textcolor{black}{The two-dimensional test cases are more representative of the properties of the HOCUS-TVD scheme as shown in Figures \ref{fig_tvd-2D}. Fig. \ref{fig:DMR-hllc-tvd} shows the density contours of the Double Mach reflection test case which indicates that the HOCUS-TVD is more dissipative than HOCUS6, shown in \ref{fig:HOCUS-VST-tvd}, as well as the MP5 scheme, shown in Fig. \ref{fig:HOCUS5_DB}. Similar results can be seen for the Viscous shock tube simulation where the primary vortex is distorted compared to the HOCUS6 method, as shown in Figs. \ref{fig:HOCUS-VST-tvd} and \ref{fig:HOCUS6-VST-tvd}.}

\item \textcolor{black}{As explained in the introduction, we need at least five points to distinguish between, and extrema and a genuine discontinuity. These results indicate that the combination of compact and MP5 schemes are superior over the combination of compact and TVD schemes at least for the present algorithm.}
\end{itemize}

\begin{figure}[H]
\centering\offinterlineskip
\subfigure[]{\includegraphics[width=0.7\textwidth]{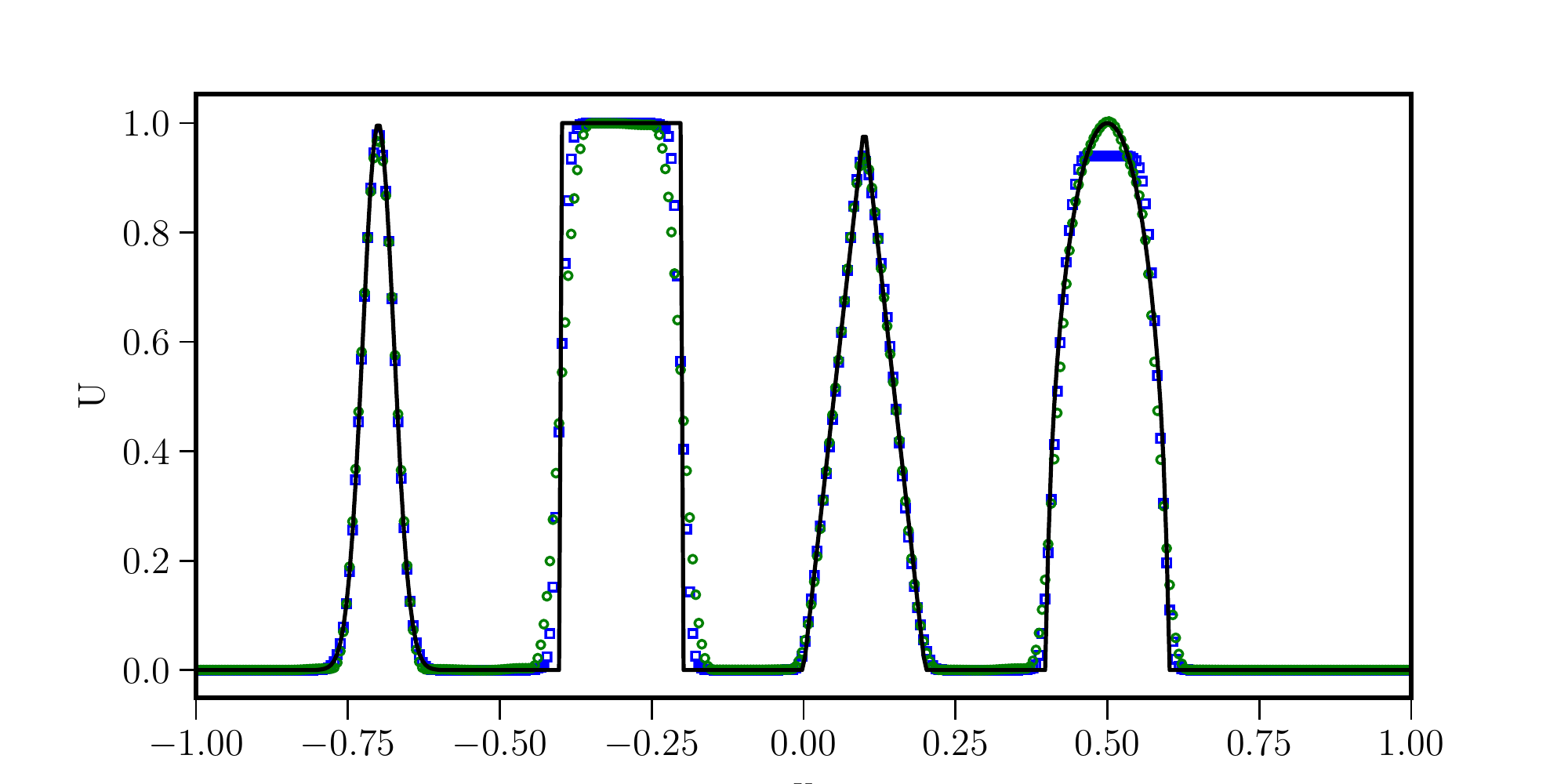}
\label{fig:Jiang-tvd}}
\subfigure[]{\includegraphics[width=0.34\textwidth]{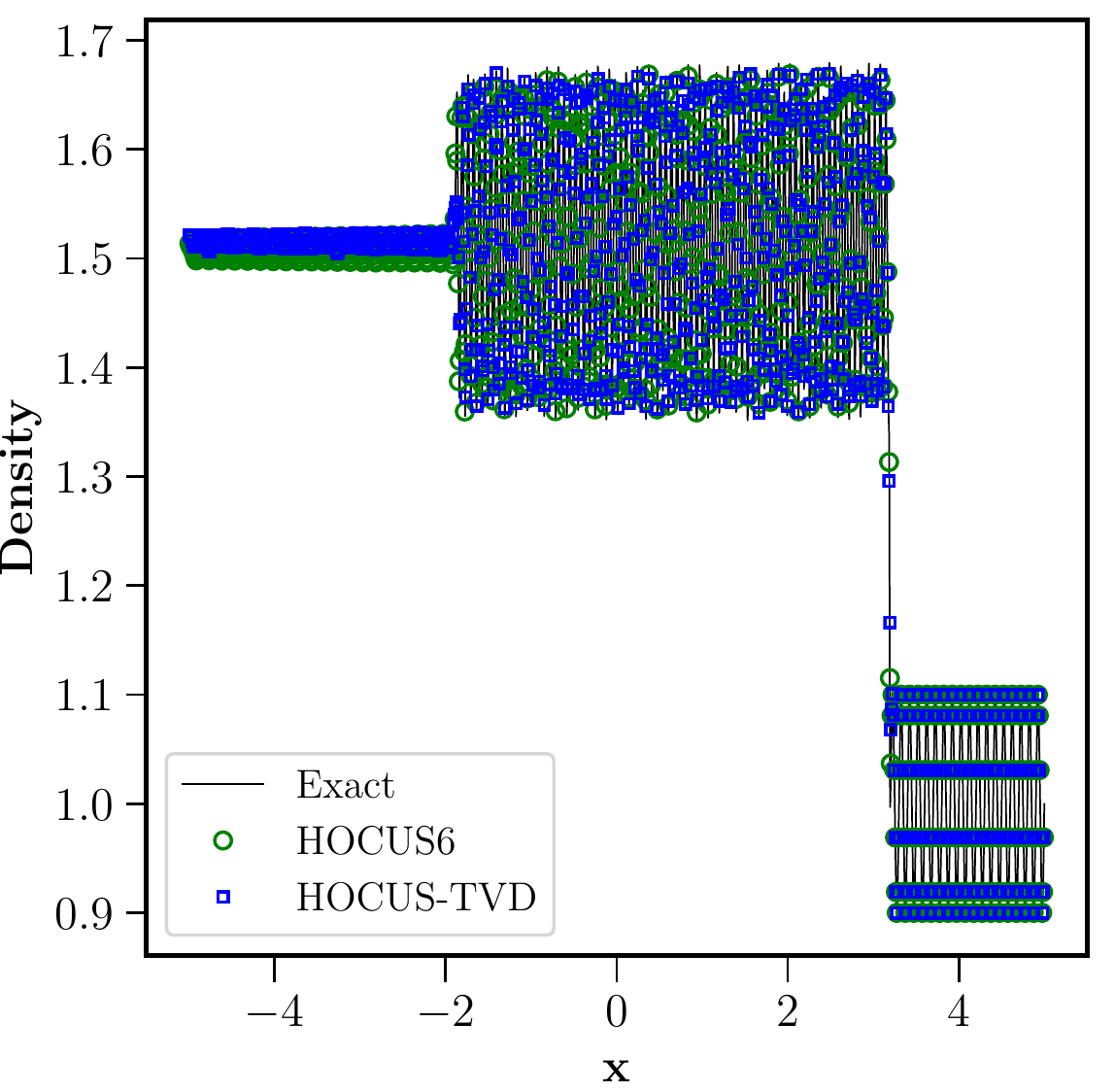}
\label{fig:tita-tvd}}
\subfigure[]{\includegraphics[width=0.34\textwidth]{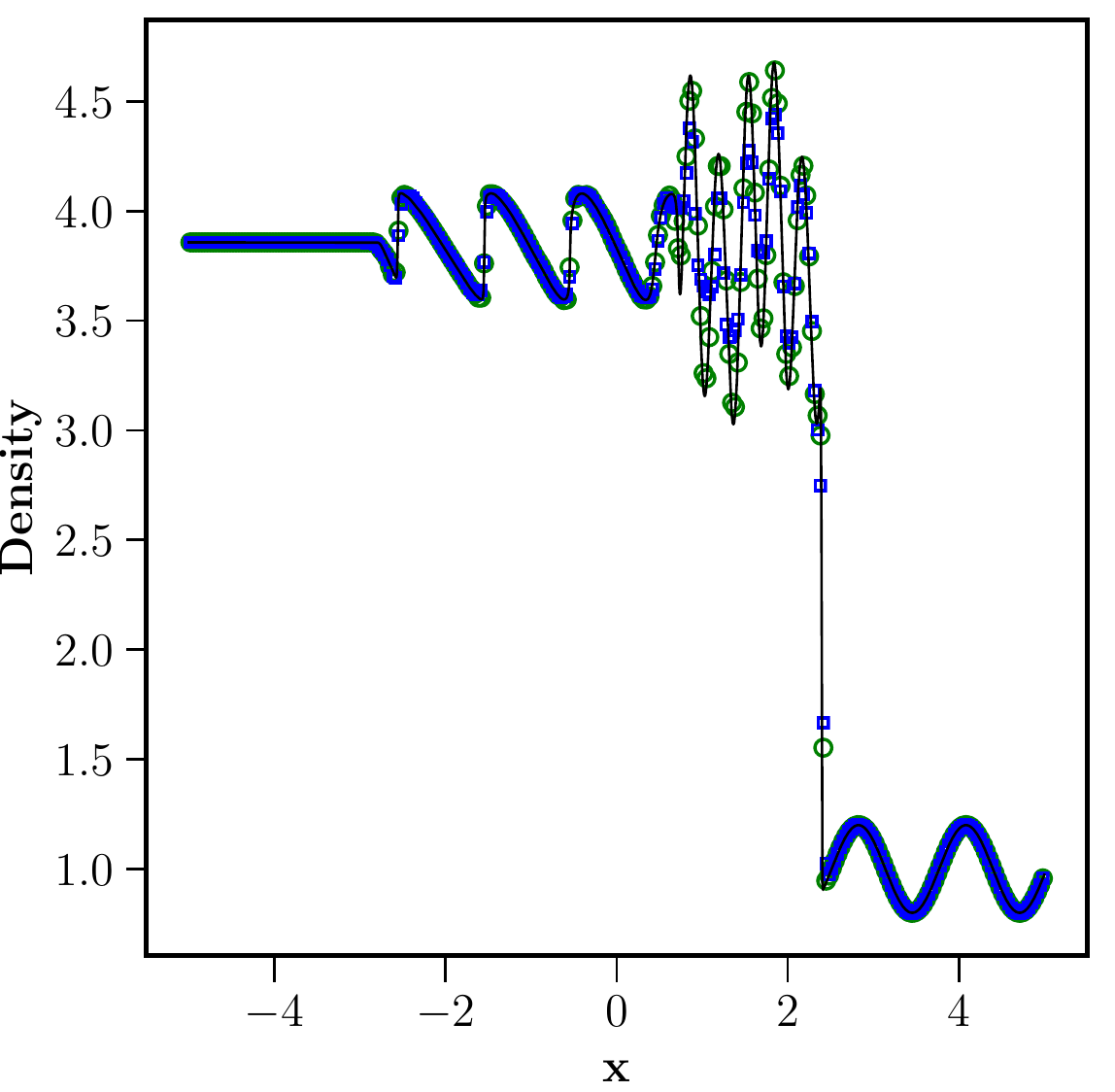}
\label{fig:shu-tvd}}
\subfigure[]{\includegraphics[width=0.34\textwidth]{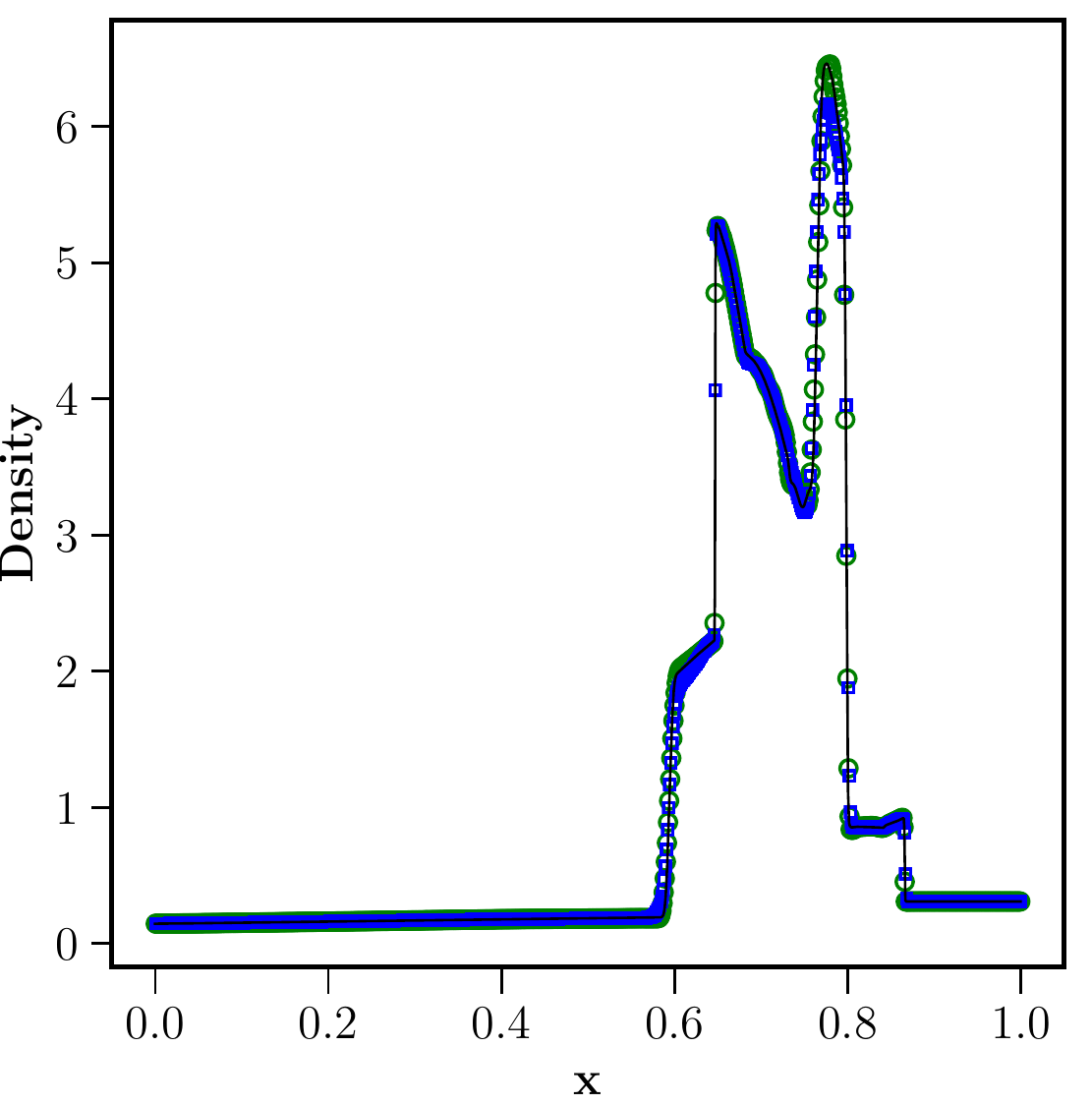}
\label{fig:tita-blast}}
%\subfigure[]{\includegraphics[width=0.44\textwidth]{HOCUS_entropy_TVD.pdf}
%\label{fig:entropy-tvd}}
\caption{\textcolor{black}{Density profiles by HOCUS-TVD and HOCUS6 schemes for one-diemnsional test cases. Figure (a) corresponds to Example \ref{ex:adv_comp}, (b) Example \ref{Titarev-Toro}, (c) Example \ref{Shu-Osher} and (d) corresponds to Example \ref{blast} respectively. solid line: reference solution; green circles: HOCUS6; blue squares: HOUCS-TVD.} }
\label{fig_tvd-1D}
\end{figure}

\begin{figure}[H]
\centering\offinterlineskip
\subfigure[]{\includegraphics[width=0.44\textwidth]{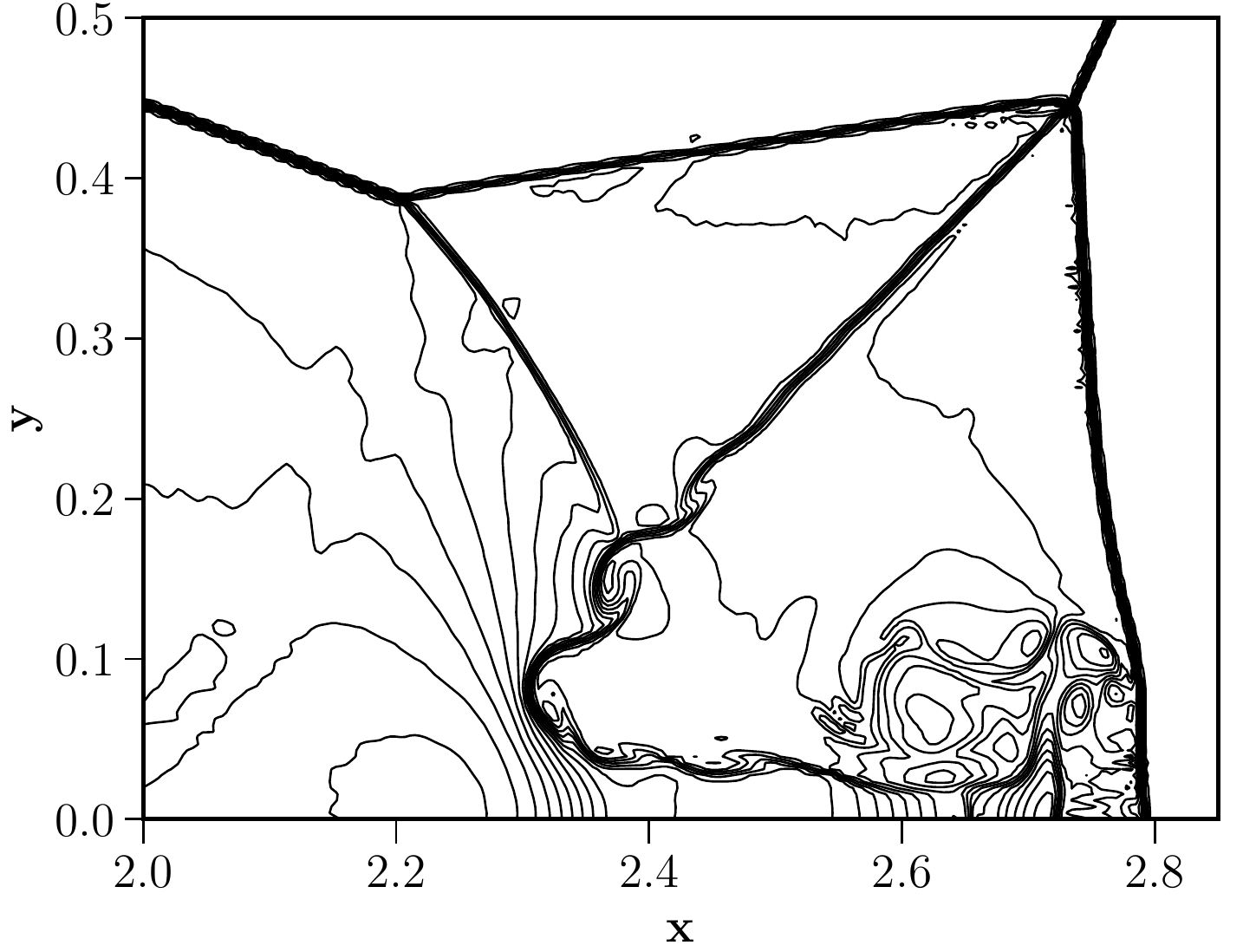}
\label{fig:DMR-hllc-tvd}}
\subfigure[]{\includegraphics[width=0.44\textwidth]{DMR-BVD-local.pdf}
\label{fig:DMR-HOCUS6-tvd}}
\subfigure[]{\includegraphics[width=0.44\textwidth]{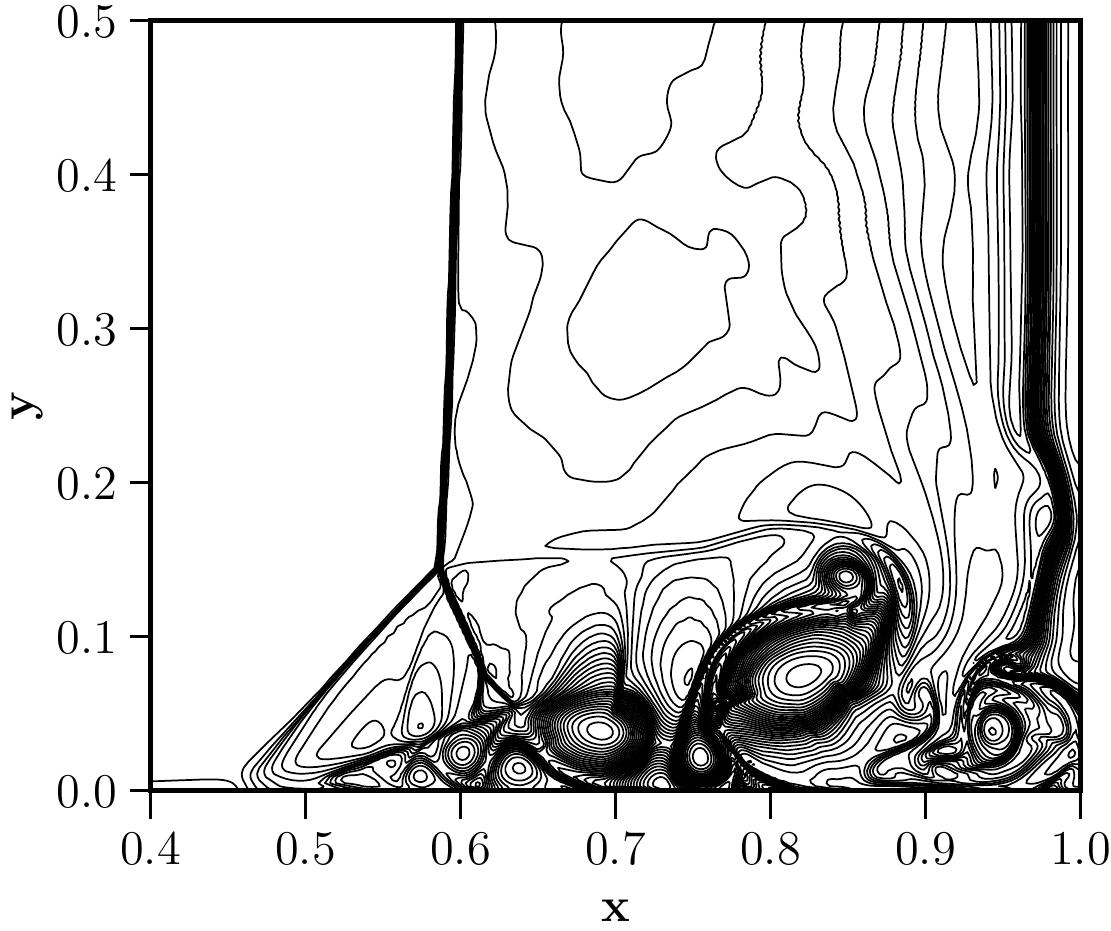}
\label{fig:HOCUS-VST-tvd}}
\subfigure[]{\includegraphics[width=0.44\textwidth]{VST_HOCUS6.pdf}
\label{fig:HOCUS6-VST-tvd}}
\caption{\textcolor{black}{Figs. (a) and (b) shows the density distribution for Example \ref{ex:dmr} and Figs. (c) and (d)  shows the density distribution for Example \ref{ex:vs} using HOCUS-TVD and HOCUS6 respectively.}}
\label{fig_tvd-2D}
\end{figure}

\subsection{\textcolor{black}{C5-T2 scheme}}\label{sec-appc}
\textcolor{black}{In this subsection, we present the hybrid approach by combining the compact reconstruction and the THINC schemes, named as C5-T2, using the two-stage boundary variation diminishing algorithm proposed by Deng et al. in \cite{deng2019fifth}. The complete algorithm is as follows.}
\begin{description}
\item[Step 1.] \textcolor{black}{Evaluate the interface values by using two different reconstruction procedures:}
\begin{enumerate}[(a)]
\item \textcolor{black}{Linear upwind compact reconstruction given by Equation (\ref{eqn:upwind-compact}) and}
\item \textcolor{black}{THINC scheme, which is as follows,}

\begin{equation}
\textcolor{black}{U_{i+1 / 2}^{L, T H I N C}=\left\{\begin{array}{l}
\bar{u}_{min}+\frac{\bar{u}_{\max }}{2}\left(1+\theta \frac{\tanh (\beta)+A}{1+A \tanh (\beta)}\right) \text { if }\left({u}_{i+1}-{u}_{i}\right)\left({u}_{i}-{u}_{i-1}\right)>0 \\
{u}_{i} \text { otherwise }
\end{array}\right.}
\end{equation}

\begin{equation}
\textcolor{black}{U_{i-1 / 2}^{R, T H I N C}=\left\{\begin{array}{ll}
\bar{u}_{\min }+\frac{\bar{u}_{\max }}{2}(1+\theta A) & \text { if }\left({u}_{i+1}-{u}_{i}\right)\left({u}_{i}-{u}_{i-1}\right)>0 \\
{u}_{i} \text { otherwise },
\end{array}\right.}
\end{equation}

where, $A=\frac{B / \cosh (\beta)-1}{\tanh (\beta)}$, $B=\exp (\theta \beta(2 C-1))$,  $C=\frac{{u}_{i}-\bar{u}_{\min }+\epsilon}{\bar{u}_{\max }+\epsilon}$, $\bar{u}_{\min }=\min \left({u}_{i-1}, {u}_{i+1}\right), {u}_{\max }=\max \left({u}_{i-1}, {u}_{i+1}\right)-{u}_{\min } \text { and } \theta=\operatorname{sgn}\left({u}_{i+1}-{u}_{i-1}\right)$

\end{enumerate}
\item[Step 2.] \textcolor{black}{Calculate the TBV values for each cell $I_{j}$ by using the compact reconstruction:}

\begin{equation}\label{Eq:TBVC5-1}
TBV_{j}^{C5}=\big|U_{j-\frac{1}{2}}^{L,C5}-U_{j-\frac{1}{2}}^{R,C5}\big|+\big|U_{j+\frac{1}{2}}^{L,C5}-U_{j+\frac{1}{2}}^{R,C5} \big|
\end{equation} 
\textcolor{black}{and the THINC scheme:}
\begin{equation}\label{Eq:TBV-THINC}
TBV_{j}^{THINC}=\big|U_{j-\frac{1}{2}}^{L,THINC}-U_{j-\frac{1}{2}}^{R,THINC}\big|+\big|U_{j+\frac{1}{2}}^{L,THINC}-U_{j+\frac{1}{2}}^{R,THINC} \big|.
\end{equation} 

\textcolor{black}{In this step the parameter $\beta$ in the THINC scheme is taken as 1.1.}

\item[Step 4.] \textcolor{black}{Now,  \textcolor{black}{all the interface values of the cells $j-1$, $j$, and $j+1$, both $R$ and $L$} are modified according to the following algorithm:}
\begin{equation}\label{eqn:BVDstep-1}
\begin{aligned} 
{\rm if } \ \ TBV_{j}^{THINC} < TBV_{j}^{C5} 
\end{aligned} 
\end{equation}

\item[Step 5.] \textcolor{black}{Denote the interface values, denoted as   ${U}_{j+\frac{1}{2}}^{L^{'},R^{'}}$, for stage one of the algorithm, S1. Calculate the TBV values for each cell $I_{j}$ by using the new interface values:}

\begin{equation}\label{Eq:TBVS-1}
TBV_{j}^{S1}=\big|U_{j-\frac{1}{2}}^{L^{'}}-U_{j-\frac{1}{2}}^{R^{'}}\big|+\big|U_{j+\frac{1}{2}}^{L^{'}}-U_{j+\frac{1}{2}}^{R^{'}} \big|
\end{equation} 
\textcolor{black}{and the THINC scheme:}
\begin{equation}\label{Eq:TBV-THINC-2}
TBV_{j}^{THINC}=\big|U_{j-\frac{1}{2}}^{L,THINC}-U_{j-\frac{1}{2}}^{R,THINC}\big|+\big|U_{j+\frac{1}{2}}^{L,THINC}-U_{j+\frac{1}{2}}^{R,THINC} \big|.
\end{equation} 

\textcolor{black}{In this step the parameter $\beta$ in the THINC scheme is taken as 1.6.}

\item[Step 6.] \textcolor{black}{Now,  \textcolor{black}{the interface values of cell  $j$, both $R$ and $L$} is modified according to the following algorithm:}
\begin{equation}\label{eqn:BVDstep-2}
\begin{aligned} 
{\rm if } \ \ TBV_{j}^{THINC} < TBV_{j}^{S1} 
\end{aligned} 
\end{equation}

\item[Step 7.] \textcolor{black}{Finally, evaluate the interface flux $\hat{F}_{j + \frac{1}{2}}$ from ${U}_{j+\frac{1}{2}}^{L,R}$ using Equation (\ref{eqn:Riemann})}
\end{description}

\textcolor{black}{We carried out numerical experiments for the test cases  \ref{Titarev-Toro}, \ref{blast}, \ref{ex:dmr} and \ref{ex:vs} using the C5-T2 algorithm. All the numerical simulations \textcolor{black}{are carried out} with a CFL of $0.2$ and the grid sizes are same as that of HOCUS6 with the following observations:}

\begin{itemize}
\item \textcolor{black}{Figure \ref{fig:blast-thinc} shows the results for the blast wave test case, \ref{fig:tita-blast}. We can see that the C5-T2 has oscillations, unlike the HCOUS6 method. Similar oscillations are also observed for the P4T2 scheme of Deng at al. \cite{deng2019fifth}. For the Titarev-Toro test case, the density profile obtained by C5-T2 is similar to that of HOCUS5 scheme, which is slightly dissipative in comparison with HOCUS6.}

\item \textcolor{black}{Fig. \ref{fig:DMR-THINC} shows the density contours of the Double Mach reflection test case. We can see the near-wall jet of C5-T2 scheme is distorted compared to that of HOCUS6. Similar results can be seen for the Viscous shock tube simulation where the primary vortex is distorted along with the vortical structures in the bottom-right corner due to the roll-up of the contact discontinuity and oscillations observed near the shock as shown in Figs. \ref{fig:VST-thinc} and \ref{fig:HOCUS6-VST-tvd}. HOCUS6 is more robust than the C5-T2 scheme.}

\begin{figure}[H]
\centering\offinterlineskip
\subfigure[]{\includegraphics[width=0.44\textwidth]{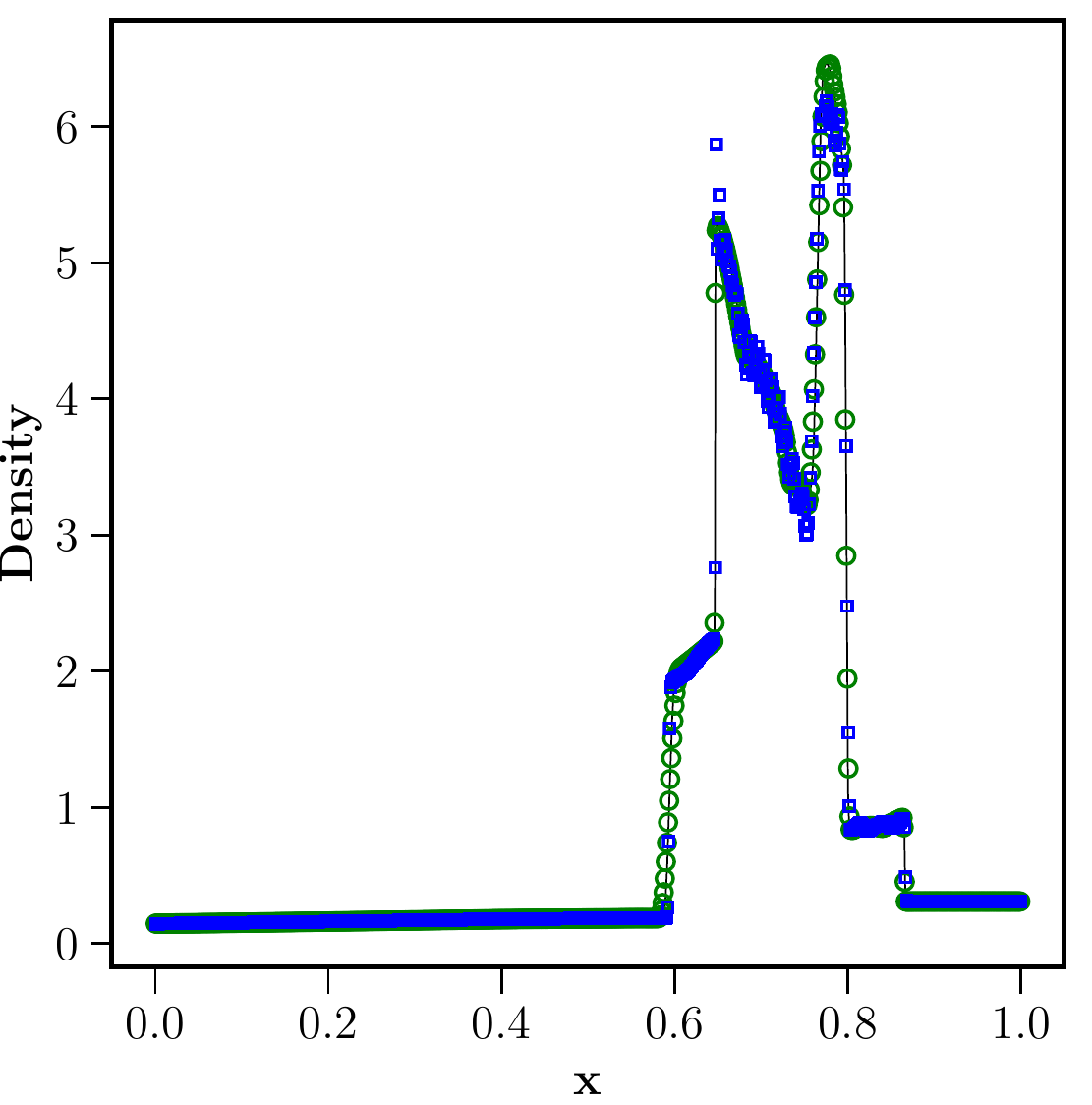}
\label{fig:blast-thinc}}
\subfigure[]{\includegraphics[width=0.44\textwidth]{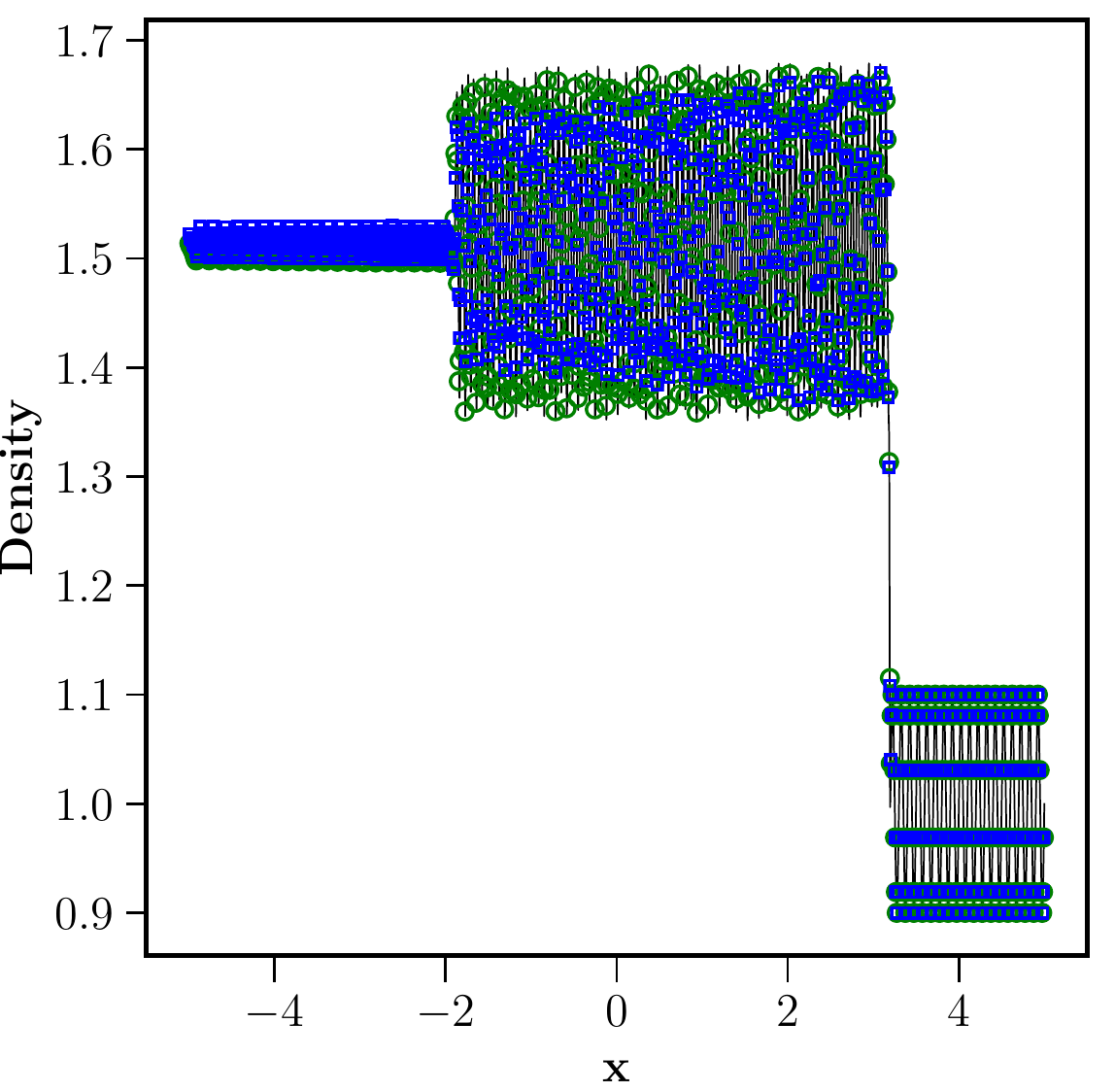}
\label{fig:tita-thinc}}
\subfigure[]{\includegraphics[width=0.44\textwidth]{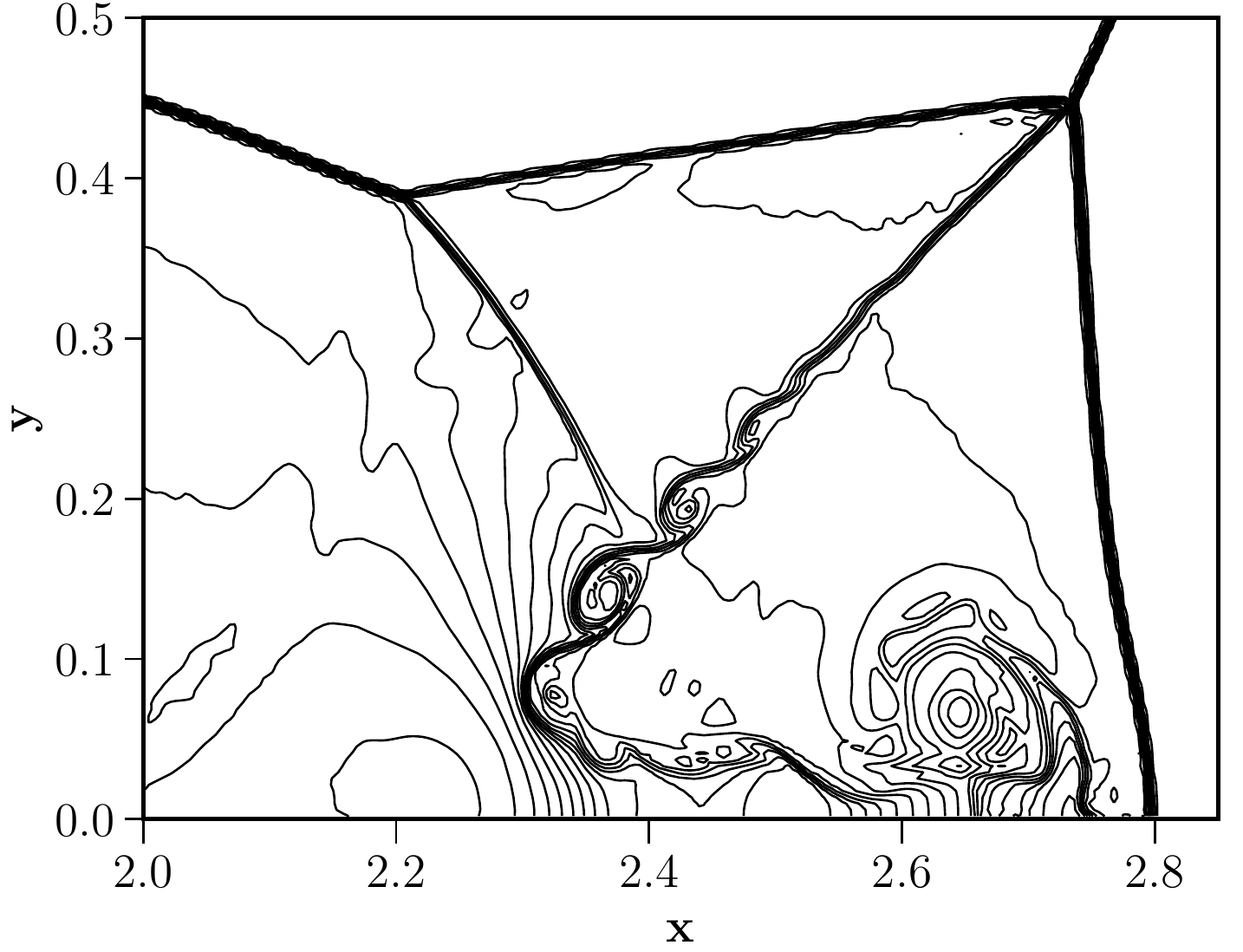}
\label{fig:DMR-THINC}}
\subfigure[]{\includegraphics[width=0.44\textwidth]{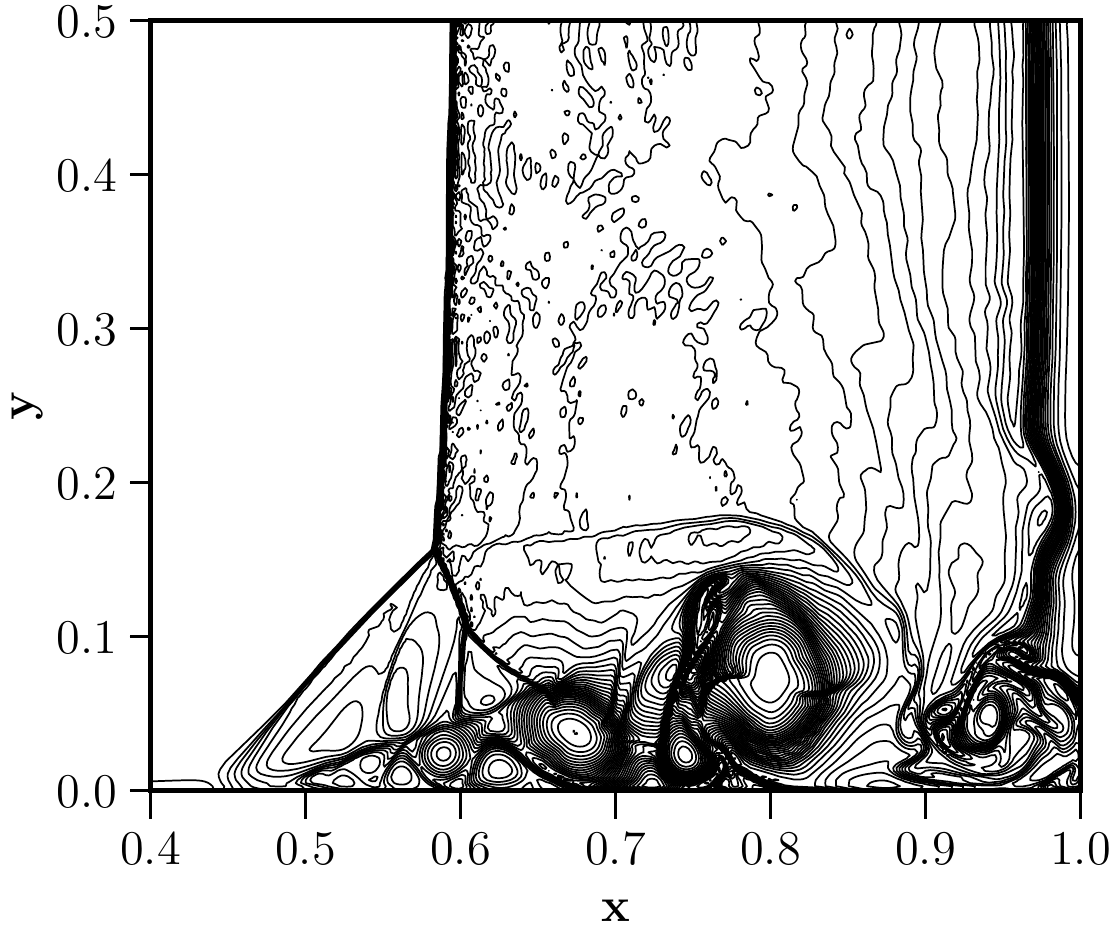}
\label{fig:VST-thinc}}
\caption{\textcolor{black}{Density profiles by C5-T2 and HOCUS6 schemes for one and two-diemnsional test cases. Figure (a) corresponds to Example \ref{blast}, (b) Example \ref{Titarev-Toro}, (c) Example \ref{ex:dmr} and (d) corresponds to Example \ref{ex:vs} respectively. solid line: reference solution; green circles: HOCUS6; blue squares: C5-T2. }}
\label{fig_thinc-12D}
\end{figure}

\item \textcolor{black}{Final test case we considered is the following two-dimensional Riemann problem taken from Schulz-Rinne et al. \cite{schulz1993numerical}. The simulation is carried out over a unit square domain $[0,1]\times [0,1]$, with the following initial data:}

\begin{equation}
\left(\rho, u, v, p_{0}\right)\\
=\left\{\begin{array}{lll}
(1.0,-0.6259,0.1,1.0) & x \leq 0.0, y \geq 0.0 \\
(0.8,0.1,0.1,1.0) & x<0.0, y<0.0 \\
(0.5197,0.1,0.1,0.4) & x>0.0, y>0.0 \\
(1.0,0.1,-0.6259,1.0) & x>0.0, y<0.0
\end{array}\right.
\end{equation}

\textcolor{black}{The small-scale complex structures generated along the slip lines due to Kelvin-Helmholtz instability serves as a benchmark to test numerical dissipation of the scheme. The numerical solutions are computed for time $t=0.25$ on a grid of size $1000\times 1000$. Non-reflective boundary conditions are applied at all the boundaries for this test case. A closer look at Figs. \ref{fig:RR-thinc} and \ref{fig:RR-HOCUS62} indicates that the small-scale structures obtained for the C5-T2 scheme are richer in comparison with the HOCUS6 scheme. Each scheme has its strengths and disadvantages, depending on the problem being addressed.}

\end{itemize}

\begin{figure}[H]
\centering\offinterlineskip
\subfigure[C5-T2]{\includegraphics[width=0.44\textwidth]{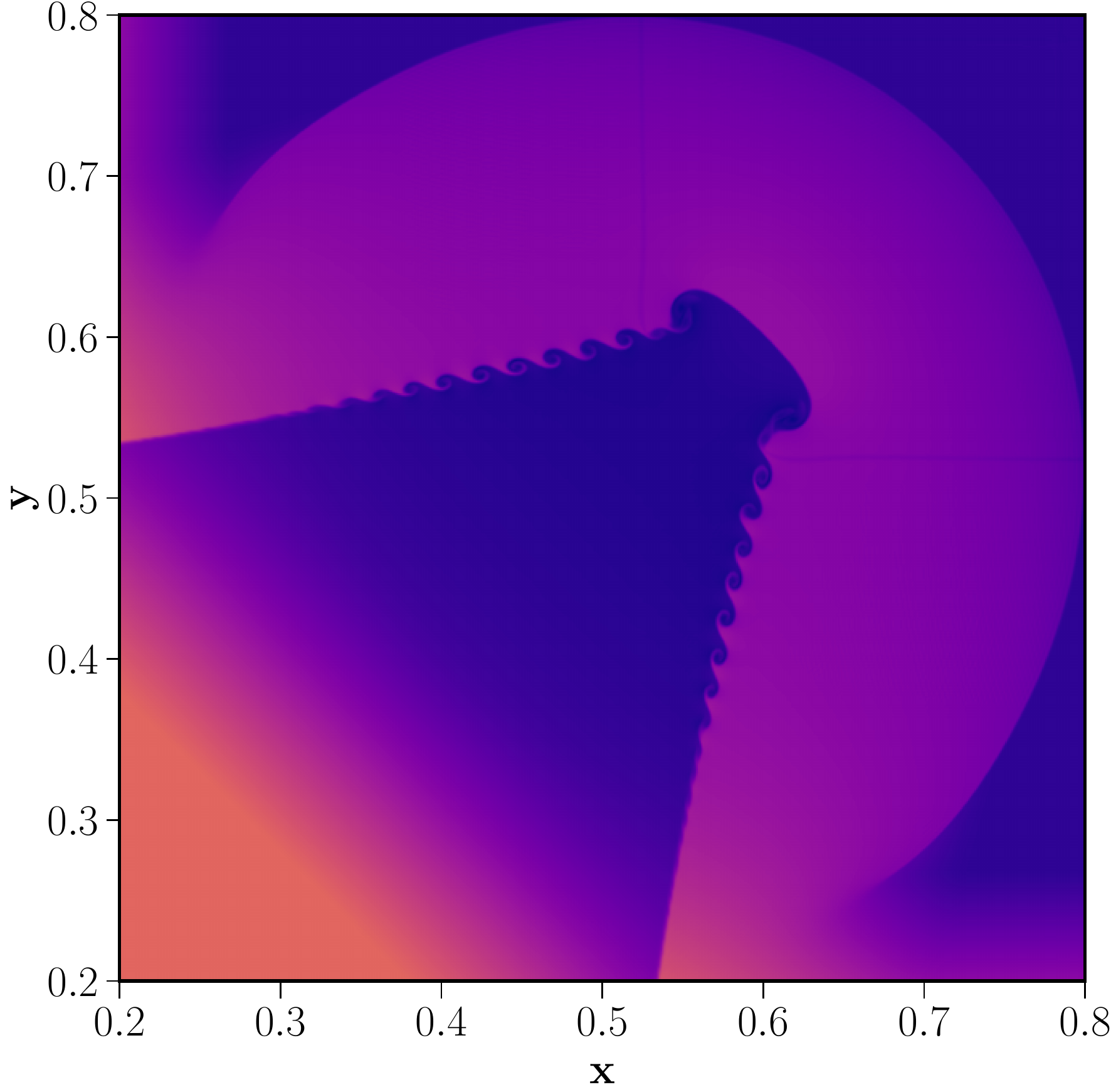}
\label{fig:RR-thinc}}
\subfigure[HOCUS6]{\includegraphics[width=0.44\textwidth]{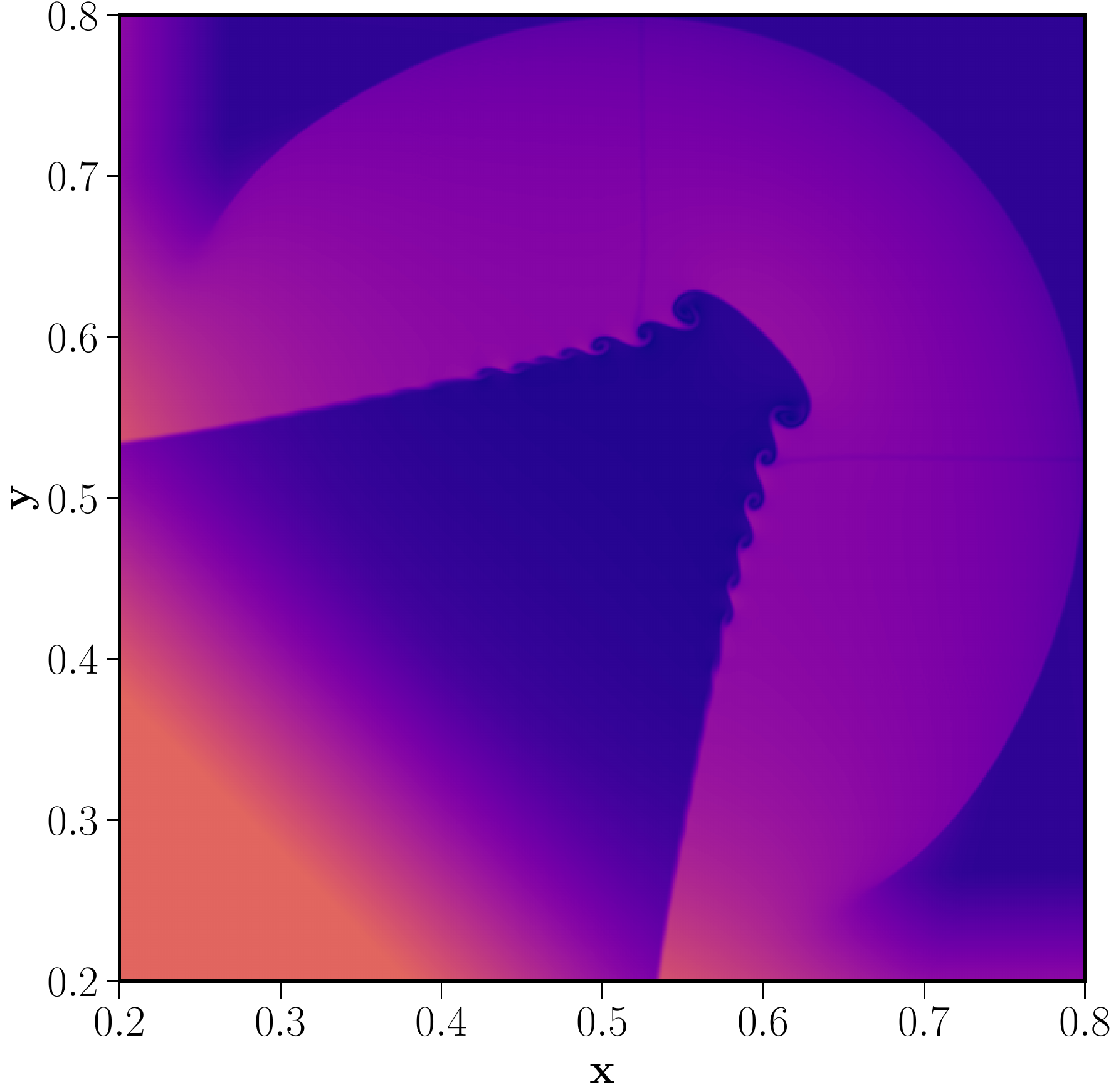}
\label{fig:RR-HOCUS62}}
\caption{Density profiles by C5-T2 and HOCUS6 schemes for Riemann problem}
\label{fig_thinc-2D-RR}s
\end{figure}

\subsection{\textcolor{black}{HOCUS-WENOZ}}
In this subsection, we present the hybrid approach by combining the compact reconstruction and the WENO-Z through boundary variation diminishing algorithm. In Fig. \ref{fig:algo1}, we can see the results for algorithm using the Equation (\ref{eqn:BVDstep}) for the advection of complex waves given in Example \ref{ex:adv_comp}. We can notice that the BVD algorithm with WENO-Z and the linear-compact scheme is diffusive and not able to preserve the wave patterns in comparison with the HOCUS6, Fig. \ref{fig:HOCUS6_AC}, which is a combination of MP5 and the linear-compact scheme. Similar results are observed for two-dimensional simulations for the Euler equations. In light of these results, we have also implemented the algorithm given in \cite{xie2017hybrid, ivan2014high} which is as follows

\begin{enumerate}[i)]
\item Compute the cell-interface value by the linear compact scheme given by Equation (\ref{eqn:upwind-compact})
\begin{equation}
\begin{cases} 
U_{{j+\frac{1}{2}}}^{L} = U_{{j+\frac{1}{2}}}^{L, C6} = \frac{1}{2}(U_{{j+\frac{1}{2}}}^{L, C5} +U_{{j+\frac{1}{2}}}^{R, C5})\\
\\
U_{{j+\frac{1}{2}}}^{R}=U_{{j+\frac{1}{2}}}^{R, C6} = \frac{1}{2}(U_{{j+\frac{1}{2}}}^{L, C5} +U_{{j+\frac{1}{2}}}^{R, C5})\\
\end{cases}
\end{equation} 

\item Compute the smoothness indicator
\begin{equation}
S = \dfrac{1-TBV_{j}^{(WZ)}}{\max(TBV_{j}^{(WZ)},10^{-20})} ,
\end{equation}  

where the total boundary variation of the target cell by WENO-Z reconstruction is given by
\begin{equation}
 TBV_{j}^{WZ}=\frac{\bigg( U^{L, WZ}_{j-\frac{1}{2}}-U^{R, WZ}_{j-\frac{1}{2}}\bigg)^4 + \bigg(U^{L, WZ}_{j+\frac{1}{2}}-U^{R, WZ}_{j+\frac{1}{2}} \bigg)^4}{\big(\hat{U}_{j}-\hat{U}_{j-1}\big)^4+\big(\hat{U}_{j}-\hat{U}_{j+1}\big)^4+10^{-20}}, 
\end{equation}  
\item Modify the admissible reconstruction via following equation
\begin{equation}\label{eqn:bvd2}
\begin{aligned} 
{\rm if } \ \ S_j < 1\times10^{6} ~~~ 
{U^{L,R}}_{j+\frac{1}{2}}&={U^{(L,R) WZ}}_{j+\frac{1}{2}}
\end{aligned} 
\end{equation}

and the threshold value $1\times10^{6}$ is used to determine the non-smooth solution as considered in \cite{ivan2014high}. The improved results are shown in Fig. \ref{fig:algo2} indicates that we can construct BVD algorithms that can effectively suppress the oscillations and also preserve the smooth solutions. Similar improved results are observed for one-dimensional Euler equations as well, but we noticed that the results are not improved for two-dimensional cases.
\end{enumerate}

\begin{figure}[H]
\begin{onehalfspacing}
\centering
\subfigure[]{\includegraphics[width=0.48\textwidth]{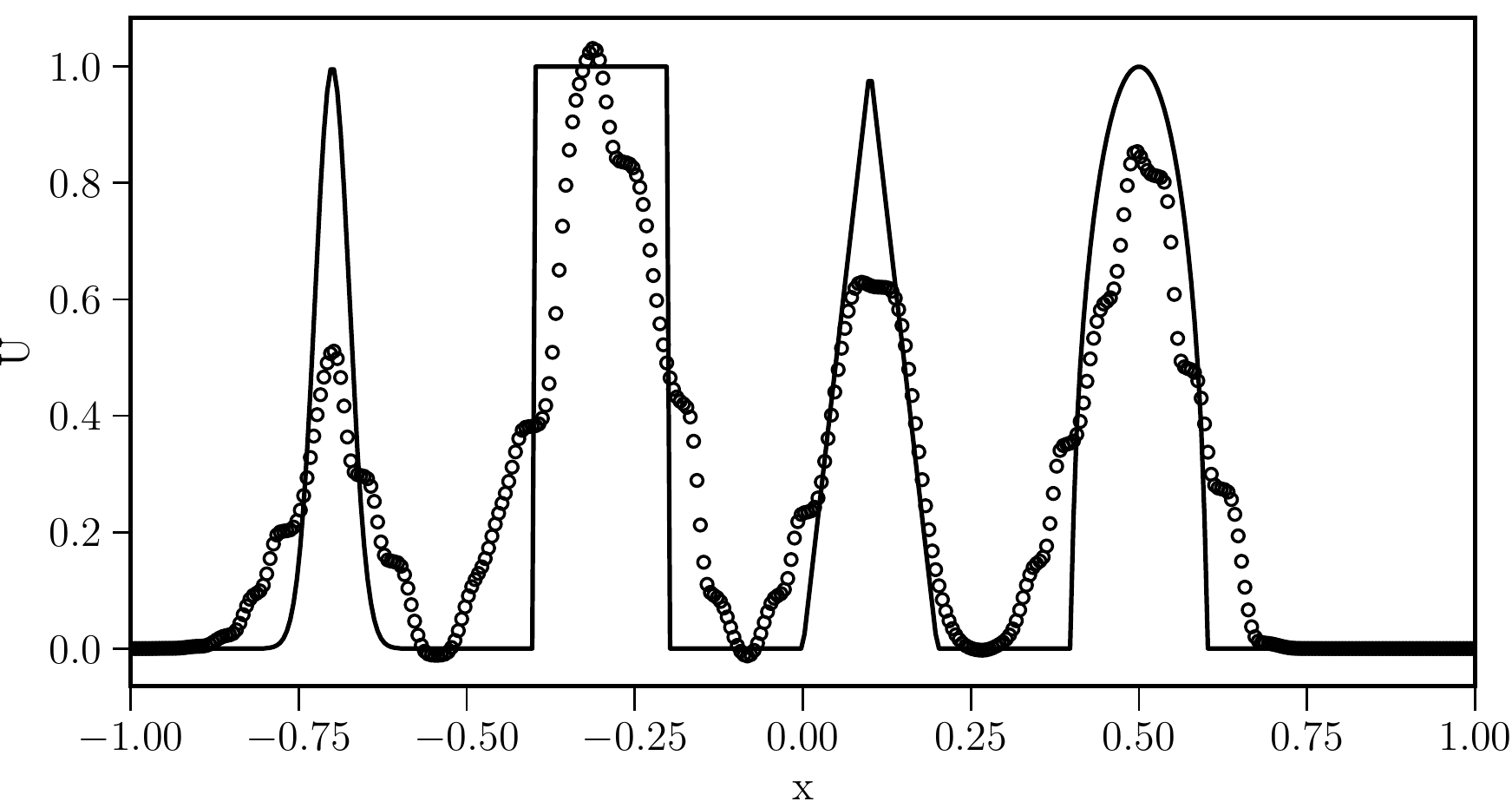}
\label{fig:algo1}}
\subfigure[]{\includegraphics[width=0.48\textwidth]{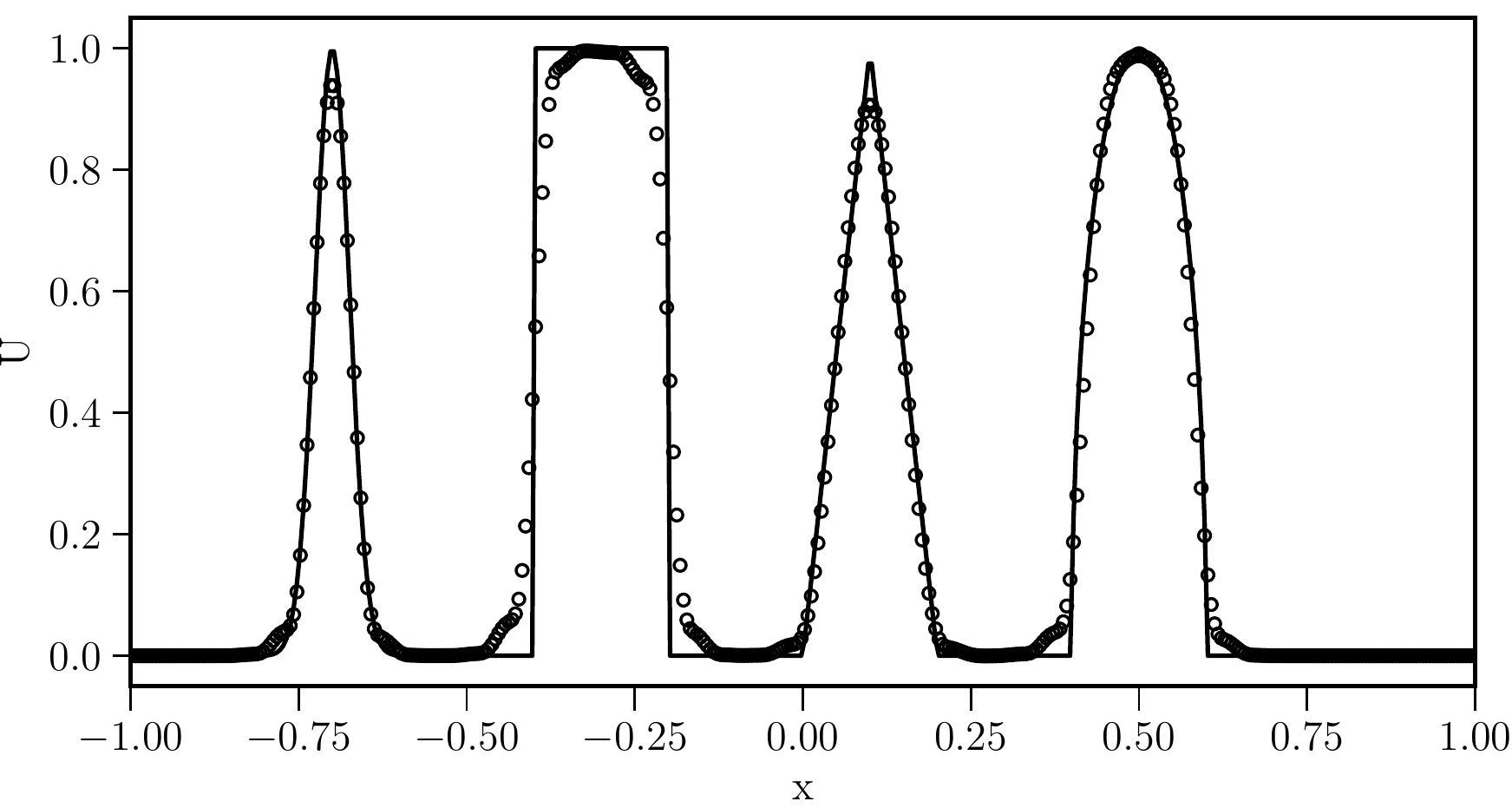}
\label{fig:algo2}}
\caption{Advection of complex waves in Example \ref{ex:adv_comp} (a) by BVD algorithm with Equation (\ref{eqn:BVDstep}) and (b) with Equation (\ref{eqn:bvd2}), by combining WENO-Z and linear central compact scheme for t = 500.}
\label{fig_appen-advec}
\end{onehalfspacing}
\end{figure}

\subsection{\textcolor{black}{Accuracy for critical points test case}}\label{sec-appd}
\textcolor{black}{In Step-4 of the HOCUS algorithm discussed in \ref{sec2.2} we have incorporated an additional criteria given by Equation \ref{eqn:BVDstep-extra} and the results are shown in Table \ref{tab:extra-cond}.}

\begin{equation}\label{eqn:BVDstep-extra}
\begin{aligned} 
{\rm if } \ \ TBV_{j}^{MP5} < TBV_{j}^{C5} \ \  \& \ \  \left({u}_{i+1}-{u}_{i}\right)\left({u}_{i}-{u}_{i-1}\right)\textcolor{magenta}{<} 0 
\end{aligned} 
\end{equation}

\textcolor{black}{
\begin{table}[H]
  \centering
  \footnotesize
 \caption{\textcolor{black} {$L_1$ errors and numerical orders of accuracy on $u_t$ + $u_x$ = 0 with $u_0$(x) = sin($\pi$x - sin($\pi$x)/$\pi$ ). $N$ is the total number of cells on a uniform mesh and $t = 2$.}}
    \begin{tabular}{ccccccc}
    \hline
    N     & HOCUS6 &       & C6 &       & HOCUS6 - Extra &  \\
    \hline
          &       & order &       & order &       & order \\
          \hline
    10    & 6.70E-03 &       & 4.02E-03 &       & 4.02E-03 &  \\
    \hline
    20    & 1.91E-04 & 5.13  & 3.40E-05 & 6.89  & 3.40E-05 & 6.89 \\
    \hline
    40    & 3.99E-06 & 5.58  & 8.84E-07 & 5.26  & 8.84E-07 & 5.26 \\
    \hline
    80    & 7.54E-08 & 5.72  & 2.06E-08 & 5.43  & 2.06E-08 & 5.43 \\
    \hline
    160   & 1.26E-09 & 5.90  & 4.69E-10 & 5.46  & 4.68E-10 & 5.46 \\
    \hline
    \end{tabular}%
  \label{tab:extra-cond}%
\end{table}}

%----------------------------------------%----------------------------------------%----------------------------------------%----------------------------------------%----------------------------------------
%----------------------------------------%----------------------------------------%----------------------------------------%----------------------------------------%----------------------------------------
%----------------------------------------%----------------------------------------%----------------------------------------%----------------------------------------%----------------------------------------
%----------------------------------------%----------------------------------------%----------------------------------------%----------------------------------------%----------------------------------------
%----------------------------------------%----------------------------------------%----------------------------------------%----------------------------------------%----------------------------------------
%----------------------------------------%----------------------------------------%----------------------------------------%----------------------------------------%----------------------------------------
%\section*{References}
\bibliographystyle{elsarticle-num}
\bibliography{bvd_ref}
\end{document}